\pdfoutput=1
\documentclass[12pt]{report}
\usepackage[titletoc]{appendix}
\usepackage{amsmath, amsthm, amsfonts}
\usepackage{a4wide}
\usepackage{amsmath, array}
\usepackage{url}
\usepackage{enumitem}
\usepackage{calligra,amsmath,amssymb}
\usepackage{scalerel}
\usepackage{booktabs}
\usepackage{multicol}
\usepackage{tikz-cd}
\usepackage{tikz}
\usetikzlibrary{arrows, decorations.text, decorations.pathmorphing, decorations.markings, positioning, shapes}
\usetikzlibrary{}
\usetikzlibrary{calc}
\usetikzlibrary{decorations.markings}
\usepackage{feynmp}
\tikzset{
	fermion/.style = {draw = black, postaction = {decorate},decoration = {markings,mark = at position .55 with {\arrow{>}}}},
	vertex/.style = {draw,shape = circle,fill = black,minimum size = 2pt,inner sep = 0pt},
	fermionbar/.style={draw=black, postaction={decorate},
		decoration={markings,mark=at position .55 with {\arrow{<}}}}}

\usepackage[all,knot,poly,color]{xy}

\usepackage[T1]{fontenc}
\usepackage[latin1]{inputenc}

\usepackage[margin=2cm]{geometry}
\usepackage[font=small,labelfont=bf,tableposition=top]{caption}
\usepackage[official,right]{eurosym}
\usepackage{rotating}

\usepackage{graphicx}

\usepackage{relsize}

\usepackage{epsfig,graphics}
\usepackage[all,knot,poly,color]{xy}

\usepackage{epigraph}

\newtheorem{thm}{Theorem}[section]
\newtheorem{cor}[thm]{Corollary}

\newtheorem{lem}[thm]{Lemma}
\newtheorem{prop}[thm]{Proposition}
\theoremstyle{definition}
\newtheorem{defn}[thm]{Definition}
\newtheorem{rem}[thm]{Remark}
\newtheorem{ex}[thm]{Example}
\newtheorem{co}[thm]{Conjecture}
\def\RR{\mathbb{R}}
\def\ZZ{\mathbb{Z}}
\def\NN{\mathbb{N}}
\def\CC{\mathbb{C}}

\newcommand{\gra}{{\alpha}} \newcommand{\grb}{{\beta}} \newcommand{\grg}{{\gamma}} \newcommand{\grd}{{\delta}}
 \newcommand{\grz}{{\zeta}}  \newcommand{\gru}{{\theta}}
\newcommand{\gri}{{\iota}}  \newcommand{\grl}{{\lambda}} \newcommand{\grm}{{\mu}}
   \newcommand{\grp}{{\pi}}
 \newcommand{\grs}{{\sigma}} \newcommand{\grt}{{\tau}} 
\newcommand{\grf}{{\phi}} \newcommand{\grx}{{\chi}} \newcommand{\grc}{{\psi}} \newcommand{\grv}{{\omega}}

\newcommand{\grF}{{\Phi}}  \newcommand{\grC}{{\Psi}} 
\newcommand{\arw}{\rightarrow} 
\newcommand{\an}{\Rightarrow} 

\renewcommand{\bibname}{References}

\begin{document}
	\begin{titlepage}
\thispagestyle{empty}
\hspace{1em}
\begin{center}

\includegraphics[width=0.7\linewidth]{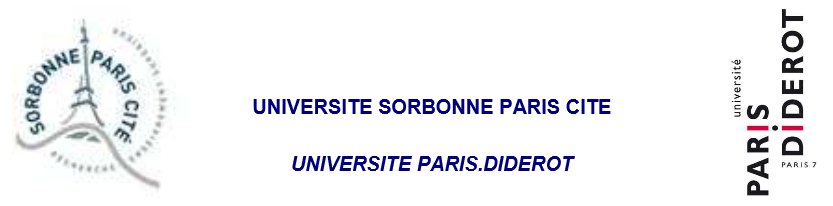}

\vspace*{0.8cm}

 {\Large \bf TH\`ESE DE DOCTORAT}

\vspace*{0.5cm}
Discipline: {\bf Math\'ematiques}\;\;
Ecole Doctorale: {\bf ED 386}\;\;
Laboratoire: {\bf IMJ-PRG}

\vspace*{0.8cm}

Pr\'esent\'ee par

\vspace*{0.5cm}

{\Large \bf EIRINI CHAVLI}

\vspace*{0.8cm}
Pour obtenir le grade de \ \\[1ex]
{\bf DOCTEUR DE L'UNIVERSIT\'E PARIS DIDEROT-PARIS 7} \ \\

\vspace*{0.8cm}
{\Large {\bf{The Brou\'e-Malle-Rouquier conjecture for the exceptional groups of rank 2}\\ }}

\vspace*{0.3cm}
{\Huge{***}}\\
\vspace*{0.008cm}
{\Large {\bf{La conjecture de Brou\'e-Malle-Rouquier pour les groupes exceptionnels de  rang 2}}}\\
  \end{center}
  \vspace*{0.5 cm}
    \begin{tabular} {r@{\ }lll}&Th\`ese dirig\'ee par:& \textbf{M. Ivan MARIN}&\textbf{Universit\'e de Picardie Jules Verne}\\ \\
    	&Rapporteurs: &\textbf {M. Pavel ETINGOF}& \textbf{Massachusetts Institute of Technology}\\
  	& &\textbf{M. Louis FUNAR} & \textbf{Universit\'e Grenoble Alpes}
  \end{tabular}
  \begin{center}
\vspace*{0.2cm} 
Soutenue le \textbf{12 mai 2016} devant le jury compos\'e de :  \\
\vspace*{0.5cm} 
\begin{tabular}{r@{\ }lll}
     & \textbf{M. Fran\c{c}ois DIGNE} & \textbf{Universit\'e de Picardie Jules Verne} &\textbf{Examinateur}\\
      & \textbf{M. Louis FUNAR} & \textbf{Universit\'e Grenoble Alpes} &\textbf{Rapporteur}\\
     & \textbf{M. Nicolas JACON} & \textbf{Universit\'e  de Reims} &\textbf{Examinateur}\\
     & \textbf{M. Gunter MALLE} & \textbf{TU Kaiserslautern} &\textbf{Examinateur}\\
      & \textbf{M. Ivan MARIN} & \textbf{Universit\'e de Picardie Jules Verne} &\textbf{Directeur de th\`ese}\\
     & \textbf{M. G\"{o}tz PFEIFFER} & \textbf{National University of Ireland,
     	Galway} &\textbf{Examinateur}
\end{tabular}

\end{center}

\end{titlepage}
	\clearpage
	\thispagestyle{empty}
	\newpage
	~
\clearpage	
\setcounter{page}{1}
	\chapter*{}
	\epigraph{\large{``A lesson without the opportunity for learners to generalise  is not a mathematics lesson.''}}{\large{J. Mason, 1996, p. 65}}
	\newpage
	~
\chapter*{Acknowledgments}

\indent

First and foremost, I would like to thank my supervisor for his support and for the trust he showed in me all these years of the Ph.D. I want to thank him for speaking to me in English until I learned French. For informing me about important conferences and seminars. For teaching me how to research properly, respecting my ideas and my disposition, and for letting me work in my way - he offered guidance but never interfered with my research. He was always fast to reply to my emails, and I need to thank him extra for the last year of the Ph.D. because of the sheer volume of my messages! He also helped me with the correction, editing, and revising of my article, even though that was not among his responsibilities. I feel incredibly grateful that he always made me feel at ease, and that I never felt intimidated to ask questions, which enhanced my self-confidence. His method of making me explain and analyze my ideas and the results that I had come to, even the obvious ones, helped me achieve greater clarity and validate my work. I feel that his students are his priority in many respects, and proof for this, among others, is the seminar he organized this year for his Ph.D. students. This seminar helped me a lot to improve the way I present my results. Furthermore, when I hit dead ends or found obstacles during my thesis, he was always calm and positive, and never stressed me or worsened my concerns. Finally, and most importantly, I would like to thank him because apart from being a teacher and a colleague, he has also been a friend and he has most successfully managed to keep all these features in balance.

I would like to thank P. Etingof and L. Funar for spending time reading, accessing and writing a report about my thesis. I would also like to thank P. Etingof for spending the time to explain to me his article about the weak version of the freeness conjecture when I met him in Corsica, and L. Funar for consenting to travel all the way to Paris to be part of the panel.  

I would like to thank all the other members of the committee, for having accepted to evaluate my thesis. Specifically, I would like to thank:
F. Digne, for sharing wholeheartedly his office with me in Amiens this year. N. Jacon, for always being positive, always having something good to say about my thesis, and for being such an encouraging person overall. G. Malle, for inviting me to Kaiserslautern and giving me the chance to be a speaker at the seminar held there. G. Pfeiffer, because his book was an inspiration to me and a point of reference when I was writing my article. Also, I would like to thank him for traveling from Ireland to be here for my defense. 

I would like to thank the people of the University of Amiens for giving me an office during my second and third year of my thesis, even though my position was in Paris. Also, I would like to thank them for granting me the ATER position this year. Specifically, I would like to thank I. Wallet and C. Calimez for offering all the necessary information, for emailing me very thoroughly and in a personal and helpful way regarding all my questions and affairs, for accommodating me with the bureaucracy, for treating me as though I were an actual member of this University, and for always being friendly and pleasant.

Furthermore, I would like to thank all the people I worked with as an ATER:  D. Chataur , Y. Palu, and K. Sorlin, who were more than eager to relieve me of my supervisions when I needed more time to work on my thesis. S. Ducay and E. Janvresse, who helped me with the schedule of the courses that I taught. Last but not least, R. Stancu for his support and collaboration, his kindness and his trust with the course of Galois theory for advanced students.

I would like to thank all the people from Paris 7 who supported me during my thesis. Specifically, E. Delos, who was in charge of the  ``ecole doctorale'', that always replied to my emails and helped me organize the presentation. O. Dudas, who helped me with my queries, offered a lot of support, and was more than willing to drive me around during a conference in Corsica, where I had a problem with my leg. M. Hindry, who accommodated me with the train routes and was very supportive of the fact that I had to commute from Amiens to Paris, when I was teaching the exercise course under his supervision. C. Lavollay, who was the administrative assistant in charge of the train fees reimbursements, and has always been very positive, helpful, indulging and pleasant and made me feel included, by inviting me to social events, such as the  ``cr\^epe parties''! J. Michel, with whom I had great discussions regarding the factorized Schur elements in GAP, and who also provided me with the reference letter I needed. 

I would like to thank M. Chlouveraki. On a professional level, she gave me valuable advice, she informed me about available post-doctorate positions, gave me a copy of her book which has been very useful, helped me with my queries about the decomposition maps, and we had constructive discussions about them.
On a more personal level, I want to thank her for the ``drinking Fridays'' I had with her and Christina. I also want to thank her for a great time we had in Italy and the fact that, to be in my presentation, she had to reschedule her holidays and leave her little son with her mother.

I would like to thank my French teachers, A. Le Corre and C. Le Dilly, for having the wish to attend my presentation. They are excellent teachers, and they have been particularly supportive to me, especially in the beginning when I couldn't speak the language at all. They encouraged me a lot, and they taught me not only the grammar or the vocabulary but also the French culture. They made me feel comfortable about speaking the language, and this motivated me to start teaching during my thesis.

I would like to thank all the friends I made in France. Specifically, I thank with all my heart Claire, who let me stay in her place when I had nowhere to go, even though she knew me only one week. Her support and hospitality during my first year in Paris were invaluable to me. I would like to thank Sandrine for being such a good and supporting friend. I thank her for all the moments we had in Paris and also for the hospitality at her mother's place at Gex, that I will always remember. She and Lorick, that I also thank a lot, always had a place for me and my husband, speaking in English when necessary. I would like to thank Julie, who was always there for me. Japanese food and ballet nights and Shrek movies are some of our activities that were an oasis for me during my thesis. I also thank her for traveling from Paris to Amiens, just to be here for my 30th birthday. I would like to thank Huafeng for the mathematical discussions and the lovely walks we used to do when I was living in Paris.

I would like to thank some Greek friends in Paris. First of all, my dear Christina, who was an inspiration for me to come to France. If I hadn't visited her the May of 2010, I wouldn't have decided to apply in Paris 7. She let me stay in her place when I first arrived in Paris, even though she was doing her Ph.D. at the time and she was very busy. She is a person that I can count on, who supported me during my thesis and encouraged me by believing in me. I will always remember all the moments we had together. Secondly, I would like to thank Eugenia, who always had a place for me to stay when visiting Paris. I also thank her for visiting me in Amiens twice, even not having a lot of money at the time. I would like to thank Evita, who helped my organize my wedding and for her willing to come to my presentation.

I would like to thank all the Ph.D. students I met in Amiens, starting with my  ``thesis' brothers''; I would like to thank Alexandre, who carefully read this manuscript and corrected some English mistakes, and who also helped me to write my introduction in French. I thank him for his support and for our discussions, which always made me feel better. I also thank Georges for his support and kindness, the great time we had at Lille, the funny times in my office and the fact that he reminds me more than once that ``c'est la vie pas le paradis''! I want to thank them both for all the time we spent at the ``Atelier d'Alex'', even if this year I had a lot of work and I couldn't always be a part of the fun.  I thank Malal, Vianney, Ruxi, Silvain, Emna, and Clara for the moments we had at the office BC01, Lucie for her help every time I needed it, Pierre for our nice conversations, Anne-Sophie for her help with my documents I had to send to Paris, Valerie for her kindness, Maxime for all his support and his sense of humor, Aktham for the delicious fruits he offered me and his kind words. I also thank Thomas for his help with French translations and his belief in me.

I would like to thank Emilie for the nice moments we had in Lille and Caen and for traveling from far away just to be a part of my defense. 

I would like to thank my ``dance friends'' Anna, Clio, Sarah, and Cassis for all the good times we had and all the girls nights we organized. I also thank Clio for the `` chocolat chaud'' times I had with her and little Pia, our pleasant discussions and her patience with my French. She is a friend I will always remember.

I want to thank all my childhood friends Eleni, Eleutheria, Dina, Dora, Roula, Tasos, Christos, Giorgos, who always support me and even if I am so far from them I feel so close. I would like to thank Nikos separately, for being such a great friend all these years. Always someone who supports me and believes in me, the person who I could say understands my fears, insecurities and always finds a way to make me feel better. I thank him for always being there for me.

I would like to thank my former students,  Olga and Efi, who went on to become my friends, for inspiring me to be the teacher I am today, for their love and support. I also thank Maria IB, my first student that is also an excellent friend now. I thank her for all the time we had together, her sense of humor, her attitude to the difficulties that I met during my thesis, which always helped me to smile and keep going.

I would like to thank Niki for her support and trust all these years and for all the happy times we have in Greece every time I go back. Moreover, I would like to thank Vasso, who is such a unique friend. A person who can be happy with my happiness and sad with my sadness. Her phone calls, her emails and messages all these years I live abroad are a great support for me. I thank her for being available anytime I need her, to calm me when I feel stressed, to give me strength when I need it. 

I would like to thank Fotini, who has been not only a great piano teacher but also a great friend. I thank her for believing that I can achieve my goals, for being supporting and caring and for all the sweet and refreshing moments I have with her and her family, especially the ``evil portions'' we make with her little son, Diogenis, every time I am in Greece! I also thank her because, even if she didn't like the idea of new technologies, she downloaded Skype and Viber, just to be in touch with me.

I would like to thank Maria ``the frog'', my best friend since we were eight years old. I thank her for supporting me during this thesis, even if she hates maths! She is like the sister I never had, always a voice in my head to say to me that I have to keep going. I thank her because even if I am in France, she always finds ways to keep in touch as if I were in Greece. I would like to thank her for the times she visited me in France and for her help during my wedding. Wherever I go, I always have her with me.

I would like to thank Maria ``piano'' for all her support and kindness. She is the friend who made all my wedding bonbonniere on her own because I didn't have the money to buy the ones I wanted; who is always there to talk to the phone when I need company; who makes me her priority, who traveled from Greece so many times to visit me, even when she didn't have the money to do so. I thank her for being here for my presentation today. I would like to thank her also for all the moments we had together, for her help with my English letters that I had to write sometimes. Most importantly, I thank her because she is one of the friends that makes me feel like I never left home. 

I thank Giorgos, who is always near my family and who traveled today from Greece in a last-minute notice. I am happy that he wanted to be a part of this experience. I thank him for his encouragement and support.

I would like to thank my parents-in-law for his love and support, for the long phone calls when I miss Greece and the feeling I have with them that I am always an important part of their family. I would like to thank Annoula for our Facebook conversations, her understanding and support and her belief in me. I also thank her for all the great moments we had in her place, together with my niece and nephew, and, of course, for the tomato salads she makes me! I thank Korina for being like the young sister I never had and for her sense of humor that makes our family meetings so unique. I thank her for traveling to Amiens to visit us, and we had such a great time. 

I would like to thank with all my heart my family, because if it weren't for them, I wouldn't have been the person I am today. I thank my grandma for being a second parent to me, for her support, love and calmness and for her ability to make me smile every time we talk. I would like to thank my mother for giving me freedom to chose my path, for trusting me and supporting my choices, even if she didn't approve or understand all of them. I thank her for her love, her kindness, and her smile all these years. I also thank her for making me feel like living in a five-star hotel every time I travel back home. I would like to thank my brother, Spyros, who traveled from Greece to be here today, even if he had a lot of work with his Ph.D. thesis. I thank him for being a friend; someone I can rely on. Most importantly, I thank him for feeling proud of me since we were kids. His belief in me made me achieve a lot all these years. I also thank him for his hospitality in Crete this summer, that made me relax and be able to finish my thesis on time. 

Last but not least, I would like to thank Angelos for all his love and support during my thesis. He accepted to leave Greece and to come to a country where he doesn't speak the language, just to be with me. I thank him for being understanding about all the weekends I had to work, about all the conferences I had to attempt and I wasn't home, and about the times I was complaining because my thesis was at a dead-end. I thank him for all the  ``don't worries'' he told me and his willing to adapt to every new situation my job will put us on. I thank him for being able to see beauty, strength, and capability in me, even if sometimes I cannot see it in myself.
	
	\chapter*{}

\textbf{Abstract:} Between 1994 and 1998, the work of M. Brou\'e, G. Malle, and R. Rouquier generalized in a natural way the definition of the Hecke algebra associated to a finite Coxeter group, for the case of an 
arbitrary complex reflection group. Attempting to also generalize the properties of the Coxeter
case, they stated a number of conjectures concerning these Hecke algebras. One specific example
of importance regarding those yet unsolved conjectures is the so-called BMR freeness conjecture. 
This conjecture is known to be true apart from 16 cases, that are almost all the exceptional groups of rank 2. These exceptional groups of rank 2 fall into three families: the tetrahedral, octahedral and icosahedral family.
We prove the validity of the BMR freeness conjecture for the exceptional groups belonging to the first two families, using a case-by-case analysis and we give a nice description of the basis, similar
to the classical case of the finite Coxeter groups. We also give a new consequence of this conjecture, by obtaining
the classification of irreducible representations of the braid group on 3 strands in dimension at most 5, recovering results of Tuba and Wenzl. \\ \\ 
\textbf{keywords:} Cyclotomic Hecke algebras, BMR freeness conjecture, Complex Reflection Groups, Braid groups, Representations.\\
\begin{center}
\vspace*{0.2cm}
{\Huge{***}}\\
\vspace*{0.008cm}
\end{center}
\textbf{R\'esum\'e:} Entre 1994 et 1998, M. Brou\'e, G. Malle et R. Rouquier ont g\'en\'eralis\'e aux groupes de r\'eflexions complexes la d\'efinition naturelle des alg\`ebres de Hecke associ\'ees aux groupes de Coxeter finis. Dans la tentative de g\'em\'eraliser certaines propri\'et\'es de ces alg\`{e}bres, ils ont annonc\'e des conjectures parmi lesquelles la conjecture importante de libert\'e de BMR. Il est connu que cette derni\`{e}re conjecture est vraie  sauf pour 16 cas qui concernent presque tous les groupes exceptionnels de rang 2.  Ces derniers se plongent dans 3 familles : t\'etra\'edrale, octa\'edrale et icosa\'edrale. Nous prouvons que la conjecture de libert\'e de BMR est vraie pour les groupes exceptionnels appartenant aux deux premi\`{e}res familles en utilisant un raisonnement cas par cas et en donnant une jolie description de la base, ce qui est similaire au cas classique d' un groupe de Coxeter fini. Nous donnons aussi une nouvelle cons\'equence de cette conjecture qui est l'obtention de la classification des repr\'esentations irr\'eductibles du groupe de tresses \`{a} 3 brins de dimension au plus 5, retrouvant ainsi des r\'esultats de Tuba et Wenzl. \\\\
\textbf{Mot cl\'es:} Alg\`ebres de Heckes cyclotomiques, Conjecture de libert\'e de BMR, Groupes de r\'eflexions complexes, Groupes de tresses, Repr\'esentations.

	\newpage
	~
	\tableofcontents

	\chapter*{Introduction}
	\addcontentsline{toc}{chapter}{Introduction}
	\indent 

\emph{Real reflection groups}, also known as finite Coxeter groups, are finite groups of matrices with
real coefficients generated by reflections (elements of order 2 whose vector space of fixed points
is a hyperplane). We often encounter finite Coxeter groups in commutative algebra, Lie theory, representation theory, singularity theory, crystallography, low-dimensional topology, and geometric group theory. They are also a powerful way to construct interesting examples in geometric and combinatorial group theory. 

All finite Coxeter groups are particular cases of \emph{complex reflection groups}. Generalizing the definition of the real reflection groups, complex reflection groups are finite groups of 
matrices with complex coefficients generated by pseudo-reflections (elements of finite order
whose vector space of fixed points is a hyperplane). Any complex reflection group can be decomposed as a direct product of the so-called irreducible ones (which means
that, considering them as subgroups of the general linear group $GL(V)$, where $V$ is a finite dimensional complex vector space, they act irreducibly on $V$). The irreducible complex reflection
groups were classified by G. C. Shephard and J. A. Todd (see \cite{shephard}); they
belong either to the infinite family $G(de, e, n)$ depending on 3 positive
integer parameters, or to the 34 exceptional groups, which are numbered
from 4 to 37 and are known as $G_4,\dots, G_{37}$, in the Shephard and Todd classification. The infinite family $G(de, e, n)$ is the group of monomial matrices whose non-zero entries are $de^{\text{th}}$
roots of unity and their product is a $d^{\text{th}}$
root of unity.

Subsequent work by a number of authors has proven that complex reflection groups have properties analogous to those of real reflection groups, such as presentations, root systems, and generic degrees. On any real reflection group, we can associate an Iwahori-Hecke algebra (a one parameter deformation of the group algebra of the real reflection group). Hecke algebras associated to reductive groups were introduced in order to decompose representations of these groups induced from parabolic subgroups.

Between 1994 and 1998, M. Brou\'e, G. Malle, and R. Rouquier  
generalized in a natural way the definition of the Iwahori-Hecke algebra to   arbitrary complex reflection groups (see \cite{bmr}). Attempting to also generalize the properties of the Coxeter case, they stated a number of conjectures concerning the Hecke algebras, which haven't been proven yet. Even without being proven, those conjectures have been used by a number of papers in the last decades as assumptions, and are still being used in various subjects, such as representation theory of finite reductive groups, Cherednik algebras, and usual braid groups (more details about these conjectures and their applications can be found in \cite{marinreport}).

One specific example of importance, regarding those yet unsolved conjectures, is the so-called freeness conjecture. In 1998, M. Brou\'e, G. Malle and R. Rouquier conjectured that the generic Hecke algebra $H$ associated to a complex reflection group $W$ is a free module of rank $|W|$ over its ring of definition $R$. They also proved that it is enough to show that $H$ is generated as $R$-module by $|W|$ elements. 

 The validity of the conjecture, even in its weak version (which states that $H$ is finitely generated as $R$-module), implies that by extending the scalars to an algebraic closure field $F$ of the field of fractions of $R$, the algebra $H\otimes_R F$ becomes isomorphic to the group algebra $FW$ (see \cite{marincubic} and \cite{marinG26}).  
G. Malle assumed the validity of the conjecture and used it to prove that the characters of $H$ take their values in a specific field (see \cite{malle1}). Moreover, he and J. Michel also used this conjecture to provide matrix models for the representations of $H$, whenever we can compute them; these matrices for the generators of $H$ have entries in the field generated by the corresponding character values (see \cite{mallem}).

The freeness conjecture is fundamental in the world of generic Hecke algebras. Once it is proved, our better knowledge of these algebras could allow the possibility of using various computer algorithms on the structure constants for the multiplication, in order to thoroughly improve our understanding in each case (see for example \S 8 in \cite{mallem} about
the determination of a canonical trace).  

The freeness conjecture has also many applications, apart from the ones connected to the properties of the generic Hecke algebra itself. Provided that the freeness conjecture is true, the category of representations of $H$ is equivalent to a category of representations of a Cherednik algebra (see \cite{Ocat}). Another application is about the algebras connected to cubic
invariants,  including  the  Kaufman  polynomial  and  the  Links-Gould  polynomial. These algebras are quotients of the generic Hecke algebra associated to $G_4$, $G_{25}$ and $G_{32}$. I. Marin used the validity of the conjecture of these cases and he proved that the generic algebra $K_n(\gra, \grb)$ introduced by P. Bellingeri and L. Funar in \cite{bfunar} is zero for
$n\geq 5$ (see theorem 1.4 in \cite{marincubic}). Furthermore, in \cite{chavli} we used the freeness conjecture for the cases of $G_4$, $G_8$ and $G_{16}$  to recover and explain a classification due to I. Tuba and H. Wenzl (see \cite{tuba}) for the irreducible representations of the braid group on 3 strands of dimension at most 5
(we explain in detail this result in Chapter 2).

The freeness
conjecture that we call here the BMR freeness conjecture, is known to be true for the finite Coxeter groups, and also for the infinite series by Ariki and Koike (see \cite{ariki} and \cite{arikii}). Considering the exceptional cases, one may divide them into two families; the family that includes the exceptional groups $G_4, \dots, G_{22}$, which are of rank 2 and the family that includes the rest of them, which are of rank at least 3 and at most 8. Among the second family we encounter 6 finite Coxeter groups for which we know the validity of the conjecture: the groups $G_{23}$, $G_{28}$,
$G_{30}$, $G_{35}$, $G_{36}$ and $G_{37}$. Thus, it remains to prove the conjecture for 28 cases: the exceptional groups of rank 2 and the exceptional groups $G_{24}$, $G_{25}$, $G_{26}$ $G_{27}$,
$G_{29}$, $G_{31}$, $G_{32}$, $G_{33}$ and $G_{34}$.

Until recently, it was
widely believed that the BMR freeness conjecture had been verified for most of
the exceptional groups of rank 2. However, there were imprecisions in the proofs, as I. Marin indicated a few years ago (for more details see the introduction of \cite{marinG26}). In the following years, his research and his work with G. Pfeiffer concluded that the exceptional complex reflection groups for which there is a complete proof for the freeness conjecture  are the groups $G_4$ (this case has also been proved in \cite{brouem} and independently in \cite{funar1995}), $G_{12}$,  $G_{22}$, $G_{23}$,  \dots, $G_{37}$ (see \cite{marincubic}, \cite{marinG26} and \cite{marinpfeiffer}). The remaining groups are almost all the exceptional groups of rank 2. The main goal of this PhD thesis is to prove the BMR freeness
 conjecture for these remaining cases. 
 
 In the first chapter we give the necessary preliminaries to accurately describe the BMR freeness conjecture. More precisely, we introduce the definition of a complex reflection group and we conclude that the study of these groups reduces to the irreducible case, which leads us to the description of the Shephard-Todd classification. Following \cite{bmr}, we also associate to a complex reflection group a complex braid group, by using tools from algebraic topology. We finish this chapter by giving in detail the description of the BMR freeness conjecture and a brief report on the recent results on the subject by I. Marin and G. Pfeiffer. 
 
 In the second chapter, as a first approach of the BMR freeness conjecture, we focus on the groups that are finite quotients of the braid group $B_3$ on 3 strands. These are the exceptional groups $G_4$, $G_8$ and $G_{16}$ and, since the case of $G_4$ has been proven earlier as we mentioned before, we prove the conjecture for the rest of the cases. This result completes the proof for the validity of the BMR conjecture for the case of the exceptional groups, whose associated complex braid group is an Artin group (see theorem \ref{braidcase}).
 
 In order to prove the validity of the conjecture we follow the idea I. Marin used in \cite{marincubic}, theorem 3.2(3) for the case of $G_4$. This approach is to presume how the basis should look like and then prove that this is in fact a basis. We recall that  $B_3$ is given by generators the braids $s_1$ and $s_2$ and the single relation $s_1s_2s_1=s_2s_1s_2$. We also recall that the center of $B_3$ is the subgroup generated by the element $z=s_1^2s_2s_1^2s_2$ (see, for example, theorem 1.24 of \cite{turaev}) and that the order of the center of the groups $G_4$, $G_8$ and $G_{16}$ is even. What we manage to prove (see theorems \ref{th} and \ref{thh2}) is that the generic Hecke algebra associated to $G_4$, $G_8$ or $G_{16}$ is of the form 
 
 $$
 H=\sum\limits_{k=0}^{p-1}u_1z^k+\sum\limits_{k=1}^{p}u_1z^{-k}+
 u_1\text{``some other elements''}u_1,$$
 where $p:=|Z(W)|/2$ ($W$ denotes the group $G_4$, $G_8$ or $G_{16}$) and  $u_1$ denotes the $R$-subalgebra of $H$ generated by the image of $s_1$ inside $H$ (for the case of $G_4$ this is a deformation of the result presented in theorem 3.2(3) in \cite{marincubic} where we ``replace'' inside the definition of $H_3$, which denotes as $A_3$ there, the element $s_2s_1^{-1}s_2$ by the element $s_2s_1^2s_2$). 
 
  Exploring the consequences of the validity of the BMR freeness conjecture for the cases of $G_4$, $G_8$ and $G_{16}$, in the last section of the second chapter we  prove that we can determine completely the irreducible representations of $B_3$ for dimension at most 5, thus recovering a classification of I. Tuba and H. Wenzl in a more general framework (see \cite{tuba}). Within our approach, we managed to answer their questions concerning the form of the matrices of these representations, as well as the nature of some polynomials that play an important role to the description of these representations. We also explained why their classification  doesn't work if the dimension of the representation is larger than 5.
 
 Observing that the center of $B_3$ plays an important role in the construction of the basis
 for the Hecke algebras of $G_4$, $G_8$ and $G_{16}$, we developed a more general
 approach to the problem, which uses the results of P. Etingof and E. Rains on the center of
 the braid group associated to a complex reflection group (see \cite{ERrank2}). In chapter 3 we explain in detail these results and we also rewrite some parts of their arguments that are sketchy there. Their approach also provides the weak version of the conjecture for the exceptional groups of rank 2. We conclude by using this weak version to prove a proposition stating that if the Hecke algebra of every exceptional group of rank 2 is torsion free as module, it is enough to prove the validity of the conjecture for the groups $G_7$, $G_{11}$ and $G_{19}$, known as the maximal exceptional groups of rank 2. Unfortunately, this torsion-free assumption does not appear to be easy to check a priori.
 
 In our last chapter, we prove the validity of the BMR freeness conjecture for the exceptional groups $G_5, \dots, G_{15}$.  As we explain in detail in chapter 3,  P. Etingof and E. Rains  connected every exceptional group of rank 2 with a finite Coxeter group of rank 3. More precisely, 
if $W$ is an exceptional group of rank 2, then $\overline{W}:=W/Z(W)$ can be considered as the group of even elements in a finite Coxeter group $C$ of rank 3. For every $x$ in $\overline{W}$, we fix a reduced word $w_x$ representing $x$ in  $\overline{W}$. P. Etingof and E. Rains corresponds $w_x$ to an element we denote by $T_{w_x}$ inside an algebra over a ring $R^C$, which they call $A_+(C)$. This algebra is in fact related with the algebra $H$. They also prove that this algebra is spanned over $R^C$ by the elements $T_{w_x}$. 

Motivated by this result, for the cases of  the exceptional groups $G_5, \dots, G_{15}$ we found a spanning set for $H$ over $R$ of $|W|$ elements. More precisely, for every $x\in \overline{W}$ we choose a specific word $\tilde w_x$ that represents $x$ in $\overline{W}$, which verifies some nice properties we give in detail in this chapter. Note that this word is not necessarily reduced. We associate this word to an element $T_{\tilde w_x}$ inside $A_+(C)$ and we denote by $v_x$ its image inside $H$ (see propositions \ref{ERSUR} and \ref{ERRSUR} and corollary \ref{corIVAN}). We set $U:=\sum\limits_{x\in \overline{W}} \sum\limits_{k=0}^{|Z(W)|-1}Rz^kv_x$. The main theorem of this section is that $H=U$. Using this approach we also managed to prove in a different way the case of $G_{12}$ that has already been verified by I. Marin and G. Pfeiffer (see theorem \ref{2case}). 

 The main part of this chapter is devoted to the proof of $H=U$. In this proof we use a lot of calculations, that we tried to make as less complicated and short as possible, in order to be fairly easy to follow. 
 
To sum up,  the freeness conjecture is still open for the groups $G_{17},\dots, G_{21}$. We are optimistic that the methodology we used for the rest of the groups of rank 2 can be applied to prove the conjecture for these cases, too. Since these groups are large we are not sure that we can provide computer-free proofs, as we did for the other cases. However, there are strong indications that with continued research, and possibly with the development of computer algorithms, we can prove these final cases.
   \begin{center}
   	\vspace*{0.2cm}
   	{\Huge{***}}\\
   	\vspace*{0.008cm}
   \end{center}
\indent 

Les \emph{groupes de r\'eflexions r\'eels}, aussi appel\'es groupes de Coxeter finis, sont des groupes finis de matrices \`a coefficients r\'eels engendr\'es par des r\'eflexions (des \'el\'ements d'ordre 2 dont l'espace des points fixes est un hyperplan). On rencontre souvent des groupes de Coxeter en alg\`ebre commutative, th\'eorie de Lie, th\'eorie des repr\'esentations, th\'eorie des singularit\'es, cristallographie, topologie en petite dimension et th\'eorie g\'eom\'etrique des groupes. Ils sont aussi un outil puissant pour construire des exemples int\'eressants en g\'eom\'etrie et en th\'eorie combinatoire des groupes.

Les groupes de Coxeter finis forment une sous-famille des \emph{groupes de r\'eflexions complexes}. G\'en\'eralisant la d\'efinition des groupes de r\'eflexions r\'eels, les groupes de r\'eflexions complexes sont des groupes finis de matrices \`a  coefficients complexes engendr\'es par des pseudo-r\'eflexions( des \'el\'ements d'ordre fini dont l'espace des points fixes est un hyperplan). Tout groupe de r\'eflexion complexe peut \^etre d\'ecompos\'e en un produit direct de groupes de r\'eflexions complexes dits irr\'eductibles (groupes qui en \'etant vu comme un sous-groupe de $GL(V)$ o\`u $V$ est un espace vectoriel complexe de dimension finie, agissent de mani\`ere  irr\'eductible sur $V$). Les groupes de r\'eflexions complexes irr\'eductibles ont \'et\'e classifi\'es par G.C. Shephard et J.A. Todd (voir \cite{shephard}); ils appartiennent soit \`a  la famille infinie $G(de,e,n)$ qui d\'epend de 3 param\`etres entiers naturels soit aux 34 groupes exceptionnels num\'erot\'es de 4 \`a  37 et appel\'es $G_4,\dots, G_{37}$ dans la classification de Shephard et Todd. Les $G(de,e,n)$ de la famille infinie sont les groupes de matrices monomiales dont les coefficients non-nuls sont des racines $de$-i\`emes de l'unit\'e et dont le produit de ces coefficients est une racine $d$-i\`eme de l'unit\'e 

Du travail provenant de nombreux auteurs a prouv\'e que les groupes de r\'eflexions complexes ont des propri\'et\'es analogues \`a celles des groupes de r\'eflexions r\'eels, telles que les pr\'esentations, syst\`emes de racines et degr\'es g\'en\'eriques. On peut associer une alg\`ebre de Iwahori-Hecke \`a tout groupe de r\'eflexion r\'eel (une d\'eformation \`a un param\`etre de l'alg\`ebre de groupe du groupe de r\'eflexion r\'eel). Les alg\`ebres de Hecke associ\'ees aux groupes r\'eductifs ont \'et\'e introduits afin de d\'ecomposer les repr\'esentations induites par des sous-groupes paraboliques de ces groupes.

De 1994 \`a 1998,  M. Brou\'e, G. Malle, et R. Rouquier ont g\'en\'eralis\'e de mani\`ere naturelle la d\'efinition de l'alg\`ebre de Iwahori-Hecke \`a un groupe de r\'eflexion complexe quelconque (voir \cite{bmr}). En tentant de g\'en\'eraliser les propri\'et\'es dans le cas des groupes de Coxeter, ils ont formul\'e plusieurs conjectures concernant les alg\`ebres de Hecke, qui n'ont pas encore \'et\'e prouv\'ees. M\^emes sans \^etre prouv\'ees, ces conjectures ont \'et\'e utilis\'ees dans de nombreux articles de ces derni\`eres d\'ecennies comme hypoth\`eses, et sont toujours utilis\'ees dans de nombreux domaines, tels que la th\'eorie des repr\'esentations des groupes r\'eductifs finis, les alg\`ebres  de Cherednik et les groupes de tresses usuels (plus de d\'etails sur ces conjectures et leurs applications peuvent \^etre trouv\'ees dans \cite{marinreport}).

Un exemple important de ces conjectures non r\'esolues est la conjecture de libert\'e. En 1998, M. Brou\'e, G. Malle et R. Rouquier ont conjectur\'e que l' alg\`ebre de Hecke g\'en\'erique $H$ associ\'ee \`a un groupe de r\'eflexion complexe $W$ est un module libre de rang $|W|$ sur son anneau de d\'efinition $R$. Ils ont aussi prouv\'e qu'il \'etait suffisant de montrer que $H$ est engendr\'e comme $R$-module par $|W|$ \'el\'ements. 

La validit\'e de la conjecture, m\^eme dans sa version faible (qui dit que $H$ est finiment engendr\'ee comme $R$-module), implique que en \'etendant les scalaires \`a une cl\^oture alg\`ebrique $F$ du corps des fractions de $R$, l'alg\`ebre $H\otimes_R F$ devient isomorphe \`a l'alg\`ebre de groupe $FW$ (voir \cite{marincubic} et \cite{marinG26}).  
G. Malle a suppos\'e la conjecture vraie et l'a utilis\'ee pour montrer que les caract\`eres de $H$ prennent leurs valeurs dans un certain corps (voir \cite{malle1}). De plus, lui et J. Michel ont utilis\'e la conjecture pour donner des mod\`eles matriciels pour les repr\'esentations de $H$, d\`es que l'on peut les calculer; les matrices pour les g\'en\'erateurs de $H$ ont des coefficients dans le corps engendr\'e par les valeurs des caract\`eres correspondants (voir \cite{mallem}).

La conjecture de libert\'e est fondamentale dans le monde des alg\`ebres de Hecke g\'en\'eriques. Une fois prouv\'ee, notre compr\'ehension meilleure de ces alg\`ebres pourrait permettre d'utiliser divers algorithmes informatiques sur les constantes structurelles, afin de profond\'ement am\'eliorer notre compr\'ehension dans chaque cas (voir par exemple \S 8 dans  \cite{mallem} \`a propos de la d\'etermination d'une trace canonique).  

La conjecture de libert\'e a aussi de nombreuses applications autres que celles li\'ees aux propri\'et\'es de l'alg\`ebre de Hecke g\'en\'erique. Si la conjecture est vraie alors la cat\'egorie des repr\'esentations de $H$ est \'equivalente  \`a une cat\'egorie de repr\'esentations d'une alg\`ebre de Cherednik (voir \cite{Ocat}). Une autre application concerne les alg\`ebres li\'ees \`a des invariants cubiques, incluant le polyn\^ome de Kaufman et le polyn\^ome de Links-Gould. Ces alg\`ebres  sont les quotients de l'alg\`ebre de Hecke g\'en\'erique associ\'ee  \`a $G_4, G_{25}$ et $G_{32}$. I. Marin a utilis\'e la validit\'e de la conjecture dans ces cas-l\`a et il a prouv\'e que l'alg\`ebre g\'en\'erique $K_n(\gra, \grb)$ introduite par P. Bellingeri et L. Funar dans \cite{bfunar} est nulle pour
$n\geq 5$ (voir th\'eor\`eme 1.4 dans \cite{marincubic}). De plus, dans \cite{chavli} l'hypoth\`ese de conjecture pour les cas $G_4$, $G_8$ et $G_{16}$ est utilis\'ee pour retrouver et expliquer une classification due  \`a I. Tuba et H. Wenzl (voir \cite{tuba}) pour les repr\'esentations irr\'eductibles du groupe de tresse  \`a 3 brins de dimension au plus 5 (on explique ce r\'esultat en d\'etail dans le chapitre 2).

La conjecture de libert\'e que l'on appelle ici conjecture de libert\'e BMR, est vraie et d\'emontr\'ee pour les groupes de Coxeter finis, et pour les s\'eries infinies de Ariki et Koike (voir \cite{ariki} et \cite{arikii}). En ce qui concerne les cas exceptionnels, on peut les s\'eparer en deux familles; la famille des groupes exceptionnels $G_4, \dots, G_{22}$, qui sont de rang 2 et la famille des autres qui sont de rang sup\'erieur ou \'egal \`a 3 et inf\'erieur ou \'egal \`a 8. Dans la seconde famille, on a 6 groupes de Coxeter finis pour lesquels nous savons que la conjecture est valide : les groupes $G_{23}$, $G_{28}$,
$G_{30}$, $G_{35}$, $G_{36}$ et $G_{37}$. Ainsi, il reste \`a prouver la conjecture dans 28 cas : les groupes exceptionnels de rang 2 et les groupes exceptionnels $G_{24}$, $G_{25}$, $G_{26}$ $G_{27}$,
$G_{29}$, $G_{31}$, $G_{32}$, $G_{33}$ et $G_{34}$.

Jusqu'\`a r\'ecemment, il \'etait pens\'e que la conjecture de libert\'e BMR avait \'et\'e v\'erifi\'ee pour la plupart des groupes exceptionnels de rang $2$. N\'eanmoins, il y avait quelques impr\'ecisions dans les preuves, comme indiqu\'e par I. Marin il y a quelques ann\'ees (pour plus de d\'etails, voir l'introduction de \cite{marinG26}). Dans les ann\'ees suivantes, sa recherche et son travail en collaboration avec G. Pfeiffer ont permis de conclure que les groupes de r\'eflexions complexes exceptionnels pour lesquels il y a une preuve compl\`ete de l'hypoth\`ese de conjecture sont les groupes $G_4$ (ce cas a aussi \'et\'e d\'emontr\'e dans \cite{brouem} et de mani\`ere ind\'ependante dans \cite{funar1995}), $G_{12}$,  $G_{22}$, $G_{23}$,  \dots, $G_{37}$ (voir \cite{marincubic}, \cite{marinG26} et \cite{marinpfeiffer}). Les groupes restants sont presque tous des groupes exceptionnels de rang 2. L'objectif principal de cette th\`ese et de d\'emontrer la conjecture pour les cas restants.

Dans le premier chapitre, nous donnons les pr\'eliminaires n\'ecessaires afin de d\'ecrire correctement la conjecture de libert\'e BMR. Plus pr\'ecis\'ement, nous introduisons la d\'efinition d'un groupe de r\'eflexion complexe et et nous ramenons  \`a l'\'etude des groupes de r\'eflexions complexes irr\'eductibles, ce qui nous am\`ene   \`a donner une description de la classification de Shephard-Todd. Comme dans \cite{bmr}, nous associons un groupe de tresse complexe  \`a un groupe de r\'eflexions complexe, en utilisant des outils de topologie alg\'ebrique. Nous terminons ce chapitre en donnant une description d\'etaill\'ee de la conjecture de libert\'e BMR et un r\'esum\'e bref des r\'esultats r\'ecents sur le sujet par I. Marin et G. Pfeiffer. 

Dans le deuxi\`eme chapitre, comme premi\`ere approche de la conjecture de libert\'e BMR, nous nous concentrons sur les groupes qui sont des quotients finis du groupe de tresse   \`a trois brins 
$B_3$. Ce sont les groupes exceptionnels $G_4$, $G_8$ et $G_{16}$ . Le cas $G_4$ ayant \'et\'e d\'emontr\'e comme nous l'avons indiqu\'e pr\'ec\'edemment, nous montrons la conjecture dans le reste des cas. Ce r\'esultat termine la preuve de la conjecture $BMR$ pour le cas des groupes exceptionnels dont le groupe de tresse complexe associ\'e est un groupe de Artin (voir th\'eor\`eme \ref{braidcase}).

Pour prouver la validit\'e de la conjecture, nous proc\'edons de mani\`ere analogue  \`a celle de I. Marin utilis\'ee dans \cite{marincubic}, th\'eor\`eme 3.2(3) pour le cas de $G_4$. Cette approche consiste  \`a estimer une base potentielle et d\'emontrer que c'est effectivement une base. Rappelons que $B_3$ a pour g\'en\'erateurs deux tresses $s_1$ et $s_2$ avec l'unique relation $s_1s_2s_1=s_2s_1s_2$. Rappelons aussi que le centre de $B_3$ est le sous-groupe engendr\'e par l'\'el\'ement $z=s_1^2s_2s_1^2s_2$ (voir, par exemple, th\'eor\`eme 1.24 de \cite{turaev}) et que le cardinal du centre de $G_4$, $G_8$ et $G_{16}$ est pair. Nous parvenons   \`a prouver (voir th\'eor\`emes \ref{th} et \ref{thh2}) que l'alg\`ebre de Hecke g\'en\'erique associ\'ee  \`a $G_4$, $G_8$ et $G_{16}$ est de la forme 

$$
H=\sum\limits_{k=0}^{p-1}u_1z^k+\sum\limits_{k=1}^{p}u_1z^{-k}+
u_1\text{``d'autres \'el\'ements''}u_1,$$
o\`u $p:=|Z(W)|/2$ ($W$ repr\'esente le groupe $G_4$, $G_8$ ou $G_{16}$) et  $u_1$ repr\'esente la sous-$R$-alg\`ebre de $H$ engendr\'ee par l'image de $s_1$ dans $H$ (pour le cas de $G_4$ c'est une d\'eformation du r\'esultat du Th\'eor\`eme 3.2(3) dans \cite{marincubic} o\`u on ``remplace'' dans la d\'efinition de $H_3$, not\'ee $A_3$ dans cet article, l'\'el\'ement $s_2s_1^{-1}s_2$ par l'\'el\'ement $s_2s_1^2s_2$). 

En \'etudiant les cons\'equences de la validit\'e de la conjecture de libert\'e BMR pour le cas de $G_4$, $G_8$ et $G_{16}$, nous d\'emontrons dans la derni\`ere section du deuxi\`eme chapitre que l'on peut d\'eterminer compl\`etement les repr\'esentations irr\'eductibles de $B_3$ de dimension au plus 5, retrouvant ainsi une classification de I. Tuba et H. Wenzl dans un cadre plus g\'en\'eral (voir \cite{tuba}). Dans le cadre de notre approche, nous sommes parvenus  \`a r\'epondre aux questions concernant la forme des matrices de ces repr\'esentations, de m\^eme que la nature de certains polyn\^omes jouant un r\^ole important dans la description de ces repr\'esentations. Nous avons aussi expliqu\'e pourquoi leur classification ne s'adaptait pas  aux repr\'esentations de dimension strictement sup\'erieure  \`a 5.

En observant que le centre de $B_3$ joue un r\^ole important dans la construction de la base pour l'alg\`ebre de Hecke de $G_4$, $G_8$ et $G_{16}$, nous avons \'elabor\'e une approche plus g\'en\'erale au probl\`eme qui utilise les r\'esultats de P. Etingof et E. Rains sur le centre du groupe de tresse associ\'e  \`a un groupe de r\'eflexion complexe (voir \cite{ERrank2}). Dans le chapitre 3, nous expliquons ces r\'esultats en d\'etail et r\'e\'ecrivons une partie de leurs arguments qui sont peu clairs l\`a. Leur approche donne aussi la version faible de la conjecture pour les groupes exceptionnels de rang 2. Nous utilisons cette version faible afin de d\'emontrer une proposition qui dit que si l'alg\`ebre de Hecke de tout groupe exceptionnel de rang 2 est sans-torsion en tant que module, il est suffisant de prouver la validit\'e de la conjecture pour les groupes $G_7$, $G_{11}$ et $G_{19}$, que l'on appelle les groupes exceptionnels de rang 2 maximaux. Malheureusement , cette hypoth\`ese d'\^etre sans-torsion ne semble pas simple  \`a v\'erifier  \`a priori.

Dans le dernier chapitre, nous prouvons que la conjecture de libert\'e BMR est v\'erifi\'ee pour les groupes exceptionnels $G_5, \dots, G_{15}$.  Comme expliqu\'e en d\'etail dans le chapitre 3,  P. Etingof et E. Rains ont li\'e chaque groupe exceptionnel de rang $2$  \`a un groupe de Coxeter fini de rang 3. Plus pr\'ecis\'ement, si $W$ est un groupe exceptionnel de rang 2 alors $\overline{W}:=W/Z(W)$ peut \^etre consid\'er\'e comme le groupe des \'el\'ements pairs d'un groupe de Coxeter fini $C$ de rang 3. Pour tout $x$ dans $\overline{W}$, nous fixons un mot r\'eduit $w_x$ repr\'esentant $x$ dans  $\overline{W}$. P. Etingof et E. Rains font correspondre $w_x$  \`a un \'el\'ement que l'on note $T_{w_x}$ dans une alg\`ebre sur un anneau $R^C$, qu'ils notent $A_+(C)$. Cette alg\`ebre est en r\'ealit\'e li\'ee   \`a l'alg\`ebre $H$. Ils prouvent \'egalement que l'alg\`ebre est engendr\'ee sur $R^C$ par les \'el\'ements $T_{w_x}$. 

En utilisant ce r\'esultat, nous avons trouv\'e un ensemble g\'en\'erateur pour $H$ sur $R$ de $|W|$ \'el\'ements pour les groupes exceptionnels $G_5, \dots, G_{15}$. Plus pr\'ecis\'ement, pour tout $x\in \overline{W}$ nous prenons un certain mot $\tilde w_x$ qui repr\'esente $x$ dans $\overline{W}$, v\'erifiants de bonnes propri\'et\'es que nous donnons en d\'etail dans ce chapitre. Remarquons que ce mot n'est pas n\'ecessairement r\'eduit. Nous associons ce mot  \`a un \'el\'ement $T_{\tilde w_x}$ dans $A_+(C)$ et notons $v_x$ son image dans $H$ (voir propositions \ref{ERSUR} et \ref{ERRSUR} et le corollaire \ref{corIVAN}). On pose $U:=\sum\limits_{x\in \overline{W}} \sum\limits_{k=0}^{|Z(W)|-1}Rz^kv_x$. Le th\'eor\`eme principal de cette section est que $H=U$. En utilisant cette approche, nous avons prouv\'e d'une mani\`ere diff\'erente le cas de $G_{12}$ qui a \'et\'e v\'erifi\'e par I. Marin et G. Pfeiffer (voir th\'eor\`eme \ref{2case}). 

Ce chapitre est principalement d\'edi\'e \`a la preuve de $H=U$. Dans la preuve, nous utilisons beaucoup de calculs, que nous tentons de rendre aussi simples et courts que possible afin de les rendre assez faciles  \`a  suivre.

Pour r\'esumer, l'hypoth\`ese de libert\'e est toujours ouverte pour les groupes $G_{17},\dots, G_{21}$. Nous sommes optimistes sur le fait  que la m\'ethodologie utilis\'ee pour le reste des groupes de rang 2 puisse \^etre appliqu\'ee afin de v\'erifier la conjecture dans ces cas-l\`a aussi. Ces groupes \'etant grands, nous ne sommes pas certains de pouvoir r\'ealiser des preuve sans calculs faits par ordinateur, comme nous l'avions fait dans les autres cas. Cependant, il y a de forts indices qui permettent d'esp\'erer qu'en poursuivant la recherche et en utilisant \'eventuellement des algorithmes informatiques, il sera possible de d\'emontrer les cas restants.
	\newpage
		~
	
	\chapter {The BMR freeness conjecture}
	In this chapter we are going to define the subject of this PhD thesis, which is to check several cases of an important
conjecture of M. Brou\'e, G. Malle and R. Rouquier, that we call the BMR freeness conjecture. We include also the necessary preliminaries to accurately describe it.

\section{Complex Reflection Groups and Complex Braid Groups}
\indent

Let $V$ be a $\CC$-vector space of finite dimension $n$. For all definitions and results we follow mostly \cite{bmr}, \cite{lehrer} and \cite{ivancourse}.
\subsection{Complex Reflection Groups}

\begin{defn} A pseudo-reflection of $V$ is a non-trivial element $s$ of $GL(V)$ of finite order, whose vector space of fixed points is a hyperplane (meaning that $dim_{\CC}\text{ker}(s-1)=n-1$). 
\end{defn}
By definition, one could easily observe that a pseudo-reflection has a unique eigenvalue not equal to 1. We give now some examples of pseudo-reflections.
\begin{ex} 
	\mbox{}
	\vspace*{-\parsep}
	\vspace*{-\baselineskip}\\
	\begin{itemize}[leftmargin=*]
		
		\item Every reflection is a pseudo-reflection of order 2.
		\item If $j$ is a third root of unity, the matrices 
		$$s_1:=\begin{bmatrix}
		j&0\\
		-j^2&1
		\end{bmatrix} \text { and }
		s_2:=\begin{bmatrix}
		1&j^2\\
		0&j^2
		\end{bmatrix}$$
		are pseudo-reflections of order 3.
		\qed
	\end{itemize}
\end{ex}
\begin{defn}
	A complex reflection group $W$ of rank $n$ is a finite subgroup of $GL(V)$, which is generated by pseudo-reflections.
\end{defn}
\begin{rem} Every pseudo-reflection is a unitary reflection with respect to some form (see lemma 1.3 of \cite{lehrer}). Therefore, every complex reflection group $W\leq GL(V)$ can be considered as a subgroup of $U_n(\CC)$, where $U_n(\CC)$ is the unitary group of degree $n$. 
	\label{unitary}
	\end{rem}
\begin{ex}
	\mbox{}
	\vspace*{-\parsep}
	\vspace*{-\baselineskip}\\
	\begin{itemize}[leftmargin=*]
		\item The group $\grm_d$ of all the $d$-th roots of unity is a complex reflection group of rank 1.
		\item The symmetric group $S_n$ is a complex reflection group of rank $n-1$, if we associate every permutation $(i$ $j)$ with a permutation matrix, obtained by permuting the rows of the $n\times n$ identity matrix according to the permutation $(i$ $j)$. 
		\item More generally, every real reflection group (also known as finite Coxeter group) $$C=\langle s_1,\dots, s_n\;|\; (s_is_j)^{m_{i,j}}=1\rangle,$$ where $m_{i,i}=1$ and $2\leq m_{i,j}<\infty $ for $i\not=j$, is a complex reflection group of rank $n$.
		\qed
	\end{itemize}
\end{ex}
\begin{ex}
	Another example of a complex reflection  group is the three-parameter family, 
	known as $G(de,e,n)$, in the notation of Shephard and Todd (see \cite{shephard}). By definition, $G(de,e,n)$ is the group of monomial $n\times n$ matrices (meaning that they have only one
	nonzero entry in each row and column) where the nonzero entries are $de$-th
	roots of unity, and the product of these entries is a $d$th root of unity. The three-parameter family $G(de,e,n)$ can also be defined via the semi-direct product 
	$$G(de,e,n):=D(de,e,n)\rtimes S_n,$$
	where $D(de,e,n)$ denotes the group of all diagonal $n\times n$ matrices, whose diagonal entries are $de$-th roots of unity and their determinant is a $d$th root of unity.
	\qed
\end{ex}
	 Let $s_i$ be the permutation matrix $(i\;i+1)$, $i=1,\dots,n-1$, $t_e:=\begin{pmatrix}
	0&\grz_e^{-1}&0\\
	\grz_e&0&0\\
	0&0& \text{Id}
	\end{pmatrix}$ and $u_d:=$diag$\{\grz_d,1,\dots,1\}$, where $\grz_m$ denotes a $m$-th root of unity. We can easily verify that these matrices are pseudo-reflections.
	The following results are based on \cite{Michelcour}, Exercise 2.9 and on Chapter 2, \S 3 in \cite{lehrer}.
	\begin{prop}
		$G(de,e,n)$ is a complex reflection group of rank $n$. Its generators are as follows: 
		\begin{itemize}[leftmargin=*]
			\item $\underline{d=1}$: The group $G(e,e,n)$ is generated by the pseudo-reflections $t_e$, $s_1$,\dots, $s_{n-1}$.
			\item $\underline{e=1}$: The group $G(d,1,n)$ is generated by the pseudo-reflections $u_d$, $s_1$,\dots, $s_{n-1}$.
			\item \underline{$d\not= 1$, $e\not=1$}: The group $G(de,e,n)$ is generated by the pseudo-reflections $u_d$, $t_{de}$, $s_1$,\dots, $s_{n-1}$.
		\end{itemize}
	\end{prop}
	\begin{rem}
		The complex reflection groups $G(de,e,n)$ give rise to the real reflection groups. More precisely:
		\begin{itemize}[leftmargin=*]
			\item The group $G(1,1,n)$ is the symmetric group $S_n$, also known as the finite Coxeter group of type $A_{n-1}$.
			\item The group $G(2,1,n)$ is the finite Coxeter group of type $B_n$.
			\item The group $G(2,2,n)$ is the finite Coxeter group of type $D_n$.
			\item The group $G(e,e,2)$ is the dihedral group of order $2e$, also known as the finite Coxeter group of type $I_2(e)$.
		\end{itemize}
		\label{cox}
	\end{rem}
\begin{defn} 
	Let $W\leq GL(V)$ be a complex reflection group. We say that $W$ is irreducible if the only subspaces of $V$ that stay stable under the action of $W$ are $\{0\}$ and $V$.
\end{defn}
\begin{ex} The group $G(de,e,n)$ is irreducible apart from the following cases:
	\begin{itemize}[leftmargin=*]
	\item The group $G(1,1,n)=S_n$, considered as a complex reflection group of rank $n$ ($S_n$ is irreducible as a complex reflection group of rank $n-1$).
	\item The group $G(2,2,2)$ (the dihedral group of order 4).
	\qed
\end{itemize}
	\label{ex1}
	\end{ex}
The following proposition (proposition 1.27 in \cite{lehrer}) states that every complex reflection group can be written as a direct product of irreducible ones.
\begin{prop}
Let $W\leq GL_n(V)$ be a complex reflection group. Then, there is a decomposition $V=V_1\oplus\dots\oplus V_m$ such that the restriction $W_i$ of $W$ in $V_i$ acts irreducibly on $V_i$ and $W=W_1\times\dots\times W_m$.
\label{irred}
\end{prop}
Therefore, the study of reflection groups reduces to the irreducible case and, as a result, it is sufficient to classify the irreducible complex reflection groups. The following classification is due to G.C. Shephard and J.A. Todd (for more details one may refer to  \cite{cohen} or \cite{shephard}), also known as the ``Shephard-Todd classification''.
\begin{thm}
	Let $W$ be an irreducible complex reflection group. Then, up to conjugacy,
	$W$ belongs to precisely one of the following classes:
	\begin{itemize}[leftmargin=*]
		\item The infinite family $G(de,e,n)$, as described in example \ref{ex1}.
		\item The 34 exceptional groups $G_n$ ($n=4,\dots,37)$.
	\end{itemize}
	\label{classif}
\end{thm}
\begin{rem}
	\mbox{}  
	\vspace*{-\parsep}
	\vspace*{-\baselineskip}\\ 
	\begin{itemize}[leftmargin=*]
		\item[(i)] G.C. Shephard and J.A. Todd numbered the irreducible unitary reflection groups from 1 to 37. However, the first three entries in their list refer
		to the three-parameter family. More precisely,  $G_1:=G(1,1,n)$, $G_2:=G(de,e,n),$ $n\geq 2$ and $G_3:=G(d,1,1)$. Hence, when we refer to the exceptional cases, we mean the groups $G_4, \dots, G_{37}$.
		\item [(ii)] Among the irreducible complex reflection groups we encounter the irreducible finite Coxeter groups. By remark \ref{cox} we have already seen the finite Coxeter groups of type $A_{n-1}$, $B_n$, $D_n$ and $I_2(e)$ as irreducible complex reflection groups inside the three-parameter family $G(de,e,n)$. The rest of the cases, which are the irreducible finite Coxeter groups of types $H_3$, $F_4$, $H_4$, $E_6$,
		$E_7$ and $E_8$ are the exceptional groups $G_{23}$, $G_{28}$, $G_{30}$, $G_{35}$, $G_{36}$ and $G_{37}$, respectively.
	\end{itemize}
	\label{remirred}
\end{rem}

	\subsection{Complex Braid groups}
	\label{seccccc}
	\indent
	
	This section is a rewriting of \cite{broue2000}, \cite{bmr} and \cite{marinkrammer} \S2.1. Let $W\leq GL(V)$ be a complex reflection group. We let $\mathcal{R}$ denote the set of pseudo-reflections of $W$, $\mathcal{A}=$\{ker$(s-1)\;|\;s\in \mathcal{R}\}$ the hyperplane arrangement associated to $\mathcal{R}$, and $X=V\setminus \cup \mathcal{A}$ the corresponding hyperplane complement. 
We assume that $\mathcal{A}$ is essential, meaning that $\cap_{H\in \mathcal{A}}H=\{0\}$. 

$W$ acts on $\mathcal{A}$ by $w\cdot\text{ker}(s-1)=\text{ker}(wsw^{-1}-1).$
This action is well-defined, since for every $w\in W$ and $s\in \mathcal{R}$, we have that $wsw^{-1}\in \mathcal{R}$ (see lemma 1.9 in \cite{lehrer} or lemma 1.4 (1) in \cite{brouebook}). Let $x\in X$ and $s\in \mathcal{R}$. We notice that $s(x)\in X$ (if $s(x)\not \in X$, then there is a $s'\in\mathcal{R}$, such that $s(x)\in$ ker$(s'-1)$. Therefore, $s^{-1}(s'(s(x)))=x$, which means that $x\in s^{-1}\cdot$ker$(s'-1)\subset \cup \mathcal{A}$. This contradicts the fact that $x$ belongs to $X$). Therefore, we have an action of $W$ on $X$, defined by $s\cdot x=s(x)$, for every $s\in \mathcal{R}$ and $x\in X$. 

Let $X/W$ be the space of orbits of the above action. For every $x\in X$ we write  $\underline{x}$  for the image of $x$  under the canonical surjection $p: X\rightarrow X/W$. By Steinberg's theorem (see \cite{steinberg}) we have that the action of $W$ on $X$ is free. Therefore it defines a Galois covering $X\rightarrow X/W$, which gives rise to the following exact sequence:
$$1\rightarrow \grp_1(X,x)\rightarrow \grp_1(X/W, \underline{x})\rightarrow W\rightarrow 1$$
(For more details about Galois coverings and their exact sequences, one may refer to \cite{berenstein}, Chapter 5, \S8, 9, 10). We define $P:=\grp_1(X,x)$ the \emph{pure complex braid group} associated to $W$ and $B:=\grp_1(X/W, \underline{x})$ the \emph{complex braid group} associated to $W$.

The canonical projection map $p: X\rightarrow X/W$ induces a natural projection map $ p_r: B \twoheadrightarrow W$, defined as follows:
Let $b\in B$ and let $\grb : \left[0,1\right]\rightarrow X$ be a path in $X$ such that $\grb(0)=x$, which lifts $b$, meaning that $b=\left[p\circ \grb\right]$.
$$
\begin{tikzcd}
& X \arrow{d}{p} \\
\left[0,1\right] \arrow{ur}{\grb} \arrow{r}{p\circ \grb} & X/W
\end{tikzcd}
$$
Then $p_r(b)$ is defined by the equality $p_r(b)(x)=\grb(1)$.

The goal of this section is to define an element $\grs\in B$ from an element $s\in \mathcal{R}$. By algebraic topology, we know that if $(X, x_0)$ and $(Y,y_0)$ are based topological spaces, then $\grp_1(X\times Y, (x_0,y_0))\simeq \grp_1(X, x_0)\times \grp_1(Y,y_0)$ (one may refer, for example, to Chapter 2, \S 7 in  \cite{massey}). Therefore, by proposition \ref{irred} we may assume that $W$ is an irreducible complex reflection group. Let $s\in \mathcal{R}$ a pseudo-reflection of order $m$ and $H=\text{ker}(s-1)$ the corresponding hyperplane. We fix an element $x\in X$ and we define a path from $x$ to $s\cdot x$ in the following way:
\begin{itemize}[leftmargin=*]
	\item We choose a point $x_0$ ``close to $H$ and far from the
	other reflecting hyperplanes'' as follows: Let $x_0^H \in H$ and $\varepsilon>0$, such that $B(x_0^H,\varepsilon)\cap \cup \mathcal{A}\subset H$. In other words, if $H\not=H'\in \mathcal{A}$, then $B(x_0^H,\varepsilon)\cap H'=\emptyset$. Let $x_0^{H^{\perp}} \in H^{\perp}$, such that $||x_0^{H^{\perp}}||<\varepsilon$. We define $x_0=x_0^H+x_0^{H^{\perp}}\in B(x_0^H,\varepsilon)$ (see figure \ref{braid} above).
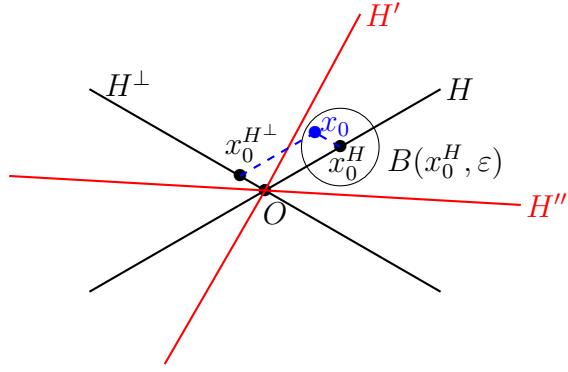
\begin{figure}[h]
	\centering
	\begin{tikzpicture}
	[x={(0.866cm,-0.5cm)}, y={(0.866cm,0.5cm)}, z={(0cm,1cm)}, scale=0.77,
	>=stealth, %
	inner sep=0pt, outer sep=2pt,%
	axis/.style={thick,->},
	wave/.style={thick,color=#1,smooth},
	polaroid/.style={fill=black!60!white, opacity=0.3},
	]
	
	\colorlet{darkgreen}{green!50!black}
	\colorlet{lightgreen}{green!80!black}
	\colorlet{darkred}{red!50!black}
	\colorlet{lightred}{red!80!black}

	\coordinate (O) at (0, 0,0);
	\draw (O) node{\textbullet};
	\draw[thick] (O) -- +(0,  3.5,0) node[right]{$H$};
	\draw[thick] (O) -- +(0,  -3.5,0);
	\draw [thick] (O)-- +(3.5, 0,   0);
	\draw [thick] (O)-- +(-3.5, 0,   0);
	\draw (-3.5,0.2,0) node [right] {$H^{\perp}$};
		\draw[thick, red] (O) -- +(2.8,  2.3,0) node[right]{$H''$};
		\draw[thick, red] (O) -- +(-2.1,  4,0)node[right]{$H'$};
		\draw[thick, red] (O) -- +(2,  -4,0);
\draw[thick, red] (O) -- +(-2.8,  -2.3,0);
\draw (0,1.5,0) node{\textbullet};
\draw (0.35,1.3,0) node{$x_0^H$};
\draw (0,1.5,0) circle (0.67cm);
\draw(-0.5,0,0) node{\textbullet};
\draw(-0.95,0.7,0) node{$x_0^{H^{\perp}}$};
\draw [thick, dashed, blue](-0.5,0,0)-- +(0,  1.5,0);
\draw [thick, dashed, blue](0,1.5,0)-- +(-0.5,  0,0);
\draw [blue](-0.4,1.8,0) node{$x_0$};
\draw [blue] (-0.5,1.5,0) node{\textbullet};
\draw (1.3,2.3,0) node{$B(x_0^H,\varepsilon)$};
\draw (0.5,-0.3,0) node{$O$};
\end{tikzpicture}
\caption{The choice of $x_0$} \label{ball}
\end{figure}
\item Let $\grg$ be an arbitrary path in $X$ from $x$ to $x_0$. Then, $\grg^{-1}$ is the path in $X$ with initial point $x_0$ and terminal point $x$, such that $\grg^{-1}(t)=\grg(1-t), \text{ for every } t\in\left[0,1\right].$
Thus, the path $s\cdot \grg^{-1}$ defined as $t\mapsto s(\grg^{-1}(t))$ is a path inside $X$, which goes from $s\cdot x_0$ to $s\cdot x$. We consider now a rotation $\grg_0$ of angle  $\gru=2\grp/m$ inside $H^{\perp}$, with initial point
$x_0$ and terminal point $s\cdot x_0$ (see Figure \ref{braid}). Since $x_0=x_0^H+x_0^{H^{\perp}}$, we can define $\grg_0: \left[0,1\right]\rightarrow H^{\perp}$ by
$\grg_0(t)=x_0^H+e^{\gru i t}x_0^{H^{\perp}}.$
\item Let $\tilde \grg$ be the path from $x$ to $s\cdot x$ defined by 
$\tilde \grg:=s\cdot \grg^{-1}*\grg_0*\grg$. 
By the choice of $x_0$ this path lies in $X$ and its homotopy class does not depend on its choice. Let $\grs_{\grg}$ denote the element that $\tilde \grg$ induces in the braid group $B$. We say that $\grs_\grg$ is \emph{a braided reflection} associated to $s$ around the image of $H$ in $X/W$.
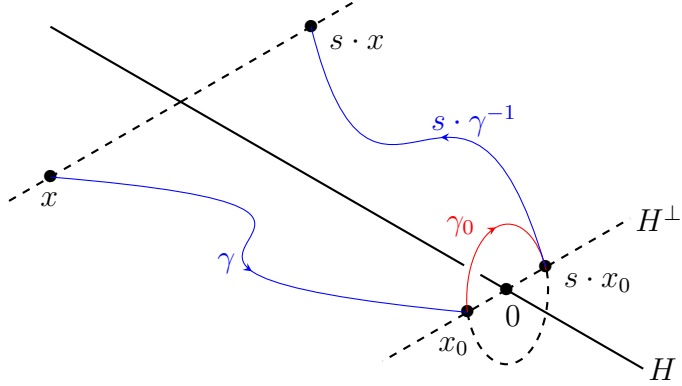
\begin{figure}[h]
	\centering
\begin{tikzpicture}
[x={(0.866cm,-0.5cm)}, y={(0.866cm,0.5cm)}, z={(0cm,1cm)}, scale=1.0,
>=stealth, %
inner sep=0pt, outer sep=2pt,%
axis/.style={thick,->},
wave/.style={thick,color=#1,smooth},
polaroid/.style={fill=black!60!white, opacity=0.3},
]

\colorlet{darkgreen}{green!50!black}
\colorlet{lightgreen}{green!80!black}
\colorlet{darkred}{red!50!black}
\colorlet{lightred}{red!80!black}

\coordinate (O) at (0, 0, 0);
\draw[thick,dashed] (O) -- +(0,  2.7, 0) ;
\draw[thick,dashed] (O) -- +(0,  -2.7, 0);
\draw(0.3,1.9,0) node [right]{$s\cdot x$};
\draw(0.4,-2.2,0) node [left] {$x$};
 \draw[thick] (-2,0,0) -- (O);
\draw[thick] (O) -- +(4.35, 0,   0) ;

\draw [thick] (4.6,0,0)-- +(2.5, 0,   0) node [right] {$H$};
\draw (0,2,0)  node{\textbullet} ;
\draw (0,-2,0)  node{\textbullet} ;
\draw[thick,dashed] (5,0,0) -- +(0,  1.9, 0) node[right]{$H^{\perp}$};
\draw (5, 0.6, 0) node{\textbullet}; 
\draw[thick,dashed] (5,0,0) -- +(0,  -1.9, 0) ;
\draw (5, -0.6, 0) node{\textbullet} ;
\draw (5.6,-1.1,0.1) node [left] {$x_0$};
\draw (5.8, 1.2,-0.1) node [left] {$s\cdot x_0$};
\draw (5,0,0) node{\textbullet};
\draw (5.5,-0.2,0) node[left]{0};

\draw[fermion, red] (5,-0.6,0) arc [start angle=190,end angle=23, x radius=0.55cm, y radius=1cm];
\draw[thick, dashed] (5,-0.6,0) arc [start angle=-162, end angle=21, x radius=0.55cm, y radius=1cm];

\path[fermion,blue,draw] (0,-2,0) .. controls (4,1,0) and (0.65,-3.7,0) .. (5,-0.6,0);
\path [fermionbar,blue, draw](0,2,0).. controls (4,-1,0) and (0.65,3.7,0) .. (5,0.6,0);

\draw [blue](2.5, -1.8,0) node{$\grg$};
\draw [blue](2.5, 2,0) node{$s\cdot \grg^{-1}$};
\draw [red](3.8, 0.5,0) node{$\grg_0$};
\end{tikzpicture}
\caption{Braided reflection associated to $s$} \label{braid}
\end{figure}
\end{itemize} 
The next lemma is proposition 2.12 (1) in \cite{broue2000}.
\begin{lem}
	Let $s$ be pseudo-reflection and $\grs_{\grg}$ a braided reflection associated to $s$, as defined above. Then, 
	$p_r(\grs_{\grg})=s$.
\end{lem}

The following proposition (see Lemma 2.12 (2) in \cite{broue2000}) states that if we have two braided reflections associated to the same pseudo-reflection $s$, they are conjugate in $P$.

\begin{prop}
	Whenever $\grg'$ is a path in $X$ from $x$ to $x_0$, if $\grt$ denotes the loop in $X$ defined by $\grt:=\grg'^{-1}*\grg$, then
	$$\tilde\grg'=s\cdot \grt*\tilde\grg*\grt^{-1}.$$
	In particular, $\grs_{\grg}$ and $\grs_{\grg'}$ are conjugate in $P$.
\end{prop}
\begin{cor}
	Let $s_1, s_2$ two pseudo-reflections, which are conjugate in $W$ and let $\grs_1$ and $\grs_2$ denote two braided reflections associated to $s_1$ and $s_2$, respectively.  Then, $\grs_1$ and $\grs_2$ are conjugate in $B$.
	\label{conbraid}
\end{cor}

Let $s\in \mathcal{R}$ and $H=\text{ker}(s-1)$ the corresponding hyperplane. Let $W_H$ be the subgroup of $W$ formed by $id_V$ and all the reflections fixing $H$. We recall that $V=H\oplus H^{\perp}$ and we set $\grf: GL(V)\rightarrow GL(H^{\perp})$ defined as $f\mapsto f|_{H^{\perp}}$. Let $w_1, w_2 \in W_H$ such that $\grf(w_1)=\grf(w_2)$. Hence, by the definition of $\grf$ we have that $w_1(h^{\perp})=w_2(h^{\perp})$, for every $h^{\perp}\in H^{\perp}$. By the definition of $W_H$ we have also that $w_1(h)=w_2(h)$, for every $h\in H$. As a result, we have that $w_1(v)=w_2(v)$, for every $v\in V$, since every
$v\in V$ is written as $h+h^{\perp}$, where $h\in H$ and $h^{\perp}\in H^{\perp}$. Therefore, the restriction of $\grf$ in $W_H$ is injective and,
hence, we have  $W_H \leq GL(H^{\perp})\simeq GL_1(\CC)\simeq \CC^{\times}$.
As a result, the group $W_H$ is cyclic, as a subgroup of $\CC^{\times}$. We denote by $e_H$ the order of the cyclic group $W_H$ and we give the following definition: 
\begin{defn}
 Let $s_H$ be the (unique) pseudo-reflection of eigenvalue $e^{2\grp i/e_H}$, which generates $W_H$. We say that $s_H$ is a \emph{distinguished pseudo-reflection}. 
	\end{defn}
We recall that $p_r$ denotes the natural projection map $B \twoheadrightarrow W$ induced by the canonical surjection $X\rightarrow X/W$.	
	\begin{defn}
		Let $s_H$ be a distinguished pseudo-reflection. A \emph{distinguished braided reflection} $\grs_H$ associated to $s_H$ is a braided reflection  around the image of $H$ in $X/W$ such that $p_r(\grs_H)=s_H$.
	\end{defn}
	The next result (proposition 2.5 (a) in \cite{broue2000}) explains the importance of the distinguished braided reflections.
	 \begin{prop}
	 	The complex braid group $B$ is generated by all the distinguished braided
	 	reflections around the images of the hyperplanes $H\in \mathcal{A}$ in $X$.
	 \end{prop}
The next result (theorem 0.1 in \cite{bessiszariski}) shows that a complex braid group has an \emph{Artin-like presentation}; that is a presentation of the form 
$$\langle \mathbf{s}\in \mathbf{S}\;|\; \mathbf{v_i}=\mathbf{w_i} \rangle_{i\in I},$$
where $\mathbf{S}$ is a finite set of distinguished braided reflections and $I$ is a finite set of relations such that, for each $i\in I$,  $\mathbf{v_i}$ and $\mathbf{w_i}$ are positive words in elements of $\mathbf{S}$.
\begin{thm}Let $W$ be a complex reflection group with associated complex braid group $B$. There exists a finite subset $\mathbf{S}=\{\mathbf{s_1},\dots, \mathbf{s_n}\}$ of $B$, such that: 
	\begin{itemize}
		\item[(i)] The elements $\mathbf{s_1}, \dots, \mathbf{s_n}$ are distinguished braided reflections and therefore, their images $s_1,\dots, s_n$ in $W$ are distinguished pseudo-reflections.
		\item[(ii)]The set  $\mathbf{S}$ generates $B$ and therefore the set $S:=\{s_1,\dots, s_n\}$ generates $W$.
		\item [(iii)] There exist a set $I$ of relations of the form $\mathbf{w_1}=\mathbf{w_2}$, where $\mathbf{w_1}$ and $\mathbf{w_2}$ are positive words of equal length in the elements of $\mathbf{S}$, such that $\langle \mathbf{S}\;|\;I\rangle$ is a presentation of $B$.
		\item[(iv)] Viewing now $I$ as a set of relations in $S$, the group
		$W$ is presented by $$\langle S\;|I\;; \forall s\in S, s^{e_s}=1\rangle,$$ where $e_s$ denotes the order of $s$ in $W$.
	\end{itemize}
	\label{Presentt}
\end{thm}
\begin{rem}
	\mbox{}  
	\vspace*{-\parsep}
	\vspace*{-\baselineskip}\\ 
	\begin{itemize}
\item[(i)]By Brieskorn's theorem (see \cite{Brieskorn}) we have the following result: Let $C$ be a finite Coxeter group with presentation 
 $$C=\langle s_1,\dots, s_n\;|\; s_i^2=1, \;(s_is_j)^{m_{i,j}}=1, \text{ for all } i\not=j\rangle.$$The corresponding braid group $A$ (known as an Artin group of finite Coxeter type), has a presentation of the form $$A=\langle s_1,\dots, s_n\;|\;  \;\underbrace{s_is_j\dots}_{m_{i,j}-\text{ times}}=\underbrace{s_js_i\dots}_{m_{i,j}-\text{ times}}, \text{ for all } i\not=j\rangle.$$
 This means that we can obtain the presentation of the finite Coxeter group from the presentation of the Artin group, if we  ``add'' to the latter the additional relations $s_i^2=1$, for every $i=1,\dots,n$. 
 Therefore, theorem \ref{Presentt} generalizes Brieskorn's result for every complex reflection group.
 \item[(ii)]Theorem \ref{Presentt} not only shows that a complex braid group has an Artin-like presentation but also implies that any complex reflection group has a \emph{Coxeter-like presentation}; that is a presentation of the form 
 $$\langle s\in S\;|\; \{v_i=w_i\}_{i\in I} , \{s^{e_s}=1\}_{s\in S}\rangle,$$
 where $S$ is a finite set of distinguished  reflections and $I$ is a finite set of relations such that, for each $i\in I$,  $v_i$ and $w_i$ are positive words with the same length in elements of $S$. The tables in Appendix 1 in \cite{broue2000}  provide  a  complete  list  of  the irreducible complex reflection groups in Shephard-Todd classification, together with a Coxeter-like presentation symbolized by a diagram.
\end{itemize}
\label{coxeterlike}
\end{rem}
We conclude this section by giving a description of the center of the complex braid group, a result that plays an important role in the sequel. For this description, we follow the arguments in the introduction of \cite{dignecenter}. Let $W\leq GL(V)$ be a complex reflection group. By proposition \ref{irred} we can write $W$ as a direct product $W_1\times\dots\times W_r$ of irreducible reflection groups. As we explained in section \ref{seccccc}, the associated complex braid group is the direct product $B_1\times\dots\times B_r$, where $B_i$ is the complex braid group associated to the irreducible complex reflection group $W_i$, $i=1,\dots, r$.
Therefore, the center of the complex braid group $Z(B)$ is the direct product  $Z(B_1)\times \dots \times Z(B_r)$ and, thus, we may assume that $W$ is  irreducible.

Since $W$ is an irreducible complex reflection group,  its center $Z(W)$ is a (finite) subgroup of $\CC^{\times}$, thus a cyclic group. Let $Z(W)=\langle \grz \rangle$, where  $\grz:=e^{2\grp i/|Z(W)|}$ and let $x$ be a basepoint of $X$. In \cite{bmr}, \S 2.24, M. Brou\'e, G. Malle and R. Rouquier defined a path $\grb$ inside $X$ by $t\mapsto \grz^tx$, where $\grz^t:=e^{2\grp it/|Z(W)|}$.  We denote by $\tilde\grb$ the homotopy class of $\grb$ in $B$. The next theorem is due to  M. Brou\'e, G. Malle and R. Rouquier (see theorem 2.2.4 in \cite{bmr}) and D. Bessis (see theorem 12.8 in \cite{bessiscenter}).

\begin{thm}
	The center of the braid group of an irreducible complex reflection
	group is cyclic, generated by  $\tilde\grb$.
	\label{centerbraid}
\end{thm}

\section{The freeness conjecture}

\subsection {Generic Hecke algebras}
\indent

In this section we follow mostly \cite{bmr} and \cite{marinG26}.
Let $W\leq GL(V)$ be a complex reflection group. Recalling the definitions and notations of the previous section, let $B$ be the complex braid group associated to $W$ and $S$ the set of the distinguished pseudo-reflections. We denote also by $e_s$ the order of $s$ in $W$.

For each $s\in S$ we choose a set of $e_s$ indeterminates $u_{s,1},\dots, u_{s,e_s}$ such that $u_{s,i}=u_{t,i}$ if $s$ and $t$ are conjugate in $W$. We denote by $R$ the Laurent polynomial ring $\ZZ[u_{s,i},u_{s,i}^{-1}]$ and we give the following definition:
\begin{defn}
	The generic Hecke algebra $H$ associated to $W$ with parameters $u_{s,1},\dots, u_{s,e_s}$ is the quotient of the group algebra $RB$ of $B$ by the ideal generated by the elements of the form 
	\begin{equation}
	(\grs-u_{s,1})(\grs-u_{s,2})\dots (\grs-u_{s,e_s}),
	\label{Hecker}
	\end{equation}
	where $s$ runs over the conjugacy classes of $S$ and $\grs$ over the set of distinguished braided reflections associated to the pseudo-reflection $s$. 
\end{defn}

\begin{rem}
	\mbox{}
	\vspace*{-\parsep}
	\vspace*{-\baselineskip}\\
	\begin{itemize}
 \item [(i)]It is enough to choose one relation of the form described in (\ref{Hecker}) per conjugacy class, since  the corresponding braided reflections are conjugate in $B$ (Corollary \ref{conbraid}).
 \item[(ii)]Let $C$ be a finite Coxeter group with Coxeter system $S$. In this case, the ring of Laurent polynomials is $R=\ZZ[u_{s,1}^{\pm},u_{s,2}^{\pm}]_{s\in S}$ and the generic Hecke algebra associated to $C$ has a presentation as follows:
 
$$\langle \grs_1,\dots, \grs_n\;|\; \underbrace{\grs_i\grs_j\dots}_{m_{i,j}-\text{ times}}=\underbrace{\grs_j\grs_i\dots}_{m_{i,j}-\text{ times}},\; (\grs_i-u_1)(\grs_i-u_2)=0\rangle,$$
where $i\not =j$ and $2\leq m_{i,j}<\infty $.
 In this case (the real case), the generic Hecke algebra is known as the \emph{Iwahori-Hecke algebra} (for more details about Iwahori-Hecke algebras one may refer, for example, to \cite{geck} \S 4.4).
\item [(iii)] Let $\grf: R\rightarrow \CC$ be the specialization morphism defined as
$u_{s,k}\mapsto e^{-2\grp \gri k/e_s}$, where $1\leq k\leq e_c$ and $\gri$ denotes an imaginary unit (a solution of the equation $x^2=-1$). Therefore, $H\otimes_{\grf}\CC=\CC B/(\grs^{e_s}-1)=\CC \big(B/(\grs^{e_s}-1)\big)$. By theorem \ref{Presentt} (iv) we have that $B/(\grs^{e_s}-1)=W$. Hence, $H\otimes_{\grf}\CC=\CC W$, meaning that $H$ is a deformation of the group algebra $RW$.
 \end{itemize} 
 \label{remi}
\end{rem}
In order to make the definition of the generic Hecke algebra clearer to the reader, we give an example of the generic Hecke algebra associated to the exceptional groups $G_4$ and $G_5$.
\begin{ex}
	\mbox{}
	\vspace*{-\parsep}
	\vspace*{-\baselineskip}\\
	\begin{itemize}[leftmargin=*]
	\item  Let $W:=G_4=\langle s_1, s_2\;|\;s_1^3=s_2^3=1, s_1s_2s_1=s_2s_1s_2\rangle$.
	Since $s_1=(s_2s_1)s_2(s_2s_1)^{-1}$, the ring of Laurent polynomials is $\ZZ[u_i^{\pm}]$, $i=1,2,3$ and the generic Hecke algebra has a presentation as follows:
	$$H_{G_4}=\langle \grs_1,\grs_2:\grs_1\grs_2\grs_1=\grs_2\grs_1\grs_2,\; (\grs_i-u_1)(\grs_i-u_2)(\grs_i-u_3)=0\rangle.$$
	
	\item Let $W:=G_{10}=\langle s_1, s_2\;|\;s_1^3=s_2^4=1, s_1s_2s_1s_2=s_2s_1s_2s_1\rangle$.
	The ring of Laurent polynomials is $\ZZ[u_i^{\pm},v_j^{\pm}]$, $i=1,2,3$, $j=1,2,3,4$ and the generic Hecke algebra has a presentation as follows:
	$$H_{G_{10}}=\langle \grs_1,\grs_2\;|\;\grs_1\grs_2\grs_1\grs_2=\grs_2\grs_1\grs_2\grs_1,\prod\limits_{i=1}^{3}(\grs_1-u_i)=\prod\limits_{j=1}^{4}(\grs_2-v_j)=0 \rangle.$$
	\qed
	\end{itemize}
\end{ex}

\subsection{The BMR freeness Conjecture: Recent work and open cases}
\indent

Let $W\leq GL(V)$ be a complex reflection group and $H$ its associated generic Hecke algebra over the ring of Laurent polynomials $R$, as defined in the previous section. We have the following conjecture due to M. Brou\'e, G. Malle and R. Rouquier (see \cite{bmr}).

\begin{co} (The BMR freeness conjecture) The generic Hecke algebra $H$ is a free module over $R$ of rank $|W|$.
	\end{co}
The next proposition (theorem 4.24 in \cite{bmr} or proposition 2.4, (1) in \cite{marinG26}) states that in order to prove the validity of the BMR conjecture, it is enough to find a spanning set of $H$ over $R$ of $|W|$ elements.
\begin{prop}
	If $H$ is generated as $R$-module by $|W|$ elements, then it is a free module over $R$ of rank $|W|$.
	\label{BMR PROP}
\end{prop}

\begin{thm}
The BMR freeness conjecture holds for the real complex reflections groups (i.e. for the Iwahori-Hecke algebras). 
\label{IwahoriBMR}
\end{thm}
We give a sketch  of the proof of this theorem in order to underline the similarities of this proof with the one of the complex cases of rank 2 we prove in Chapter 4. For more details one may refer for example to \cite{geck}, lemma 4.4.3. 

\begin{proof}[Sketch of the proof] Let $C$ be a finite Coxeter group with Coxeter system $S=\{s_1,\dots, s_n\}$ and $\mathcal{H}$ the associated Iwahori-Hecke algebra over the ring $R=\ZZ[u_{s,1}^{\pm},u_{s,2}^{\pm}]_{s\in S}$, as defined in remark \ref{remi} (ii). By Matsumoto's lemma (see \cite{Matsumoto}) there is a natural map (not a group homomorphism) from  $C$ to the corresponding braid group, taking any element of $C$ represented by some reduced word in the generators to the same word in the generators of the braid group. Therefore, for every element $c\in C$ represented by a reduced word  $s_{i_1}\dots s_{i_r}$  there is a well-defined element in $\mathcal{H}$ we denote by $T_c$, such that $T_c=T_{s_{i_1}}\dots T{s_{i_r}}$. In particular, $T_1=1_{\mathcal{H}}$. 
	
	We will prove that a spanning set of $\mathcal{H}$ over $R$ is $\{T_c,\;c\in C\}$ and, therefore, by proposition \ref{BMR PROP} we prove also the validity of the BMR freeness conjecture in the real case. Let $U$ be the $R$-submodule of $\mathcal{}H$ generated by  $\{T_c,\;c\in C\}$. We must prove that $\mathcal{H}=U$. Since $1_{\mathcal{H}}\in U$, it will be sufficient to show that
	$U$ is a left ideal of $\mathcal{H}$. For this purpose, one may
	check that $U$ is invariant under left multiplication by all $T_s$, $s\in S$. This is a straightforward consequence of the fact that 
	$$T_sT_w=\begin{cases}T_{sw}, &\text { if } \ell(sw)=\ell(w)+1\\
	(u_{s,1}+u_{s,2})T_w-u_{s,1}u_{s,2}T_{sw},&\text { if } \ell(sw)=\ell(w)-1
	\end{cases},$$
	where $\ell(sw)$ and $\ell(w)$ denote the length of the words $w$ and $sw$, respectively.
\end{proof}
We go back now to the case of an arbitrary complex reflection group $W$. In proposition \ref{irred} we saw that $W$ is a direct product of some irreducible complex reflection groups. Therefore, we restrict ourselves to proving the validity of the conjecture for the irreducible complex reflection groups. Due to the classification of Shephard and Todd (theorem \ref{classif}), we have to prove the conjecture for the three-parameter family $G(de,e,n)$ and for the exceptional groups $G_4, \dots, G_{37}$.
Thanks to S. Ariki and S. Ariki and K. Koike (see \cite{ariki} and \cite{arikii}) we have the following theorem:
\begin{thm}
	The BMR freeness conjecture holds for the infinite family $G(de,e,n)$.
\end{thm}
 As a consequence of the above result, one has to concentrate on the exceptional groups, which are divided into two families: The first family includes the groups $G_4, \dots, G_{22}$, which are of rank 2 and the second one includes the rest of them, which are of rank at least 3 and at most 8. We recall that among these groups we encounter 6 finite Coxeter groups (remark \ref{remirred} (ii)), for which we know the validity of the conjecture: the groups $G_{23}$, $G_{28}$,
$G_{30}$, $G_{35}$, $G_{36}$ and $G_{37}$. Thus, it remains to prove the conjecture for 28 cases: the exceptional groups of rank 2 and the exceptional groups $G_{24}$, $G_{25}$, $G_{26}$ $G_{27}$,
$G_{29}$, $G_{31}$, $G_{32}$, $G_{33}$ and $G_{34}$.

Among these 28 cases, we encounter 6 groups whose associated complex braid group is an Artin group: The groups $G_4$, $G_8$ and $G_{16}$ related to the Artin  group of Coxeter type $A_2$, and the groups $G_{25}$, $G_{26}$ and $G_{32}$ related to the Artin group of Coxeter type $A_3$, $B_3$ and $A_4$, respectively.  
The next theorem is thanks to the results of I. Marin (see \cite{marincubic} and \cite{marinG26}). 
\begin{thm}
	The BMR freeness conjecture holds for the exceptional groups $G_4$, $G_{25}$, $G_{26}$ and $G_{32}$.
	\label{braidcase}
\end{thm}
The case of $G_4$ has been also proven independently by M. Brou\'e and G. Malle and B. Berceanu and L. Funar (see \cite{brouem} and  \cite{funar1995}).
Exploring the rest of the cases, we notice that we encounter 9 groups generated by reflections (i.e. pseudo-reflections of order 2): These groups are the exceptional groups $G_{12}$, $G_{13}$, $G_{22}$ of rank 2 and the exceptional groups $G_{24}$, $G_{27}$,
$G_{29}$, $G_{31}$, $G_{33}$ and $G_{34}$ of rank at least 3 and at most 6. I. Marin and G. Pfeiffer proved the following result by using computer algorithms (see \cite{marinpfeiffer}).
\begin{thm}
	The BMR freeness conjecture holds for the exceptional groups $G_{12}$,  $G_{22}$, $G_{24}$, $G_{27}$,
	$G_{29}$, $G_{31}$, $G_{33}$ and $G_{34}$.
	\label{2case}
\end{thm}
To sum up, the BMR freeness conjecture is still open for the exceptional groups of rank 2, apart from the cases of $G_4$, $G_{12}$ and $G_{22}$ for which we know the validity of the conjecture (theorems \ref{braidcase} and \ref{2case}). The following chapters are devoted to the proof of 11 of the 16 remaining cases, including also another proof of the case of $G_{12}$.

\underline{}

\chapter{The freeness conjecture for the finite quotients of $B_3$}

 In this chapter we prove that  the quotients of the group algebra of the braid group on 3 strands by a generic quartic and quintic relation, respectively have finite rank. This is the special case of the BMR freeness conjecture for the generic Hecke algebra of the groups $G_8$ and $G_{16}$. This result completes the proof of this conjecture in the case of the exceptional groups whose associated complex braid group is an Artin group (see theorem \ref{braidcase}). Exploring the consequences of this case, we prove that we can determine completely the irreducible representations of the braid group on 3 strands of dimension at most 5, thus recovering a classification of Tuba and Wenzl in a more general framework. This chapter is based on the author's article (see \cite{chavli}).

\section{The finite quotients of $B_n$}
\label{s}
\indent

Let $B_n$ be the braid group on $n$ strands,  defined via the following presentation:
$$\langle s_1, \dots, s_{n-1}\;|\;s_is_{i+1}s_i=s_{i+1}s_is_{i+1}, \; s_is_j=s_js_i\rangle,$$
where in the first group of relations $1\leq i\leq n-2$, and in the second one $|i-j|\geq 2$.
There is a Coxeter's classification of the finite quotients of $B_n$ by the additional relation $s_i^k=1$ (for more details one may refer to \textsection10 in \cite{Coxeter}); these quotients are finite if and only if  $\frac{1}{k}+\frac{1}{n}>\frac{1}{2}$. If we exclude the obvious cases $n= 2$ and $k=2$, which lead to the
cyclic groups and to the symmetric groups respectively, there is only a finite number of such groups, which
are irreducible complex reflection groups: these are the exceptional groups $G_4, G_8$ and $G_{16}$, for $n=3$ and $k=3,4,5$ and the exceptional groups  $G_{25}$, $G_{32}$ for $n=4,5$ and $k=3$.

 The BMR freeness conjecture is known  for the cases of $G_4$, $G_{25}$ and $G_{32}$, as we explained in the previous chapter (see theorem \ref{braidcase}).
Therefore, it remains to prove the validity of the conjecture for the groups $G_8$ and $G_{16}$, which are finite quotients of $B_3$ together with the groups 
 $S_3$ (if we consider the case of the symmetric group, as well) and $G_4$. Moreover, these exceptional groups belong to the class of complex reflection groups of rank two.
 
Let $B_3$ be the braid group on 3 strands, given by generators the braids $s_1$ and $s_2$ and the single relation $s_1s_2s_1=s_2s_1s_2$ that we call braid relation.
For every $k=2,\dots,5$ we denote by  $R_k$ the Laurent polynomial ring $\ZZ[a_{k-1},...,a_1, a_0, a_0^{-1}]$. Let $H_k$ denote the quotient of the group algebra $R_kB_3$ by the relations \begin{equation}s_i^k=a_{k-1}s_i^{k-1}+...+a_1s_i+a_0,\label{one} \end{equation}
for $i=1,2$.
 \begin{defn}For $k=2, 3,4$ and 5 we call the algebra $H_k$ the quadratic, cubic, quartic and quintic Hecke algebra, respectively.
\end{defn}
We identify $s_i$ with their images in $H_k$. We multiply $(\ref{one})$ by $s_i^{-k}$ and since $a_0$ is invertible in $R_k$ we have:
\begin{equation}s_i^{-k}=-a_0^{-1}a_1s_i^{-k+1}-a_0^{-1}a_2s_i^{-k+2}-...-a_0^{-1}a_{k-1}s_i^{-1}+a_0^{-1}, \label{two}\end{equation}
for $i=1,2$.
If we multiply ($\ref{two}$) with a suitable power of $s_i$ we can expand $s_i^{-n}$ as a linear combination of 
 $s_i^{-n+1},...,s_i^{-n+(k-1)}, s_i^{-n+k}$, for every $n\in \NN$.
Moreover, comparing (\ref{one}) and (\ref{two}), we can define an automorphism $\grF$  of $H_k$ as $\ZZ$-algebra, where $$\begin{array}{lcl}
s_i &\mapsto &s_i^{-1},\text{ for } i=1,2\\
a _j&\mapsto&-a_0^{-1}a_{k-j},\text{ for } j=1,...,k-1\\
a_0&\mapsto&a_0^{-1}
\end{array}$$

We will prove now an easy lemma that plays an important role in the sequel. This lemma is in fact a generalization of  lemma 2.1 of \cite{marincubic}. 
\begin{lem}For every $m \in \ZZ$ we have  
$s_2s_1^ms_2^{-1}=s_1^{-1}s_2^ms_1$ and $s_2^{-1}s_1^ms_2=s_1s_2^ms_1^{-1}$.
\label{lem1}
\end{lem}
\begin{proof}
By using the braid relation we have that $(s_1s_2)s_1(s_1s_2)^{-1}=s_2$. Therefore, for every $m \in \ZZ$ we have $(s_1s_2)s_1^m(s_1s_2)^{-1}=s_2^m$, that gives us the first equality. Similarly, we prove the second one.
\end{proof}
 If we assume $m$ of lemma \ref{lem1} to be positive we have  $s_1s_2s_1^n=s_2^ns_1s_2$ and $s_1^ns_2s_1=s_2s_1s_2^n$, where $n\in \NN$. Taking inverses, we also get $ s_1^{-n}s_2^{-1}s_1^{-1}=s_2^{-1}s_1^{-1}s_2^{-n}$ and $ s_1^{-1}s_2^{-1}s_1^{-n}=s_1^{-n}s_2^{-1}s_1^{-1}$. We call all the above relations \emph{the generalized braid relations.}

We will denote by $u_i$ the $R_k$-subalgebra of $H_k$ generated by $s_i$ (or equivalently by $s_i^{-1}$) and by $u_i^{\times}$ the group of units of $u_i$ and we let $\grv=s_2s_1^2s_2$. Since the center of $B_3$ is the subgroup generated by the element $z=s_1^2\grv$ (see, for example, theorem 1.24 of \cite{turaev}), for all $x\in u_1$ and $m\in \ZZ$ we have that $x\grv^m=\grv^mx$. We will see later that $\grv$ plays an important role in the description of $H_k$.

Let $W_k$ denote the quotient group $B_3/\langle s_i^k \rangle$, 
 $k=2, 3, 4$ and 5 and let $r_k<\infty$ denote its order.
Our goal now is to prove that $H_k$ is a free $R_k$-module of rank $r_k$, a statement that holds for $H_2$ since $W_2=S_3$ is a Coxeter group (see 
theorem \ref{IwahoriBMR}). We also know that this holds for the cubic Hecke algebra $H_3$ (see theorem 3.2 $(3)$ in \cite{marincubic}).
For the remaining cases, we will use the following proposition.
\begin{prop}Let $k\in\{4,5\}$. If $H_k$ is generated as a module over $R_k$ by $r_k$ elements, then $H_k$ is a free $R_k$-module of rank $r_k$ and, therefore, the BMR freeness conjecture holds for the exceptional groups $G_8$ and $G_{16}$.
\label{rp}
\end{prop}
\begin{proof} Let $\tilde H_k$ denote the generic Hecke algebra of 
	 $G_8$ and $G_{16}$, for $k=4$ and 5, respectively. 
We know that $\tilde{H_k}$ is a free $\tilde{R_k}$-module of rank $r_k$ if and only if  $H_k$ is a free $R_k$-module of rank $r_k$ (see lemma 2.3 in \cite{marinG26}). The result then follows from proposition 2.4(1) in \cite{marinG26}.
\end{proof}
\section{The BMR freeness conjecture for $G_8$ and $G_{16}$}
\indent

In proposition \ref{rp} we saw that in order to prove the BMR freeness conjecture for $G_8$ and $G_{16}$ we only 
need to find a spanning set of $H_k$, $k=4,5$ of $r_k$ elements. For this purpose we follow the idea I. Marin used in \cite{marincubic}, theorem 3.2 $(3)$, in order to find a spanning set for the cubic Hecke algebra. 

More precisely, for every $k\in\{4,5\}$, let $w_1$ denote the subgroup of $W_k$ generated by $s_1$ and let $\mathcal{J}$ denote a system of representatives of the double cosets
$w_1 \backslash W_k / w_1$. We have $W_k=\bigsqcup\limits_{w\in \mathcal{J}}w_1\cdot w \cdot w_1$. For every $w\in \mathcal{J}$ we fix a factorization $f_w$ of $w$ into a product of the generators $s_1$ and $s_2$ of $W_k$ and we define an element $T_{f_w}$ inside $H_k$ as follows:
 Let $f_w=x_1^{a_1}x_2^{a_2}\dots x_r^{a_r},$
where $x_i\in \{s_1,s_2\}$ and $a_i\in \ZZ$. We define the element $T_{f_w}$ to be the product $x_1^{a_1}x_2^{a_2}\dots x_r^{a_r}$ inside $H_k$ (recall that we identify $s_i$ with their images in $H_k$).
For $f_w=1$, we define $T_1:=1_{H_k}$. We notice that the element $T_{f_w}$ depends on the factorization $f_w$, meaning that if we choose a different factorization $f'_w$ of $w$ we may have $T_{f_w}\not=T_{f'_w}$. We set $U:=\sum\limits_{w\in \mathcal{J}}u_1\cdot T_{f_w} \cdot u_1$. The main result of this chapter is that $U$ is generated as $R_k$-module by $r_k$ elements and that $H_k=U$. We prove this result by using a case-by-case analysis.
\subsection{The quartic Hecke algebra $H_4$}
\label{sec2}
\indent
 Our ring of definition is $R_4=\ZZ[a, b,c,d, d^{-1}]$ and therefore, relation (\ref{one}) becomes  $s_i^4=as_i^3+bs_i^2+cs_i+d$, for $i=1,2$. We set
 
 $$\begin{array}{lcl}
 U'&=&u_1u_2u_1+u_1s_2s_1^{-1}s_2u_1+
 u_1s_2^{-1}s_1s_2^{-1}u_1+u_1s_2^{-1}s_1^{-2}s_2^{-1} \\  U&=&U'+u_1s_2s_1^{-2}s_2u_1+u_1s_2^{-2}s_1^{-2}s_2^{-2}u_1.\end{array}$$
 It is obvious that $U$ is a $u_1$-bimodule and that $U'$ is a $u_1$-sub-bimodule of $U$. Before proving our main theorem (theorem \ref{th}) we  need a few preliminaries results. 
 \begin{lem}For every $m \in \ZZ$ we have 
 	\begin{itemize}[leftmargin=0.8cm]
 		\item[(i)] $s_2s_1^{m}s_2\in U$.
 		\item [(ii)]$s_2^{-1}s_1^{m}s_2^{-1}\in U'$.
 		\item [(iii)]$s_2^{-2}s_1^{m}s_2^{-1}\in U'$.
 	\end{itemize}
 	\label{lem2}
 \end{lem}
 \begin{proof}By using the relations (\ref{one}) and (\ref{two}) we can assume that $m\in\{0,1,-1,-2\}$. Hence, we only have to prove (iii), since (i) and (ii)  follow from the definition of $U$ and $U'$ and the braid relation.  For (iii), we can assume that $m\in\{-2, 1\}$, since the case where $m=-1$  is obvious by using the generalized braid relations. We have: 
 	$s_2^{-2}s_1^{-2}s_2^{-1}=s_1^{-1}(s_1s_2^{-2}s_1^{-1})s_1^{-1}s_2^{-1}
 	=s_1^{-1}s_2^{-1}s_1^{-2}(s_2s_1^{-1}s_2^{-1})
 	=s_1^{-1}(s_2^{-1}s_1^{-3}s_2^{-1})s_1.$ The result then follows from $(ii)$.
 	For the element $s_2^{-2}s_1s_2^{-1}$, we expand $s_2^{-2}$ as  a linear combination of $s_2^{-1}, 1, s_2, s_2^2$ and by using the definition of $U'$ and lemma $\ref{lem1}$,  we only have to check that $s_2^2s_1s_2^{-1}\in U'$. Indeed, we have: $s_2^2s_1s_2^{-1}= s_2(s_2s_1s_2^{-1})=(s_2s_1^{-1}s_2)s_1
 	\in U'.$
 \end{proof}
 \begin{prop}$u_2u_1u_2\subset U.$
 	\label{prop1}
 \end{prop}
 \begin{proof}We need to prove that every element of the form $s_2^{\gra}s_1^{\grb}s_2^{\grg}$ belongs to $U$, for $\gra,\grb,\grg\in\{-2,-1,0,1\}$. However, when $\gra\grb\grg=0$ the result is obvious. Therefore, we can assume $\gra,\grb,\grg\in\{-2,-1,1\}$. We have the following cases:
 	\begin{itemize}[leftmargin=*]
 		\item\underline{$\gra=1$}:
 		The cases where $\grg\in\{-1, 1\}$ follow from lemmas $\ref{lem1}$ and \ref{lem2}(i). Hence, we need to prove that $s_2s_1^{\grb}s_2^{-2}\in U$. 
 		For $\grb=-1$ we use lemma $\ref{lem1}$ and we have $s_2s_1^{-1}s_2^{-2}=(s_2s_1^{-1}s_2^{-1})s_2^{-1}=s_1^{-1}(s_2^{-1}s_1s_2^{-1})\in U.$
 		For $\grb=1$ we  expand $s_2^{-2}$ as  a linear combination of $s_2^{-1}, 1, s_2, s_2^2$ and the result follows from the cases where $\grg\in\{-1,0,1\}$ and the generalized braid relations. It remains 
 		to prove that $s_2s_1^{-2}s_2^{-2}\in U$. By expanding now $s_1^{-2}$ as  a linear combination of $s_1^{-1}, 1, s_1, s_1^2$  we only need to prove that $s_2s_1^2s_2^{-2}\in U$ (the rest of the cases correspond to $b=-1$, $b=0$ and $b=1$).
 		We use  lemma $\ref{lem1}$ and we have:
 		$s_2s_1^2s_2^{-2}=(s_2s_1^2s_2^{-1})s_2^{-1}=s_1^{-1}s_2(s_2s_1s_2^{-1})=
 		s_1^{-1}(s_2s_1^{-1}s_2)s_1 \in U.$
 		\item \underline{$\gra=-1$}: Exactly as in the case where $\gra=1$, we only have to prove that $s_2^{-1}s_1^{\grb}s_2^{-2}\in U$. For $\grb=-1$ the result is obvious by using the generalized braid relations. For $\grb=-2$ we have:
 		$s_2^{-1}s_1^{-2}s_2^{-2}=(s_2^{-1}s_1^{-2}s_2)s_2^{-3}=s_1s_2^{-1}(s_2^{-1}s_
 		1^{-1}s_2^{-3})
 		=s_1(s_2^{-1}s_1^{-3}s_2^{-1})s_1^{-1}$. However,  by lemma \ref{lem2}(ii) we have that the element $s_2^{-1}s_1^{-3}s_2^{-1}$ is inside $U'$ and, hence, inside $U$. It remains to prove that
 		$s_2^{-1}s_1s_2^{-2}\in U$.  For this purpose, we expand $s_2^{-2}$ as  a linear combination of $s_2^{-1}, 1, s_2, s_2^2$ and by the definition of $U$ and lemma \ref{lem1} we only need to prove that $s_2^{-1}s_1s_2^{2}\in U$.
 		Indeed, using lemma \ref{lem1} again we have:
 		$s_2^{-1}s_1s_2^2=(s_2^{-1}s_1s_2)s_2
 		=s_1(s_2s_1^{-1}s_2)
 		\in U.$  \item
 		\underline{$\gra=-2$}:
 		We can assume that $\grg\in\{1, -2\}$, since the case where $\grg=-1$ follows immediately from lemma \ref{lem2}$(iii)$.
 		For $\grg=1$ we use lemma $\ref{lem1}$ and we have $s_2^{-2}s_1^{\grb}s_2=s_2^{-1}(s_2^{-1}s_1^{\grb}s_2)=(s_2^{-1}s_1s_2^{\grb})s_1^{-1}$. The latter is an element in $U$, as we proved in the case where $\gra=-1$. For $\grg=-2$ we only need to prove the cases where $\grb=\{-1,1\}$, since the case where $\grb=-2$ follows from the definition of $U$. We use the generalized braid relations and we have  $s_2^{-2}s_1^{-1}s_2^{-2}=(s_2^{-2}s_1^{-1}s_2^{-1})s_2^{-1}=s_1^{-1}(s_2^{-1}s_1^{-2}s_2^{-1})\in U$. Moreover,  $s_2^{-2}s_1s_2^{-2}=s_1(s_1^{-1}s_2^{-2}s_1)s_2^{-2}= s_1(s_2s_1^{-2}s_2^{-3})$. The result follows from the case where $\gra=1$, if we expand $s_2^{-3}$ as a linear combination of $s_2^{-2}$, $s_2^{-1}$, 1 and $s_2$.
 		\qedhere
 	\end{itemize}
 \end{proof}
 We can now prove the main theorem of this section.
 \begin{thm}
 	\mbox{}  
 	\vspace*{-\parsep}
 	\vspace*{-\baselineskip}\\ 
 	\begin{itemize}[leftmargin=0.8cm]
 		\item [(i)]
 		$U=u_1u_2u_1+u_1s_2s_1^{-1}s_2u_1
 		+u_1s_2^{-1}s_1s_2^{-1}u_1
 		+u_1\grv
 		+u_1\grv^{-1}
 		+u_1\grv^{-2}$.
 		\item[(ii)]$H_4=U$.
 	\end{itemize}
 	\label{th}
 \end{thm}
 \begin{proof}\mbox{}  
 	\vspace*{-\parsep}
 	\vspace*{-\baselineskip}\\ 
 	\begin{itemize}[leftmargin=0.6cm]
 		\item [(i)]
 		We recall that $\grv=s_2s_1^2s_2$.
 		We must prove that the RHS, which is by definition $U'+u_1\grv
 		+u_1\grv^{-2}$, is  equal to $U$. For this purpose we will ``replace''  inside the definition of $U$ the elements $s_2s_1^{-2}s_2$ and $s_2^{-2}s_1^{-2}s_2^{-2}$ with the elements $\grv$ and $\grv^{-2}$ modulo $U'$, by proving that $s_2s_1^{-2}s_2\in u_1^{\times}\grv +U'$ and $s_2^{-2}s_1^{-2}s_2^{-2}\in u_1^{\times}\grv^{-2} +U'$. 
 		
 		For the element $s_2s_1^{-2}s_2$, we expand $s_1^{-2}$ as a linear combination of $s_1^{-1}, 1, s_1, s_1^2$, where the coefficient of $s_1^2$ is invertible. The result then follows from  the definition of $U'$ and the braid relation.
 		For the element $s_2^{-2}s_1^{-2}s_2^{-2}$ we apply lemma \ref{lem1} and the generalized braid relations and we have: 
 		$s_2^{-2}s_1^{-2}s_2^{-2}=s_2^{-2}s_1^{-1}(s_1^{-1}s_2^{-2}s_1)s_1^{-1}=s_2^{-1}(s_2^{-1}s_1^{-1}s_2)s_1^{-1}(s_1^{-1}s_2^{-1}s_1)s_1^{-2}=s_2^{-1}s_1(s_2^{-1}s_1^{-2}s_2)s_1^{-1}s_2^{-1}s_1^{-2}\in s_2^{-1}s_1^2s_2^{-2}s_1^{-2}s_2^{-1}u_1$.
 		We expand $s_1^2$ as a linear combination of $s_1$, 1, $s_1^{-1}$, $s_1^{-2}$, where the coefficient of $s_1^{-2}$ is invertible and by the generalized braid relations and the fact that $s_1^{-2}\grv^{-2}=\grv^{-2}s_1^{-2}=s_2^{-1}s_1^{-2}s_2^{-2}s_1^{-2}s_2^{-1}s_1^{-2}$ we have that $$s_2^{-2}s_1^{-2}s_2^{-2}
 		\in s_2^{-1}s_1s_2^{-2}s_1^{-2}s_2^{-1}u_1+s_2^{-3}s_1^{-2}s_2^{-1}u_1+u_1s_2^{-1}s_1^{-3}s_2^{-1}u_1+
 		u_1^{\times}\grv^{-2}.$$
 		Therefore, by lemma \ref{lem2}(ii) it is enough to prove that the elements $s_2^{-1}s_1s_2^{-2}s_1^{-2}s_2^{-1}$ and $s_2^{-3}s_1^{-2}s_2^{-1}$ belong to $U'$. However, the latter is an element in $U'$,  if we expand $s_2^{-3}$ as a linear combination of $s_2^{-2}, s_2^{-1}, 1, s_2$ and use lemma \ref{lem2}(iii), the definition of $U'$ and lemma \ref{lem1}. Moreover,
 		$s_2^{-1}s_1s_2^{-2}s_1^{-2}s_2^{-1}=s_2^{-2}(s_2s_1s_2^{-1})\grv^{-1}=
 		s_2^{-2}s_1^{-1}s_2s_1\grv^{-1}=s_2^{-2}s_1^{-1}s_2\grv^{-1}s_1^{-1}=(s_2^{-2}s_1^{-4}s_2^{-1})s_1^{-1}
 		\in U'$, by lemma \ref{lem2}(iii).
 		
 		\item [(ii)]
 		Since $1\in U$, it will be sufficient to show that $U$ is a left ideal of $H_4$. We know that $U$ is a $u_1$-sub-bimodule of $H_4$. Therefore, we only need to prove that $s_2U\subset U$. Since $U$ is equal to the RHS of (i) we have that
 		$$s_2U\subset s_2u_1u_2u_1+s_2u_1s_2s_1^{-1}s_2u_1+s_2u_1s_2^{-1}s_1s_2^{-1}u_1+ s_2u_1\grv+s_2u_1\grv^{-1}+s_2u_1\grv^{-2}.$$ However, $s_2u_1u_2u_1+s_2u_1\grv+s_2u_1\grv^{-1}+s_2u_1\grv^{-2}=s_2u_1u_2u_1+s_2\grv u_1+s_2\grv^{-1}u_1+s_2\grv^{-2}u_1=s_2u_1u_2u_1+s_2^3s_1^2s_2u_1+s_1^{-2}s_2^{-1}u_1+s_1^{-2}s_2^{-2}s_1^{-1}s_2^{-1}
 		\subset u_1u_2u_1u_2u_1$. Furthermore,
 		by lemma \ref{lem1} we have $s_2u_1s_2^{-1}=s_1^{-1}u_2s_1$. Hence, $
 		s_2u_1s_2s_1^{-1}s_2u_1=(s_2u_1s_2^{-1})s_2^2s_1^{-1}s_2u_1=s_1^{-1}u_2(s_1s_2^2s_1^{-1})s_2u_1=
 		s_1^{-1}u_2s_1^2s_2^2u_1\subset u_1u_2u_1u_2u_1.$ Moreover, by using \ref{lem1} again we have that $(s_2u_1s_2^{-1})s_1s_2^{-1}u_1=s_1^{-1}u_2s_1^2s_2^{-1}u_1\subset u_1u_2u_1u_2u_1$ . Therefore, $$ s_2u_1u_2u_1+s_2u_1s_2s_1^{-1}s_2u_1+s_2u_1s_2^{-1}s_1s_2^{-1}u_1+ s_2u_1\grv+s_2u_1\grv^{-1}+s_2u_1\grv^{-2}\subset u_1u_2u_1u_2u_1.$$
 		The result follows directly from proposition \ref{prop1}.
 		\qedhere
 	\end{itemize}
 \end{proof}
\begin{cor}$H_4$ is a free $R_4$-module of rank $r_4=96$ and, therefore, the BMR freeness conjecture holds for the exceptional group $G_8$.
	\label{G8}
\end{cor}
\begin{proof}By proposition \ref{rp} it will be sufficient to show that $H_4$ is generated as $R_4$-module by $r_4$ elements. By Theorem  \ref{th} and the fact that $u_1u_2u_1=u_1(R_4+R_4s_2+R_4s_2^{-1}+R_4s_2^2)u_1=u_1+u_1s_2u_1+u_1s_2^{-1}u_1+u_1s_2^2u_1$ we have that $H_4$ 
is generated as left $u_1$-module by 24 elements. Since $u_1$ is generated by 4 elements as a $R_4$-module, we have that $H_4$ is generated over $R_4$ by $r_4=96$ elements.
\end{proof}
\subsection{The quintic Hecke algebra $H_5$}
\indent

Our ring of definition is $R_5=\ZZ[a, b,c, d,e,e^{-1}]$ and therefore, relation (\ref{one}) becomes $s_i^5=as_i^4+bs_i^3+cs_i^2+ds_i+e$, for $i=1,2$.
We recall that $\grv=s_2s_1^2s_2$ and we set 
$$\small{\begin{array}{lcl}U'&=&u_1u_2u_1+u_1\grv+u_1\grv^{-1}+
	u_1s_2^{-1}s_1^2s_2^{-1}u_1+u_1s_2s_1^{-2}s_2u_1+u_1s_2^2s_1^2s_2^{2}u_1+
	u_1s_2^{-2}s_1^{-2}s_2^{-2}u_1+\\&&+u_1s_2s_1^{-2}s_2^{2}u_1+
	u_1s_2^{-1}s_1^2s_2^{-2}u_1+u_1s_2^{-1}s_1s_2^{-1}u_1
	+u_1s_2s_1^{-1}s_2u_1+u_1s_2^{-2}s_1^{-2}s_2^{2}u_1+
	u_1s_2^{2}s_1^2s_2^{-2}u_1+\\&&+u_1s_2^{2}s_1^{-2}s_2^{2}u_1+
	u_1s_2^{-2}s_1^2s_2^{-2}u_1+u_1s_2^{-2}s_1s_2^{-1}u_1+u_1s_2^{-1}s_1s_2^{-2}u_1\\ \\
	U''&=&U'+u_1\grv^2+u_1\grv^{-2}+
	u_1s_2^{-2}s_1^2s_2^{-1}s_1s_2^{-1}u_1+u_1s_2^{2}s_1^{-2}s_2s_1^{-1}s_2u_1+
	u_1s_2s_1^{-2}s_2^{2}s_1^{-2}s_2^{2}u_1+\\&&+u_1s_2^{-1}s_1^2s_2^{-2}s_1^2s_2^{-2}u_1\\ \\
	U'''&=&U''+u_1\grv^3+u_1\grv^{-3}\\ \\
	U''''&=&U'''+u_1\grv^4+u_1\grv^{-4}\\ \\
	U&=&U''''+u_1\grv^5+u_1\grv^{-5}.
	\end{array}}$$

It is obvious that $U$ is a $u_1$-bi-module and that $U', U'', U'''$ and $U''''$ are $u_1-$ sub-bi-modules of $U$. Again, our goal is to 
prove that $H_5=U$ (theorem \ref{thh2}). As we explained in the proof of theorem \ref{th}, since $1\in U$ and $s_1U\subset U$ (by the definition of $U$), it is enough to prove that $s_2U\subset U$.
We notice that $$U=\sum_{k=1}^5u_1\grv^{\pm k}+\underbrace{u_1u_2u_1+u_1\text{``some elements of length 3''}u_1}_{\in U'}+\underbrace{u_1\text{``some elements of length 5''}u_1}_{\in U''}.$$ 
By the definition of $U'$ and $U''$ we have that $u_1\grv^{\pm1}\subset U'$ and $u_1\grv^{\pm2}\subset U''$. Therefore, in order to prove that $s_2U\subset U$ we only need to prove that $s_2u_1\grv^{\pm k}$ ($k=3,4,5$), $s_2U'$ and $s_2U''$ are subsets of $U$. 

The rest of this section is devoted to this proof (see proposition \ref{p2}, lemma \ref{oo}$(ii)$, proposition \ref{xx}$(i),(ii)$ and theorem \ref{thh2}). The reason we define also $U'''$ and $U''''$ is because, in order to prove that $s_2u_1\grv^k$ and $s_2u_1\grv^{-k}$ ($k=3,4,5$) are subsets of $U$, we want to ``replace'' inside the definition of $U$ the elements $\grv^k$ and $\grv^{-k}$ by some other elements modulo  $U'', U'''$ and $U''''$, respectively (see lemmas \ref{cc}, \ref{ll} and \ref{lol}).

Recalling that $\grF$ is the automorphism of $H_5$ as defined in section \ref{s}, we have the following lemma:
\begin{lem}The $u_1$-bi-modules $U', U'', U''', U''''$ and $U$ are stable under $\grF$.
	\label{r1}
\end{lem}
\begin{proof}We notice that $U',U'', U''', U''''$ and $U$ are of the form $$u_1s_2^{-2}s_1s_2^{-1}u_1+u_1s_2^{-1}s_1s_2^{-2}u_1+\sum u_1\grs u_1+\sum u_1\grs^{-1}u_1,$$ for some $\grs\in B_3$ satisfying $\grs^{-1}=\grF(\grs)$ and $\grs=\grF(\grs^{-1})$. Therefore, 
	we restrict ourselves to proving that  the elements $\grF(s_2^{-2}s_1s_2^{-1})=s_2^2s_1^{-1}s_2$ and $\grF(s_2^{-1}s_1s_2^{-2})=s_2s_1^{-1}s_2^2$ belong to $U'$. We expand $s_2^2$ as a linear combination of $s_2,1, s_2^{-1}, s_2^{-2}$ and $s_2^{-3}$ and by the definition of $U'$ and lemma \ref{lem1} we have to prove that the elements $s_2^{k}s_1^{-1}s_2$ and $s_2s_1^{-1}s_2^{k}$  are elements in $U'$, for $k=-3, -2$. Indeed, by using lemma \ref{lem1} we have: $s_2^{k}s_1^{-1}s_2=s_2^{k+1}(s_2^{-1}s_1^{-1}s_2)=(s_2^{k+1}s_1s_2^{-1})s_1^{-1}\in U'$ and  $s_2s_1^{-1}s_2^{k}=(s_2s_1^{-1}s_2^{-1})s_2^{k+1}=s_1^{-1}(s_2^{-1}s_1s_2^{k+1})\in U'$.
\end{proof}
From now on, we will use lemma \ref{lem1} without mentioning it.
\begin{prop} $u_2u_1u_2\subset U'.$
	\label{p1}
\end{prop}
\begin{proof}
	We have to prove that every element of the form $s_2^{\gra}s_1^{\grb}s_2^{\grg}$ belongs to $U'$, for $\gra,\grb,\grg\in\{-2,-1,0,1,2\}$. However, when $\gra\grb\grg=0$ the result is obvious. Therefore, we can assume that $\gra,\grb,\grg\in \{-2,-1,1,2\}.$ We continue the proof as in the proof of proposition \ref{prop1}, which is by distinguishing cases for $\gra$. However, by using lemma \ref{r1} we can assume that $\gra\in\{1,2\}$. We have:
	\begin{itemize}[leftmargin=*] 
		\item\underline{$\gra=1$}: 
		\begin{itemize}[leftmargin=*] 
			\item \underline{$\grg\in\{-1,1\}$}: The result follows from lemma \ref{lem1}, the braid relation and the definition of  $U'$.
			\item \underline{$\grg=-2$}:  $s_2s_1^{\grb}s_2^{-2}=(s_2s_1^{\grb}s_2^{-1})s_2^{-1}=s_1^{-1}(s_2^{\grb}s_1s_2^{-1})$. For $\grb\in\{1,-1,-2\}$ the result follows from lemma \ref{lem1} and the definition of $U'$. For $\grb=2$, we have $s_1^{-1}s_2^2s_1s_2^{-1}=s_1^{-1}s_2(s_2s_1s_2^{-1})=s_1^{-1}(s_2s_1^{-1}s_2)s_1\in U'$.
			\item \underline{$\grg=2$}:
			We need to prove that the element $s_2s_1^{\grb}s_2^{2}$ is inside $U'$. For $\grb\in\{-2,1\}$ the result is obvious by using the definition of $U'$ and the generalized braid relations. For $\grb=-1$ we have $s_2s_1^{-1}s_2^2=\grF(s_2^{-1}s_1s_2^{-2})\in\grF(U')\stackrel{\ref{r1}}{=}U'$. For $\grb=2$ we have $s_2s_1^2s_2^2= s_1^{-1}(s_1s_2s_1^2)s_2^2=s_1^{-1}s_2(s_2s_1s_2^3)=s_1^{-1}(s_2s_1^3s_2)s_1$. The result then follows from the case where $\grg=1$, if we expand $s_1^3$ as a linear combination of $s_1^2, s_1, 1, s_1^{-1}, s_1^{-2}$.
		\end{itemize}
		\item \underline {$\gra=2$}:
		\begin{itemize}[leftmargin=*] 
			\item \underline{$\grg=-1$}:  $s_2^2s_1^{\grb}s_2^{-1}=s_2(s_2s_1^{\grb}s_2^{-1})=(s_2s_1^{-1}s_2^{\grb})s_1\in U'$ (case where $\gra=1$).
			\item \underline{$\grg=2$}: We only have to prove the cases where $\grb\in\{-1,1\}$, since the cases where $\grb\in \nolinebreak \{2,-2\}$ follow from the definition of $U'$. We have $s_2^2s_1s_2^2=(s_2^2s_1s_2)s_2=s_1\grv\in U'$. Moreover,
			$s_2^2s_1^{-1}s_2^2=s_1^{-1}(s_1s_2^2s_1^{-1})s_2^2=s_1^{-1}\grF(s_2s_1^{-2}s_2^{-3})$. The result follows from the case where $\gra=1$ and lemma
			\ref{r1}, if we expand $s_2^{-3}$ as  a linear combination of $s_2^{-2}, s_2^{-1}, 1, s_2, s_2^{2}$.
			\item \underline{$\grg=1$}: We  have to check the cases where $\grb\in\{-2,-1,2\}$, since the case where $\grb=1$ is a direct result from the generalized braid relations. However,   $s_2^2s_1^{-1}s_2=\grF(s_2^{-2}s_1s_2^{-1})\in\grF(U')\stackrel{\ref{r1}}{=}U'$. Hence, it remains to prove the cases where $\grb\in \{-2,2\}$. 
			We have $s_2^2s_1^{-2}s_2=s_2^3(s_2^{-1}s_1^{-2}s_2)=s_1(s_1^{-1}s_2^3s_1)s_2^{-2}s_1^{-1}=
			s_1(s_2s_1^3s_2^{-3})s_1^{-1}$. The latter is an element in $U'$, if we expand $s_1^3$ and $s_2^{-3}$ as  linear combinations of $s_1^2, s_1, 1, s_1^{-1}, s_1^{-2}$ and $s_2^{-2}, s_2^{-1}, 1, s_2, s_2^{2}$, respectively and use the case where $\gra=1$. Moreover, 
			$s_2^2s_1^2s_2=
			s_2^2s_1(s_1s_2s_1)s_1^{-1}=(s_2^2s_1s_2)s_1s_2s_1^{-1}=s_1(s_2s_1^3s_2)s_1.$ The result follows again from the case where $\gra=1$, if we expand $s_1^3$ as a linear combination of $s_1^2, s_1, 1, s_1^{-1}, s_1^{-2}$. 
			
			\item\underline{$\grg=-2$}: We need to prove that $s_2^2s_1^{\grb}s_2^{-2}\in U'$. For $\grb=2$ the result follows from the definition of $U'$. For $\grb\in\{1,-1\}$ we have: $s_2^2s_1s_2^{-2}=s_2^2(s_1s_2^{-2}s_1^{-1})s_1=(s_2s_1^{-2}s_2)s_1\in U'$.
			$s_2^2s_1^{-1}s_2^{-2}=s_2(s_2s_1^{-1}s_2^{-1})s_2^{-1}=(s_2s_1^{-1}s_2^{-1})s_1s_2^{-1}=
			s_1^{-1}(s_2^{-1}s_1^2s_2^{-1})\in U'$. 
			It remains to prove the case where $\grb=-2$. We recall that $\grv =s_2s_1^2s_2$ and we have:
			$ s_2^2s_1^{-2}s_2^{-2}=s_1^{-1}(s_1s_2^2s_1^{-1})s_1^{-1}s_2^{-2}=s_1^{-1}s_2^{-2}\grv s_1^{-1}s_2^{-2}=s_1^{-1}s_2^{-2} s_1^{-1}\grv s_2^{-2}=
			s_1^{-1}s_2^{-2}s_1^{-1}(s_2s_1^2s_2^{-1})= s_1^{-1}(s_2^{-2}s_1^{-2}s_2^2)s_1. $ The result follows from the definition of $U'$.
			\qedhere
		\end{itemize}
	\end{itemize}
\end{proof}
From now on, in order to make it easier for the reader to follow the calculations, we will underline the elements belonging to $u_1u_2u_1u_2u_1$ and we will use immediately the fact that these elements belong to $U'$ (see proposition \ref{p1}).
\begin{lem}
	\mbox{}  
	\vspace*{-\parsep}
	\vspace*{-\baselineskip}\\ 
	\begin{itemize}[leftmargin=0.8cm]
		\item [(i)]$s _2u_1s_2u_1s_2u_1\subset\grv^2u_1+u_1u_2u_1u_2u_1\subset U''$.  
		\item [(ii)]$s_2\grv^2u_1=s_1s_2s_1^4s_2s_1^3s_2u_1\subset U''.$
	\end{itemize}  
	\label{oo}
\end{lem}
\begin {proof}
We recall that $\grv=s_2s_1^2s_2$.
\begin{itemize}[leftmargin=0.8cm]
	\item[(i)] The fact that $\grv^2u_1+u_1u_2u_1u_2u_1\subset U''$ follows directly from the definition of $U''$ and  proposition \ref{p1}. For the rest of the proof, we use the definition of $u_1$ and we have that  $s_2u_1s_2u_1s_2u_1=s_2u_1s_2(R_5+R_5s_1^{-1}+R_5s_1+R_5
	s_1^{2}+R_5s_1^3)s_2u_1\subset \underline{s_2u_1s_2^2u_1}+s_2u_1s_2s_1^{-1}s_2u_1
	+\underline{s_2u_1(s_2s_1s_2)u_1}+
	s_2u_1\grv+s_2u_1s_2s_1^{3}s_2u_1$.
	However, $s_2u_1\grv=\underline{s_2\grv u_1}$ and $s_2u_1s_2s_1^{-1}s_2u_1=s_2u_1(s_1s_2s_1^{-1})s_2u_1= (s_2u_1s_2^{-1})s_1s_2^2u_1=\underline{s_1^{-1}u_2s_1^2s_2^2u_1}$. Therefore, it is enough to prove that $s_2u_1s_2s_1^3s_2u_1\subset\grv^2u_1+u_1u_2u_1u_2u_1$. For this purpose, we use again the definition of $u_1$ and we have:
	
	$\small{\begin{array}[t]{lcl}
		s_2u_1s_2s_1^3s_2u_1
		&\subset&
		s_2(R_5+R_5s_1+R_5s_1^{-1}+R_5s_1^2+R_5s_1^3)s_2s_1^{3}s_2u_1\\
		&\subset& \underline{s_2^2s_1^3s_2u_1}+
		\underline{s_2(s_1s_2s_1^{3})s_2u_1}+
		s_2(s_1^{-1}s_2s_1)s_1^{2}s_2u_1+\grv s_1^{3}s_2u_1+s_2s_1^{2}(s_1s_2s_1^{3})s_2u_1\\
		&\subset& \underline{s_2^2s_1(s_2^{-1}s_1^2s_2)u_1}+\underline{s_1^{3}\grv s_2u_1}+s_2s_1^2s_2^2(s_2s_1s_2^2)u_1+u_1u_2u_1u_2u_1\\
		&\subset&\grv^2u_1+u_1u_2u_1u_2u_1.
		\end{array}}$
	\item[(ii)] We have that $s_2\grv^2=s_1(s_1^{-1}s_2^2s_1)(s_1s_2s_1)s_1^{-1}\grv=s_1s_2s_1^4(s_1^{-1}s_2s_1)s_1^{-2}\grv=s_1s_2s_1^4s_2s_1s_2^{-1}s_1^{-2}\grv=s_1s_2s_1^4s_2s_1s_2^{-1}\grv s_1^{-2}=s_1s_2s_1^4s_2s_1^3s_2s_1^{-2}$. Therefore, $s_2\grv^2u_1
	\subset
	u_1s_2u_1s_2u_1s_2u_1.$	The fact that $u_1s_2u_1s_2u_1s_2u_1\subset U''$ follows immediately from (i).
	\qedhere
\end{itemize}
\end{proof}
\begin{prop}\mbox{}  
	\vspace*{-\parsep}
	\vspace*{-\baselineskip}\\ 
	\begin{itemize}[leftmargin=0.8cm]
		\item [(i)]$u_2u_1s_2^{-1}s_1s_2^{-1}\subset u_1\grv^{-2}+R_5s_2^{-2}s_1^2s_2^{-1}s_1s_2^{-1}+u_1u_2u_1u_2u_1\subset U''.$
		\item[(ii)] $u_2u_1s_2s_1^{-1}s_2\subset u_1 \grv^{2}+R_5s_2^{2}s_1^{-2}s_2s_1^{-1}s_2+u_1u_2u_1u_2u_1\subset U''.$
	\end{itemize}
	\label{l2}
\end{prop}
\begin{proof}
	We restrict ourselves to proving $(i)$, since $(ii)$ follows from $(i)$ by applying $\grF$ (see lemma \ref{r1}). By the definition of $U''$ and by proposition \ref{p1} we have that  $u_1\grv^{-2}+R_5s_2^{-2}s_1^2s_2^{-1}s_1s_2^{-1}+u_1u_2u_1u_2u_1\subset U''$. Therefore, it remains to prove that  $u_2u_1s_2^{-1}s_1s_2^{-1}\subset u_1\grv^{-2}+R_5s_2^{-2}s_1^2s_2^{-1}s_1s_2^{-1}+u_1u_2u_1u_2u_1.$
	
	By he definition of $u_1$ we have
	$u_2u_1s_2^{-1}s_1s_2^{-1}=u_2(R_5+R_5s_1+R_5s_1^{-1}+R_5s_1^{-2}
	+R_5s_1^{2})s_2^{-1}s_1s_2^{-1}
	\subset \underline{u_2s_1s_2^{-1}}+
	u_2s_1s_2^{-1}s_1s_2^{-1}+\underline{u_2(s_1^{-1}s_2^{-1}s_1)s_2^{-1}}+u_2s_1^{-2}s_2^{-1}s_1s_2^{-1}
	+u_2s_1^{2}s_2^{-1}s_1s_2^{-1}$.
	However, $u_2s_1s_2^{-1}s_1s_2^{-1}=u_2(s_2s_1s_2^{-1})s_1s_2^{-1}=\underline{u_2s_1^{-1}(s_2s_1^2s_2^{-1})}$.
	Therefore, we only have to prove that $u_2s_1^{-2}s_2^{-1}s_1s_2^{-1}$ and 
	$u_2s_1^{2}s_2^{-1}s_1s_2^{-1}$ are subsets of $ u_1\grv^{-2}+R_5s_2^{-2}s_1^2s_2^{-1}s_1s_2^{-1}+u_1u_2u_1u_2u_1$. We have: \\ 
	$\small{\begin{array}[t]{lcl}
		u_2s_1^{-2}s_2^{-1}s_1s_2^{-1}
		&\subset&(R_5+R_5s_2+R_5s_2^{-1}+
		R_5s_2^2+R_5s_2^3)s_1^{-2}s_2^{-1}s_1s_2^{-1}\\
		&\subset&\underline{R_5s_1^{-2}s_2^{-1}s_1s_2^{-1}}+\underline{
			R_5(s_2s_1^{-2}s_2^{-1})s_1s_2^{-1}}+R_5\grv^{-1}s_1s_2^{-1}+
		R_5s_2(s_2s_1^{-2}s_2^{-1})s_1s_2^{-1}+\\&&+
		R_5s_2^2(s_2s_1^{-2}s_2^{-1})s_1s_2^{-1}\\
		
		&\subset&\underline{R_5s_1\grv^{-1}s_2^{-1}}+R_5(s_2s_1^{-1}s_2^{-1})s_2^{-1}s_1^2s_2^{-1}+
		R_5s_2(s_2s_1^{-1}s_2^{-1})s_2^{-1}s_1^2s_2^{-1}+u_1u_2u_1u_2u_1\\
		&\subset&R_5s_1^{-1}s_2^{-1}s_1s_2^{-1}s_1^2s_2^{-1}+R_5(s_2s_1^{-1}s_2^{-1})s_1s_2^{-1}s_1^2s_2^{-1}+u_1u_2u_1u_2u_1\\
		&\subset&\grF(u_1s_2u_1s_2u_1s_2)+u_1u_2u_1u_2u_1.
		\end{array}}$\\ \\
	However, by lemma \ref{oo}(i) we have that $\grF(u_1s_2u_1s_2u_1s_2)\subset \grF(\grv^2u_1+u_1u_2u_1u_2u_1)=\grv^{-2}u_1+u_1u_2u_1u_2u_1$. Therefore, 
	$ u_2s_1^{-2}s_2^{-1}s_1s_2^{-1}\subset 	\grv^{-2}u_1+u_1u_2u_1u_2u_1$.	
	By using analogous calculations, we have: \\
	$\small{\begin{array}[t]{lcl}
		u_2s_1^{2}s_2^{-1}s_1s_2^{-1}	&\subset&(R_5+R_5s_2
		+R_5s_2^{-1}+R_5s_2^2+R_5s_2^{-2})s_1^{2}s_2^{-1}s_1s_2^{-1}\\
		&\subset& \underline{R_5s_1^2s_
			2^{-1}s_1s_2^{-1}}+\underline{R_5(s_2s_1^2s_2^{-1})s_1s_2^{-1}}+
		R_5s_2^{-1}s_1^3(s_1^{-1}s_2^{-1}s_1)s_2^{-1}+
		R_5s_2(s_2s_1^2s_2^{-1})s_1s_2^{-1}+\\&&+R_5s_2^{-2}s_1^2s_2^{-1}s_1s_2^{-1}\\
		&\subset&
		\underline{R_5(s_2^{-1}s_1^3s_2)s_1^{-1}s_2^{-2}}+
		R_5s_2s_1^{-1}s_2^{2}s_1^2s_2^{-1}+R_5s_2^{-2}s_1^2s_2^{-1}s_1s_2^{-1}+u_1u_2u_1u_2u_1.
		\end{array}}$\\ \\
	It is enough  to prove that $s_2s_1^{-1}s_2^{2}s_1^2s_2^{-1}\subset u_1u_2u_1u_2u_1$. Indeed, we have that $s_2s_1^{-1}s_2^2s_1^2s_2^{-1}=s_1^{-1}(s_1s_2s_1^{-1})s_2(s_2s_1^2s_2^{-1})=\underline{s_1^{-1}s_2^{-1}(s_1s_2^2s_1^{-1})s_2^2s_1}$.

	\qedhere
\end{proof}
We can now prove a lemma that helps us to ``replace''  inside the definition of $U'''$ the element $\grv^3$ with the element $s_2s_1^3s_2^2s_1^2s_2^2$ modulo $U''$.
\begin{lem}$s_2s_1^3s_2^2s_1^2s_2^2\in u_1s_2u_1s_2s_1^3s_2u_1+u_1s_2^2s_1^3s_2s_1^{-1}s_2u_1+
	u_1u_2u_1u_2u_1+ u_1^{\times}\grv^3\subset u_1^{\times}\grv^3+\nolinebreak U''.$
	\label{cc}
\end{lem}
\begin{proof}
	The fact that 	$u_1s_2u_1s_2s_1^3s_2u_1+u_1s_2^2s_1^3s_2s_1^{-1}s_2u_1+
	u_1u_2u_1u_2u_1+ u_1^{\times}\grv^3$ is a subset of $u_1^{\times}\grv^3+U''$ follows  from lemma \ref{oo}$(i)$ and propositions \ref{l2}(ii) and \ref{p1}. 
	For the rest of the proof, 	we notice that we have	
	$s_2s_1^3s_2^2s_1^2s_2^2=s_2s_1^2(s_1s_2^2s_1^{-1})s_1^2(s_1s_2^2s_1^{-1})s_1=s_2s_1^2s_2^{-2}\grv s_1^2s_2^{-1}s_1(s_1s_2s_1^{-1})s_1^2=
	s_2s_1^2s_2^{-2}s_1^2\grv s_2^{-1}s_1s_2^{-1}(s_1s_2s_1^{-1})s_1^3=s_2s_1^2s_2^{-2}s_1^2s_2s_1^3s_2^{-2}s_1s_2s_1^3=s_2s_1^2s_2^{-3}\bold{\boldsymbol{\grv} s_1^3s_2^{-2}}s_1s_2s_1^3$.
	However, $\bold{\boldsymbol{\grv} s_1^3s_2^{-2}}=s_1^3\grv s_2^{-2}=s_1^3(s_2s_1^2s_2^{-1})=s_1^2s_2^2s_1$ and, hence, $s_2s_1^3s_2^2s_1^2s_2^2=s_2s_1^2s_2^{-3}s_1^2s_2^2s_1^2s_2s_1^3$.
	
	Our goal now is to prove that  $s_2s_1^2s_2^{-3}s_1^2s_2^2s_1^2s_2s_1^3\in  u_1s_2u_1s_2s_1^3s_2u_1+u_1s_2^2s_1^3s_2s_1^{-1}s_2u_1+
	u_1u_2u_1u_2u_1+ u_1^{\times}\grv^3$. For this purpose we expand $s_2^{-3}$ as a linear combination of $s_2^{-2}$, $s_2^{-1}$, 1, $s_2$ and $s_2^2$, where the coefficient of $s_2^2$ is invertible, and we have that  
	$s_2s_1^2s_2^{-3}s_1^2s_2^2s_1^2s_2s_1^3\in s_2s_1^2s_2^{-2}s_1^2s_2^2s_1^2s_2u_1+s_2s_1^2s_2^{-1}s_1^2s_2^2s_1^2s_2u_1+
	s_2s_1^4s_2^2s_1^2s_2u_1+s_2\grv^2u_1+u_1^{\times}\grv^3$.
	However, by lemma \ref{oo}(ii) we have that $s_2\grv^2u_1\subset u_1s_2u_1s_2s_1^3s_2u_1$. Moreover, $s_2s_1^4s_2^2s_1^2s_2u_1=s_2s_1^5(s_1^{-1}s_2^2s_1)(s_1s_2s_1)u_1\subset u_1s_2u_1s_2s_1^3s_2u_1$. It remains to prove that the elements
	$s_2s_1^2s_2^{-2}s_1^2s_2^2s_1^2s_2$ and $s_2s_1^2s_2^{-1}s_1^2s_2^2s_1^2s_2$
	are inside $u_1s_2u_1s_2s_1^3s_2u_1+u_1s_2^2s_1^3s_2s_1^{-1}s_2u_1+
	u_1u_2u_1u_2u_1$.
	
	On one hand,  we have
	$s_2s_1^2s_2^{-2}s_1^2s_2^2s_1^2s_2=s_2s_1^3(s_1^{-1}s_2^{-2}s_1)s_1s_2\grv=s_2s_1^3s_2s_1^{-1}(s_1^{-1}s_2^{-1}s_1)s_2\grv=s_2s_1^3s_2^2(s_2^{-1}s_1^{-1}s_2)s_1^{-1}\grv=s_2s_1^3s_2^2s_1s_2^{-1}s_1^{-2}\grv=s_2s_1^3s_2^2s_1s_2^{-1}\grv s_1^{-2}=s_2s_1^3s_2^2s_1^3s_2s_1^{-2}$, meaning that the element $s_2s_1^2s_2^{-2}s_1^2s_2^2s_1^2s_2$ is inside $s_2s_1^3s_2^2u_1s_2u_1$. On the other hand, 
	$s_2s_1^2s_2^{-1}s_1^2s_2^2s_1^2s_2
	=s_2s_1^2(s_2^{-1}s_1^2s_2)\grv=s_2s_1^3s_2^2s_1^{-1}\grv=s_2s_1^3s_2^2\grv s_1^{-1}=s_2s_1^3s_2^3s_1^2s_2s_1^{-1}$ and, if we expand $s_2^3$ as a linear combination of $s_2^2$, $s_2$, 1, $s_2^{-1}$ and $s_2^{-2}$, we have that $s_2s_1^3s_2^3s_1^2s_2s_1^{-1}\in s_2s_1^3s_2^{2}s_1^2s_2u_1+s_2s_1^3\grv u_1+
	\underline{s_2s_1^5s_2u_1}+\underline{
		(s_2s_1^3s_2^{-1})s_1^2s_2u_1}+\underline{(s_2s_1^3s_2^{-1})(s_2^{-1}s_1^2s_2)u_1}\subset 
	s_2s_1^3s_2^{2}u_1s_2u_1+\underline{s_2\grv u_1}+u_1u_2u_1u_2u_1$, meaning that the element $s_2s_1^2s_2^{-1}s_1^2s_2^2s_1^2s_2$ is inside 	$s_2s_1^3s_2^{2}u_1s_2u_1+u_1u_2u_1u_2u_1$. As a result, in order to finish the proof, it will be sufficient to show that $s_2s_1^3s_2^{2}u_1s_2u_1$ is a subset of
	$u_1s_2u_1s_2s_1^3s_2u_1+u_1s_2^2s_1^3s_2s_1^{-1}s_2u_1+
	u_1u_2u_1u_2u_1$. Indeed, we have: 
	$$\small{\begin{array}[t]{lcl}
		s_2s_1^3s_2^{2}u_1s_2u_1		&\subset&s_2s_1^3s_2^2(R_5s_1^2+R_5s_1+R_5+R_5s_1^{-1}+R_5s_1^{-2})s_2u_1\\
		&\subset&s_2s_1^3s_2^2s_1^2s_2u_1+\underline{s_2s_1^3(s_2^2s_1s_2)u_1}+\underline{s_2s_1^3s_2^3u_1}
		+s_2s_1^2(s_1s_2^2s_1^{-1})s_2u_1+
		s_2s_1^2(s_1s_2^2s_1^{-1})s_1^{-1}s_2u_1\\
		&\subset&s_2s_1^4(s_1^{-1}s_2^2s_1)(s_1s_2s_1)u_1+\underline{
			(s_2s_1^2s_2^{-1})s_1^2s_2^2u_1}+
		(s_2s_1^2s_2^{-1})s_1^2s_2s_1^{-1}s_2u_1+u_1u_2u_1u_2u_1\\
		&\subset& u_1s_2u_1s_2s_1^3s_2u_1+u_1s_2^2s_1^3s_2s_1^{-1}s_2u_1+u_1u_2u_1u_2u_1.
		\end{array}}$$
	
\end{proof} 
\begin{prop}
	\mbox{} 
	\vspace*{-\parsep}
	\vspace*{-\baselineskip}\\
	\begin{itemize}[leftmargin=0.8cm]
		\item [(i)]$s_2u_1u_2u_1u_2 \subset U'''$.
		\item [(ii)]$s_2^{-1}u_1u_2u_1u_2 \subset U'''$.
	\end{itemize}
	\label{p2}
\end{prop}
\begin{proof}By lemma \ref{r1}, we only have to prove $(i)$, since $(ii)$ is a consequence of $(i)$ up to applying $\grF$. We know that $u_2u_1u_2 \subset U'$ (proposition \ref{p1}) hence it is enough to prove that $s_2U'\subset U'''$. Set 
	$$\small{\begin{array}{lcl}V&=&u_1u_2u_1+\grv u_1+\grv^{-1}u_1+u_1s_2^{-1}s_1^2s_2^{-1}u_1+
		u_1s_2^{-1}s_1s_2^{-1}u_1+u_1s_2s_1^{-1}s_2u_1+u_1s_2^{-2}s_1s_2^
		{-1}u_1+\\&&+u_1s_2^{-1}s_1^2s_2^{-2}u_1+u_1s_2^{-1}s_1s_2^{-2}u_1+u_1s_2s_1^{-2}s_2u_1+
		u_1s_2^{-2}s_1^{-2}s_2^{-2}u_1+
		u_1s_2^{-2}s_1^{-2}s_2^{2}u_1.
		\end{array}}$$
	We notice that 
	$$\begin{array}{lcl}
	U'&=&V+u_1s_2s_1^{-2}s_2^{2}u_1+
	u_1s_2^{2}s_1^{-2}s_2^{2}u_1+u_1s_2^{2}s_1^{2}s_2^{2}u_1+
	u_1s_2^{-2}s_1^{2}s_2^{-2}u_1+u_1s_2^{2}s_1^{2}s_2^{-2}u_1.
	\end{array}$$
	Therefore, in order to prove that $s_2U'\subset U'''$, we will prove first that $s_2V\subset U'''$ and then we will check the other five cases separately. 
	We have: 
	$$\small{\begin{array}{lcl}s_2V&\subset&
		\underline{s_2u_1u_2u_1}+\underline{s_2\grv u_1}+\underline{s_2\grv^{-1} u_1}+\underline{(s_2u_1s_2^{-1})u_1u_2u_1}+s_2u_1s_2u_1s_2+s_2u_1s_2^{-2}s_1s_2^{-1}+\\&&+
		s_2u_1s_2^{-2}s_1^{-2}s_2^{-2}u_1+
		s_2u_1s_2^{-2}s_1^{-2}s_2^{2}u_1+U'''
		\end{array}}$$
	However, by proposition \ref{oo}(i) we have that $s_2u_1s_2u_1s_2\subset U''\subset U'''$. It remains to prove that $A:=s_2u_1s_2^{-2}s_1s_2^{-1}+
	s_2u_1s_2^{-2}s_1^{-2}s_2^{-2}u_1+
	s_2u_1s_2^{-2}s_1^{-2}s_2^{2}u_1$ is a subset of $U'''$. We have: 
	$$\small{\begin{array}{lcl}
		A&=&s_2u_1s_2^{-2}s_1s_2^{-1}+s_2u_1s_2^{-2}s_1^{-2}s_2^{-2}u_1+
		s_2u_1s_2^{-2}s_1^{-2}s_2^{2}u_1\\
		&=&(s_2u_1s_2^{-1})s_2^{-1}s_1s_2^{-1}u_1+
		(s_2u_1s_2^{-1})s_2^{-1}s_1^{-2}s_2^{-2}u_1+
		(s_2u_1s_2^{-1})s_2^{-1}s_1^{-2}s_2^2u_1\\
		&=&s_1^{-1}u_2(s_1s_2^{-1}s_1^{-1})s_1^2s_2^{-1}u_1+s_1^{-1}u_2(s_1s_2^{-1}s_1^{-1})s_1^{-1}s_2^{-2}u_1+
		s_1^{-1}u_2s_1(s_2^{-1}s_1^{-2}s_2)s_2u_1\\
		&=&s_1^{-1}u_2s_1^{-1}(s_2s_1^2s_2^{-1})u_1+s_1^{-1}u_2s_1^{-1}(s_2s_1^{-1}s_2^{-1})s_2^{-1}u_1+s_1^{-1}u_2s_1^2s_2^{-1}(s_2^{-1}s_1^{-1}s_2)u_1\\
		&\subset&u_1(u_2u_1s_2^{-1}s_1s_2^{-1})u_1.
		\end{array}}$$
	By proposition \ref{l2} we have then $A\subset U'''$ and, hence, we proved that
	\begin{equation}s_2V\subset U'''\label{7}\end{equation}
	In order to finish the proof that $s_2U'\subset U''$,  it will be sufficient to prove that $u_1s_2s_1^{-2}s_2^{2}u_1$, 
	$u_1s_2^{2}s_1^{-2}s_2^{2}u_1$, $u_1s_2^{2}s_1^{2}s_2^{2}u_1$,
	$u_1s_2^{-2}s_1^{2}s_2^{-2}u_1$ and $u_1s_2^{2}s_1^{2}s_2^{-2}u_1$ are subsets of $U'''$. 
	\begin{itemize}[leftmargin=0.8cm]
		\item[C1.] We will prove that $s_2u_1s_2s_1^{-2}s_2^2u_1\subset U'''$. We expand $s_2^2$ as a linear combination of $s_2$, 1 $s_2^{-1}$, $s_2^{-2}
		$ and $s_2^{-3}$ and we have that 		
		$s_2u_1s_2s_1^{-2}s_2^2u_1\subset s_2u_1s_2s_1^{-2}s_2u_1+\underline{s_2u_1s_2u_1}+\underline{s_2u_1(s_2s_1^{-2}s_2^{-1})u_1}+
		s_2u_1(s_2s_1^{-2}s_2^{-1})s_2^{-1}u_1+
		s_2u_1s_2s_1^{-2}s_2^{-3}u_1
		\subset s_2u_1s_2s_1^{-2}s_2^{-3}u_1+s_2V+U'''$ and, hence, 
		by relation (\ref{7}) we have that $s_2u_1s_2s_1^{-2}s_2^2u_1\subset s_2u_1s_2s_1^{-2}s_2^{-3}u_1+U'''$. Therefore, it will be sufficient to prove that $s_2u_1s_2s_1^{-2}s_2^{-3}u_1\subset U'''$. We use the definition of $u_1$ and we have:
		
		$\small{\begin{array}[t]{lcl}
			s_2u_1s_2s_1^{-2}s_2^{-3}u_1
			&\subset&s_2(R_5+R_5s_1+R_5s_1^{-1}+R_5s_1^2+R_5s_1^{3})s_2s_1^{-2}s_2^{-3}u_1\\
			&\subset&\underline{s_2^2s_1^{-2}s_2^{-3}u_1}+\underline{(s_2s_1s_2)s_1^{-2}s_2^{-3}u_1}+
			s_2s_1^{-1}s_2s_1^{-2}s_2^{-3}u_1+\grv s_1^{-2}s_2^{-3}u_1+\\&&+s_2s_1^{3}s_2s_1^{-2}s_2^{-3}u_1\\
			&\subset&s_1^{-1}(s_1s_2s_1^{-1})s_2s_1^{-2}s_2^{-3}u_1+\underline{s_1^{-2}\grv s_2^{-3}u_1}+
			s_2s_1^{2}(s_1s_2s_1^{-1})s_1^{-1}s_2^{-3}u_1+U'''\\
			&\subset&s_1^{-1}s_2^{-1}(s_1s_2^2s_1^{-1})s_1^{-1}s_2^{-3}u_1+
			(s_2s_1^{2}s_2^{-1})s_1s_2s_1^{-1}s_2^{-3}u_1+U'''\\
			&\subset&s_1^{-1}s_2^{-2}s_1^2(s_2s_1^{-1}s_2^{-1})s_2^{-2}u_1+
			s_1^{-1}s_2^2s_1(s_1s_2s_1^{-1})s_2^{-3}u_1+
			U'''\\
			&\subset&s_1^{-1}s_2^{-3}(s_2s_1s_2^{-1})s_1s_2^{-2}u_1+s_1^{-1}s_2(s_2s_1s_2^{-1})s_1s_2^{-2}u_1+U'''\\
			&\subset&s_1^{-1}s_2^{-3}s_1^{-1}(s_2s_1^2s_2^{-1})s_2^{-1}u_1+s_1^{-1}s_2s_1^{-1}(s_2s_1^2s_2^{-1})s_2^{-1}u_1+U'''\\
			&\subset&
			s_1^{-1}s_2^{-3}s_1^{-2}s_2(s_2s_1s_2^{-1})u_1+s_1^{-1}s_2s_1^{-2}s_2(s_2s_1s_2^{-1})u_1+U'''\\
			&\subset&u_1(u_2u_1s_2s_1^{-1}s_2)u_1+U'''.
			\end{array}}$
		
		The result follows from proposition \ref{l2}(ii).
		\item[C2.] We will prove that $s_2u_1s_2^2s_1^{-2}s_2^2u_1\subset U'''$. For this purpose, we expand $u_1$ as $R_5+R_5s_1+R_5s_1^4+R_5s_1^2+R_5s_1^{-2}$ and we have that $s_2u_1s_2^2s_1^{-2}s_2^2u_1\subset 
		\underline{s_2^3s_1^{-2}s_2^2u_1}+\underline{
			(	s_2s_1s_2^2)s_1^{-2}s_2^2u_1}+s_2s_1^4s_2^2s_1^{-2}s_2^2u_1
		+s_2s_1^2s_2^2s_1^{-2}s_2^2u_1+s_2s_1^{-2}s_2^2s_1^{-2}s_2^2u_1$.
		By the definition of $U'''$ we have that $s_2s_1^{-2}s_2^2s_1^{-2}s_2^2u_1\subset U'''$. Therefore, it remains to  prove that $s_2s_1^4s_2^2s_1^{-2}s_2^2u_1
		+s_2s_1^2s_2^2s_1^{-2}s_2^2u_1\subset U'''$. We notice that

		$\small{\begin{array}{lcl}s_2s_1^4s_2^2s_1^{-2}s_2^2u_1
			+s_2s_1^2s_2^2s_1^{-2}s_2^2u_1&\subset&
			s_2s_1^3(s_1s_2^2s_1^{-1})s_1^{-1}s_2^2u_1
			+\grv (s_2s_1^{-2}s_2^{-1})s_2^3u_1\\&\subset&(s_2s_1^3s_2^{-1})s_1(s_1s_2s_1^{-1})s_2^2u_1+
			\grv s_1^{-2}(s_1s_2^{-2}s_1^{-1})s_1^2s_2^3u_1\\
			&\subset& s_1^{-1}s_2^3s_1^2s_2^{-1}s_1s_2^3u_1+\underline{s_1^{-2}\grv s_2^{-1}s_1^{-2}s_2s_1^2s_2^3u_1}
			\end{array}}$
		
		Therefore, we have to prove that the element $s_2^3s_1^2s_2^{-1}s_1s_2^3$ is inside $U'''$. For this purpose, we expand $s_2^3$ as a linear combination  of 
		$s_2^2$, $s_2$, 1 $s_2^{-1}$ and $s_2^{-2}$ and we have:

		$\small{\begin{array}[t]{lcl}
			s_2^3s_1^2s_2^{-1}s_1s_2^3
			&\in&
			R_5s_2^3s_1^{3}(s_1^{-1}s_2^{-1}s_1)s_2^2+\underline{R_5s_2^3s_1^2(s_2^{-1}s_1s_2)}+\underline
			{R_5s_2^3s_1^2s_2^{-1}}+R_5s_2^3s_1^2s_2^{-1}s_1s_2^{-1}+\\&&+
			R_5s_1^{-1}(s_1s_2^2s_1^{-1})s_1(s_2s_1^2s_2^{-1})s_1s_2^{-2}\\
			&\in&u_2u_1s_2s_1^{-1}s_2+u_2u_1s_2^{-1}s_1s_2^{-1}u_1+u_1s_2^{-1}s_1^2s_2^3s_1^2s_2^{-2}+U'''.
			
			\end{array}}$
		
		However, by proposition 	\ref{l2} we have that 	$u_2u_1s_2s_1^{-1}s_2$ and $u_2u_1s_2^{-1}s_1s_2^{-1}$ are subsets of $U'''$. Therefore, we only need to prove that the element $s_2^{-1}s_1^2s_2^3s_1^2s_2^{-2}$ is inside $U'''$. We expand $s_2^3$ as a linear combination  of 
		$s_2^2$, $s_2$, 1 $s_2^{-1}$ and $s_2^{-2}$ and we have that $s_2^{-1}s_1^2s_2^3s_1^2s_2^{-2}\in \grF(s_2V)+\underline{s_2^{-1}(s_1^2s_2s_1)s_1s_2^{-2}}+      \grF(s_2u_1s_2s_1^{-2}s_2^2)+R_5s_2^{-1}s_1^{2}s_2^{-2}s_1^{2}s_2^{-2}$. However, by the definition of $U'''$ we have that $s_2^{-1}s_1^{2}s_2^{-2}s_1^{2}s_2^{-2}\in U'''$.
		Moreover, by relation (\ref{7}) and by the previous case (case C1) we have that $\grF(s_2V)+\grF(s_2u_1s_2s_1^{-2}s_2^2)\subset \grF(U''')\stackrel{\ref{r1}}\subset U'''.$ 
		
		\item[C3.] We will prove that $s_2u_1s_2^2s_1^{2}s_2^2\subset U'''$. For this purpose, we expand $u_1$ as $R_5+R_5s_1+R_5s_1^{-1}+
		R_5s_1^2+R_5s_1^3$ and we have $s_2u_1s_2^2s_1^{2}s_2^2\subset \underline{s_2^3s_1^2s_2^2u_1}+\underline{(s_2s_1s_2^2)s_1^2s_2^2u_1}+ s_2s_1^{-1}s_2^2s_1^{2}s_2^2u_1+s_2s_1^{2}s_2^2s_1^2s_2^2u_1+s_2s_1^3s_2^2s_1^{2}s_2^2$. However, be lemma \ref{cc} we have that $s_2s_1^3s_2^2s_1^{2}s_2^2\subset u_1\grv^3+U''\subset U'''$. Therefore, it remains to  prove that $s_2s_1^{-1}s_2^2s_1^{2}s_2^2u_1+s_2s_1^{2}s_2^2s_1^2s_2^2u_1\subset U'''$. We have: 
		
		$\small{\begin{array}[t]{lcl}s_2s_1^{-1}s_2^2s_1^{2}s_2^2u_1+s_2s_1^{2}s_2^2s_1^2s_2^2u_1&=&
			s_2^2(s_2^{-1}s_1^{-1}s_2)s_2s_1^{2}s_2^2u_1+s_1^{-1}s_1\grv^2s_2\\
			&=&s_2^2s_1(s_2^{-1}s_1^{-1}s_2)s_1^{2}s_2^2u_1+
			s_1^{-1}\grv^2s_1s_2\\
			&=&s_2^2s_1^2(s_2^{-1}s_1s_2)s_2u_1+s_1^{-1}s_2s_1^2s_2^2s_1^2(s_2s_1s_2)\\
			&\subset& u_2u_1s_2s_1^{-1}s_2u_1+u_1s_2s_1^2s_2^2s_1^3s_2u_1.
			
			\end{array}}$.
		
		By lemma \ref{l2}(ii) it will be sufficient to prove that $s_2s_1^2s_2^2s_1^3s_2\in U'''$. We expand $s_2^3$ as a linear combination  of 
		$s_2^2$, $s_2$, 1 $s_2^{-1}$ and $s_2^{-2}$ and we have: 
		
		$\small{\begin{array}[t]{lcl}
			s_2s_1^2s_2^2s_1^3s_2&\in& R_5\grv^2+\underline{R_5s_2s_1^2(s_2^2s_1s_2)}+\underline{R_5s_2s_1^2s_2^2}+
			R_5s_2s_1(s_1s_2^2s_1^{-1})s_2+R_5s_2s_1^2s_2^2s_1^{-2}s_2\\
			
			&\in&\underline{R_5(s_2s_1s_2^{-1})s_1^2s_2^2}+R_5s_1^{-1}(s_1s_2s_1^{2})s_2^2s_1^{-2}s_2+U'''\\
			&\in&u_1s_2^2(s_1s_2^3s_1^{-1})s_1^{-1}s_2+U'''\\
			&\in& u_1u_2u_1s_2s_1^{-1}s_2u_1+U'''.
			\end{array}}$
		
		The result follows from proposition \ref{l2}(ii).
		\item[C4.] We will prove that $s_2u_1s_2^{-2}s_1^{2}s_2^{-2}u_1\subset U'''$. Since $s_2u_1s_2^{-2}s_1^{2}s_2^{-2}u_1=(s_2u_1s_2^{-1})s_2^{-1}s_1^{2}s_2^{-2}u_1=s_1^{-1}u_2s_1s_2^{-1}s_1^2s_2^{-2}u_1$, it will be sufficient to prove that $u_2s_1s_2^{-1}s_1^2s_2^{-2}\subset U'''$. We expand $u_2$ as $R_5+R_5s_2+R_5s_2^{-1}+R_5s_2^2+R_5s_2^3$ and we have:  $u_2s_1s_2^{-1}s_1^2s_2^{-2}\subset \underline{R_5s_1s_2^{-1}s_1^2s_2^{-2}}+\underline{R_5(s_2s_1s_2^{-1})s_1^2s_2^{-2}}+\grF(u_1s_2u_1s_2s_1^{-2}s_2^{2})+R_5s_2^2s_1s_2^{-1}s_1^2s_2^{-2}+
		R_5s_2^3s_1s_2^{-1}s_1^2s_2^{-2}$.
		By the first case (case C1) we have that $\grF(u_1s_2u_1s_2s_1^{-2}s_2^{2})
		\subset u_1\grF(U''')u_1\stackrel{\ref{r1}}{\subset} U'''$. It remains to prove that the elements $s_2^2s_1s_2^{-1}s_1^2s_2^{-2}$ and $s_2^3s_1s_2^{-1}s_1^2s_2^{-2}$ are inside $U'''$. We have:
		$s_2^2s_1s_2^{-1}s_1^2s_2^{-2}=s_2(s_2s_1s_2^{-1})s_1^2s_2^{-2}=s_2s_1^{-1}(s_2s_1^3s_2^{-1})s_2^{-1}=s_2s_1^{-2}s_2^2(s_2s_1s_2^{-1})=s_1^{-1}(s_1s_2s_1^{-1})s_1^{-1}s_2^2s_1^{-1}s_2s_1=s_1^{-1}s_2^{-1}(s_1s_2s_1^{-1})s_2^2s_1^{-1}s_2s_1=\underline{s_1^{-1}s_2^{-2}(s_1s_2^3s_1^{-1})s_2s_1}$. Moreover,
		$s_2^3s_1s_2^{-1}s_1^2s_2^{-2}=s_2^2(s_2s_1s_2^{-1})s_1(s_1s_2^{-2}s_1^{-1})s_1=s_2^2s_1^{-2}(s_1s_2s_1^2)s_2^{-1}s_1^{-2}s_2s_1\in s_2^2s_1^{-2}s_2^2s_1^{-1}s_2u_1$.
		We expand $s_1^{-2}$ as a linear combination  of 
		$s_1^{-1}$, $1$, $s_1$, $s_2^{2}$ and $s_2^{3}$ and we have:
		
		$\small{\begin{array}[t]{lcl}s_2^2s_1^{-2}s_2^2s_1^{-1}s_2&\in& R_5s_2^2s_1^{-1}s_2^2s_1^{-1}s_2+\underline{R_5s_2^4s_1^{-1}s_2}+
			\underline{R_5s_2(s_2s_1s_2^2)s_1^{-1}s_2}+
			R_5s_2^2s_1^{2}s_2^2s_1^{-1}s_2+\\&&+
			R_5s_2^2s_1^{3}s_2^2s_1^{-1}s_2\\
			&	\in& R_5s_2^3(s_2^{-1}s_1^{-1}s_2)s_2s_1^{-1}s_2+R_5
			s_2^2s_1(s_1s_2^{2}s_1^{-1})s_2
			+R_5s_2^2s_1^2(s_1s_2^2s_1^{-1})s_2+U'''\\
			&\in& R_5s_2^3s_1(s_2^{-1}s_1^{-1}s_2)s_1^{-1}s_2+
			R_5s_2(s_2s_1s_2^{-1})s_1^{2}s_2^{2}
			+R_5s_2(s_2s_1^2s_2^{-1})s_1^2s_2^2+U'''\\
			&\in&\underline{R_5s_2^3s_1^2(s_2^{-1}s_1^{-2}s_2)}+R_5s_2s_1^{-1}s_2s_1^{3}s_2^{2}+
			R_5s_2s_1^{-1}s_2^2s_1^3s_2^2+U'''
			\end{array}}$
		
		Therefore, it remains to prove that $B:=R_5s_2s_1^{-1}s_2s_1^{3}s_2^{2}+
		R_5s_2s_1^{-1}s_2^2s_1^3s_2^2\subset U'''$. We expand  $s_1^3$ as a linear combination  of 
		$s_1^2$, $s_1$, 1 $s_1^{-1}$ and $s_1^{-2}$ and we have that $B\subset R_5s_2s_1^{-1}s_2(R_5s_1^2+R_5s_1+R_5+R_5s_1^{-1}+R_5s_1^{-2})s_2^{2}+
		R_5s_2s_1^{-1}s_2^2(R_5s_1^2+R_5s_1+R_5+R_5s_1^{-1}+R_5s_1^{-2})s_2^2$. 
		By cases C1, C2 and C3 we have: 
		
		$\small{\begin{array}[t]{lcl}
			B&\subset&
			R_5s_2s_1^{-1}\grv s_2+\underline{R_5s_2s_1^{-1}(s_2s_1s_2^2)}+\underline{R_5s_2s_1^{-1}s_2^3u_1}+
			R_5s_2s_1^{-1}s_2s_1^{-1}s_2^2+\underline{R_5s_2s_1^{-1}(s_2^2s_1s_2)s_2}+\\&&+\underline{R_5s_2s_1^{-1}s_2^4}+
			R_5s_2s_1^{-1}s_2^2s_1^{-1}s_2^2+U'''\\
			&\subset&R_5s_2\grv s_1^{-1}s_2+
			R_5s_2s_1^{-1}s_2s_1^{-1}s_2^2+
			R_5s_2s_1^{-1}s_2^2s_1^{-1}s_2^2+U'''\\
			&\subset&u_2u_1s_2s_1^{-1}s_2+R_5s_2^2(s_2^{-1}s_1^{-1}s_2)s_1^{-1}s_2^2+
			R_5s_2s_1^{-2}(s_1s_2^2s_1^{-1})s_2^2+U'''\\
			&\stackrel{\ref{l2}}{\subset}&R_5s_2^2s_1(s_2^{-1}s_1^{-2}s_2)s_2+
			R_5(s_2s_1^{-2}s_2^{-1})s_1^2s_2^3+U'''\\
			&\subset&R_5s_2^2s_1^{2}s_2^{-1}(s_2^{-1}s_1^{-1}s_2)+U''' \\&\subset& u_1u_2u_1s_2^{-1}s_1s_2^{-1}+U'''
			.
			\end{array}}$
		
		The result follows from proposition \ref{l2}(ii).
		\item[C5.] We 
		will prove that $s_2u_1s_2^2s_1^2s_2^{-2}u_1\subset U'''$. For this purpose, we use straight-forward calculations and we have
		$s_2u_1s_2^2s_1^2s_2^{-2}=(s_2u_1s_2^{-1})s_2^2(s_2s_1^2s_2^{-1})s_2^{-1}
		=s_1^{-1}u_2(s_1s_2^2s_1^{-1})s_2(s_2s_1s_2^{-1})=
		s_1^{-1}u_2s_1(s_1s_2^2s_1^{-1})s_2s_1=
		s_1^{-1}u_2(s_2s_1s_2^{-1})s_1^2s_2^2s_1=s_1^{-2}(s_1u_2s_1^{-1})s_2s_1^3s_2^2s_1=s_1^{-2}s_2^{-1}u_1s_2^2s_1^3s_2^2s_1$, meaning that $s_2u_1s_2^2s_1^2s_2^{-2}u_1
		\subset u_1s_2^{-1}u_1s_2^2s_1^3s_2^2u_1$.
		Hence, we have to prove that $s_2^{-1}u_1s_2^2s_1^3s_2^2\subset U'''$.
		For this purpose, we expand  $s_1^3$ as a linear combination  of 
		$s_1^2$, $s_1$, 1 $s_1^{-1}$ and $s_1^{-2}$ and we have that $s_2^{-1}u_1s_2^2s_1^3s_2^2\subset \grF(s_2V+s_2u_1s_2^{-2}s_1^2s_2^{-2})+s_2^{-1}u_1s_2^2s_1s_2^2+
		s_2^{-1}u_1s_2^2s_1^{-1}s_2^2$.
		By relation (\ref{7}) and case C4 we have that $\grF(s_2V+s_2u_1s_2^{-2}s_1^2s_2^{-2})\subset \grF(U''')\stackrel{\ref{r1}}{\subset} U'''$. Moreover, 
		$s_2^{-1}u_1s_2^2s_1s_2^2=s_2^{-1}u_1(s_2^2s_1s_2)s_2=s_2^{-1}u_1\grv=\underline{s_2^{-1}\grv u_1}.$ It remains to prove that $s_2^{-1}u_1s_2^2s_1^{-1}s_2^2\subset U'''$. We have:
		$s_2^{-1}u_1s_2^2s_1^{-1}s_2^2=(s_2^{-1}u_1s_2)s_2s_1^{-1}s_2^2=s_1u_2(s_1^{-1}s_2s_1)s_1^{-2}s_2^2=s_1u_2s_1(s_2^{-1}s_1^{-2}s_2)s_2\subset u_1u_2u_1s_2^{-1}s_1s_2^{-1}u_1$.
		The result follows from proposition \ref{l2}(i).
		\qedhere
	\end{itemize}
\end{proof}
From now on we will double-underline the elements of the form $u_1s_2^{\pm}u_1u_2u_1u_2u_1$ and  we will use the fact that they are elements of $U'''$ (proposition \ref{p2}) without mentioning it.

We can now prove the following lemma that helps us to ``replace''  inside the definition of $U''''$ the element $\grv^4$ by the element $s_2^{-2}s_1^2s_2^2s_1^3s_2^2$ modulo $U'''$.  
\begin{lem}  $s_2^{-2}s_1^2s_2^2s_1^3s_2^2\in u_1\grv^3+u_1^{\times}\grv^4+u_1s_2u_1u_2u_1u_2u_1\subset U''''.$
	\label{ll}
\end{lem}
\begin{proof} In this proof we will double-underline only the elements of the form $u_1s_2u_1u_2u_1u_2u_1$ (and not of the form $u_1s_2^{-1}u_1u_2u_1u_2u_1$ ). The fact that $u_1\grv^3+u_1^{\times}\grv^4+u_1s_2u_1u_2u_1u_2u_1$ is a subset of $ U''''$ follows from the definition of $U''''$ and proposition \ref{p2}. As a result, we restrict ourselves to proving that $s_2^{-2}s_1^2s_2^2s_1^3s_2^2\in u_1\grv^3+u_1^{\times}\grv^4+u_1s_2u_1u_2u_1u_2u_1$.
	We first notice that
	$$\small{\begin{array}{lcl}
		s_2^{-2}s_1^2s_2^2s_1^3s_2^2&=&s_1(s_1^{-1}s_2^{-2}s_1)s_2^{-2}(s_2^2s_1s_2)s_2s_1^2(s_1s_2^2s_1^{-1})s_1^{-1}s_1^2\\
		&=&
		s_1s_2s_1^{-2}s_2^{-3}s_1\grv s_1^2s_2^{-1}s_1(s_1s_2s_1^{-1})s_1^2\\
		&=&
		s_1s_2s_1^{-2}s_2^{-3}s_1^3(s_2s_1^3s_2^{-1})s_1s_2s_1^2\\
		&=&
		s_1s_2s_1^{-2}s_2^{-3}s_1^2s_2^3s_1^2s_2s_1^2\\
		&\in& u_1s_2s_1^{-2}s_2^{-3}s_1^2s_2^3s_1^2s_2u_1.
		\end{array}}$$
	We expand $s_2^{-3}$ as a linear combination of $s_2^{-2}$, $s_2^{-1}$, 1, $s_2$ and $s_2^2$, where the coefficient of $s_2^2$ is invertible, and we have:
	
	$$\small{\begin{array}[t]{lcl}
		s_2s_1^{-2}s_2^{-3}s_1^2s_2^3s_1^2s_2&\in& R_5s_2s_1^{-2}s_2^{-2}s_1^2s_2^3s_1^2s_2+R_5s_2s_1^{-2}s_2^{-1}s_1^2s_2^3s_1^2s_2+\underline{\underline{R_5s_2s_2^3s_1^2s_2}}+\\&&+R_5s_2s_1^{-2}s_2s_1^2s_2^3s_1^2s_2+u_1^{\times}s_2s_1^{-2}s_2^{2}s_1^2s_2^3s_1^2s_2u_1^{\times}.
		\end{array}}$$
	However, we notice that $s_2s_1^{-2}s_2^{-2}s_1^2s_2^3s_1^2s_2=s_2s_1^{-1}(s_1^{-1}s_2^{-2}s_1)s_1s_2^3s_1^2s_2=s_2s_1^{-1}s_2^2\grv^{-1}s_1s_2^3s_1^2s_2=s_2s_1^{-1}s_2^2s_1\grv^{-1}s_2^3s_1^2s_2=s_2s_1^{-1}s_2^2s_1(s_2^{-1}s_1^{-2}s_2)\grv=s_2s_1^{-1}s_2^2s_1^2s_2^{-2}\grv s_1^{-1} =\underline{\underline{s_2s_1^{-1}s_2^2s_1^2(s_2^{-1}s_1^2s_2)}}s_1^{-1}$. Moreover,  $s_2s_1^{-2}s_2^{-1}s_1^2s_2^3s_1^2s_2=s_2s_1^{-2}(s_2^{-1}s_1^2s_2)s_2^2s_1^2s_2=s_2s_1^{-1}s_2^2(s_1^{-1}s_2^2s_1)s_2(s_2^{-1}s_1s_2)=\underline{\underline{s_2s_1^{-1}s_2^3(s_1^3s_2s_1)s_1^{-2}}}.$ We also have $s_2s_1^{-2}s_2s_1^2s_2^3s_1^2s_2=s_2s_1^{-3}(s_1s_2s_1^2)s_2^3s_1^2s_2=s_2s_1^{-3}s_2(s_2s_1s_2^4)s_1^2s_2\in s_2s_1^{-3}(s_2u_1s_2u_1s_2u_1).$ However, by lemma \ref{oo}(i) we have that $s_2s_1^{-3}(s_2u_1s_2u_1s_2u_1)\subset s_2s_1^{-3}(\grv^2u_1+u_1u_2u_1u_2u_1)\subset s_2\grv^2 u_1+ \underline{\underline{u_1s_2u_1u_2u_1u_2u_1}}$. By lemma \ref{oo}(ii) we  also have $s_2\grv^2 u_1\subset \underline{\underline{u_1s_2u_1u_2u_1u_2u_1}}$.
	
	It remains to prove that $s_2s_1^{-2}s_2^{2}s_1^2s_2^3s_1^2s_2\in u_1\grv^3+u_1^{\times}\grv^4+u_1s_2u_1u_2u_1u_2u_1$. We have:
	$$\small{\begin{array}[t]{lcl}
		s_2s_1^{-2}s_2^{2}s_1^2s_2^3s_1^2s_2&=& 
		s_2(-de^{-1}s_1^{-1}-ce^{-1}-e^{-1}bs_1-e^{-1}as_1^2+e^{-1}s_1^3)s_2^{2}s_1^2s_2^3s_1^2s_2\\

		&\in&
		R_5s_2s_1^{-1}s_2^2s_1^2s_2^3s_1^2s_2+R_5s_2^3s_1^2s_2^3s_1^2s_2
		+\underline{\underline{R_5(s_2s_1s_2^2)s_1^2s_2^3s_1^2s_2u_1}}+R_5s_2s_1^2s_2^2s_1^2s_2^3s_1^2s_2+\\&&+
		u_1^{\times}s_2s_1^3s_2^2s_1^2s_2^2\grv.
		\end{array}}$$
	We first notice that we have  $s_2s_1^{-1}s_2^2s_1^2s_2^3s_1^2s_2=s_2(s_1^{-1}s_2^2s_1)s_1s_2^3s_1^2s_2=s_2^2s_1^2(s_2^{-1}s_1s_2)s_2^2s_1^2s_2=s_1(s_1^{-1}s_2^2s_1)s_1^2s_2(s_1^{-1}s_2^2s_1)(s_1s_2s_1)s_1^{-1}=s_1s_2s_1^2(s_2^{-1}s_1^2s_2)s_2s_1^3s_2s_1^{-1}=\underline{\underline{s_1s_2s_1^3s_2^3(s_2^{-1}s_1^{-1}s_2)(s_1^3s_2s_1)s_1^{-2}}}.$
	Moreover, we have that $s_2^3s_1^2s_2^3s_1^2s_2=s_1(s_1^{-1}s_2^3s_1)s_1s_2^3s_1^2s_2=s_1s_2s_1^3(s_2^{-1}s_1s_2)s_2^2s_1(s_1s_2s_1)s_1^{-1}=\underline{\underline{s_1s_2s_1^4s_2(s_1^{-1}s_2^2s_1)s_2s_1s_2s_1^{-1}}}.$ Using analogous calculations, $s_2s_1^2s_2^2s_1^2s_2^3s_1^2s_2=s_1^{-1}(s_1s_2s_1^2)s_2^2s_1^2s_2^3s_1^2s_2=
	s_1^{-1}s_2(s_2s_1s_2^3)s_1^2s_2^3s_1^2s_2=s_1^{-1}s_2s_1^2(s_1s_2s_1^3)s_2^3s_1^2s_2=s_1^{-1}s_2s_1^2s_2^2(s_2s_1s_2^4)s_1^2s_2\in u_1\grv (s_2u_1s_2u_1s_2u_1).$ However, by lemma \ref{oo}(i) we have that $u_1\grv (s_2u_1s_2u_1s_2u_1)\subset u_1\grv(\grv^2u_1+u_1u_2u_1u_2u_1)\subset u_1\grv^3+\underline{\underline{u_1\grv u_2u_1u_2u_1}}.$ 
	
	In order to finish the proof, it remains to prove that $s_2s_1^3s_2^2s_1^2s_2^2\grv\in u_1^{\times}\grv^4+u_1s_2u_1u_2u_1u_2u_1$. 
	We use lemma \ref{cc} and we have: \\\\
	$\small{\begin{array}[t]{lcl}
		s_2s_1^3s_2^2s_1^2s_2^2\grv&\in&
		u_1^{\times}(u_1s_2u_1s_2s_1^3s_2+u_1s_2^2s_1^3s_2s_1^{-1}s_2+
		u_1u_2u_1u_2u_1+u_1^{\times}\grv^3)\grv\\
		&\in&
		u_1s_2u_1s_2s_1^4(s_1^{-1}s_2^2s_1)s_1s_2+u_1s_2^2s_1^3s_2(s_1^{-1}s_2s_1)s_1^{-1}\grv +u_1^{\times}\grv^4\\

		&\in&
		u_1s_2u_1s_2s_1^4s_2s_1^2(s_2^{-1}s_1s_2)+u_1s_2^2s_1^3s_2^2s_1s_2^{-1}\grv +u_1^{\times}\grv^4\\

		&\in&
		u_1s_2u_1s_2s_1^4s_2s_1^3s_2u_1+u_1s_2^2s_1^3s_2^2s_1^3s_2+u_1^{\times}\grv^4\\
		
		&\in&
		u_1s_2u_1(s_2u_1s_2u_1s_2u_1)+u_1(s_1^{-1}s_2^2s_1)s_1^2s_2^2s_1^3s_2+
		u_1^{\times}\grv^4\\
		&\stackrel{\ref{oo}(i)}{\in}& u_1s_2u_1(\grv^2u_1+u_1u_2u_1u_2u_1)+
		u_1s_2s_1^2(s_2^{-1}s_1^2s_2)s_2s_1^3s_2+u_1^{\times}\grv^4\\

		&\in&
		u_1s_2\grv^2u_1+u_1s_2u_1u_2u_1u_2u_1+ u_1s_2s_1^3s_2^2(s_1^{-1}s_2s_1)s_1^2s_2+u_1^{\times}\grv^4\\
		&\stackrel{\ref{oo}(ii)}{\in}&u_1s_2u_1u_2u_1u_2u_1+\underline{\underline{u_1s_2s_1^3s_2^3s_1(s_2^{-1}s_1^2s_2)}}+u_1^{\times}\grv^4
		.
		\end{array}}$\\
	
\end{proof}
\begin{prop}\mbox{}  
	\vspace*{-\parsep}
	\vspace*{-\baselineskip}\\ 
	\begin{itemize}[leftmargin=0.8cm]
		\item [(i)]$s_2u_1u_2u_1s_2s_1^{-1}s_2\subset U''''.$
		\item[(ii)]$s_2u_1u_2u_1s_2^{-1}s_1s_2^{-1}\subset U''''.$
		\item[(iii)]$s_2u_1u_2u_1u_2\grv \subset U''''$.
	\end{itemize}
	\label{xx}
\end{prop}
\begin{proof}\mbox{} 
	\vspace*{-\parsep}
	\vspace*{-\baselineskip}\\
	\begin{itemize}[leftmargin=0.8cm]

		\item[(i)] By proposition \ref{l2}(ii) we have  $	s_2u_1u_2u_1s_2s_1^{-1}s_2\subset s_2u_1(u_1\grv^{2}+R_5s_2^2s_1^{-2}s_2s_1^{-1}s_2+u_1u_2u_1u_2u_1)$ and, hence, by lemma \ref{oo}(ii) we have  $s_2u_1u_2u_1s_2s_1^{-1}s_2\subset s_2u_1s_2^2s_1^{-2}s_2s_1^{-1}s_2+\underline{\underline{s_2u_1u_2u_1u_2u_1}}+U''''.$ As a result, we must prove that $s_2u_1s_2^2s_1^{-2}s_2s_1^{-1}s_2\subset U''''$. For this purpose, we expand $u_1$ as $R_5+R_5s_1+R_5s_1^{-1}+R_5s_1^2+R_5s_1^3$ and we have:
		
		$\begin{array}{lcl}		s_2u_1s_2^2s_1^{-2}s_2s_1^{-1}s_2&\subset& u_2u_1s_2s_1^{-1}s_2+\underline{\underline{
				R_5(s_2s_1s_2^2)s_1^{-2}s_2s_1^{-1}s_2}}+R_5s_2s_1^{-1}s_2^2s_1^{-2}s_2s_1^{-1}s_2
		+\\&&+R_5s_2s_1^2s_2^2s_1^{-2}s_2s_1^{-1}s_2+R_5s_2s_1^3s_2^2s_1^{-2}s_2s_1^{-1}s_2.
		\end{array}$
		
		By proposition \ref{l2}(ii) we have that $u_2u_1s_2s_1^{-1}s_2\subset U''''$. Moreover, $s_2s_1^2s_2^2s_1^{-2}s_2s_1^{-1}s_2=s_1^{-1}(s_1s_2s_1^2)s_2^2s_1^{-3}(s_1s_2s_1^{-1})s_2=s_1^{-1}s_2(s_2s_1s_2^3)s_1^{-3}s_2^{-1}s_1s_2^2=
		\underline{\underline{s_1^{-1}s_2s_1^3(s_2s_1^{-2}s_2^{-1})s_1s_2^2}}$. We also notice that $s_2s_1^3s_2^2s_1^{-2}s_2s_1^{-1}s_2=s_1^{-1}(s_1s_2s_1^3)s_2^2s_1^{-3}(s_1s_2s_1^{-1})s_2=s_1^{-1}s_2^3(s_1s_2^3s_1^{-1})s_1^{-2}s_2^{-1}s_1s_2^2=s_1^{-1}s_2^2s_1^3(s_2s_1^{-2}s_2^{-1})s_1s_2^2=s_1^{-2}(s_1s_2s_1^{-1})s_1(s_2s_1^2s_2^{-1})s_2^{-1}s_1^2s_2^2=s_1^{-2}s_2^{-1}(s_1s_2^3s_1^{-1})s_2^{-1}(s_2s_1^2s_2^{-1})s_1^2s_2^2
		\in u_1s_2^{-2}s_1^2s_2^2s_1^3s_2^2\stackrel{\ref{ll}}{\subset} U''''.$
		
		It remains to prove that the element $s_2s_1^{-1}s_2^2s_1^{-2}s_2s_1^{-1}s_2$ is inside $U''''$.
		We expand $s_2^2$ as a linear combination of $s_2$, 1, $s_2^{-1}$,  $s_2^{-2}$ and  $s_2^{-3}$ and we have
		
		$\small{\begin{array}[t]{lcl}
			s_2s_1^{-1}s_2^2s_1^{-2}s_2s_1^{-1}s_2&\in &	s_2s_1^{-1}(R_5s_2+R_5+R_5s_2^{-1}+R_5s_2^{-2}+R_5s_2^{-3})s_1^{-2}s_2s_1^{-1}s_2\\

			&\in&R_5s_2s_1^{-2}(s_1s_2s_1^{-1})s_1^{-1}(s_2s_1^{-1}s_2^{-1})s_2^2+
			\underline{\underline{R_5s_2s_1^{-3}s_2s_1^{-1}s_2}}+\\&&+\underline{\underline{
					R_5(s_2s_1^{-1}s_2^{-1})s_1^{-2}s_2s_1^{-1}s_2}}+\underline{\underline{
					R_5(s_2s_1^{-1}s_2^{-1})(s_2^{-1}s_1^{-2}s_2)s_1^{-1}s_2}}+\\&&+R_5(s_2s_1^{-1}s_2^{-1})s_2^{-1}(s_2^{-1}s_1^{-2}s_2)s_1^{-1}s_2\\

			&\in&
			\underline{\underline{R_5s_2s_1^{-2}s_2^{-1}s_1(s_2s_1^{-2}s_2^{-1})s_1s_2^2}}+
			R_5s_1^{-1}s_2^{-1}s_1^2(s_1^{-1}s_2^{-1}s_1)s_2^{-2}s_1^{-2}s_2
			+U''''\\
			&\in&R_5s_1^{-1}s_2^{-1}s_1^2s_2(s_1^{-1}s_2^{-3}s_1)s_1^{-3}s_2+
			U''''\\

			&\in&\underline{\underline{R_5s_1^{-1}s_2^{-1}s_1^2s_2^2s_1^{-3}(s_2^{-1}s_1^{-3}s_2)}}+U''''.
			\end{array}}$
		\item[(ii)]By proposition \ref{l2} (i)  we have that $s_2u_1(u_2u_1s_2^{-1}s_1s_2^{-1})\subset s_2u_1(u_1\grv^{-2}+R_5s_2^{-2}s_1^{2}s_2^{-1}s_1s_2^{-1}+u_1u_2u_1u_2u_1)
		\subset\underline{\underline{s_2\grv^{-2}u_1}}+
		s_2u_1s_2^2s_1^{-2}s_2s_1^{-1}s_2+\underline{\underline{s_2u_1u_2u_1u_2u_1}}$. Therefore, it remains to prove that  $s_2u_1s_2^2s_1^{-2}s_2s_1^{-1}s_2 \subset U''''$. We expand $u_1$ as $R_5+R_5s_1+R_5s_1^{-1}+R_5s_1^2+R_5s_1^{3}$ and we have:
		
		$\small{\begin{array}[t]{lcl}

			s_2u_1s_2^2s_1^{-2}s_2s_1^{-1}s_2	&\subset&
			\underline{\underline{R_5s_2^{-1}s_1^2s_2^{-1}s_1s_2^{-1}}}+\underline{\underline{R_5s_1^{-2}(s_1^2s_2s_1)s_2^{-2}s_1^2s_2^{-1}s_1s_2^{-1}}}+\\&&
			+R_5(s_2s_1^{-1}s_2^{-1})(s_2^{-1}s_1^2s_2)s_2^{-2}s_1s_2^{-1}+R_5(s_2s_1^2s_2^{-1})s_2^{-2}(s_2s_1^2s_2^{-1})s_1s_2^{-1}+\\&&+R_5(s_2s_1^3s_2^{-1})s_2^{-1}s_1^2s_2^{-1}(s_1s_2s_1^{-1})s_1\\
			&\subset&\underline{\underline{R_5s_1^{-1}s_2^{-1}s_1^2s_2^2(s_1^{-1}s_2^{-2}s_1)s_2^{-1}}}+
			\underline{\underline{R_5s_1^{-1}s_2^2(s_1s_2^{-2}s_1^{-1})s_2(s_2s_1s_2^{-1})}}+\\
			&&+R_5s_1^{-1}s_2^2(s_2s_1s_2^{-1})s_1(s_1s_2^{-2}s_1^{-1})s_2s_1+U''''\\
			&\subset&R_5s_1^{-1}s_2^2s_1^{-1}(s_2s_1^{2}s_2^{-1})s_1^{-2}s_2^2s_1+U''''\\
			&\subset&R_5s_1^{-1}s_2^2s_1^{-3}(s_1s_2^2s_1^{-1})s_2^2s_1+U''''\\
			&\subset& \underline{\underline{R_5s_1^{-1}s_2(s_2s_1^{-3}s_2^{-1})s_1^2s_2^3s_1}}+U''''.
			\end{array}}$
		\item[(iii)] We notice that $s_2u_1u_2u_1u_2\grv=s_2u_1u_2u_1u_2s_1^2s_2
		s_2u_1u_2u_1u_2(s_2^{-1}s_1^2s_2)
		\subset s_2u_1u_2u_1u_2s_1s_2^2u_1$. We expand $\bold{u_2}$ as $R_5+R_5s_2+R_5s_2^{-1}+R_5s_2^{2}+R_5s_2^{-2}$ and we have:

		$\small{\begin{array}[t]{lcl}

			s_2u_1u_2u_1\bold{u_2}s_1s_2^2u_1&\subset&\underline{\underline{s_2u_1u_2u_1s_2^2u_1}}+\underline{\underline{s_2u_1u_2u_1(s_2s_1s_2^2)u_1}}+
			s_2u_1u_2u_1(s_2^{-1}s_1s_2)s_2u_1+\\&&+
			s_2u_1u_2u_1s_2(s_2s_1s_2^2)u_1+s_2u_1u_2u_1s_2^{-1}(s_2^{-1}s_1s_2)s_2u_1\\
			&\subset&s_2u_1(u_2u_1s_2s_1^{-1}s_2u_1)+\underline{\underline{s_2u_1u_2u_1\grv u_1}}+s_2u_1u_2u_1(s_2^{-1}s_1s_2)s_1^{-1}s_2u_1+U''''\\
			&\stackrel{(i)}{\subset}&s_2u_1u_2u_1s_2s_1^{-2}s_2u_1+U''''.
			\end{array}}$
		
		It remains to prove that $s_2u_1u_2u_1s_2s_1^{-2}s_2u_1\subset U''''$. For this purpose, we expand $s_1^{-2}$ as a linear combination of $s_1^{-1}$, 1, $s_1$, $s_1^2$, $s_1^3$ and we have: $s_2u_1u_2u_1s_2s_1^{-2}s_2u_1	\subset s_2u_1u_2u_1s_2s_1^{-1}s_2u_1+\underline{\underline{s_2u_1u_2u_1s_2^2u_1}}+
		\underline{\underline{s_2u_1u_2u_1(s_2s_1s_2)u_1}}+s_2u_1u_2u_1\grv u_1+s_2u_1u_2u_1(s_1^{-1}s_2s_1)s_1^2s_2u_1\stackrel{(i)}{\subset}
		\underline{\underline{s_2u_1u_2\grv u_1}}+s_2u_1u_2u_1s_2s_1(s_2^{-1}s_1^2s_2)u_1+U''''
		\subset s_2u_1u_2u_1\grv s_2u_1+U''''
		\subset s_2u_1u_2\grv u_1s_2u_1+U''''$.
		
		However, 
		$\small{\begin{array}[t]{lcl}
			s_2u_1u_2\grv u_1s_2u_1&\subset&s_2u_1u_2s_1^2s_2(R_5+R_5s_1+
			R_5s_1^{-1}+R_5s_1^2+R_5s_1^3)s_2u_1\\
			&\subset&\underline{\underline{s_2u_1u_2s_1^2s_2^2u_1}}+\underline{\underline{s_2u_1u_2(s_1^2s_2s_1)s_2u_1}}
			+s_2u_1(u_2u_1s_2s_1^{-1}s_2u_1)+\\&&+
			s_2u_1u_2s_1^2\grv u_1+s_2u_1(s_1^{-1}u_2s_1)(s_1s_2s_1^3)s_2u_1\\
			&\stackrel{(i)}{\subset}&\underline{\underline{s_2u_1u_2\grv u_1}}+s_2u_1s_2u_1(s_2^2s_1s_2)s_2u_1+U''''\\
			&\subset&s_2u_1s_2u_1\grv u_1+U''''.
			\end{array}}$
		
		The result follows from the fact that $s_2u_1s_2u_1\grv u_1=\underline{\underline{s_2u_1s_2\grv u_1}}$.
		\qedhere
	\end{itemize}
\end{proof}
We can now prove the following lemma that helps us to ``replace''  inside the definition of $U$ the elements $\grv^5$ and $\grv^{-5}$ by the elements $s_2^{-2}s_1^2s_2^3s_1^2s_2^3$ and $s_2^{-2}s_1^{2}s_2^{-2}s_1^2s_2^{-2}$ modulo $U''''$, respectively.
\begin{lem}\mbox{}  
	\vspace*{-\parsep}
	\vspace*{-\baselineskip}\\ 
	\begin{itemize}[leftmargin=0.8cm]
		\item [(i)]$s_2^{-2}s_1^2s_2^3s_1^2s_2^3\in u_1^{\times}\grv^5+ U''''.$
		\item [(ii)]$s_2^{-2}s_1^{2}s_2^{-2}s_1^2s_2^{-2}\in u_1^{\times}\grv^{-5}+U''''.$
		\item  [(iii)]$s_2^{-2}u_1s_2^{-2}s_1^{2}s_2^{-2}\subset U.$
	\end{itemize}
	\label{lol}
\end{lem}
\begin{proof}\mbox{}  
	\vspace*{-\parsep}
	\vspace*{-\baselineskip}\\ 
	\begin{itemize}[leftmargin=0.6cm]
		\item[(i)]  $s_2^{-2}s_1^2s_2^3s_1^2s_2^3=s_2^{-2}s_1^2s_2^3s_1(s_1s_2^3s_1^{-1})s_1
		=s_2^{-2}s_1^2s_2^2(s_2s_1s_2^{-1})s_1^2(s_1s_2s_1^{-1})s_1^2		=s_2^{-2}s_1^2s_2^2s_1^{-1}(s_2s_1^3s_2^{-1})s_1s_2s_1^2
		=s_2^{-2}s_1^2s_2^2s_1^{-2}s_2^2\grv s_1^2$. 
		We expand $s_1^{-2}$ as a linear combination of $s_1^{-1}$, 1, $s_1$ $s_1^{2}$ and $s_1^3$, where the coefficient of $s_1^3$ is invertible, and we have:
		
		$\small{\begin{array}[t]{lcl}
			s_2^{-2}s_1^2s_2^2s_1^{-2}s_2^2\grv s_1^2	
			&\in&s_2^{-2}s_1(s_1s_2^2s_1^{-1})s_2^2\grv u_1+ s_2^{-1}(s_2^{-1}s_1^4s_2)s_2^2\grv u_1+\\&&+s_2^{-2}s_1^2(s_2^2s_1s_2)s_2\grv u_1+s_2^{-1}(s_2^{-1}s_1^2s_2)s_2s_1^2s_2^2\grv s_2\grv u_1+
			u_1^{\times}s_2^{-2}s_1^2s_2^2s_1^{3}s_2^2\grv s_1^2\\
			&\stackrel{\ref{ll}}{\in}&\underbrace{s_1(s_1^{-1}s_2^{-2}s_1)s_2^{-1}s_1^2s_2^3\grv u_1+s_1(s_1^{-1}s_2^{-1}s_1)s_2^{4}s_1^{-1}s_2^2 s_1\grv u_1}_{\in u_1s_2u_1u_2u_1u_2\grv u_1}
			+s_2^{-2}s_1^3\grv^2u_1+\\&&+
			(s_2^{-1}s_1s_2)s_2(s_1^{-1}s_2s_1)s_1s_2^2\grv u_1+u_1^{\times}(u_1\grv^3+u_1^{\times}\grv^4+u_1s_2u_1u_2u_1u_2u_1)\grv s_1^2\\
			&\in&\underline{\underline{s_2^{-2}\grv^2u_1}}+s_1s_2s_1^{-1}s_2^2s_1(s_2^{-1}s_1s_2)s_2\grv  u_1+u_1\grv^3+u_1^{\times}\grv^5+u_1s_2u_1u_2u_1u_2u_1\grv u_1\\
			&\in&s_1s_2s_1^{-1}(u_2s_1^2s_2s_1^{-1}s_2)\grv u_1+u_1^{\times}\grv^5+u_1s_2u_1u_2u_1u_2u_1\grv u_1+U''''.
			\end{array}}$
		
		By proposition \ref{l2}(ii) we have that $s_2s_1^{-1}(u_2s_1^2s_2s_1^{-1}s_2)\grv u_1\subset s_2s_1^{-1}(u_1\grv^2+R_5s_2^2s_1^{-2}s_2s_1^{-1}s_2+u_1u_2u_1u_2u_1)\grv u_1	\stackrel{\ref{ll}}{\in}s_2\grv^2u_1+u_1s_2u_1u_2u_1u_2u_1\grv u_1$
		and, hence, the element $s_2^{-2}s_1^2s_2^2s_1^{-2}s_2^2\grv s_1^2 $ is inside $s_2\grv^2u_1+u_1^{\times}\grv^5+u_1s_2u_1u_2u_1u_2u_1\grv u_1+U''''$. We  notice that $u_1s_2u_1u_2u_1u_2u_1\grv u_1=u_1s_2u_1u_2u_1u_2\grv u_1$ and, hence, by lemma \ref{oo}(ii) and proposition \ref{xx}(iii)   we have that the element $s_2^{-2}s_1^2s_2^2s_1^{-2}s_2^2\grv s_1^2$ is inside $u_1^{\times}\grv^5+U''''$.

		\item[
		(ii)]$\small{\begin{array}[t]{lcl}
			s_2^{-2}s_1^{2}s_2^{-2}s_1^2s_2^{-2}&=&s_2^{-2}(as_1+b+cs_1^{-1}+ds_1^{-2}+es_1^{-3})s_2^{-2}s_1^2s_2^{-2}
			\\
			&\in&\underline{\underline{s_1(s_1^{-1}s_2^{-2}s_1)s_2^{-2}s_1^2s_2^{-2}}}+
			\underline{R_5s_2^{-4}s_1^2s_2}+\underline{\underline{R_5s_2^{-1}(s_2^{-1}s_1^{-1}s_2^{-2})s_1^2s_2^{-2}}}+\\&&+R_5s_2^{-1}(s_2^{-1}s_1^{-2}s_2)s_2^{-3}s_1^2s_2^{-2}+R_5
			s_2^{-2}s_1^{-2}(s_1^{-1}s_2^{-2}s_1)s_1s_2^{-2}\\
			&\in&R_5s_2^{-1}s_1s_2^{-2}(s_1^{-1}s_2^{-3}s_1)s_1s_2^{-2}+R_5
			s_2^{-1}(s_2^{-1}s_1^{-2}s_2)s_1^{-2}(s_2^{-1}s_1s_2)s_2^{-3}+U''''\\ 
			&\in&R_5s_2^{-1}s_1^2(s_1^{-1}s_2^{-1}s_1^{-3})s_2^{-1}s_1s_2^{-2}+R_5s_2^{-1}s_1s_2^{-1}(s_2^{-1}s_1^{-2}s_2)s_1^{-1}s_2^{-3}+U''''\\
			&\in&\underline{\underline{R_5s_2^{-1}s_1^2s_2^{-2}(s_2^{-1}s_1^{-1}s_2^{-2})s_1s_2^{-2}}}+R_5s_2^{-1}s_1(s_2^{-1}s_1s_2)s_2^{-2}s_1^{-2}s_2^{-3}+U''''
			\\
			&\in&R_5(s_2^{-1}s_1^2s_2)s_1^{-1}s_2^{-2}s_1^{-2}s_2^{-3}+U''''\\
			&\in& \grF(u_1s_2^{-2}s_1^2s_2^3s_1^2s_2^3)+U''''.
			\end{array}}$
		
		The result then follows from (i) and Lemma \ref{r1}.
		\item[(iii)] We expand $u_1$ as $R_5+R_5s_1+R_5s_1^{-1}+R_5s_1^2+R_5s_1^{-2}$ and we have that
		
		$\small{\begin{array}{lcl}s_2^{-2}u_1s_2^{-2}s_1^2s_2^{-2}&\subset& \underline{R_5s_2^{-4}s_1^2s_2^{-2}}+\underline{\underline{R_5s_2^{-1}(s_2^{-1}s_1s_2)s_2^{-3}s_1^2s_2^{-2}}}+\underline{\underline{R_5(s_2^{-2}s_1^{-1}s_2^{-1})s_2^{-1}s_1^2s_2^{-2}}}+\\&&+R_5s_2^{-2}s_1^2s_2^{-2}s_1^2s_2^{-2}+R_5s_2^{-2}s_1^{-2}s_2^{-2}s_1^2s_2^{-2}.
			\end{array}}$
		
		The element $s_2^{-2}s_1^2s_2^{-2}s_1^2s_2^{-2}$ is inside $U''''$, by (ii). Moreover, $s_2^{-2}s_1^{-2}s_2^{-2}s_1^2s_2^{-2}=
		s_2^{-1}\grv^{-1}(s_2^{-1}s_1^2s_2)s_2^{-3}
		=s_2^{-1}\grv^{-1}s_1s_2^2s_1^{-1}s_2^{-3}
		=s_2^{-1}s_1\grv^{-1}s_2^2s_1^{-1}s_2^{-3}=\underline{\underline{s_2^{-1}s_1^2(s_1^{-1}s_2^{-1}s_1^{-2})s_2s_1^{-1}s_2^{-3}}}$.
		
		\qedhere
	\end{itemize}
\end{proof}
We can now prove the main theorem of this section.
\begin{thm}
	\mbox{}  
	\vspace*{-\parsep}
	\vspace*{-\baselineskip}\\ 
	\begin{itemize}[leftmargin=0.8cm]
		\item [(i)]	$U=U''''+u_1\grv^{-5}$.
		\item[(ii)]$H_5=U$.
		\label{thh2}
	\end{itemize}
\end{thm}
\begin{proof}
	\mbox{}  
	\vspace*{-\parsep}
	\vspace*{-\baselineskip}\\ 
	\begin{itemize}[leftmargin=0.8cm]
		\item [(i)]	
		By definition, $U=U''''+u_1\grv^5+u_1\grv^{-5}$. Hence, it is enough to prove that $\grv^{-5}\in u_1^{\times}\grv^{5}+U''''.$
		By lemma \ref{lol}(ii) we have that $\grv^{-5}\in u_1^{\times}s_2^{-2}s_1^{2}s_2^{-2}s_1^2s_2^{-2}+U''''$. We expand $\bold
		{s_2^{-2}}$ as a linear combination of $s_2^{-1}$, 1, $s_2$, $s_2^2$ and $s_2^3$, where the coefficient of $s_2^3$ is invertible, and we have: 
		
		$\small{\begin{array}{lcl}
			s_2^{-2}s_1^{2} s_2^{-2}s_1^2\bold {s_2^{-2}}&\in& R_5s_2^{-1}(s_2^{-1}s_1^2s_2)s_2^{-3}s_1^2s_2^{-1}+\underline{R_5s_2^{-2}s_1^2s_2^{-2}s_1^2}+
			R_5(s_1^{-1}s_2^{-2}s_1)s_1s_2^{-1}(s_2^{-1}s_1^2s_2)+\\&&+R_5s_1(s_1^{-1}s_2^{-2}s_1)s_1s_2^{-1}(s_2^{-1}s_1^2s_2)s_2+R_5^{\times}s_2^{-2}s_1^{2}s_2^{-2}s_1^2s_2^3\\
			&\in&u_1s_2^{-1}s_1s_2^2(s_1^{-1}s_2^3s_1)s_1s_2^{-1}+u_1s_2s_1^{-1}(s_1^{-1}s_2^{-1}s_1)(s_2^{-1}s_1s_2)s_2s_1^{-1}+\\&&+u_1s_2s_1^{-1}(s_1^{-1}s_2^{-1}s_1)(s_2^{-1}s_1s_2)s_2s_1^{-1}s_2+
			u_1^{\times}s_2^{-2}s_1^{2}s_2^{-2}s_1^2s_2^3+U''''\\
			&\in&u_1\grF(s_2s_1^{-1}u_2u_1s_2s_1^{-1}s_2)+
			\underline{\underline{u_1s_2s_1^{-1}s_2(s_1^{-1}s_2^{-1}s_1)s_2s_1^{-1}s_2
					s_1^{-1}}}+\\&&+u_1s_2s_1^{-1}s_2(s_1^{-1}s_2^{-1}s_1)s_2s_1^{-1}s_2
			s_1^{-1}s_2+u_1^{\times}s_2^{-2}s_1^{2}s_2^{-2}s_1^2s_2^3+U''''\\
			&\stackrel{\ref{l2}}{\in}&u_1\grF\big(s_2s_1^{-1}\grv^2u_1+s_2s_1^{-1}s_2^2s_1^{-2}s_2s_1^{-1}s_2+
			s_2s_1^{-1}u_1u_2u_1u_2u_1\big)+\\&&+
			u_1s_2s_1^{-1}s_2^2s_1^{-2}s_2s_1^{-1}s_2+u_1^{\times}s_2^{-2}s_1^{2}s_2^{-2}s_1^2s_2^3+U''''
			\end{array}}$
		
		However, by lemma \ref{oo}(ii) proposition \ref{xx}(i)  we have that $\grF\big(s_2s_1^{-1}\grv^2u_1+s_2s_1^{-1}s_2^2s_1^{-2}s_2s_1^{-1}s_2+\underline{\underline{s_2s_1^{-1}u_1u_2u_1u_2u_1}}\big)\subset \grF(U'''')\stackrel{\ref{r1}}{\subset}U''''$. Therefore, it will be sufficient to prove that the element $s_2^{-2}s_1^{2}s_2^{-2}s_1^2s_2^3$ is inside  $u_1^{\times}\grv^5+U''''$. We expand $\bold
		{s_2^{-2}}$ as a linear combination of $s_2^{-1}$, 1, $s_2$, $s_2^2$ and $s_2^3$, where the coefficient of $s_2^3$ is invertible, and we have: 
		$s_2^{-2}s_1^{2}\bold {s_2^{-2}}s_1^2s_2^3 \in u_1s_2^{-3}(s_2s_1^2s_2^{-1})s_1^2s_2^3+\underline{u_1s_2^{-2}s_1^4s_2^3}+
		\underline{\underline{u_1s_2^{-2}(s_1^2s_2s_1)s_1s_2^3}}+u_1s_2^{-1}(s_2^{-1}s_1^2s_2)s_2s_1^2s_2^3+u_1^{\times}s_2^{-2}s_1^2s_2^3s_1^2s_2^3$.
		By lemma \ref{lol}(i) we have that $u_1^{\times}s_2^{-2}s_1^2s_2^3s_1^2s_2^3\subset u_1^{\times}\grv^5+U''''.$
		Therefore, $s_2^{-2}s_1^{2}s_2^{-2}s_1^2s_2^3 \in\underline{\underline{u_1(s_1s_2^{-3}s_1^{-1})s_2^2s_1^3s_2^3}}+u_1s_2^{-1}s_1s_2^2s_1^{-1}s_2s_1^2s_2^3+u_1^{\times}\grv^5+U''''$. It remains to prove that the element $s_2^{-1}s_1s_2^2s_1^{-1}s_2s_1^2s_2^3$ is inside $U''''$. Indeed, 
		$s_2^{-1}s_1s_2^2s_1^{-1}s_2s_1^2s_2^3=(s_2^{-1}s_1s_2)s_2(s_1^{-1}s_2s_1)s_1s_2^3=s_1s_2(s_1^{-1}s_2^2s_1)(s_2^{-1}s_1s_2)s_2^2=s_1s_2^2s_1^2(s_2^{-1}s_1s_2)s_1^{-1}s_2^2=s_1s_2^3(s_2^{-1}s_1^3s_2)s_1^{-2}s_2^2=\underline{\underline{s_1^2(s_1^{-1}s_2^3s_1)s_2^3s_1^{-3}s_2^2}}$.

		\item[(ii)] As we explained in the beginning of this section, since $1\in U$ it will be sufficient to prove that $U$ is invariant under left multiplication by $s_2$. We use the fact that $U$ is equal to the RHS of (i) and by the definition of $U''''$ we have:
		
		$\small{\begin{array}{lcl} U&=&U'+\sum\limits_{k=2}^4\grv^{k}u_1+\sum\limits_{k=2}^5\grv^{-k}u_1+
			u_1s_2^{-2}s_1^2s_2^{-1}s_1s_2^{-1}u_1+u_1s_2^{2}s_1^{-2}s_2s_1^{-1}s_2u_1+
			u_1s_2s_1^{-2}s_2^{2}s_1^{-2}s_2^{2}u_1+\\&&+u_1s_2^{-1}s_1^2s_2^{-2}s_1^2s_2^{-2}u_1.
			\end{array}}$
		
		On one hand, $s_2(U'+\grv^2u_1+
		u_1s_2^{-2}s_1^2s_2^{-1}s_1s_2^{-1}u_1+u_1s_2^{2}s_1^{-2}s_2s_1^{-1}s_2u_1)\subset U$ (proposition \ref{p2}, lemma \ref{oo}$(ii)$ and proposition \ref{xx}$(i),(ii)$). 
		On the other hand, $$\sum\limits_{k=2}^{5}s_2\grv^{-k}u_1=\sum\limits_{k=2}^{5}s_1^{-2}s_2^{-1}\grv^{-k+1}u_1\subset \grF(\sum\limits_{k=2}^{5}u_1s_2\grv^{k-1}u_1)\stackrel{\ref{p1}\text{ and }\ref{oo}(ii)}{\subset}u_1\grF(U+s_2\grv^{3}+s_2\grv^4)u_1
		$$
		Therefore, by lemma \ref{r1} we only need to prove that 
		$$s_2(\grv^3u_1+\grv^4 u_1+
		u_1s_2s_1^{-2}s_2^{2}s_1^{-2}s_2^{2}u_1+
		u_1s_2^{-1}s_1^2s_2^{-2}s_1^2s_2^{-2}u_1)\subset U.$$
		We first notice that $s_2\grv^4u_1=s_2\grv^3\grv u_1=s_2\grv^3u_1\grv$. Therefore, in order to prove that $s_2(\grv^3u_1+\grv^4u_1)\subset U$, it will be sufficient to prove that $s_2\grv^3u_1\subset u_1s_2u_1u_2u_1u_2u_1$ (propositions \ref{p2} and \ref{xx}(iii)).
		Indeed, we have: $s_2\grv^3 u_1=s_2\grv^2\grv u_1
		\stackrel{\ref{oo}(ii)}{=}s_1s_2s_1^4s_2s_1^3s_2\grv u_1
		=s_1s_2s_1^4s_2s_1^4(s_1^{-1}s_2^2s_1)s_1s_2u_1
		=s_1s_2s_1^4s_2s_1^4s_2s_1^2(s_2^{-1}s_1s_2)u_1
		\subset u_1s_2u_1(s_2u_1s_2u_1s_2u_1)$.
		However, by lemma \ref{oo}(i)  we have that $u_1s_2u_1(s_2u_1s_2u_1s_2u_1)\subset u_1s_2u_1(\grv^2u_1+u_1u_2u_1u_2u_1)$.
		The result follows from lemma \ref{oo}(ii).

		It remains to prove that $s_2
		u_1s_2^{-1}s_1^2s_2^{-2}s_1^2s_2^{-2}u_1$ and $s_2u_1s_2s_1^{-2}s_2^{2}s_1^{-2}s_2^{2}u_1$ are subsets of  $U.$ We have: 
		
		$\small{\begin{array}{lcl}
			s_2u_1s_2^{-1}s_1^2s_2^{-2}s_1^2s_2^{-2}&=&s_2(R_5+R_5s_1+R_5s_1^{-1}+
			R_5s_1^2+R_5s_1^{-2})s_2^{-1}s_1^2s_2^{-2}s_1^2s_2^{-2}\\
			&\subset&\underline{R_5s_1^2s_2^{-2}s_1^2s_2^{-2}}+\underline{\underline{
					R_5(s_2s_1s_2^{-1})s_1^2s_2^{-2}s_1^2s_2^{-2}}}+\underline{\underline{
					R_5(s_2s_1^{-1}s_2^{-1})s_1^2s_2^{-2}s_1^2s_2^{-2}}}+\\&&
			+R_5(s_2s_1^{2}s_2^{-1})s_1^2s_2^{-2}s_1^2s_2^{-2}+
			R_5(s_2s_1^{-2}s_2^{-1})s_1^2s_2^{-2}s_1^2s_2^{-2}
			\\
			&\subset&u_1s_2^2s_1^3s_2^{-2}s_1^2s_2^{-2}+u_1s_2^{-2}u_1s_2^{-2}s_1^2s_2^{-2}+U.
			\end{array}}$
		
		However, by lemma \ref{lol}(iii) we have that $u_1s_2^{-2}u_1s_2^{-2}s_1^2s_2^{-2}\subset U$. Therefore, it remains to prove that the element $s_2^2s_1^3s_2^{-2}s_1^2s_2^{-2}$ is inside $ U$. For this purpose, we expand 
		$s_1^3$ as a linear combination of $s_1^2$, $s_1$, 1, $s_1^{-1}$ and $s_1^{-2}$ and  we have:
		
		$\small{\begin{array}{lcl}s_2^2s_1^3s_2^{-2}s_1^2s_2^{-2}&\in&
			R_5s_2^2s_1^2s_2^{-2}s_1^2s_2^{-2}+\underline{\underline{
					R_5s_2^2(s_1s_2^{-2}s_1^{-1})s_1^3s_2^{-2}}}+\underline{u_1s_2^{-2}}+
			\underline{\underline{R_5s_1^{-1}(s_1s_2^2s_1^{-1})s_2^{-2}s_1^2s_2^{-2}}}+\\&&+
			R_5s_2^{2}s_1^{-2}s_2^{-2}s_1^{2}s_2^{-2}.
			\end{array}}$
		
		However, $s_2^2s_1^2s_2^{-2}s_1^2s_2^{-2}=s_2(s_2s_1^2s_2^{-1})s_2^{-1}s_1(s_1s_2^{-2}s_1^{-1})s_1=s_2s_1^{-1}s_2(s_2s_1s_2^{-1})s_1s_2^{-1}s_1^{-2}s_2s_1=s_2s_1^{-2}(s_1s_2s_1^{-1})(s_2s_1^2s_2^{-1})s_1^{-2}s_2s_1=s_2s_1^{-2}s_2^{-1}(s_1s_2s_1^{-1})s_2^2s_1^{-1}s_2s_1=\underline{\underline{u_1s_2s_1^{-2}s_2^{-2}(s_1s_2^3s_1^{-1})s_2s_1}}$. 
		Moreover, we expand $s_1^2$ as a linear combination of $s_1$, 1, $s_1^{-1}$, $s_1^{-2}$ and $s_1^{-3}$ and we have:

		$\small{\begin{array}{lcl}

			s_2^{2}s_1^{-2}s_2^{-2}s_1^{2}s_2^{-2}

			&\in &\underline{R_5s_2^2s_1^{-2}s_2^{-4}}+\underline{\underline{
					R_5s_1^{-1}(s_1s_2^2s_1^{-1})(s_1^{-1}s_2^{-2}s_1)s_2^{-2}}}+R_5s_2^2s_1^{-2}s_2^{-1}(s_2^{-1}s_1^{-1}s_2^{-2})+\\&&+R_5s_2(s_2s_1^{-2}s_2^{-1})\grv^{-1}s_2^{-1}+
			\grF(s_2^{-2}s_1^2s_2^2s_1^3s_2^2)\\

			&\stackrel{\ref{ll}}{\in}&R_5s_2^2s_1^{-2}\grv^{-1}s_1^{-1}+
			R_5s_2s_1^{-1}s_2^{-2}s_1\grv^{-1}s_2^{-1}+\grF(U)+U\\

			&\stackrel{\ref{r1}}{\subset}&\underline{s_2^2\grv^{-1}u_1}+
			R_5s_2s_1^{-1}s_2^{-2}\grv^{-1}s_1s_2^{-1}+U\\
			&\subset&s_2s_1^{-1}u_2u_1s_2^{-1}s_1s_2^{-1}+U.
			\end{array}}$
		
		Therefore, by proposition \ref{xx}(ii) we have that the element $s_2^{2}s_1^{-2}s_2^{-2}s_1^{2}s_2^{-2}$ is inside $U$ and, hence, $s_2
		u_1s_2^{-1}s_1^2s_2^{-2}s_1^2s_2^{-2}u_1\subset U$.
		
		In order to finish the proof that $H_5=U$ it remains to prove that $s_2u_1s_2s_1^{-2}s_2^2s_1^{-2}s_2^2\subset U$. For this purpose we expand $u_1$ as 
		$R_5+R_5s_1+R_5s_1^2+R_5s_1^3+
		R_5s_1^4$ and we have:
		
		$\small{\begin{array}[t]{lcl}
			s_2u_1s_2s_1^{-2}s_2^2s_1^{-2}s_2^2
			&\subset&\grF(s_2^{-2}u_1s_2^{-2}s_1^2s_2^{-2})+\underline{\underline{
					R_5(s_2s_1s_2)s_1^{-2}s_2^2s_1^{-2}s_2^2}}+R_5\grv s_1^{-2}s_2^2s_1^{-2}s_2^2+\\&&
			+R_5s_2s_1^3s_2s_1^{-2}s_2^2s_1^{-2}s_2^{2}	+		R_5s_2s_1^4s_2s_1^{-2}s_2^2s_1^{-2}s_2^2.
			\end{array}}$
		
		However, by lemma \ref{lol}(iii) we have that $\grF(s_2^{-2}u_1s_2^{-2}s_1^2s_2^{-2})$ is a subset of $\grF(U)$ and, hence, by lemma \ref{r1}, a subset of $U$. Moreover, $\grv s_1^{-2}s_2^2s_1^{-2}s_2^2=\underline{\underline{R_5s_1^{-2}\grv s_2^2s_1^{-2}s_2^2}}$. It remains to prove that $C:=R_5s_2s_1^3s_2s_1^{-2}s_2^2s_1^{-2}s_2^{2}	+		R_5s_2s_1^4s_2s_1^{-2}s_2^2s_1^{-2}s_2^2$ is a subset of $U$. We have:

		$\small{\begin{array}[t]{lcl}
			C
			&=&
			R_5s_2s_1^2(s_1s_2s_1^{-1})s_1^{-1}s_2(s_2s_1^{-2}s_2^{-1})s_2^3+
			R_5s_1^{-1}(s_1s_2s_1^4)s_2s_1^{-2}s_2(s_2s_1^{-2}s_2^{-1})s_2^3\\
			&=&
			R_5(s_2s_1^2s_2^{-1})(s_1s_2s_1^{-1})s_2s_1^{-1}s_2^{-2}s_1s_2^3+R_5s_1^{-1}s_2^4(s_1s_2s_1^{-1})s_1^{-1}(s_2s_1^{-1}s_2^{-1})s_2^{-1}s_1s_2^3\\
			&=&R_5s_1^{-1}s_2(s_2s_1s_2^{-1})(s_1s_2^2s_1^{-1})s_2^{-2}s_1s_2^3+R_5s_1^{-1}s_2^2\grv s_1^{-2}s_2^{-1}s_1s_2^{-1}s_1s_2^3\\
			&=&R_5s_1^{-1}s_2s_1^{-1}(s_2s_1s_2^{-1})s_1^2s_2^{-1}s_1s_2^3+R_5s_1^{-1}s_2^2s_1^{-2}(s_2s_1^3s_2^{-1})s_1s_2^3+U\\
			&=&\underline{\underline{R_5s_1^{-1}s_2s_1^{-2}(s_2s_1^3s_2^{-1})s_1s_2^3}}+u_1s_2^2s_1^{-3}s_2^3s_1^2s_2^3+U. 
			
			\end{array}}$
		
		We expand $s_1^{-3}$ as a linear combination of $s_1^{-2}$, $s_1^{-1}$, 1, $s_1$, $s_1^2$ and $s_1^3$ and we have that 
		
		$\small{\begin{array}{lcl}u_1s_2^2s_1^{-3}s_2^3s_1^2s_2^3
			&\subset& u_1s_2^2s_1^{-2}s_2^3s_1^2s_2^3+\underline{\underline{u_1(s_1s_2^2s_1^{-1})s_2^3s_1^2s_2^3}}+\underline{u_1s_2^5s_1^2s_2^3}+
			\underline{\underline{u_1s_2(s_2s_1s_2^3)s_1^2s_2^3}}+\\&&+
			u_1s_2^2s_1^2s_2^3s_1^2s_2^3
			.
			\end{array}}$
		
		Hence, in order to finish the proof that $H_5=U$ we have to prove that $u_1s_2^2s_1^{-2}s_2^3s_1^2s_2^3$ and $u_1s_2^2s_1^2s_2^3s_1^2s_2^3$ are subsets of $U$. We have:

		$\small{\begin{array}{lcl}
			u_1s_2^2s_1^{-2}s_2^3s_1^2s_2^3&=&u_1s_2^2s_1^{-2}(as_2^2+bs_2+c+ds_2^{-1}+es_2^{-2})s_1^2s_2^3\\
			&\subset&u_1s_2^2s_1^{-2}s_2^2s_1^2s_2^3+
			u_1s_2^2s_1^{-2}s_2^2s_1^2s_2^3+\underline{u_1s_2^5}+\underline{\underline{
					u_1s_2(s_2s_1^{-2}s_2^{-1})s_1^2s_2^3}}+u_1s_2^2s_1^{-2}s_2^{-2}s_1^2s_2^3
			\end{array}}$
		
		However, we have that	$u_1s_2^2s_1^{-2}s_2^2s_1^2s_2^3=
		\underline{\underline{u_1(s_1^{-1}s_2^2s_1)s_2(s_2^{-1}s_1^{-3}s_2)s_1^2s_2^3}}$. Moreover,  we have
		$u_1s_2^2s_1^{-2}s_2^{-2}s_1^2s_2^3=u_1(s_1s_2^2s_1^{-1})(s_1^{-1}s_2^{-2}s_1)s_1s_2^3=	u_1s_2^{-1}s_1^2s_2^3\grv^{-1}s_1s_2^3=u_1s_2^{-1}s_1^2s_2^{3}s_1(s_2^{-1}s_1^{-2}s_2)s_2=u_1s_2^{-1}s_1^2s_2^3s_1^2s_2^{-1}(s_2^{-1}s_1^{-1}s_2)\subset\grF(s_2u_2u_1s_2s_1^{-1}s_2)$. By proposition \ref{xx}(i) and lemma \ref{r1} we have that $u_1s_2^2s_1^{-2}s_2^{-2}s_1^2s_2^3\subset U$. It remains to prove that $u_1s_2^2s_1^{-2}s_2^2s_1^2s_2^3\subset U$. We notice that $u_1s_2^2s_1^{-2}s_2^2s_1^2s_2^3=
		u_1s_2^3(s_2^{-1}s_1^{-2}s_2)s_2s_1^2s_2^3
		=u_1(s_1^{-1}s_2^3s_1)s_2^{-2}s_1^{-1}s_2s_1^2s_2^3
		=u_1s_2s_1^3s_2^{-3}s_1^{-1}s_2s_1^2s_2^3$. We expand $s_2^3$ as a linear combination of $s_2^2$, $s_2$, 1, $s_2^{-1}$ and $s_2^{-2}$ and we have:
		
		$\small{\begin{array}{lcl}
			
			u_1s_2s_1^3s_2^{-3}s_1^{-1}s_2s_1^2s_2^3

			&\subset&u_1s_2s_1^3s_2^{-2}(s_2^{-1}s_1^{-1}s_2)s_1^2s_2^2+u_1s_2s_1^3s_2^{-3}s_1^{-1}\grv+
			\underline{\underline{u_1s_2s_1^3s_2^{-3}s_1^{-1}s_2s_1^2}}+\\&&+
			\underline{\underline{u_1s_2s_1^3s_2^{-3}s_1^{-1}(s_2s_1^2s_2^{-1})}}+u_1s_2s_1^3s_2^{-3}s_1^{-1}(s_2s_1^2s_2^{-1})s_2^{-1}\\
			&\subset&u_1s_2s_1^3s_2^{-2}s_1(s_2^{-1}s_1s_2)s_2+
			u_1s_2s_1^3s_2^{-3}s_1^{-2}s_2(s_2s_1s_2^{-1})+\underline{\underline{u_1s_2s_1^3s_2^{-3}\grv s_1^{-1}}}+U\\
			&\subset&
			u_1(s_2u_1u_2u_1s_2s_1^{-1}s_2)u_1+U.
			\end{array}}$\\
		The result follows from \ref{xx}(i) and \ref{r1}.
		
		Using analogous calculations we will prove that $u_1s_2^2s_1^2s_2^3s_1^2s_2^3$ is a subset of $U$. We have:
		
		$\small{\begin{array}{lcl}
			u_1s_2^2s_1^2s_2^3s_1^2s_2^3&=&u_1s_2^2s_1^2(as_2^2+bs_2+c+ds_2^{-1}+es_2^{-2})s_1^2s_2^3\\
			&\subset&u_1s_2^2s_1^2s_2^2s_1^2s_2^3+u_1s_2\grv s_1^2s_2^3+\underline{u_1s_2^2s_1^4s_2^3}
			+\underline{\underline{u_1s_2(s_2s_1^2s_2^{-1})s_1^2s_2^3}}+
			u_1s_2^2s_1^2s_2^{-2}s_1^2s_2^3.
			\end{array}}$
		
		However, $u_1s_2\grv s_1^2s_2^3=\underline{\underline{u_1s_2s_1^2\grv
				s_2^3}}$. Therefore, it remains to prove that $u_1s_2^2s_1^2s_2^2s_1^2s_2^3$ and 	$u_1s_2^2s_1^2s_2^{-2}s_1^2s_2^3$ are subsets of $U$. We have:
		
		$\small{\begin{array}{lcl}
			u_1s_2^2s_1^2s_2^2s_1^2s_2^3&=&u_1s_2^2s_1^2s_2^2s_1^2(as_2^2+bs_2+c+ds_2^{-1}+es_2^{-2})\\
			&\subset&u_1s_2\grv^2s_2+u_1s_2\grv^2+\underline{u_1s_2^2s_1^2s_2^2s_1^2}+
			u_1s_2\grv s_2s_1^2s_2^{-1}+u_1s_2\grv s_2s_1^2s_2^{-2}\\
			\end{array}}$
		
		By lemma \ref{oo}(ii) we have that $u_1s_2\grv^2s_2+u_1s_2\grv^2\subset u_1s_2s_1^4s_2s_1^3s_2u_1s_2+U$. However, by lemma \ref{oo}(i) we have  $u_1s_2s_1^4(s_2s_1^3s_2u_1s_2)\subset u_1s_2s_1^4(\grv^2u_1+u_1u_2u_1u_2u_1)\subset u_1s_2\grv^2u_1+\underline{\underline{u_1s_2u_1u_2u_1u_2u_1}}$. Using another time lemma \ref{oo}(ii) we have that $u_1s_2s_1^4(\grv^2u_1+u_1u_2u_1u_2u_1)\subset U$. It remains to prove that $D:=u_1s_2\grv s_2s_1^2s_2^{-1}+u_1s_2\grv s_2s_1^2s_2^{-2}$ is a subset of  $U$. We have:

		$\small{\begin{array}{lcl}
			D &=& u_1s_2\grv(s_2s_1^2s_2^{-1})+u_1s_2\grv(s_2s_1^2s_2^{-1})s_2^{-1}\\
			&=&u_1s_2\grv s_1^{-1}s_2^2s_1+u_1s_2\grv s_1^{-1}s_2^2s_1s_2^{-1}\\
			&=&\underline{\underline{u_1s_2s_1^{-1}\grv s_2^2s_1}}+u_1s_2s_1^{-1}\grv s_2^2s_1s_2^{-1}.
			
			\end{array}}$
		
		However, we have $u_1s_2s_1^{-1}\grv s_2^2s_1s_2^{-1}=u_1s_2s_1^{-1}\grv s_2(s_2s_1s_2^{-1})=u_1s_2s_1^{-1}s_2s_1(s_1s_2s_1^{-1})s_2s_1=\underline{\underline{u_1s_2s_1^{-1}(s_2s_1s_2^{-1})s_1s_2^2s_1}},$ meaning that $D \subset U$.
		
		Using analogous calculations we will prove that $u_1s_2^2s_1^2s_2^{-2}s_1^2s_2^3\subset U$. We have:

		$\small{\begin{array}{lcl}
			u_1s_2^2s_1^2s_2^{-2}s_1^2s_2^3
			&=&
			u_1s_2(s_2s_1^2s_2^{-1})s_2^{-1}s_1^2s_2^3\\
			
			&=&u_1s_2s_1^{-1}s_2^2s_1s_2^{-1}s_1^2(as_2^2+bs_2+c+ds_2^{-1}+es_2^{-2})\\
			&\subset&u_1(s_1s_2s_1^{-1})s_2(s_2s_1s_2^{-1})s_1^2s_2^2+
			\underline{\underline{u_1s_2s_1^{-1}s_2^2s_1(s_2^{-1}s_1^2s_2)}}+
			
			\underline{
				\underline{u_1s_2s_1^{-1}s_2^2s_1s_2^{-1}s_1^2}}+\\&&+
			u_1s_2s_1^{-1}s_2(s_2s_1s_2^{-1})s_1^2s_2^{-1}+u_1s_2s_1^{-1}s_2(s_2s_1s_2^{-1})s_1^2s_2^{-2}\\
			
			&\subset&u_1s_2^{-1}(s_1s_2^2s_1^{-1})s_2s_1^3s_2^2
			+\underline{\underline{u_1s_2s_1^{-1}s_2s_1^{-1}(s_2s_1^3s_2^{-1})}}+\\&&+u_1s_2s_1^{-2}(s_1s_2s_1^{-1})(s_2s_1^3s_2^{-1})s_2^{-1}+U\\ 
			&\subset&u_1s_2^{-2}s_1^2s_2^2s_1^3s_2^2+u_1s_2s_1^{-2}s_2^{-1}(s_1s_2s_1^{-1})s_2^2(s_2s_1s_2^{-1})+U\\
			&\stackrel{
				\ref{ll}}{\subset}&\underline{\underline{u_1s_2s_1^{-2}s_2^{-2}(s_1s_2^3s_1^{-1})s_2s_1}}+U.
			\end{array}}$
		
		\qedhere
	\end{itemize}
\end{proof}

\begin{cor}$H_5$ is a free $R_5$-module of rank $r_5=600$ and, therefore, the BMR freeness conjecture holds for the exceptional group $G_{16}$.
\end{cor}
\begin{proof}
By proposition \ref{rp} it will be sufficient to show that $H_5$ is generated as $R_5$-module by $r_5$ elements. By theorem \ref{thh2}, the definition of $U''''$ and the fact that $u_1u_2u_1=u_1+u_1s_2u_1+u_1s_2^{-1}u_1+u_1s_2^2u_1+u_1s_2^{-2}u_1$ we have that $H_5$ is spanned as left $u_1$-module by 120 elements. Since $u_1$ is spanned by 5 elements as a $R_5$-module, we have that $H_5$ is spanned over $R$ by $r_5=600$ elements. 
\end{proof}
\section{An application: The irreducible representations of $B_3$ of dimension at most 5}
\indent

In 1999 I. Tuba and H. Wenzl classified the irreducible representations of the braid group $B_3$ of dimension $k$ at most 5 over an algebraically closed field $K$ of any characteristic (see \cite{tuba}) and, therefore, of $PSL_2(\ZZ)$, since the quotient group $B_3$ modulo its center is isomorphic to $PSL_2(\ZZ)$. Recalling that $B_3$ is given by generators $s_1$ and $s_2$ that satisfy the relation $s_1s_2s_1=s_2s_1s_2$, we assume that $s_1\mapsto A, s_2\mapsto B$ is an irreducible representation of $B_3$, where $A$ and $B$ are invertible $k\times k$ matrices over $K$ satisfying $ABA=BAB$. I. Tuba and H. Wenzl proved that $A$ and $B$ can be chosen to be in \emph{ordered triangular form}\footnote{Two $k\times k$ matrices are in ordered triangular form if one of them is an upper triangular matrix with eigenvalue $\grl_i$ as $i$-th diagonal entry, and the other is a lower triangular matrix with eigenvalue $\grl_{k+1-i}$
	as $i$ -th diagonal entry.} with coefficients completely determined by the eigenvalues (for $k\leq3$) or by the eigenvalues and by the choice of a $k$-th root of det$A$ (for $k>3$). Moreover, they proved that such irreducible representations exist if and only if the eigenvalues do not annihilate some polynomials $P_k$ in the eigenvalues and the choice of the $k$th root of det$A$, which they determined explicitly. 

At this point, a number of questions arise: what is the reason we do not expect their methods to work for any dimension beyond 5 (see remark 2.11, 3 in \cite{tuba} )? Why are the matrices in this neat form? In \cite{tuba}, remark 2.11, 4 there is an explanation for the nature of the polynomials $P_k$. However, there is no argument connected with the nature of $P_k$ that explains the reason why these polynomials provide  a  necessary condition for a representation of this form to be irreducible. In this section we answer these questions by proving in a different way this classification of the irreducible representations of the braid group $B_3$ as a consequence of the BMR freeness conjecture for the generic Hecke algebra of the finite quotients of the braid group $B_3$, we proved in the previous section. The fact that there is a connexion between the classification of irreducible representations of dimension at most 5 and the finite quotients of $B_3$ has already been suspected by I. Tuba and H. Wenzl (see  remark 2.11, 5 in \cite{tuba}).

\subsection{Some preliminaries}
\indent

We set $\tilde{R_k}=\ZZ[u_1^{\pm1},...,u_k^{\pm1}]$, $k=2,3,4,5$. Let $\tilde{H_k}$ denote the quotient of the group algebra $\tilde{R_k}B_3$ by the relations $(s_i-u_1)...(s_i-u_k)$.  In the previous sections we proved that $H_k$ is a free $R_k$-module of rank $r_k$. Hence, $\tilde{H_k}$ is a free $\tilde{R_k}$-module of rank $r_k$  (Lemma 2.3 in \cite{marinG26}). We now assume that $\tilde{H_k}$ has a unique symmetrizing trace $t_k: \tilde{H_k} \arw \tilde{R_k}$ (i.e. a trace function such that the bilinear form $(h, h')\mapsto t_k(hh')$ is non-degenerate), having nice properties (see \cite{bmm}, theorem 2.1): for
example, $t_k(1)=1$, which means that $t_k$ specializes to the canonical symmetrizing form on $\CC W_k$. 

Let $\grm_{\infty}$  be the group of all roots of unity in $\CC$. We recall that $W_k$ is the finite quotient group $B_3/\langle s_i^k \rangle$, $k=2, 3, 4$ and 5 and we let $K_k$ be the \emph{field of definition} of $W_k$, i.e. the number field contained in $\mathbb{Q}(\grm_{\infty})$, which is generated by the traces of all elements of $W_k$.
We denote by $\grm(K_k)$ the group of all roots of unity of $K_k$ and, for every integer $m>1$, we set $\grz_m:=$exp$(2
\grp i/m)$. 

Let $\mathbf{v}=(v_1,...,v_k)$ be a set of $k$ indeterminates such that, for every $i\in\{1,...,k\}$, we have $v_i^{|\grm(K_k)|}=\grz_k^{-i}u_i$. By extension of scalars we obtain a $\CC(\mathbf{v})$-algebra
$\CC(\mathbf{v})\tilde{H_k}:=\tilde{H_k}\otimes_{\tilde{R_k}}\CC(\mathbf{v})$, which is split semisimple (see \cite{malle1}, theorem 5.2). Since the algebra $\CC(\mathbf{v})\tilde{H_k}$ is split, by Tits' deformation theorem (see theorem 7.4.6 in \cite{geck}), the specialization $v_i \mapsto 1$ induces a bijection
Irr$(\CC(\mathbf{v})\tilde{H_k})\rightarrow$ Irr$(W_k)$.
\begin{ex}
	Let $W_4:=G_{8}$. The field of definition of $G_{8}$ is $K_4:=\mathbb{Q}(i)$ (One can find this field in Appendix A, table A.1 in \cite{brouebook}). Since $|\grm(K_4)|=4$ we set $v_1^{4}:=i^{-1}u_1=-iu_1$, $v_2^4:=i^{-2}u_2=-u_2$, $v_3^4:=i^{-3}u_3=iu_3$ and $v_4^4:=i^{-4}u_4=u_4$. The theorem of G. Malle we mentioned above states that the algebra $\CC(v_1,v_2,v_3,v_4)\tilde{H_4}$ is split semisimple and, hence, its irreducible characters are in bijection with
	the irreducible characters of $G_8$.
\end{ex}
Let $\varrho: B_3 \arw GL_n(\CC)$ be an irreducible representation of $B_3$ of dimension $k\leq 5$. We set $A:=\varrho(s_1)$ and $B:=\varrho(s_2)$. The matrices $A$ and $B$ are similar since $s_1$ and $s_2$ are conjugate $( s_2=(s_1s_2)s_1(s_1s_2)^{-1})$. Hence, by  Cayley-Hamilton theorem of linear algebra, there is a monic polynomial $m(X)=X^k+m_{n-1}X^{n-1}+...+m_1X+m_0\in \CC[X]$ of degree $k$ such that $m(A)=m(B)=0$. Let $R^k_K$ denotes  the integral closure of
$R_k$ in $K_k$. We fix $\gru: R^k_K\arw \mathbb{C}$ a \emph{specialization} of $R^k_K$, defined by $u_i\mapsto \grl_i$, where $\grl_i$ are the eigenvalues of $A$ (and $B$). We notice that $\gru$ is well-defined, since $m_0=$det$A\in \CC^{\times}$. Therefore, in order to determine $\varrho$ it will be sufficient to  describe the irreducible $\CC\tilde{ H_k}:=\tilde{H_k}\otimes_{\gru}\CC$-modules of dimension $k$.

When the algebra $\CC\tilde H_k$ is semisimple, we can use again Tits' deformation theorem and we have a canonical bijection between  the set of irreducible characters of $\CC\tilde{H_k}$ and  the set of irreducible characters of $\CC(\mathbf{v})\tilde{H_k}$,  which are in bijection with
the irreducible characters of $W_k$. However, this is not always the case. In order to determine the irreducible representations of  $\CC\tilde H_k$ in the general case (when we don't know a priori that  $\CC\tilde H_k$ is semisimple) we use a different approach.

Let $R_0^{+}\big(\CC(\mathbf{v})\tilde{H_k}\big)$ (respectively $R_0^{+}(\CC\tilde{ H_k})$) denote the subset of the \emph{Grothendieck group} of the category of finite dimensional $\CC(\mathbf{v})\tilde{H_k}$ (respectively $\CC\tilde{H_k}$)-modules consisting of elements $[V]$, where $V$ is a $\CC(\mathbf{v})\tilde{H_k}$ (respectively $\CC\tilde{H_k}$)-module (for more details, one may refer to \textsection 7.3 in \cite{geck}). By theorem 7.4.3 in \cite {geck} we obtain a well-defined decomposition map $$d_{\gru}: R_0^{+}\big(\CC(\mathbf{v})\tilde{H_k}) \arw R_0^{+}(\CC\tilde{ H_k}).$$
The corresponding \emph{decomposition matrix} is the Irr$\big(\CC(\mathbf{v})\tilde{H_k}\big)\times$ Irr$(\CC\tilde{ H_k})$ matrix $(d_{\grx\grf})$ with non-negative integer entries such that $d_{\gru}([V_{\grx}])=\sum\limits_{\grf}d_{\grx\grf}[V'_{\grf}]$, where $V_{\grx}$ is an irreducible $\CC(\mathbf{v})\tilde{H_k}$-module with character $\grx$ and $V_{\grf}$ is an irreducible $\CC \tilde{H_k}$-module with character $\grf$. 
This matrix records in which way the irreducible representations of
the semisimple algebra $\CC(\mathbf{v})\tilde{H_k}$ break up into irreducible representations 
of $\CC\tilde{ H_k}$. 

The form of the decomposition matrix is controlled by the \emph{Schur elements}, denoted as $s_{\grx}$, $\grx \in$ Irr$\big(\CC(\mathbf{v})\tilde{H_k}\big)$ with respect to the symmetric form $t_k$. The Schur elements belong to $R^k_K$ (see \cite{geck}, Proposition
7.3.9) and they depend only on the symmetrizing form $t_k$  and the isomorphism class of the representation. M. Chlouveraki has shown  that the Schur elements are products of cyclotomic polynomials over $K_k$ evaluated on monomials of degree 0 (see theorem 4.2.5 in \cite{chlouverakibook}). In the following section we are going to use these elements in order to determine the irreducible representations of  $\CC\tilde H_k$ (more details about the definition and the properties of the Schur elements, one may refer to \S 7.2 in \cite{geck}).

We say that the $\CC(\mathbf{v})\tilde{H_k}$-modules $V_{\grx}, V_{\grc}$ \emph{belong to the same block} if the corresponding characters $\grx, \grc$  label the rows of the same block in the decomposition matrix $(d_{\grx\grf})$ (by definition, this means that there is a $\grf\in\text{Irr}(\CC\tilde H_k)$ such that $d_{\grx,\grf}\not=0\not=d_{\grc,\grf}$).
If an irreducible $\CC(\mathbf{v})\tilde{H_k}$-module  is alone in its block, then we call it a \emph{module of defect 0}. Motivated by the idea of M. Chlouveraki and H. Miyachy in \cite{chlouveraki} \textsection 3.1 we use the following criteria in order to determine whether two modules belong to the same block:

\begin{itemize}
\item We have $\gru(s_{\grx})\not =0$ if and only if $V_{\grx}$ is a module of defect 0 (see \cite{maller}, Lemma 2.6). 

This criterium  together with theorem 7.5.11 in \cite{geck}, states that $V_{\grx}$ is a module of defect 0 if and only if the decomposition matrix is of the form
$$\begin{pmatrix} *& \dots&*&0 &*&\dots&*\\
	\vdots& \dots& \vdots&\vdots &\vdots& \dots&\vdots\\
	 *& \dots&*&0 &*&\dots&*\\
	 0& \dots& 0&1&0&\dots&0\\
	 *&  \dots&*&0 &*&\dots&*\\
	 	\vdots&  \dots&\vdots &\vdots& \vdots& \dots&\vdots\\
	 	*& \dots&*&0 &*&\dots&*
\end{pmatrix}$$

\item If $V_{\grx}, V_{\grc}$ are in  the same block, then $\gru(\grv_{\grx}(z_0))=\gru(\grv_{\grc}(z_0))$ (see \cite{geck}, Lemma 7.5.10), where $\grv_{\grx}, \grv_{\grc}$ are the corresponding \emph{central characters}\footnote{If $z$  lies in the center of $\CC(\mathbf{v})\tilde{H_k}$ then Schur's lemma implies that $z$ acts as  scalars in $V_{\grx}$ and $V_{\grc}$. We denote these scalars as $\grv_{\grx}(z)$ and $\grv_{\grc}(z)$ and we call the associated $\CC(\mathbf{v})$-homomorphisms $\grv_{\grx},\grv_{\grc}: Z\big(\CC(\mathbf{v})\tilde{H_k}\big)\rightarrow \CC(\mathbf{v})$ central characters (for more details, see \cite{geck} page 227).} and $z_0$ is the central element $(s_1s_2)^3$.
\end{itemize}
\subsection{The irreducible representations of $B_3$}
\indent

We recall that in order to describe the irreducible representations of $B_3$ of dimension at most 5, it is enough to describe the irreducible $\CC\tilde{H_k}$-modules of dimension $k$. Let $S$ be an irreducible $\CC \tilde{H_k}$-module of dimension $k$ and $s\in S$ with $s\not=0$. The morphism $f_{s}: \CC \tilde{H_k}\arw S$ defined by $h\mapsto hs$ is surjective since $S$ is irreducible. Hence, by the definition of the Grothendieck group we have that $d_{\gru}\big([\CC(\mathbf{v})\tilde{H_k}]\big)=[\CC \tilde{H_k}]=[$kerf$_{s}]+[S]$. However, since $\CC(\mathbf{v})\tilde{H_k}$ is semisimple we have  $\CC(\mathbf{v})H_k=M_1\oplus...\oplus M_r$, where the $M_i$ are (up to isomorphism) all the simple $\CC(\mathbf{v})\tilde{H_k}$-modules (with redundancies). Therefore, we have $\sum_{i=1}^{r}d_{\gru}([M_i])=[$kerf$_{s}]+[S].$ 
Hence, there is a simple  $\CC(\mathbf{v})\tilde{H_k}$-module $M$ such that 
\begin{equation}d_{\gru}([M])=[S]+[J],\label{eqqq}\end{equation} where $J$ is a $\CC \tilde{H_k}$-module.
\begin{rem}
\mbox{}  
\vspace*{-\parsep}
\vspace*{-\baselineskip}\\ 
\begin{itemize}
	\item[(i)] The $\CC(\mathbf{v})\tilde{H_k}$-module $M$ is of dimension at least $k$.
	\item[(ii)] If $J$ is of dimension 1, there is a $\CC(\mathbf{v})\tilde{H_k}$-module $N$ of dimension 1, such that $d_{\gru}([N])=[J]$. This result comes from the fact that the 1 dimensional $\CC \tilde H_k$ modules are of the form $(\grl_i)$ and, by definition, $\grl_i=\gru(u_i)$.
		\end{itemize}	
		\label{brrrr}
	\end{rem}
The irreducible  $\CC(\mathbf{v})\tilde{H_k}$-modules are known (see \cite{mallem} or \cite{brouem} \textsection 5B and \textsection 5D, for $n=3$ and $n=4$, respectively). Therefore, we
can determine $S$ by using (\ref{eqqq}) and a case-by-case analysis. 
\begin{itemize}

\item \underline{$k=2$} : Since $\tilde{H_2}$ is the generic Hecke algebra of $\mathfrak{S}_3$, which is a Coxeter group, the irreducible representations of $\CC \tilde{H_2}$ are well-known; we have two irreducible representations of dimension 1 and one of dimension 2. By $(\ref{eqqq})$ and remark \ref{brrrr} (i),  $M$ must be the irreducible $\CC(\mathbf{v})\tilde{H_k}$-module of dimension 2 and $(\ref{eqqq})$ becomes $[S]=d_{\gru}([M])$. Hence, we have: 
$$A=\begin{bmatrix}
\begin{array}{rr}
\grl_1&\grl_1\\
0&\grl_2
\end{array}
\end{bmatrix},\;
B=\begin{bmatrix}
\begin{array}{rr}
\grl_2&0\\
-\grl_2&\grl_1\\
\end{array}
\end{bmatrix}$$
Moreover, $[S]=d_{\gru}([M])$ is irreducible and $M$ is the only irreducible $\CC(\mathbf{v})\tilde{H_k}$-module of dimension 2. As a result, $M$ has to be alone in its block i.e. $\gru(s_{\grx})\not=0$, where $\grx$ is the character that corresponds to $M$. Therefore, an irreducible representation of $B_3$ of dimension 2 can be described by the explicit matrices $A$ and $B$ we have above, depending only on a choice of $\grl_1, \grl_2$ such that $\gru(s_{\grx})=\grl_1^2-\grl_1\grl_2+\grl_2^2\not=0$.

\item \underline{$k=3$} : Since the algebra $\CC(\mathbf{v})\tilde{H_3}$ is split semisimple, we have a bijection between the set 
 Irr$(\CC(\mathbf{v})\tilde{H_3})$ and the set Irr$(W_3)$, as we explained in the previous section.
We refer to J. Michel's version of CHEVIE package of GAP3 (see \cite{michelgap}) in order to find the irreducible characters of $W_3$. 
We type:
\begin{verbatim}
gap> W_3:=ComplexReflectionGroup(4);
gap> CharNames(W_3);
[ "phi{1,0}", "phi{1,4}", "phi{1,8}", "phi{2,5}", "phi{2,3}", "phi{2,1}",
  "phi{3,2}" ]
\end{verbatim}
We have 7 irreducible characters $\grf_{i,j}$, where $i$ is the dimension of the representation and $j$ the valuation of its fake degree (see \cite{malle1} \textsection 6A). Since $S$ is of dimension 3,  the equation $(\ref{eqqq})$ becomes $[S]=d_{\gru}([M])$, where $M$ is the irreducible $\CC(\mathbf{v})\tilde{H_3}$-module  that corresponds to the character $\grf_{3,2}$ (see remark \ref{brrrr}(i)). However, we have explicit matrix models for this representation (see \cite{brouem}, \textsection 5D or we can refer to CHEVIE package of GAP3 again) and since $[S]=d_{\gru}([M])$ we have:
$$A=\begin{bmatrix}
\grl_3&0&0\\
\grl_1\grl_3+\grl_2^2& \grl_2& 0\\
\grl_2& 1&\grl_1
\end{bmatrix},\;
B=\begin{bmatrix}
\grl_1& -1&\grl_2\\
0&\grl_2&-\grl_1\grl_3-\grl_2^2\\
0&0&\grl_3
\end{bmatrix}.$$
$M$ is the only  irreducible $\CC(\mathbf{v})\tilde{H_3}$-module of dimension 3, therefore, as in the previous case where $k=2$, we must have that $\gru(s_{\grf_{3,2}})\not=0$. The Schur element $s_{\grf_{3,2}}$ has been determined in \cite{malle2} and the condition $\gru(s_{\grf_{3,2}})\not=0$ becomes \begin{equation}\gru(s_{\grf_{3,2}})=\frac{(\grl_1^2+\grl_2\grl_3)(\grl_2^2+\grl_1\grl_3)(\grl_3^2+\grl_1\grl_2)}{(\grl_1\grl_2\grl_3)^2}\not=0.\label{tt1}\end{equation}
To sum up, an irreducible representation of $B_3$ of dimension 3 can be described by the explicit matrices $A$ and $B$ we gave above, depending only on a choice of $\grl_1, \grl_2, \grl_3$ such that (\ref{tt1}) is satisfied. 

\item \underline{$k=4$} : We use again the program GAP3 package CHEVIE in order to find the irreducible characters  of $W_4$. In this case we have 16 irreducible characters among which 2 of dimension 4; the characters $\grf_{4,5}$ and $\grf_{4,3}$ (we follow again the notations in GAP3, as in the case where $k=3$). Hence, by remark \ref{brrrr}(i) and relation $(\ref{eqqq})$, we have $[S]=d_{\gru}([M])$, where $M$ is the irreducible $\CC(\mathbf{v})\tilde{H_4}$-module  that corresponds either to the character $\grf_{4,5}$ or to the character $\grf_{4,3}$. We have again explicit matrix models for these representations (see \cite{brouem}, \textsection 5B, where we multiply the matrices described there by a scalar $t$ and we set $u_1=t, u_2=tu, u_3=tv$ and $u_4=tw$):
$$A=\begin{bmatrix}
\grl_1&0&0&0\\ \\
\frac{\grl_1^2}{\grl_2}&\grl_2& 0& 0\\\\
\frac{\grl_1^3}{r}&\frac{\grl_1\grl_2\grl_3-\grl_1r}{r}& \grl_3& 0\\\\
-\grl_2&\grl_2\gra&\frac{r\gra}{\grl_1^2}&\grl_4
\end{bmatrix},\;
B=\begin{bmatrix}
\grl_4&\grl_3\gra&\frac{\grl_2\grl_3\gra}{\grl_1}&-\frac{\grl_2\grl_3^2}{r}\\\\
0&\grl_3&\frac{\grl_2\grl_3-r}{\grl_1}&\frac{\grl_1^2\grl_3}{r}\\\\
0&0&\grl_2&\frac{\grl_1^3}{r}\\\\
0&0&0&\grl_1
\end{bmatrix},$$
where $r:=\pm\sqrt{\grl_1\grl_2\grl_3\grl_4}$ and $\gra:=\frac{r-\grl_2\grl_3-\grl_1\grl_4}{\grl_1^2}$.

Since $d_{\gru}([M])$ is irreducible either $M$ is of defect 0 or it is in the same block as the other irreducible module of dimension 4  i.e. $\gru(\grv_{\grf_{4,5}}(z_0))=\gru(\grv_{\grf_{4,3}}(z_0))$. We use the program GAP3 package CHEVIE in order to calculate  these central characters. More precisely, we have 16 representations where the last 2 are of dimension 4. These representations will be noted in GAP3 as $\verb+R[15]+$ and $\verb+R[16]+$. Since $z_0=(s_1s_2)^3$  we need to calculate the matrices $R[i](s_1s_2s_1s_2s_1s_2), i=15, 16$. These are the matrices $\verb+Product(R[15]{[1,2,1,2,1,2]})+$ and $\verb+Product(R[16]{[1,2,1,2,1,2]})+$, in GAP3 notation, as we can see below:
\begin{verbatim}
gap> R:=Representations(H_4);;
gap> Product(R[15]{[1,2,1,2,1,2]});
[ [ u_1^3/2u_2^3/2u_3^3/2u_4^3/2, 0, 0, 0 ], 
  [ 0, u_1^3/2u_2^3/2u_3^3/2u_4^3/2, 0, 0 ],
  [ 0, 0, u_1^3/2u_2^3/2u_3^3/2u_4^3/2, 0 ],
  [ 0, 0, 0, u_1^3/2u_2^3/2u_3^3/2u_4^3/2] ]
gap> Product(R[16]{[1,2,1,2,1,2]});
[ [ -u_1^3/2u_2^3/2u_3^3/2u_4^3/2, 0, 0, 0 ], 
  [ 0, -u_1^3/2u_2^3/2u_3^3/2u_4^3/2, 0, 0 ],
  [ 0, 0, -u_1^3/2u_2^3/2u_3^3/2u_4^3/2, 0 ],
  [ 0, 0, 0, -u_1^3/2u_2^3/2u_3^3/2u_4^3/2 ] ]
\end{verbatim}
We have that $\gru(\grv_{4,5}(z_0))=-\gru(\grv_{4,3}(z_0))$, which means that 
$M$ is of defect zero i.e. $\gru(s_{\grf_{4,i}})\not=0$, where $i=3$ or 5.  The Schur elements $s_{\grf_{4,i}}$ have been determined in \cite{malle2} \textsection 5.10, hence we must have
\begin{equation}\gru(s_{\grf_{4,i}})=\frac{-2r\prod\limits_{p=1}^4(r+\grl_p^2)\prod\limits_{r,l}(r+\grl_r\grl_l+\grl_s\grl_t)}{(\grl_1\grl_2\grl_3\grl_4)^4}\not=0,
\text {where } \{r,l,s,t\}=\{1,2,3,4\}
\label{tt2}
\end{equation}
Therefore, an irreducible representation of $B_3$ of dimension 4 can be described by the explicit matrices $A$ and $B$ depending only on a choice of $\grl_1, \grl_2, \grl_3, \grl_4$ and a square root of $\grl_1\grl_2\grl_3\grl_4$ such that (\ref{tt2}) is satisfied.

\item \underline{$k=5$} :In this case, compared to the previous ones, we have two possibilities for $S$. The reason is that we have characters of dimension 5 and  dimension 6, as well. Therefore, by remark \ref{brrrr}(i) and (ii) and (\ref{eqqq})  we either have $d_{\gru}([M])=[S]$, where $M$ is one  irreducible $\CC(\mathbf{v})\tilde{H_5}$-module of dimension 5 or  $d_{\gru}([N])=[S]+d_{\gru}([N'])$, where $N, N'$ are some irreducible $\CC(\mathbf{v})\tilde{H_5}$-modules of dimension 6 and 1, respectively.

In order to exclude the latter case, it is enough to show that $N$ and $N'$ are not in the same block. Therefore, at this point, we may assume that $\gru(\grv_{\grx}(z_0))\not=\gru(\grv_{\grc}(z_0))$, for every irreducible character $\grx, \grc$ of $W_5$ of dimension 6 and 1, respectively. We use GAP3 in order to calculate the central characters, as we did in the case where $k=4$ and we have:
$\gru(\grv_{\grc}(z_0))=\grl_i^6$, $i\in \{1,...,5\}$ and  $\gru(\grv_{\grx}(z_0))=-x^2yztw$, where $\{x,y,z,t,w\}=\{\grl_1, \grl_2, \grl_3, \grl_4, \grl_5\}$. We notice that $\gru(\grv_{\grx}(z_0))=-\grl_j$det$A$, $j\in \{1,...,5\}$.  Therefore, the assumption  $\gru(\grv_{\grx}(z_0))\not=\gru(\grv_{\grc}(z_0))$ becomes det$A\not=-\grl_i^6\grl_j^{-1}$, $i,j \in\{1,2,3,4,5\}$, where $i, j$ are not necessarily distinct.

By this assumption we have that $d_{\gru}([M])=[S]$, where $M$ is some irreducible $\CC(\mathbf{v})\tilde{H_5}$-module of dimension 5. We have again explicit matrix models for these representation (see \cite{mallem} or the CHEVIE package of GAP3))
We notice that these matrices depend only on the choice of eigenvalues and of a fifth root of det$A$.

Since $d_{\gru}([M])$ is irreducible either $M$ is of defect 0 or it is in the same block with another irreducible module of dimension 5. However, since the central characters of the irreducible modules of dimension 5 are distinct fifth roots of  $(u_1u_2u_3u_4u_5)^{6}$, we can exclude the latter case.
Hence, $M$ is of defect zero i.e. $\gru(s_{\grf})\not=0$, where $\grf$ is the character that corresponds to $M$. The Schur elements  have been determined in \cite{malle2} (see also Appendix A.3 in \cite{chlouverakibook}) and one can also find them in CHEVIE package of GAP3; they are 
$$\frac{5\prod\limits_{i=1}^5(r+u_i)(r-\grz_3u_i)(r-\grz_3^2u_i)\prod\limits_{i\not=j}(r^2+u_iu_j)}{(u_1u_2u_3u_4u_5)^7},$$
where $r$ is a 5th root of $u_1u_2u_3u_4u_5$. However, due to the assumption det$A\not=-\grl_i^5$, $i \in\{1,2,3,4,5\}$ (case where $i=j$), we have that $\gru(r)+\grl_i \not =0$. Therefore, the condition $\gru(s_{\grf})\not=0$ becomes
\begin{equation}\prod\limits_{i=1}^5(\tilde{r}^2+\grl_i\tilde{r}+\grl_i^2)\prod\limits_{i\not=j}(\tilde{r}^2+\grl_i\grl_j)\not=0,
\label{oua}
\end{equation}
where $\tilde{r}$ is a fifth root of det$A$. 

Therefore, an irreducible representation of $B_3$ of dimension 5 can be described by the explicit matrices $A$ and $B$, that one can find for example in CHEVIE package of GAP3, depending only on a choice of $\grl_1, \grl_2, \grl_3, \grl_4, \grl_5$ and a fifth root of det$A$ such that (\ref{oua}) is satisfied.
\end{itemize}
\begin{rem}
\begin{enumerate}
\mbox{}  
\vspace*{-\parsep}
\vspace*{-\baselineskip}\\
\item We can generalize our results for a representation of $B_3$ over a field of positive characteristic, using similar arguments. However, the cases where $k=4$ and $k=5$ need some extra analysis; For $k=4$ we have two irreducible $\CC(\mathbf{v})\tilde{H_4}$-modules of dimension 4, which are not in the same block if we are in any characteristic but 2. However, when we are in characteristic 2, these two modules coincide and, therefore, we obtain an irreducible module of $B_3$ which is of defect 0, hence we arrive to the same result as in characteristic 0. We have exactly the same argument for the case where $k=5$ and we are over a field of characteristic 5.

\item The irreducible representations of $B_3$ of dimension at most 5 have been classified in \cite{tuba}. Using a new framework, we arrived to the same results. The matrices $A$ and $B$ described by Tuba and Wenzl are the same (up to equivalence) with the matrices we provide in this paper. For example, in the case where $k=3$, we have given explicit matrices $A$ and $B$. If we take the matrices  $DAD^{-1}$ and $DBD^{-1}$, where $D$ is the invertible matrix $$D=\begin{bmatrix}
-\grl_1\grl_2-\grl_3^2&\grl_1(\grl_3-\grl_1)&(\grl_2-\grl_3)(\grl_3-\grl_1)\\
(\grl_2-\grl_1)(\grl_3^2+\grl_1\grl_2)& \grl_1(2\grl_2\grl_1-\grl_1^2+2\grl_1\grl_3-\grl_3\grl_2)& (\grl_1-\grl_3)(\grl_2^2+\grl_1\grl_3)\\
0& \grl_1(\grl_1-\grl_3)&-\grl_3\grl_1(\grl_1+\grl_2)
\end{bmatrix},$$ we just obtain the matrices determined in \cite{tuba} (The matrix $D$ is invertible since det$D=\grl_1(\grl_1^2+\grl_2\grl_3)(\grl_3^2+\grl_1\grl_2)^2\not=0$, due to (\ref{eqqq})).
\end{enumerate}
\end{rem}

	\newpage
	~
\chapter{The approach of Etingof-Rains: The exceptional groups of rank 2}

In this chapter we explain in detail the arguments Etingof and Rains used in order to prove a weak version of the BMR freeness conjecture for the exceptional groups of rank 2. This weak version states that the  generic Hecke algebra $H$ associated to an exceptional group of rank 2 is finitely generated as $R$-module (where $R$ is the Laurent polynomial ring over which we define $H$). Their approach also explains the appearance of the center of $B_3$ in the description of the basis of the generic Hecke algebra associated to $G_4$, $G_8$ and $G_{16}$ (see theorem 3.2 (3) in \cite{marincubic} and theorems \ref{th} and \ref{thh2}).

\section{The exceptional groups of rank 2}
\indent

In this section we follow mostly Chapter 6 of \cite{lehrer}. Let $W$ be an irreducible complex reflection group of rank 2, meaning that $W$ is one of the groups $G_4,$\dots, $G_{22}$. We know that $W\leq U_2(\CC)$, where $U_2(\CC)$ denotes the unitary group of degree 2 (see remark \ref{unitary}). For every $w\in W$ we choose $\grl_w\in \CC$ such that $\grl_w^2=\text{det}(w)$ and we set
 $$\widehat W:=\{\pm \grl_w^{-1}w;w\in W\}.$$
 
By definition,  $\widehat W\leq SU_2(\CC)$, where  $SU_2(\CC)$ denotes the special unitary group of degree 2. Therefore, $\widehat W/Z(\widehat W)\leq SU_2(\CC)/Z(SU_2(\CC))\simeq S^3/\{\pm1\}$, where $S^3$ denotes the group of quaternions of norm 1. Since $S^3/\{\pm1\}$ is isomorphic to the three dimensional rotation group $SO_3(\RR)$, we can consider $\widehat W/Z(\widehat W)$ as a subgroup of the latter. We can also assume that $W/Z(W)\leq SO_3(\RR)$, since by construction  $W/Z(W)\simeq \widehat W/Z(\widehat W)$. 

Any finite subgroup of $SO_3(\RR)$ is either a cyclic group, a dihedral group $D_n$ or the group of rotations of a Platonic solid, which is the tetrahedral, octahedral or  icosahedral group (for the classification of the finite subgroups of $SO_3(\RR)$ one may refer to theorem 5.13 of \cite{lehrer}). Since $W$ is irreducible we may exclude the case of the cyclic group. The case of the dihedral group falls into the case of the infinite family $G(de,e,2)$. Hence, it remains to examine the case of the tetrahedral, octahedral and  icosahedral group. 

Using the Shephard-Todd notation we have the following three families (see table \ref{t1} below), according to whether $W/Z(W)$ is the tetrahedral, octahedral or icosahedral group (for more details one may refer to Chapter 6 of \cite{lehrer}); the first family, known as \emph{the tetrahedral family}, includes the groups $G_4,\dots, G_7$, the second one, known as the \emph{octahedral family} includes the groups $G_8,\dots, G_{15}$ and the last one, known as the \emph{icosahedral family}, includes the rest of them, which are the groups $G_{16},\dots, G_{22}$. In each family, there is a maximal group of order $|W/Z(W)|^2$ and all the other groups are its subgroups. These groups belonging to the three families are known as  \emph{the exceptional groups of rank 2}.
\begin{table}[!ht]
	\begin{center}
		\small
		\caption{
			\bf{The three families}}
			\label{t1}
		\scalebox{0.98}
		{\begin{tabular}{|c|c|c|}
				\hline
			$W/Z(W)$& $W$& Maximal Group \\ 
				\hline
			 $\begin{array}[t]{lcl} \\\text{Tetrahedral Group }
				\mathcal{T}\simeq Alt(4)\\ \\
				\langle a,b,c\;|\;a^2=b^3=c^3=1, abc=1\rangle\end{array}$
				&
			$\begin{array}[t]{lcl}\\G_4,\dots,G_7\end{array}$
				&
				$\begin{array}[t]{lcl}\\ G_7=
				\langle a,b,c\;|\;a^2=b^3=c^3=1, abc=\text{central }\rangle
				\end{array}$\\
				\hline \hline
					
					$\begin{array}[t]{lcl} \\\text{Octahedral Group }
					\mathcal{O}\simeq Sym(4)\\\\
					\langle a,b,c\;|\;a^2=b^3=c^4=1, abc=1\rangle\end{array}$
					&$\begin{array}[t]{lcl}\\
					G_8,\dots,G_{15}
					\end{array}$&
					$\begin{array}[t]{lcl}\\ G_{11}=
					\langle a,b,c\;|\;a^2=b^3=c^4=1, abc=\text{central }\rangle
					\end{array}$\\
					\hline \hline
						
						$\begin{array}[t]{lcl} \\\text{Icosahedral Group }
						\mathcal{I}\simeq Alt(5)\\\\
						\langle a,b,c\;|\;a^2=b^3=c^5=1, abc=1\rangle\end{array}$
						&$\begin{array}[t]{lcl}\\
						G_{16},\dots,G_{22}
						\end{array}$&
						$\begin{array}[t]{lcl}\\ G_{19}=
						\langle a,b,c\;|\;a^2=b^3=c^5=1, abc=\text{central }\rangle
						\end{array}$\\
						\hline
				\end{tabular}}
			\end{center}
			
			\end{table}

We know that for every exceptional group of rank 2 we have  a Coxeter-like presentation (see remark \ref{coxeterlike}(ii)); that is a presentation of the form 
$$\langle s\in S\;|\; \{v_i=w_i\}_{i\in I} , \{s^{e_s}=1\}_{s\in S}\rangle,$$
where $S$ is a finite set of distinguished  reflections and $I$ is a finite set of relations such that, for each $i\in I$,  $v_i$ and $w_i$ are positive words with the same length in elements of $S$. We also know that for the associated complex braid group $B$ we have an Artin-like presentation (see theorem \ref{Presentt}); that is a presentation of the form 
$$\langle \mathbf{s}\in \mathbf{S}\;|\; \mathbf{v_i}=\mathbf{w_i} \rangle_{i\in I},$$
where $\mathbf{S}$ is a finite set of distinguished braided reflections and $I$ is a finite set of relations such that, for each $i\in I$,  $\mathbf{v_i}$ and $\mathbf{w_i}$ are positive words in elements of $\mathbf{S}$.
We call these presentations \emph{the BMR presentations}, due to M. Brou\'e, G. Malle and R. Rouquier.

In 2006 P. Etingof and E. Rains gave different presentations of $W$ and $B$, based on the BMR presentations associated to the maximal groups $G_7$, $G_{11}$ and $G_{19}$ (see \textsection 6.1 of \cite{ERrank2}). We call these presentations \emph{the ER presentations}.
In tables \ref{t3} and \ref{t2} we give the two representations for every $W$ and $B$, as well as the isomorphisms between the BMR and ER presentations. In Appendix A we also prove that we  indeed have such isomorphisms. Notice that for the maximal groups, the ER presentations coincide with the BMR presentations. 

In the following section we will explain how Etingof and Rains used the ER presentation in order to prove a weak version of the BMR conjecture for the exceptional groups of rank 2.

\section{A weak version of the BMR freeness conjecture}

\indent

Let $W$ be an exceptional group of rank 2 with associated complex braid group $B$ and let $H$ denote the generic Hecke algebra associated to $W$, defined over the  ring  $R=\ZZ\left[u_{s,i}^{\pm}\right]$, where $s$ runs over the conjugacy classes of distinguished reflections and $1\leq i \leq e_s$, where $e_s$ denotes the order of the pseudo-reflection $s$ in $W$. For the rest of this section we follow the notations of \textsection 2.2 of \cite{marinG26}. 

For $k$ a unitary ring we set 
$R_{ k}:=R\otimes_{\ZZ}\calligra{k}$ and $H_{k}:=H\otimes_{R}R_{ k}$. We denote by $\tilde u_{s,i}$ the images of $u_{s,i}$ inside $R_k$.
By definition, $H_{k}$ is the quotient of the group algebra $R_{k}B$ of $B$ by the ideal generated by  $P_{s}(\grs)$, where $s$ runs over the conjugacy classes of distinguished reflections, $\grs$ over the set of distinguished braided reflections associated to $s$ and $P_{s}\left[X\right]$ are the monic polynomials $(X-\tilde u_{s,1})\dots(X-\tilde u_{s,e_s})$ inside $R_{k}\left[X\right]$. Notice that if $s$ and $t$ are conjugate in $W$ the polynomials $P_s(X)$ and $P_t(X)$ coincide.

Let $Z(B)$ the center of $B$. By theorem \ref{centerbraid} we know that $Z(B)$ is cyclic and, therefore, we set $Z(B):=\langle z\rangle$. We also set $\bar B:=B/\langle z\rangle$ and $R_{ k}^+:=R_{ k}\left[x,x^{-1}\right]$. Let $f$ be a set-theoretic section of the natural projection $\grp: B\arw \bar B$, meaning that $f: \bar B \arw B$  is a map such that $\grp\circ f=id_{\bar B}$. The next proposition (Proposition 2.10 in \cite{marinG26}) states that $H_{k}$ inherits a structure of $R_{k}^+$-module. More precisely, $H_{ k}$ is isomorphic to the quotient of the group algebra $R_{ k}^+\bar B$ of $\bar B$ by some relations defined by the polynomials $P_{\mathbf{s}}\left[X\right]$ and by $f$. 
\begin{prop}
	\mbox{}  
	\vspace*{-\parsep}
	\vspace*{-\baselineskip}\\ 
	\begin{itemize}
		\item[(i)] $R_{ k}B$ admits a $R_{ k}^+$-module structure defined by $x\cdot b=zb$, for every $b\in B$. Moreover, as $R_{k}^+$-modules, $R_{ k}B\simeq R_{k}^+\bar B$.
		\item [(ii)] Under the isomorphism described in (i), the defining ideal of $H_{k}$ generated by the $P_s(\grs)$ is mapped to the ideal of $ R_{ k}^+\bar B$ generated by the $Q_s(\bar \grs)$, where $Q_s(X)=x^{c_{\grs}\text{deg}P_s}\cdot
		P_s(x^{-c_{\grs}}\cdot X)\in R_{ k}^+[X]$, the $c_{\grs}\in \ZZ$ being defined by $f(\bar\grs)=z^{c_{\grs}}\grs$.
			
	\end{itemize}
	\label{prIVAN}

\end{prop}

\begin{proof}
	This proof is a rewriting of the proof of Proposition 2.10 in \cite{marinG26}. 
	For every $b\in B$ we have
	$\grp(f(\bar b)) =\bar b=\grp(b)$ $\an f(\bar b)b^{-1}\in \text{ker}\grp=Z(B)$. We also have that $Z(B)\simeq \ZZ$. Therefore, there is a uniquely defined 1-cocycle $c: B \arw \ZZ$, $b\mapsto c_b$ such that $\forall b\in B$, $f(\bar b)=z^{c_{b}}b$.
	We have:
		$$\begin{array}[t]{lcl}
		R_kB&=&\underset{b\in B}{\oplus}R_kb\\
		&=&\underset{b\in B}{\oplus}R_kz^{-c_b}f(\bar b)\\
		&=&\underset{d\in \bar B}{\oplus}\underset{\bar b=d}{\oplus}R_kz^{-c_b}f(d)\\
		&=&\underset{d\in \bar B}{\oplus}\underset{n\in \ZZ}{\oplus}R_kz^nf(d)\\
		&=&\underset{d\in \bar B}{\oplus}R_k^{+}f(d).
		\end{array}$$
	However, $\underset{d\in \bar B}{\oplus}R_k^{+}f(d)\simeq \underset{d\in \bar B}{\oplus}R_k^{+}d=R_k^{+}\bar B$ 	and, therefore, 	$R_kB\simeq R_k^{+}\bar B$, which proves (i).
	
	For (ii), let $I$ be the $R_k$-submodule of $R_kB$ spanned by the $bP_{s}(\grs)b'$, for $b,b' \in B$. We notice that $bP_{s}(\grs)b'=z^{-c_b}f(\bar b)
		P_{s}(\grs)z^{-c_{b'}}f(\bar b')$. Hence, as $R_k^{+}$-module $I$ is spanned by the $f(\bar b)
		P_{s}(\grs)f(\bar b')$. Let $Q_s(X)=x^{c_s\text{deg}P_s}
		P_s(x^{-c_s}X)$. We have:
		
		$$Q_s(f(\bar \grs))=x^{c_s\text{deg}P_s}\cdot P_s(x^{-c_s}\cdot f(\bar \grs))=z^{c_s\text{deg}P_s}P_s(z^{-c_s}f(\bar \grs))= z^{c_s\text{deg}P_s}P_s(\grs).$$
		Therefore, $P_s(\grs)=z^{-c_s\text{deg}P_s}Q_s(f(\bar \grs))$.
		However, $I$ is identified inside $R_k^+\bar B$ by the $R_k^+$-submodule generated by the $f(\bar b)
		P_{s}(\grs)f(\bar b')$. As a result, $I$ is identified inside $R_k^+\bar B$ by the ideal of $R_k^+\bar B$ generated by $f(\bar b)Q_s(f(\bar \grs))f(\bar b')$, which is isomorphic to the ideal of $R_k^{+}\bar B$ generated by $Q_s(\bar \grs)$.
		
\end{proof}
\begin{cor}
	Let $f: \bar B \arw B$ be a set-theoretic section of the natural projection $B\arw \bar B$. There is an isomorphism $\grF_{f}$ between the $R_k^+$-modules  $R_k^+\bar B/Q_s(\bar \grs)$ and  $H_k=R_kB/P_s(\grs)$.
	\label{corIVAN}
\end{cor}
We now give an example to make all the above notations clearer to the reader.
	\begin{ex} Let $W:=G_4=\langle s,t\;|\;s^3=t^3=1, sts=tst\rangle$. Since $s,t$ are conjugate the generic Hecke algebra of $G_4$ is defined over the ring  $R=\ZZ[u_{s,i}^{\pm}]$, where $i=1,2,3$ and it has a presentation as follows:
	 $$H=\langle \grs,\grt\;|\; \grs\grt\grs=\grt\grs\grt,\;\prod\limits_{i=1}^{3}(\grs-u_{s,i})=\prod\limits_{i=1}^{3}(\grt-u_{s,i})=0\rangle.$$
	 Keeping the above notations, we let $P_s(X)=P_t(X)=(X-\tilde u_{s,1})(X-\tilde u_{s,2})(X-\tilde u_{s,3})$.
	We fix a set-theoretic section  $f: \bar B \arw B$.   We set $Q_s(X)=(X-x^{c_{\grs}}\cdot\tilde u_{s,1})(X-x^{c_{\grs}}\cdot\tilde u_{s,2})(X-x^{c_{\grs}}\cdot\tilde u_{s,3})$ and $Q_t(X)=(X-x^{c_{\grt}}\cdot\tilde u_{s,1})(X-x^{c_{\grt}}\cdot\tilde u_{s,2})(X-x^{c_{\grt}}\cdot\tilde u_{s,3})$, where $c_s, c_t \in \ZZ$, defined by  $f(\bar \grs)=z^{a_{\grs}}\grs$ and $f(\bar \grt)=z^{c_{\grt}}\grt$.
	The corollary \ref{corIVAN} states that, as $R_{\calligra k}^+$-modules
	$$H_{ k}\simeq \langle \bar \grs,\bar \grt\;|\;(\bar{\grs}\bar{\grt})^3=1,\; \bar{\grs}\bar{\grt}\bar{\grs}=\bar{\grt}\bar{\grs}\bar{\grt},\;\prod\limits_{i=1}^{3}(\bar \grs-x^{c_{\grs}}\cdot\tilde u_{s,i})=\prod\limits_{i=1}^{3}(\bar \grt-x^{c_{\grt}}\cdot\tilde u_{s,i})=0\rangle. \makeatletter\displaymath@qed$$
\end{ex}
We set $\overline{W}:=W/Z(W)$. We note that $\overline{W}$ is the group of even elements in a finite Coxeter group $C$ of rank 3 (of type $A_3$, $B_3$ and $H_3$ for the tetrahedral, octahedral and icosahedral cases, respectively), with Coxeter system $y_1,y_2,y_3$ and Coxeter matrix $(m_{ij})$. We set $\tilde \ZZ:=\ZZ\left[e^{\frac{2\pi i}{m_{ij}}}\right]$ and let $A(C)$ be the $\tilde\ZZ$-algebra presented as follows:
\begin{itemize}[leftmargin=*]
	\item \underline{Generators}: $Y_1, Y_2, Y_3$, $t_{ij,k}$, where $i,j\in\{1,2,3\}$, $i\not=j$ and $k\in\ZZ/m_{ij}\ZZ$.
	\item \underline{Relations}: $Y_i^2=1$, $t_{ij,k}^{-1}=t_{ji,-k}$, $\prod\limits_{k=1}^{m_{ij}}(Y_iY_j-t_{ij,k})=0$, $t_{ij,k}Y_r=Y_rt_{ij,k}$, $t_{ij,k}t_{i'j',k'}=t_{i'j',k'}t_{ij,k}$.
\end{itemize}
This construction of $A(C)$ is more general and can be done for every Coxeter group $C$ (for more details one may refer to \textsection 2 of \cite{ERcoxeter}).
Let $R^C=\tilde \ZZ\left[t_{ij,k}^{\pm}\right]=\tilde \ZZ\left[t_{ij,k}\right]$. The subalgebra $A_+(C)$ generated by $Y_iY_j$, $i\not=j$ becomes an $R^C$ algebra and can be presented as follows:
\begin{itemize}[leftmargin=*]
	\item \underline{Generators}: $A_{ij}:=Y_iY_j$, where $i,j\in\{1,2,3\}$, $i\not=j$.
	\item \underline{Relations}: $A_{ij}^{-1}=A_{ji}$, $\prod\limits_{k=1}^{m_{ij}}(A_{ij}-t_{ij,k})=0$, $A_{ij}A_{jl}A_{li}=1$, for $\#\{i,j,l\}=3$.
\end{itemize}
If $w$ is a word in letters $y_i$ we let $T_w$ denote the corresponding element of $A(C)$. For every $x\in \overline{W}$ let us choose a reduced word $w_x$ that represents $x$ in $\overline{W}$. We notice that $T_{w_x}$ is an element in $A_+(C)$, since $w_x$ is reduced and $\overline{W}$ is the group of even elements in $C$.
The following theorem is theorem 2.3(ii) in \cite{ERcoxeter}, which is generalized there for every finite Coxeter group.
\begin{thm}
 The algebra $A_+(C)$ is generated as $R^C$-module by the elements $T_{w_x}$, $x\in \overline{W}$.
 \label{thmER}
 \end{thm}
We recall that $C$ is a Coxeter group of type $A_3$, $B_3$ or $H_3$ respectively for each family. Therefore, the Coxeter matrix is of the form $\left( \begin{array}{ccc}
 1 & m & 2 \\
 m & 1 & 3 \\
 2 & 3 & 1 \end{array} \right)$,
 where $m=3$ for the tetrahedral family,  $m=4$ for the octahedral family and $m=5$ for the icosahedral family, respectively. Hence,
 we can present $A_+(C)$ as follows: $$\left\langle
 	\begin{array}{l|cl}
 	&(A_{13}-t_{13,1})(A_{13}-t_{13,2})=0&\\
 	A_{13}, A_{32}, A_{21}&(A_{32}-t_{32,1})(A_{32}-t_{32,2})(A_{32}-t_{32,3})=0,& A_{13}A_{32}A_{21}=1\\
 	&(A_{21}-t_{21,1})(A_{21}-t_{21,2})\dots(A_{21}-t_{21,m})=0&
 	\end{array}
 \right\rangle$$
 
 For every exceptional group of rank 2 we consider the ER presentation of  $B$ (see table \ref{t2}). We notice that apart from the cases of $G_{13}$ and $G_{15}$, this is a presentation in 3 generators $\gra$, $\grb$ and $\grg$ and for the cases of $G_{13}$ and $G_{15}$ this is a presentation of 4 generators $\gra$, $\grb$, $\grg$ and $\grd$. We also notice that in every case the group $\bar B$ can be presented as follows:
 $$\langle\bar \gra, \bar \grb, \bar \grg\;|\;\bar \gra ^{k_{\gra}}=\bar \grb ^{k_{\grb}}=\bar \grg ^{k_{\grg}}=1,\; \bar \gra \bar \grb \bar \grg=1	\rangle,$$
 where $k_{\gra}\in\{0,2\}$, $k_{\grb}\in\{0,3\}$ and the values of $k_{\grg}$ depend on the family in which the group belongs; for the tetrahedral family $k_{\grg}=0$, for the octahedral family $k_{\grg}\in\{0,4\}$ and for the icosahedral family $k_{\grg}\in\{0,5\}$.
 
 In the next two propositions, which are actually a rewriting with correction of proposition 2.13 in \cite{marinG26}, we relate the algebra $A_+(C)$ with the algebra $H_{\tilde \ZZ}$.
 \begin{prop}
 Let $W$ be an exceptional group of rank 2, apart from 	$G_{13}$ and $G_{15}$. There is a ring morphism $\gru: R^C\twoheadrightarrow R_{\tilde \ZZ}^+$ inducing  $ \grC: A_+(C)\otimes_{\gru} R_{\tilde \ZZ}^+ \twoheadrightarrow R_{\tilde \ZZ}^+\bar B/Q_{s}(\bar \grs)$ through $A_{13}\mapsto \bar \gra$, $A_{32}\mapsto \bar \grb$, $A_{21}\mapsto \bar \grg$.
 	\label{ERSUR}
 \end{prop}
\begin{proof}
We define $\gru$ by distinguishing the following cases:
\begin{itemize}[leftmargin=*]
\item[-]When $k_{\gra}=2$ (cases of $G_4$, $G_5$, $G_8$, $G_{10}$, $G_{16}$, $G_{18}$, $G_{20}$) we have $(\bar\gra-1)(\bar\gra+1)=0$. Therefore, we may define $\gru(t_{13,1}):=1$ and $\gru(t_{13,2}):=-1$.

\item[-] When $k_{\gra}=0$ (cases of $G_6$, $G_7$, $G_9$, $G_{11}$, $G_{12}$, $G_{14}$, $G_{17}$, $G_{19}$, $G_{21}$, $G_{22}$) we notice that  $\gra$  is  a distinguished braided reflection associated to a distinguished reflection $a$ of order 2 (see in tables \ref{t3} and \ref{t2} the images of $a$ and $\gra$ in BMR presentations). Keeping the  notations of proposition \ref{prIVAN}, we recall that $\bar \gra$ annihilates the polynomial $Q_{a}(X)=(X-x^{c_{\gra}}\cdot \tilde u_{a,1})(X-x^{c_{\gra}} \cdot\tilde u_{a,2})$. Therefore, we may define $\gru(t_{13,1}):=x^{c_{\gra}} \cdot\tilde u_{a,1}$ and $\gru(t_{13,2}):=x^{c_{\gra}}\cdot \tilde u_{a,2}$.

\item[-]When $k_{\grb}=3$ (cases of $G_4$, $G_6$, $G_8$, $G_9$, $G_{12}$, $G_{16}$, $G_{17}$, $G_{22}$) we have $(\bar\grb-1)(\bar\grb-j)(\bar\grb-j^2)=0$,  where $j$ is a third root of unity. Therefore, we define $\gru(t_{32,1}):=1$, $\gru(t_{32,2}):=j$
	and $\gru(t_{32,3}):=j^2$.
	
	\item[-] When $k_{\grb}=0$ (cases of $G_5$, $G_7$, $G_{10}$, $G_{11}$, $G_{14}$, $G_{18}$, $G_{19}$, $G_{20}$, $G_{21}$) we notice that $\grb$ is a distinguished braided reflection associated to a distinguished reflection $b$ of order 3 (see in tables \ref{t3} and \ref{t2} the images of $b$ and $\grb$ in BMR presentations). Therefore, similarly to the case where $k_{\gra}=0$, we may define $\gru(t_{32,1}):=x^{c_{\grb}} \cdot\tilde u_{b,1}$, $\gru(t_{32,2}):=x^{c_{\grb}} \cdot\tilde u_{b,2}$ and $\gru(t_{32,3}):=x^{c_{\grb}}\cdot \tilde u_{b,3}$.
	\item[-] When $k_{\grg}=m$ (cases of $G_{12}$, $G_{14}$, $G_{20}$, $G_{21}$, $G_{22}$) we have $(\bar\grg-1)(\bar\grg-\grz_m)(\bar\grg-\grz_m^2)\dots(\bar\grg-\grz_m^{m-1})=0$, where $\grz_m$ is a $m$-th root of unity. Therefore, we define 
$\gru(t_{21,1}):=1$, $\gru(t_{21,2}):=\grz_m$, $\gru(t_{21,3}):=\grz_m^2$,\dots,
	$\gru(t_{21,m}):=\grz_m^{m-1}$.
	\item[-] When $k_{\grg}=0$ (cases of $G_4$, $G_5$, $G_6$, $G_7$, $G_8$, $G_9$, $G_{10}$, $G_{11}$, $G_{16}$, $G_{17}$, $G_{18}$, $G_{19}$) we notice that $\grg$ is a distinguished braided reflection  associated to a distinguished reflection $c$ of order $m$ (see in tables \ref{t3} and \ref{t2} the images of $\grg$ and $c$ in BMR presentations). Therefore, similarly to the case where $k_{\gra}=0$, we may define $\gru(t_{21,1}):=x^{c_{\grg}}\cdot \tilde u_{c,1}$, $\gru(t_{21,2}):=x^{c_{\grg}}\cdot\tilde u_{c,2}$, \dots, $\gru(t_{21,m}):=x^{c_{\grg}} \cdot\tilde u_{c,m}$.
\end{itemize}
We now consider the algebra $A_+(C)\otimes_{\gru}R_{\tilde \ZZ}^+$ and we define 
$\begin{array}[t]{rcl}\grC:A_+(C)\otimes_{\gru} R_{\tilde \ZZ}^+ &\twoheadrightarrow& R_{\tilde \ZZ}^+\bar B/Q_s(\bar s)\\
A_{13}&\mapsto& \bar \gra\\
A_{32}&\mapsto& \bar \grb\\
A_{21}&\mapsto& \bar \grg.
\end{array}$

\end{proof} 

We will now deal with the cases of $G_{13}$ and $G_{15}$. For these cases, we  replace $A_+(C)$ with a specialized algebra $\tilde A_+(C)$ and we use a similar technique.
More precisely, both $G_{13}$ and $G_{15}$ belong to the octahedral family. We consider the ER presentation of the complex braid group associated to these groups (see table \ref{t2}) and we notice that $\bar B$ can be presented as follows:
$$\langle\bar \gra, \bar \grb, \bar \grg\;|\;\bar \grb^{k_{\grb}}=\bar \grg^{4}=1, \bar \gra \bar \grb \bar \grg=1	\rangle,$$
where $k_{\grb}\in\{0,3\}$. 
We set $\tilde R^C:=\tilde \ZZ\left[t_{13,1}, t_{13,2}, t_{32,1},  t_{32,2},  t_{32,3}, \sqrt{t_{21,1}},\sqrt{t_{21,3}},\right]$
and we define  $$ \begin{array}[t]{rlc}\grf: R^C &\longrightarrow\tilde R^C&\\
t_{21,1}&\mapsto \;\;\sqrt{t_{21,1}}&\\
t_{21,2}&\mapsto -\sqrt{t_{21,1}}&\\
t_{21,3}&\mapsto \;\;\sqrt{t_{21,3}}&\\
t_{21,4}&\mapsto -\sqrt{t_{21,3}}&.
\end{array}$$
Let $\tilde A_+(C)$ denote the  $\tilde R^C$ algebra $A_+(C)\otimes_{\grf}\tilde R^C$ and let $\tilde A_{ij}$ denote the image of $A_{ij}$ inside the algebra $\tilde A_+(C)$. The latter can be presented as follows:
$$\left\langle
	\begin{array}{l|cl}
	&(\tilde A_{13}-t_{13,1})(\tilde A_{13}-t_{13,2})=0&\\
	\tilde A_{13}, \tilde A_{32}, \tilde A_{21}&(\tilde A_{32}-t_{32,1})(\tilde A_{32}-t_{32,2})(\tilde A_{32}-t_{32,3})=0,& \tilde A_{13}\tilde A_{32}\tilde A_{21}=1\\
	&(\tilde A_{21}^2-t_{21,1})(\tilde A_{21}^2-t_{21,3})=0&
	\end{array}
\right\rangle$$

\begin{prop}
	Let $W$ be the exceptional group $G_{13}$ or $G_{15}$. There is a ring morphism $\gru: \tilde R^C\twoheadrightarrow R_{\tilde \ZZ}^+$ inducing  $ \grC: \tilde  A_+(C)\otimes_{\gru} R_{\tilde \ZZ}^+ \twoheadrightarrow R_{\tilde \ZZ}^+\bar B/Q_{s}(\bar s)$ through $\tilde A_{13}\mapsto \bar \gra$, $\tilde A_{32}\mapsto \bar \grb$, $\tilde A_{21}\mapsto \bar \grg$.
	\label{ERRSUR}
\end{prop}
\begin{proof}
We use the same arguments we used for the rest of the exceptional groups of rank 2 (see the proof of proposition \ref{ERSUR}). More precisely, since $(\bar{\grg}^2-1)(\bar{\grg}^2+1)=0$, we may define $\gru(t_{21,1})=1$ and $\gru(t_{21,3})=-1$. Moreover, we notice that $\gra$ is a distinguished braided reflection associated to a distinguished reflection a of order 2 (see in tables \ref{t3} and \ref{t2} the images of $\gra$ and $a$ in BMR presentations). 
Keeping the  notations of proposition \ref{prIVAN}, we recall that $\bar \gra$ annihilates the polynomial $Q_{a}(X)=(X-x^{c_{\gra}}\cdot \tilde u_{a,1})(X-x^{c_{\gra}} \cdot\tilde u_{a,2})$. Therefore, we may define $\gru$ such that $\gru(t_{13,1}):=x^{a_s} \tilde u_{s_1}$ and $\gru(t_{13,2}):=x^{a_s} \tilde u_{s_2}$.

   In the case of $G_{13}$ we have $(\bar\grb-1)(\bar\grb-j)(\bar\grb-j^2)=0$, where $j$ is a third root of unity. Therefore, we may define $\gru(t_{32,1}):=1$, $\gru(t_{32,2}):=j$
   and $\gru(t_{32,3}):=j^2$.
 In the case of $G_{15}$ we notice that $\grb$ is a distinguished braided reflection  associated to a distinguished reflection $b$ of order 3 (see in tables \ref{t3} and \ref{t2} the images of $\grb$ and $b$ in BMR presentations). We recall again that $\bar \grb$ annihilates the polynomial $Q_{b}(X)=(X-x^{c_{\grb}}\cdot \tilde u_{b,1})(X-x^{c_{\grb}} \cdot\tilde u_{b,2})(X-x^{c_{\grb}} \cdot\tilde u_{b,3})$. Therefore, we may define $\gru(t_{32,1}):=x^{c_{\grb}} \cdot\tilde u_{b,1}$, $\gru(t_{32,2}):=x^{c_{\grb}} \cdot\tilde u_{b,2}$ and $\gru(t_{32,3}):=x^{c_{\grb}}\cdot \tilde u_{b,3}$.
 
  We now consider the algebra $\tilde A_+(C)\otimes_{\gru}R_{\tilde \ZZ}^+$ and we define 
  $\begin{array}[t]{rcl}\grC:\tilde A_+(C)\otimes_{\gru} R_{\tilde \ZZ}^+ &\twoheadrightarrow& R_{\tilde \ZZ}^+\bar B/Q_s(\bar \grs)\\
  \tilde A_{13}&\mapsto& \bar \gra\\
  \tilde A_{32}&\mapsto& \bar \grb\\
  \tilde A_{21}&\mapsto& \bar \grg.
  \end{array}$
  
\end{proof} 
\begin{defn}
For every exceptional group of rank 2 we call the surjection $\grC$ as described in propositions \ref{ERSUR} and \ref{ERRSUR} the \emph{ER surjection} associated to $W$.
\end{defn}

\begin{prop}
	$H_{\tilde \ZZ}$ is generated as $R_{\tilde \ZZ}^+$-module by $|\overline{W}|$ elements.
\end{prop}
\begin{proof}
The result follows from corollary \ref{prIVAN}, theorem \ref{thmER} and from propositions \ref{ERSUR} and \ref{ERRSUR}.
\end{proof}
The following corollary is Theorem 2.14 of \cite{marinG26} and it is the main result in \cite{ERrank2}.
\begin{cor}(A weak version of the BMR freeness conjecture)
	$H$ is finitely generated over $R$.
	\label{corER}
\end{cor}
A first consequence of the weak version of the conjecture is the following:
\begin{prop}Let $F$ denote the field of fractions of $R$ and $\bar F$ an algebraic closure. Then, $H\otimes_R \bar F\simeq \bar F W$.
	\label{F}
\end{prop}
\begin{proof}
	The result follows directly from proposition 2.4(2) of \cite{marinG26} and from corollary \ref{corER}.
\end{proof}
Another consequence of the weak version of the conjecture is the following proposition, which states that if the generic Hecke algebra of an exceptional group of rank 2 is torsion-free as $R$-module, it will be sufficient to prove the BMR freeness conjecture for the maximal groups (see table \ref{t1}).
\begin{prop}
Let $W_0$ be the maximal group $G_7$, $G_{11}$ or $G_{19}$ and	let $W$ be an exceptional group of rank 2, whose associated Hecke algebra $H$ is torsion-free as $R$-module. If the BMR freeness conjecture holds for $W_0$, then it holds for $W$, as well.
\end{prop}
\begin{proof}
	Let $R_0$ and $R$ be the rings over which we define the Hecke algebras $H_0$ and $H$ associated to $W_0$ and $W$, respectively. There is a specialization $\gru :R_0\rightarrow R$, that maps some of the parameters of $R_0$ to roots of unity (see tables 4.6, 4.9 and 4.12 in \cite{malle2}). We set $A:=H_0\otimes_{\gru}R$. Due to hypothesis that the BMR freeness conjecture holds for $W_0$, we have that $A$ is a free $R$-module of rank $|W_0|$. 
	
	In proposition 4.2 in \cite{malle2}, G. Malle found a subalgebra $\bar A$ of $A$, such that $H\twoheadrightarrow \bar A$. A presentation of $\bar A$ is given in Appendix A of \cite{chlouverakibook}. He also noticed that if $m=|W_0|/|W|$, then there is an element $\grs\in A$ such that $A=\oplus_{i=0}^{m-1}\grs^i\bar A$. We highlight here that for all these results G. Malle does not use the validity of the BMR freeness conjecture. Since $A$ is a free $R$-module of rank $|W_0|$, we also have that $\bar A$ is a free module of rank $|W|$. 
	
	Our goal is to prove that as $R$-modules, $\bar A \simeq H$. Let $F$ denotes the field of fractions of $R$ and $\bar F$ an algebraic closure. By proposition \ref{F} we have that $\bar A\otimes_R \bar F \simeq H\otimes_R \bar F$. We have the following commutative diagram:
	$$
	\xymatrix{
		H \ar[r]^{\phi_2} \ar@{->>}[d]^{\psi_2} \ar@{.>}[dr]& H \otimes_R \overline{F}  \ar[d]^{\psi_1}_{\rotatebox{90}{$\simeq$}}\\
		\overline{A} \ar[r]^{\phi_1}& \overline{A} \otimes_R \overline{F}
	}
	$$
	
	Let $h\in$ ker$\grc_2$. Then, $\grf_1(\grc_2(h))=\grc_1(\grf_2(h))$ and, hence, $\grf_2(h)=(\grc_1^{-1}\circ \grf_1)(0)=\grc_1(0)=0$, which means that ker$\grc_2\subset$ ker$\grf_2$. 
	 Let $h\in$ ker$\grf_2$. Then, $\grf_1(\grc_2(h))=\grc_1(\grf_2(h))=0$, which means that $\grc_2(h)\in \bar A^{\text{tor}}$. However, $\bar A$ is a free $R$-module, therefore $A^{\text{tor}}=0$. As a result, ker$\grf_2\subset$ ker$\grc_2$ and, hence, $\grf_2=0$, since $H$ is torsion-free as $R$-module. Therefore, $\grc_2$ is an isomorphism and, as a result, $H$ is a free $R$-module of rank $|W|$.
	 
\end{proof}

We now give an example that makes clearer to the reader the nice properties of the algebra $\bar A$ we mentioned in the proof of the above proposition.
\begin{ex}
	Let $W=G_5$, which is an exceptional group that belongs to the tetrahedral family. In this family the maximal group is the group $W_0=G_7$ (see table \ref{t1}). Let $R_0$ be the Laurent polynomial ring $\ZZ\left[x_i^{\pm},y_j^{\pm}, z_l^{\pm}\right]$ over which we define the generic Hecke algebra $H_0$ associated to $G_7$, where $i=1,2$ and $j,l=1,2,3$.
	$H_0$ can be presented as follows:
	
$$H_0 =\langle s,t,u\;|\; stu=tus=ust, \prod\limits_{i=1}^2(s-x_i)=\prod\limits_{j=1}^3(t-y_j)=\prod\limits_{l=1}^3(u-z_l)=0\rangle.$$ We assume that $H_0$ is a $R_0$-module of rank $|W_0|$.
	Let $R$ be the Laurent polynomial ring $\ZZ\left[\bar y_j^{\pm}, \bar z_l^{\pm}\right]$ over which we define the generic Hecke algebra $H$ associated to $G_5$, where $j,l=1,2,3$. $H$ can be presented as follows:
	
	$$H =\langle S,T\;|\; STST=TSTS, \prod\limits_{i=1}^3(S-\bar x_i)=\prod\limits_{j=1}^3(T-\bar y_j)=0\rangle.$$
	Let $\gru: R_0\rightarrow R$ be the specialization defined by $x_1\mapsto 1, x_2\mapsto -1$. The $R$ algebra $A:=H_0\otimes_{\gru}R$ is presented as follows:
	$$A=\langle s,t,u\;|\; stu=tus=ust, s^2=1,\prod\limits_{j=1}^3(t-y_j)=\prod\limits_{l=1}^3(u-z_l)=0\rangle.$$
	Let $\bar A=\langle t,u\rangle\leq A$. We define $\grf: H \rightarrow A$, $S\mapsto t$, $T\mapsto u$. Since  $stu$ generates the center of the complex braid group associated to $G_7$ we have $tutu=tus^2tu=(tus)stu=ust(stu)=us(stu)t=utut$, meaning that $\grf$ is surjective. We also notice that $stu=ust\an us=stut^{-1}$ and $ust=tus\an t=su^{-1}tus\an ts=su^{-1}tu$. Therefore,
	$A=\bar A+s\bar A$. 
	
	In order to prove that this sum is direct, we consider an algebraic closure $\bar F$ of the field of fractions of $F$. By Proposition \ref{F} we have that $A\otimes_R\bar F =\bar A\otimes_R\bar F\oplus s\bar A\otimes_R\bar F$, since  $A\otimes_R\bar F$ is a free $R$-module of rank $|W_0|=144$ and $\bar A\otimes_R\bar F$ is a free $R$-module of rank $|W|=72$. 
	
	Let $x\in \bar A \cap sA$. Then, $x\otimes 1_{\bar F}\in \bar A\otimes_R\bar F\cap s\bar A\otimes_R\bar F$, which means that $x\in A^{\text{tor}}={0_A}$, since $A$ is a free $R$-module. Therefore, $A=\bar A\oplus s\bar A$, meaning that $\bar A$ is a free $R$-module of rank $|W|$.
	
	\qed
\end{ex}

\chapter{The BMR freeness conjecture for the tetrahedral and octahedral families}

In this chapter we prove the BMR freeness conjecture for the exceptional groups belonging to the first two families, using a case-by-case analysis. As one may notice, in this proof we ``guess'' how the basis should look
like and then we establish that this is in fact a basis through a series of 
computations. This method allowed us not
only to prove the validity of the conjecture, but also to give a nice description of the basis, similar
to the classical case of the finite Coxeter groups. In the first part of this chapter we explain how we arrived to suspect there exists a basis of this nice form, using the basis of the algebra $A_+(C)$ that P. Etingof and E. Rains also used in order to prove the weak version of the conjecture for the exceptional groups of rank 2, a result we explored in detail in the previous chapter.
\section{Finding the basis}
\indent

Let $W$ denote an exceptional group of rank 2. In the previous chapter we explained that $\overline{W}:=W/Z(W)$ can be considered as the group of even elements in a finite Coxeter group $C$ of rank 3 with Coxeter system $\{y_1, y_2, y_3\}$. For every $x$ in $\overline{W}$, we fix a reduced word $w_x$ in letters $y_1,y_2$ and $y_3$ representing $x$ in  $\overline{W}$ and we set $a_{ij}:=y_iy_j$, $i,j\in\{1,2,3\}$ with $i\not=j$.  Since $\overline{W}$ is the group of even elements in $C$, the reduced word $w_x$ is a word in letters $a_{ij}$ and it corresponds to an element in $A_+(C)$ denoted by $T_{w_x}$. In particular,  $T_{w_1}=1_{A_+(C)}$.

Keeping the notations of the previous chapter, we recall that $H_{\tilde \ZZ}$ is generated as $R_{\tilde \ZZ}^+$-module by the elements $\grC(T_{\grv_x})$, $x\in \overline{W}$, where $\grC$ is the ER-surjection associated to $W$ (see \ref{ERSUR} and \ref{ERRSUR}). Therefore, we have that $H_{\tilde \ZZ}$  is spanned over $R_{\tilde \ZZ}^+$ by $|\overline{W}|$ elements. Motivated by this idea, we will explain in general how we found a spanning set of $H$ over $R$ of $|W|$ elements.

For every $x\in \overline{W}$ we fix a reduced word $w_x$ in letters $y_1,y_2$ and $y_3$ that represents $x$ in $\overline{W}$. From the reduced word $w_x$ one can obtain a word $\bar w_x$ that also represents $x$ in $\overline{W}$, defined as follows:
$$\bar w_x= \begin{cases}
w_x  &\text{for }x=1\\
w_x(y_1y_1)^{n_1}(y_2y_2)^{n_2}(y_3y_3)^{n_3}&\text{for }x\not=1
\end{cases},$$
 where $n_i \in \ZZ_{\geq 0}$ and $(y_iy_i)^{n_i}$ is a shorter notation for the word $\underbrace{(y_iy_i)\dots(y_iy_i)}_{n_i-\text{ times }}$. We notice that if we choose $n_1=n_2=n_3=0$, the word $\bar w_{x}$ coincides with the word $w_{x}$.
 
 Moving some of the pairs   $(y_iy_i)^{n_i}$ somewhere inside $\bar w_x$ and using the braid relations between the generators $y_i$ of the Coxeter group $C$ one can obtain a word $\tilde w_x$, which also represents $x$ in $\bar W$, such that:
\begin{itemize}
	\item $\ell(\tilde w_x)=\ell(\bar w_x)$, where $\ell(w)$ denotes the length of the word $w$.
	\item Let $m$ be an odd number. Whenever in the word $\tilde w_x$ there is a letter $y_i$ at the $m$th-position from left to right, then in the $(m+1)$th-position there is a letter $y_j$, $j\not=i$.
		\item $\tilde w_x=w_x$ if and only if $\bar w_x=w_x$. In particular, 	$\tilde w_1=w_1$.
\end{itemize} 
\begin{defn}
A word $\tilde w_x$ as described above is called \emph{a base word} associated  to $\bar w_x$. 
\end{defn}

Let $\tilde w_x$ be a base word. We recall that $a_{ij}:=y_iy_j$, $i,j\in\{1,2,3\}$ with $i\not=j$. By the definition of $\tilde w_x$ and the fact that $\overline{W}$ is the group of even elements in  $C$, the word $\tilde w_x$ can be considered as a word in letters $a_{ij}$.
Let $T_{\tilde w_x}$ denote the corresponding element in $A_+(C)$. In particular,  $T_{\tilde w_1}=T_{w_1}=1_{A_+(C)}$. 

Let $B$ be the complex braid group associated to $W$ and let  $\bar B$ denote the quotient $B/Z(B)$. For every $b\in B$ we denote by $\bar b$ the image of $b$ under the natural projection $B\rightarrow \bar B$. In the previous chapter we explained that $\bar B$ is generated by the elements 
$\bar \gra,$ $\bar \grb$ and  $\bar \grg$, where $\gra$, $\grb$ and $\grg$ are generators of $B$ in ER presentation (see table \ref{t2} in Appendix B).

 By the definition of the ER-surjection the element $\grC(T_{\tilde w_x})$ is
a product of $\bar \gra$, $\bar \grb$ and $\bar \grg$ (see \ref{ERSUR} and \ref{ERRSUR}). We use the group isomorphism $\grf_2$ we describe in table \ref{t2}, Appendix B  to write the elements $\gra$, $\grb$ and $\grg$ in BMR presentation and we set $\grs_{\gra}:=\grf_2(\gra)$, $\grs_{\grb}:=\grf_2(\grb)$ and $\grs_{\grg}:=\grf_2(\grg)$. Therefore, we can also consider the element $\grC(T_{\tilde w_x})$ as being
a product of $\bar\grs_{\gra} $, $\bar\grs_{\grb} $ and $\bar\grs_{\grg}$. We denote this element by $\bar v_x$.

We now explain  how we arrived to guess a spanning set of $|W|$ elements for the generic Hecke algebra associated to every exceptional group belonging to the first two families. 
\begin{itemize}

\item  Let $W$ be an exceptional group of rank 2, which belongs either to the tetrahedral or octahedral family. For every $x\in \overline{W}$ we choose a specific reduced word $w_x$, specific non-negative integers $(n_i)_{\substack{1\leq i\leq 3}}$ and a specific  base word $\tilde w_x$ associated to the word $\bar w_x$, which is determined by $w_x$ and $(n_i)$.

\item For every $x\in \overline {W}$, let $x_1^{m_1}x_2^{m_2}\dots x_r^{m_r}$ be the corresponding factorization of $\bar v_x$ into a product of $\bar \grs_{\gra}$, $\bar \grs_{\grb}$ and $\bar \grs_{\grg}$ (meaning that $x_i \in\{\bar \grs_{\gra}, \bar \grs_{\grb},\bar \grs_{\grg}\}$ and $m_i\in \ZZ$). 
Let $f_0: \bar B\rightarrow B$ be a set theoretic section such that 
$f_0(x_1^{m_1}x_2^{m_2}\dots x_r^{m_r})=f_0(x_1)^{m_1}f_0(x_2)^{m_2}\dots f_0(x_r)^{m_r}$, $f_0(\bar \grs_{\gra})=\grs_{\gra}$,  $f_0(\bar \grs_{\grb})=\grs_{\grb}$ and $f_0(\bar \grs_{\grg})=\grs_{\grg}$. 
 Keeping the notations of the previous chapter, we use corollary \ref{corIVAN} and we obtain an isomorphism $\grf_{
 	f_0}$ between the $R_{\ZZ}^+$-modules $R_{\tilde \ZZ}^+\bar B/Q_{s}(\bar \grs)$ and $H_{\tilde \ZZ}$. We set $v_x:=\grf_{f_0}(\bar v_x)$.
\item  We set $U:=\sum\limits_{x\in \overline{W}} \sum\limits_{k=0}^{|Z(W)|-1}Rz^kv_x$, where $R$ is the Laurent polynomial ring over which we define the generic Hecke algebra $H$ associated to $W$.

\end{itemize}
	The main theorem of this chapter is the following. Notice that the second part of this theorem follows directly from proposition \ref{BMR PROP}.
\begin{thm}
	$H=U$ and, therefore, the BMR freeness conjecture holds for all the groups belonging to the tetrahedral and octahedral family. 
	\label{main}
\end{thm}

\begin{rem}
	\mbox{}
	\vspace*{-\parsep}
	\vspace*{-\baselineskip}\\
	\begin{itemize}	
	\item[1.] The choice of $w_x$, the non-negative integers $(n_i)_{\substack{1\leq i\leq 3}}$ and $\tilde w_x$  is a product of experimentation, to provide a simple and robust proof for theorem \ref{main}. We tried more combinations that lead to more complicated and bloated proofs, or others where we couldn't arrive to a conclusion.
	\item[2.] The case of $G_{12}$ has already been proven (see theorem \ref{2case}). However, using this approach we managed to give a proof in a different way. We recall that $G_{12}$ has a presentation of the form $$\langle s,t,u\;|\;s^2=t^2=u^2=1, stus=tust=ustu\rangle.$$ Since $s$, $t$ and $u$ are conjugate (notice that $(stu)s(stu)^{-1}=u$, $(ust)u(ust)^{-1}=t$), the Laurent polynomial ring over which we define $H$ is $\ZZ[u_{s,i}^{\pm}]_{1\leq i\leq 2}$. In this chapter we prove the BMR freeness conjecture for $H$ defined over the ring 
	$\ZZ[u_{s,i}^{\pm},u_{t,j}^{\pm},u_{u,l}^{\pm}]_{\substack{1\leq i,j,l\leq 2 }}$, using the method we described above for the rest of the exceptional groups belonging to the first two families.
	\end{itemize}
\end{rem}

The rest of this chapter is devoted to the proof of theorem \ref{main}, using a case-by-case analysis. We finish this section by giving an example of the way we find $U$ for the exceptional group $G_{15}$.

\begin{ex}
	Let $W:=G_{15}$, an exceptional group that belongs to the octahedral family. In this family $C$ is the Coxeter group of type $B_3$. For every $x_i\in \overline{W}$, $i=1,\dots, 24$ we fix a reduced word $w_{x_i}$ in letters $y_1, y_2, y_3$ that represents $x_i$ in $\overline{W}$:
	
	\scalebox{0.9}{%
		\vbox{%
	\begin{multicols}{3}
		\begin{itemize}[leftmargin=*]
			\item[1.] $w_{x_1}=1$
			\item[2.] $w_{x_2}=y_3y_2$
			\item[3.] $w_{x_3}=y_2y_3$
			\item[4.] $w_{x_4}=y_2y_1y_2y_1$
			\item[5.] $w_{x_5}=y_3y_1y_2y_1$  
			\item[6.] $w_{x_6}=y_2y_3y_2y_1y_2y_1$
			\item [7.]  $w_{x_7}=y_1y_3$
			\item [8.] $w_{x_8}=y_3y_2y_1y_3$
			\item[9.] $w_{x_9}=y_2y_1$
			\item[10.] $w_{x_{10}}=y_2y_1y_2y_3$
			\item[11.] $w_{x_{11}}=y_3y_1y_2y_3$
			\item[12.] $w_{x_{12}}=y_2y_3y_2y_1y_2y_3$
			\item[13.] $w_{x_{13}}=y_3y_2y_1y_2$
			\item[14.] $w_{x_{14}}=y_2y_3y_1y_2$
			\item[15.] $w_{x_{15}}=y_1y_2$
			\item [16.]  $w_{x_{16}}=y_2y_1y_2y_3y_2y_1y_2y_1$
			\item [17.]  $w_{x_{17}}=y_1y_2y_1y_3y_2y_1$
			\item[18.] $w_{x_{18}}=y_3y_2y_1y_2y_3y_1y_2y_1$
			\item[19.] $w_{x_{19}}=y_1y_3y_2y_1y_2y_3$
			\item[20.] $w_{x_{20}}=y_2y_3y_1y_2y_1y_3$
			\item[21.] $w_{x_{21}}=y_1y_2y_1y_3$
			\item[22.] $w_{x_{22}}=y_3y_2y_1y_2y_3y_2y_1y_2$
			\item[23.] $w_{x_{23}}=y_2y_1y_2y_3y_1y_2$
			\item [24.]  $w_{x_{24}}=y_1y_2y_3y_2y_1y_2$
			\end{itemize}	\end{multicols}}}
		
We choose the words $\bar w_{x_i}$, $i=1,\dots, 24$ as follows:

\scalebox{0.9}{%
	\vbox{%
\begin{multicols}{3}
\begin{itemize}[leftmargin=*]
	
		\item[1.]  $\bar w_{x_1}=w_{x_1}$
		\item[2.] $\bar w_{x_2}=w_{x_2}$
		\item[3.] $\bar w_{x_3}= w_{x_3}$
		\item[4.] $\bar w_{x_4}= w_{x_4}$
		\item[5.] $\bar w_{x_5}= w_{x_5}\mathbf{y_2y_2}$
		\item[6.] $\bar w_{x_6}= w_{x_6}$
		\item [7.]  $\bar w_{x_7}= w_{x_7}$
		\item [8.]  $\bar w_{x_8}= w_{x_5}$
		\item[9.] $\bar w_{x_9}= w_{x_9}\mathbf{y_3y_3}$
		\item[10.] $\bar w_{x_{10}}=w_{x_{10}}\mathbf{y_1y_1}$
		\item[11.] $\bar w_{x_{11}}=w_{x_{11}}\mathbf{y_1y_1}\mathbf{y_2y_2}$
		\item[12.] $\bar w_{x_{12}}=w_{x_{12}}\mathbf{y_1y_1}$
		\item[13.] $\bar w_{x_{13}}=w_{x_{13}}\mathbf{y_1y_1}$
		\item[14.] $\bar w_{x_{14}}=w_{x_{14}}\mathbf{y_1y_1y_2y_2}$
		\item[15.] $\bar w_{x_{15}}=w_{x_{15}}\mathbf{y_1y_1y_2y_2y_3y_3}$
		\item [16.]  $\bar w_{x_{16}}=w_{x_{16}}\mathbf{y_1y_1}$
		\item [17.] $\bar w_{x_{17}}=w_{x_{17}}\mathbf{y_2y_2y_3y_3}$
		\item[18.] $\bar w_{x_{18}}=w_{x_{18}}\mathbf{y_3y_3}$
		\item[19.] $\bar w_{x_{19}}=w_{x_{19}}\mathbf{y_1y_1}$
		\item[20.] $\bar w_{x_{20}}=w_{x_{20}}\mathbf{y_1y_1}\mathbf{y_3y_3}$
		\item[21.] $\bar w_{x_{21}}=w_{x_{21}}\mathbf{y_1y_1y_2y_2y_3y_3}$
		\item[22.] $\bar w_{x_{22}}=w_{x_{22}}\mathbf{(y_1y_1)^2}$
		\item[23.] $\bar w_{x_{23}}=w_{x_{23}}\mathbf{(y_1y_1)^2}\mathbf{y_2y_2y_3y_3}$
		\item [24.]  $\bar w_{x_{24}}=w_{x_{24}}\mathbf{(y_1y_1)^2}\mathbf{y_2y_2y_3y_3}$
		
	\end{itemize}	\end{multicols}}}
	
We now choose a base word $\tilde w_{x_i}$ for every $i=1,\dots, 24$.  We give as indicative example the base word $\tilde w_{x_{24}}$.
$$\begin{array}{lcl}
\bar w_{x_{24}}&=&w_{x_{24}}\mathbf{(y_1y_1)^2}\mathbf{y_2y_2y_3y_3}\\
&=&y_1y_2y_3(y_2y_1y_2\mathbf{y_1)y_1}\mathbf{y_1y_1y_2y_2y_3y_3}\\
&=&y_1y_2y_3y_1y_2y_1y_2y_1\mathbf{y_1y_1y_2y_2y_3y_3}\\
&=&y_1y_2(y_3y_1)y_2y_1y_2y_1\mathbf{y_1y_1y_2y_2y_3y_3}\\
	&=&y_1y_2y_1y_3y_2y_1y_2y_1\mathbf{y_1y_1y_2y_2y_3y_3}\\
		&=&(y_1y_2y_1\mathbf{y_2)y_2}y_3y_2y_1y_2y_1\mathbf{y_1y_1y_3y_3}\\
		&=&y_2y_1y_2y_1y_2y_3y_2y_1y_2y_1\mathbf{y_1y_1y_3y_3}\\
			&=&y_2(y_1\mathbf{y_3)y_3}y_2y_1y_2\mathbf{y_1y_1}y_3y_2y_1y_2y_1\\
			&=&y_2y_3y_1y_3y_2y_1y_2y_1y_1y_3y_2y_1y_2y_1.
\end{array}$$
We choose $\tilde w_{x_{24}}=y_2y_3y_1y_3y_2y_1y_2y_1y_1y_3y_2y_1y_2y_1=a_{23}a_{13}a_{21}a_{21}a_{13}a_{21}a_{21}$ and, hence, $T_{\tilde w_{x_{24}}}=A_{23}A_{13}A_{21}A_{21}A_{13}A_{21}A_{21}$.
The rest of the base words are chosen as follows. Note that for $i=1,\dots,4$ and for $i=6,\dots, 8$ the base word $\tilde w_{x_i}$ is chosen to be $\bar w_{x_i}$, since $\bar w_{x_i}=w_{x_i}$.

\scalebox{0.9}{%
	\vbox{%
	\begin{multicols}{3}
		
		\begin{itemize}[leftmargin=*]
			
		\item[1.] 
			
				$\tilde w_{x_1}=1$
			\item[2.] 
			$\tilde w_{x_2}=y_3y_2$
			\item[3.] 
			$\tilde w_{x_3}=y_2y_3$
			\item[4.] 
			$\tilde w_{x_4}=y_2y_1y_2y_1$
			\item[5.]  $\tilde w_{x_5}=y_3y_2y_2y_1y_2y_1$ 
			\item[6.] $\tilde w_{x_6}=y_2y_3y_2y_1y_2y_1$
			\item [7.]  $\tilde w_{x_7}=y_1y_3$
			\item [8.]$\tilde w_{x_8}=y_3y_2y_1y_3$
			\item[9.] $\tilde w_{x_9}=y_2y_3y_3y_1$
			\item[10.] $\tilde w_{x_{10}}=y_2y_1y_2y_1y_1y_3$
			\item[11.] $\tilde w_{x_{11}}=y_3y_2y_2y_1y_2y_1y_1y_3$
			\item[12.] $\tilde w_{x_{12}}=y_2y_3y_2y_1y_2y_1y_1y_3$
			\item[13.] $\tilde w_{x_{13}}=y_1y_3y_2y_1y_2y_1$
			\item[14.] $\tilde w_{x_{14}}=y_3y_2y_1y_3y_2y_1y_2y_1$
			\item[15.] $\tilde w_{x_{15}}=y_2y_3y_1y_3y_2y_1y_2y_1$
			\item [16.] 
				$\tilde w_{x_{16}}=y_2y_1y_2y_1y_1y_3y_2y_1y_2y_1$
			\item [17.]  $\tilde w_{x_{17}}=y_3y_2y_2y_1y_2y_3y_2y_1y_2y_1$
			\item[18.] $\tilde w_{x_{18}}=y_2y_3y_2y_3y_1y_2y_3y_1y_2y_1$
			\item[19.] $\tilde w_{x_{19}}=y_1y_3y_2y_1y_2y_1y_1y_3$
			\item[20.] $\tilde w_{x_{20}}=y_3y_2y_1y_3y_2y_1y_2y_1y_1y_3$
			\item[21.] $\tilde w_{x_{21}}=y_2y_3y_1y_3y_2y_1y_2y_1y_1y_3$
			\item[22.] $\tilde w_{x_{22}}=y_1y_3y_2y_1y_2y_1y_1y_3y_2y_1y_2y_1$
			\item[23.] $\tilde w_{x_{23}}=y_3y_2y_1y_3y_2y_1y_2y_1y_1y_3y_2y_1y_2y_1$
			\item [24.] $\tilde w_{x_{24}}=y_2y_3y_1y_3y_2y_1y_2y_1y_1y_3y_2y_1y_2y_1$
			\end{itemize}
		\end{multicols}}}
		
The corresponding elements $T_{\tilde{w_{x_i}}}$ in $A_+(C)$ are as follows:

\scalebox{0.93}{%
	\vbox{%
\begin{multicols}{3}
	
	\begin{itemize}[leftmargin=*]
		
		\item[1.] $T_{\tilde w_{x_1}}=1$
		\item[2.] $T_{\tilde w_{x_2}}=A_{32}$
		\item[3.] $T_{\tilde w_{x_3}}=A_{23}$
		\item[4.] $T_{\tilde w_{x_4}}=A_{21}A_{21}$
		\item[5.] $T_{\tilde w_{x_5}}=A_{32}A_{21}A_{21}$ 
		\item[6.] $T_{\tilde w_{x_6}}=A_{23}A_{21}A_{21}$
		\item [7.]  $T_{\tilde w_{x_7}}=A_{13}$
		\item [8.] $T_{\tilde w_{x_8}}=A_{32}A_{13}$
		\item[9.] $T_{\tilde w_{x_9}}=A_{23}A_{31}$
		\item[10.] $T_{\tilde w_{x_{10}}}=A_{21}A_{21}A_{13}$
		\item[11.] $T_{\tilde w_{x_{11}}}=A_{32}A_{21}A_{21}A_{13}$
		\item[12.] $T_{\tilde w_{x_{12}}}=A_{23}A_{21}A_{21}A_{13}$
		\item[13.] $T_{\tilde w_{x_{13}}}=A_{13}A_{21}A_{21}$
		\item[14.] $T_{\tilde w_{x_{14}}}=A_{32}A_{13}A_{21}A_{21}$
		\item[15.] $T_{\tilde w_{x_{15}}}=A_{23}A_{13}A_{21}A_{21}$
		\item [16.]  $T_{\tilde w_{x_{16}}}=A_{21}A_{21}A_{13}A_{21}A_{21}$
		\item [17.]  $T_{\tilde w_{x_{17}}}=A_{32}A_{21}A_{23}A_{21}A_{21}$
		\item[18.] $T_{\tilde w_{x_{18}}}=A_{23}A_{23}A_{12}A_{31}A_{21}$
		\item[19.] $T_{\tilde w_{x_{19}}}=A_{13}A_{21}A_{21}A_{13}$
		\item[20.] $T_{\tilde w_{x_{20}}}=A_{32}A_{13}A_{21}A_{21}A_{13}$
		\item[21.] $T_{\tilde w_{x_{21}}}=A_{32}A_{13}A_{21}A_{21}A_{13}$
		\item[22.] $T_{\tilde w_{x_{22}}}=A_{13}A_{21}A_{21}A_{13}A_{21}A_{21}$
		\item[23.] $T_{\tilde w_{x_{23}}}=A_{32}A_{13}A_{21}A_{21}A_{13}A_{21}A_{21}$
		\item [24.] $T_{\tilde w_{x_{24}}}=A_{23}A_{13}A_{21}A_{21}A_{13}A_{21}A_{21}$
	\end{itemize}
\end{multicols}	}}

We recall that $\grC(A_{13})=\bar \gra$, $\grC(A_{32})=\bar \grb$ and $\grC(A_{21})=\bar \grg$, or equivalently,  in BMR presentation $\grC(A_{13})=\bar t$, $\grC(A_{32})=\bar u$ and $\grC(A_{21})=\bar u^{-1}\bar t^{-1}$ (see table \ref{t2} in Appendix B). We also notice that
$\grC(A_{21}A_{21})=\bar \grg^2$ and, therefore, in BMR presentation,  $\grC(A_{21}A_{21})=\bar s$. 

Using again the example of $x_{24}$ we notice that
$\bar v_{x_{24}}=\grC(w_{x_{24}})=\grC(A_{23}A_{13}A_{21}A_{21}A_{13}A_{21}A_{21})=\bar u^{-1}\bar t \bar s \bar t \bar s$ and, hence, $v_{x_{24}}=u^{-1}tsts$.
Similarly, we the elements $v_{x_i}$, $i=1,\dots, 24$ are as follows:

\scalebox{0.9}{%
	\vbox{%
		\begin{multicols}{3}
			
			\begin{itemize}[leftmargin=*]
				
				\item[1.] $v_{x_1}=1$
				\item[2.] $v_{x_2}=u$
				\item[3.] $v_{x_3}=u^{-1}$
				\item[4.] $v_{x_4}=s$
				\item[5.] $v_{x_5}=us$
				\item[6.] $v_{x_6}=u^{-1}s$
				\item [7.]  $v_{x_7}=t$
				\item [8.] $v_{x_8}=ut$
				\item[9.] $v_{x_9}=u^{-1}t$
				\item[10.] $v_{x_{10}}=st$
				\item[11.] $v_{x_{11}}=ust$
				\item[12.] $v_{x_{12}}=u^{-1}st$
				\item[13.] $v_{x_{13}}=ts$
				\item[14.] $v_{x_{14}}=uts$
				\item[15.]$v_{x_{15}}=u^{-1}ts$
				\item [16.]  $v_{x_{16}}=sts$
				\item [17.]  $v_{x_{17}}=usts$
				\item[18.] $v_{x_{18}}=u^{-1}sts$
				\item[19.] $v_{x_{19}}=tst$
				\item[20.]  $v_{x_{20}}=utst$
				\item[21.] $v_{x_{21}}=u^{-1}tst$
				\item[22.] $v_{x_{22}}=tsts$
				\item[23.]$v_{x_{23}}=utsts$
				\item [24.] $v_{x_{24}}=u^{-1}tsts$
			\end{itemize}
		\end{multicols}	}}
		
	We denote by $u_1:=R+Rs$ the subalgebra of $H$ generated by $s$, $u_2:=R+Rt$ the subalgebra of $H$ generated by $t$ and $u_3:=R+Ru+Ru^{-1}$ the subalgebra of $H$ generated by $u$, one may notice that $\sum\limits_{i=1}^{24}Rv_i=u_3u_2u_1u_2u_1$. We recall that $|Z(G_{15})|=12$ and we set $U=\sum\limits_{k=0}^{11}z^ku_3u_2u_1u_2u_1$, where $z=stutu$ the generator of the center of the corresponding complex braid group.
\end{ex}

\section{The Tetrahedral family}
\indent

In this family we encounter the exceptional groups $G_4$, $G_5$, $G_6$ and $G_7$. We know that the BMR freeness conjecture holds for $G_4$ (see theorem \ref{braidcase}). We prove the conjecture for the rest of the  groups belonging in this family, using a case-by-case analysis. Keeping the notations of Chapter 3, let $P_s(X)$ denote the polynomials defining $H$ over $R$. If we  expand the relations $P_s(\grs)=0$, where $\grs$ is a distinguished braided reflection associated to $s$, we obtain equivalent relations of the form \begin{equation}\grs^n=a_{n-1}\grs^{n-1}+...+a_1\grs+a_0,\label{ones} \end{equation}
where $n$ is the order of $s$, $a_i\in R$, for every $i\in\{1, \dots n-1\}$ and $a_0 \in R^{\times}$.
 We multiply $(\ref{ones})$ by $\grs_i^{-n}$ and since $a_0$ is invertible in $R$ we have:
\begin{equation}\grs_i^{-n}=-a_0^{-1}a_1\grs^{-n+1}-a_0^{-1}a_2\grs^{-n+2}-...-a_0^{-1}a_{n-1}\grs^{-1}+a_0^{-1} \label{twos}\end{equation}
We multiply ($\ref{ones}$) with a suitable power of $\grs$. Then, for every $m\in \NN$ we have :\begin{equation}\grs^{m}\in R\grs^{m-1}+\dots+R\grs^{m-(n-1)}+R^{\times}\grs^{m-n}\label{ooo}\end{equation}
Similarly, we multiply ($\ref{twos}$) with a suitable power of $\grs$. Then, for every $m\in \NN$, we have: \begin{equation}\grs^{-m}\in  R\grs^{-m+1}+\dots+R\grs^{-m+(n-1)}+R^{\times}\grs^{-m+n}.\label{oooo}\end{equation}
For the rest of this section we use directly (\ref{ooo}) and (\ref{oooo}). 
\subsection{The case of $G_5$}

\indent

Let $R=\ZZ[u_{s,i}^{\pm},u_{t,j}^{\pm}]_{\substack{1\leq i,j\leq 3}}$ and let $H_{G_5}:=\langle s,t\;|\; stst=tsts,\prod\limits_{i=1}^{3}(s-u_{s,i})=\prod\limits_{j=1}^{3}(t-u_{t,i})=0\rangle$ be the generic Hecke algebra associated to $G_5$. Let $u_1$ be the subalgebra of $H_{G_5}$ generated by $s$ and $u_2$ the subalgebra of $H_{G_5}$ generated by $t$. 
We recall that $z:=(st)^2=(ts)^2$  generates the center of the associated complex braid group and that $|Z(G_5)|=6$. We set  $U=\sum\limits_{k=0}^5(z^ku_1u_2+z^kt^{-1}su_2)$. By the definition of $U$ we have the following remark:
\begin{rem} $Uu_2\subset U$.
	\label{r55}
\end{rem}
To make it easier for the reader to follow the calculations, we will underline the elements that  belong to $U$ by definition. Moreover, we will use directly remark \ref{r55}; this means that every time we have a power of $t$ at the end of an element we may ignore it. To remind that to the reader, we put a parenthesis around the part of the element we consider.

Our goal is to prove that $H_{G_5}=U$ (theorem \ref{thm55}). In order to do so, we  first need to prove some preliminary results.
\begin{lem}
	\mbox{}
	\vspace*{-\parsep}
	\vspace*{-\baselineskip}\\
	\begin{itemize}[leftmargin=0.8cm]
		\item [(i)] For every $k\in\{1,\dots, 4\}$, $z^ktu_1\subset U$.
		\item[(ii)]For every $k\in\{1,\dots,5\}$, $z^kt^{-1}u_1\subset U$.
		\item[(iii)]For every $k\in\{1,\dots,4\}$, $z^ku_2u_1\subset U$.
		\end{itemize}
		\label{ts55}
\end{lem}
\begin{proof}
	By definition, $u_2=R+Rt+Rt^{-1}$. Therefore, (iii) follows from (i) and (ii).
	\begin{itemize}[leftmargin=0.8cm]
		\item[(i)]
	$z^ktu_1=z^kt(R+Rs+Rs^{-1})\subset \underline{z^ku_2}+Rz^k(ts)^2s^{-1}t^{-1}+Rz^kts^{-1}\subset U+\underline{z^{k+1}u_1u_2}+Rz^kts^{-1}$.
	However, $z^kts^{-1}\in z^k(R+Rt^{-1}+Rt^{-2})s^{-1}\subset U+\underline{z^ku_1}+Rz^k(st)^{-2}st+Rz^kt^{-1}(st)^{-2}st\subset U+\underline{z^{k-1}u_1u_2}+\underline{z^{k-1}t^{-1}su_2}\subset U$.
		\item[(ii) ]
		$z^kt^{-1}u_1=z^kt^{-1}(R+Rs+Rs^{-1})\subset \underline{z^ku_2}+\underline{z^kt^{-1}su_2}+zR^k(st)^{-2}st\subset U+\underline{z^{k-1}u_1u_2}\subset U.$
		\qedhere
		\end{itemize}
		
\end{proof}

From now on, we will double-underline the elements described in lemma \ref{ts55} and we will use directly the fact that these elements are inside $U$.
The following proposition leads us to the main theorem of this section.
\begin{prop}
	$u_1U\subset U$.
	\label{su55}
\end{prop}
\begin{proof}
	Since $u_1=R+Rs+Rs^2$, it is enough to prove that $sU\subset U$. By the definition of $U$ and by remark \ref{r55}, we can restrict ourselves to proving that  $z^kst^{-1}s\in U$, for every $k\in\{0,\dots,5\}$. We distinguish the following cases:
	\begin{itemize}[leftmargin=*]
		\item \underline{$k\in\{0,\dots,3\}$}:
		$\small{\begin{array}[t]{lcl}
		z^kst^{-1}s&\in&z^ks(R+Rt+Rt^2)s\\
		&\in&\underline{z^ku_1}+Rz^k(st)^2t^{-1}+Rz^k(st)^2t^{-1}s^{-1}ts\\
		&\in&U+\underline{z^{k+1}u_2}+Rz^{k+1}t^{-1}(R+Rs+Rs^2)ts\\
		&\in&U+\underline{z^{k+1}u_1}+Rz^{k+1}t^{-1}(st)^2t^{-1}+Rz^{k+1}t^{-1}s(st)^2t^{-1}\\
		&\in&U+\underline{z^{k+2}u_2}+\underline{z^{k+2}t^{-1}su_2}.
		\end{array}}$
		\item \underline{$k\in\{4,5\}$}:
		$\small{\begin{array}[t]{lcl}z^kst^{-1}s&\in&z^k(R+Rs^{-1}+Rs^{-2})t^{-1}(R+Rs^{-1}+Rs^{-2})\\
		&\in&\underline{\underline{z^kt^{-1}u_1}}+\underline{z^ku_1u_2}+Rz^ks^{-1}t^{-1}s^{-1}+Rz^ks^{-1}t^{-1}s^{-2}+
		Rz^ks^{-2}t^{-1}s^{-1}+\\&&+Rz^ks^{-2}t^{-1}s^{-2}\\
		&\in&U+Rz^k(ts)^{-2}t+Rz^k(ts)^{-2}ts^{-1}+Rz^ks^{-1}(ts)^{-2}t+
		Rz^ks^{-1}(ts)^{-2}ts^{-1}\\
		&\in&U+\underline{z^{k-1}u_2}+\underline{z^{k-1}ts^{-1}u_2}+\underline{z^{k-1}u_1u_2}+Rz^{k-1}s^{-1}(R+Rt^{-1}+Rt^{-2})s^{-1}\\
		&\in&U+\underline{z^{k-1}u_1}+Rz^{k-1}(ts)^{-2}t+Rz^{k-1}(ts)^{-2}ts(st)^{-2}st\\
		&\in&U+\underline{z^{k-2}u_2}+\underline{\underline{(z^{k-3}tu_1)t}}.
		\end{array}}$\\
		\qedhere
	\end{itemize}

\end{proof}
	We can now prove the main theorem of this section.

\begin{thm} $H_{G_5}=U$.
	\label{thm55}
\end{thm}
\begin{proof}
	Since $1\in U$, it is enough to prove that $U$ is a left-sided ideal of $H_{G_5}$. For this purpose, one may check  that $sU$ and $tU$ are subsets of $U$. However, by proposition \ref{su55} we restrict ourselves to proving that $tU\subset U$. By the definition of $U$ we have that $tU\subset \sum\limits_{k=0}^5(z^ktu_1u_2+\underline{z^ku_1u_2})\subset U+\sum\limits_{k=0}^5z^ktu_1u_2.$ Therefore, by remark \ref{r55} it will be sufficient to prove  that, for every $k\in\{0,\dots,5\}$, $z^ktu_1\subset U$. However, this holds for every $k\in\{1,\dots,4\}$, by lemma \ref{ts55}(iii).
	For $k=0$ we have: $tu_1=(ts)^2s^{-1}t^{-1}u_1\subset u_1(\underline{\underline{zt^{-1}u_1}})\subset u_1U\stackrel{\ref{su55}}{\subset}U$. It remains to prove the case where $k=5$. We have:
	$$\small{\begin{array}{lcl}
	z^5tu_1&=&z^5t(R+Rs^{-1}+Rs^{-2})\\
	&\subset&\underline{z^5u_2}+Rz^5t(ts)^{-2}tst+Rz^5(R+Rt^{-1}+Rt^{-2})s^{-2}\\
	&\subset&U+\underline{\underline{(z^4u_2u_1)t}}+\underline{z^5u_1}+\underline{\underline{z^5t^{-1}u_1}}+Rz^5t^{-1}(st)^{-2}sts^{-1}\\
	&\subset&U+Rz^4t^{-1}(R+Rs^{-1}+Rs^{-2})ts^{-1}\\
	&\subset&U+\underline{z^4u_1}+Rz^4(st)^{-2}st^2s^{-1}+Rz^4(st)^{-2}sts^{-1}ts^{-1}\\
	&\subset&U+u_1\underline{\underline{z^3u_2u_1}}+Rz^3sts^{-1}(R+Rt^{-1}+Rt^{-2})s^{-1}\\
	&\subset&U+u_1U+u_1\underline{\underline{z^3tu_1}}+Rz^3st(ts)^{-2}t+Rz^3st(ts)^{-2}ts(st)^{-2}st\\
	&\subset&U+u_1U+\underline{z^2u_1u_2}+u_1\underline{\underline{(zu_2u_1)t}}\\
	&\subset&U+u_1U.
	\end{array}}$$
	The result follows from proposition \ref{su55}.
\end{proof}

\begin{cor}
	The BMR freeness conjecture holds for the generic Hecke algebra $H_{G_5}$.
\end{cor}
\begin{proof}
By theorem \ref{thm55} we have that $H_{G_5}=U$. Therefore,	the result follows from proposition \ref{BMR PROP} since, by definition, $U$ is generated as $R$-module by $|G_5|=72$ elements.
\end{proof}

\subsection{The case of $G_6$}
	\indent
	
	Let $R=\ZZ[u_{s,i}^{\pm},u_{t,j}^{\pm}]_{\substack{1\leq i\leq 2 \\1\leq j\leq 3}}$ and let $H_{G_6}:=\langle s,t\;|\; ststst=tststs,\prod\limits_{i=1}^{2}(s-u_{s,i})=\prod\limits_{j=1}^{3}(t-u_{t,j})=0\rangle$ be the generic Hecke algebra associated to $G_6$. Let $u_1$ be the subalgebra of $H_{G_5}$ generated by $s$ and $u_2$ the subalgebra of $H_{G_6}$ generated by $t$. We recall that $z:=(st)^3=(ts)^3$  generates the center of the associated complex braid group and that $|Z(G_6)|=4.$ We set $U=\sum\limits_{k=0}^3z^ku_2u_1u_2$. By the definition of $U$ we have the following remark:
	\begin{rem} $u_2Uu_2\subset U$.
		\label{r66}
	\end{rem}
	
	Our goal is to prove that $H_{G_6}=U$ (theorem \ref{thm66}). In order to do so, we first need to prove some preliminary results.
	\begin{lem}
		\mbox{}
		\vspace*{-\parsep}
		\vspace*{-\baselineskip}\\
		\begin{itemize}[leftmargin=0.8cm]
			\item [(i)] For every $k\in\{0,1,2\}$, $z^ku_1tu_1\subset U$.
			\item[(ii)]For every $k\in\{1,2,3\}$, $z^ku_1t^{-1}u_1\subset U$.
			\item[(iii)]For every $k\in\{1,2\}$, $z^ku_1u_2u_1\subset U$.
		\end{itemize}
		\label{sts66}
	\end{lem}
	\begin{proof}
		By definition, $u_2=R+Rt+Rt^{-1}$. Therefore, we only need to prove (i) and (ii). 
		\begin{itemize} [leftmargin=0.8cm]
			\item[(i)]$z^ku_1tu_1=z^k(R+Rs)t(R+Rs)\subset z^ku_2u_1u_2+z^ksts\subset U+z^k(st)^3t^{-1}s^{-1}t^{-1}\subset U+z^{k+1}u_2u_1u_2$.  The result follows from the definition of $U$.
			\item[(ii)] $z^ku_1t^{-1}u_1=z^k(R+Rs^{-1})t^{-1}(R+Rs^{-1})\subset z^ku_2u_1u_2+z^ks^{-1}t^{-1}s^{-1}\subset U+z^k(ts)^{-3}tst\subset U+z^{k-1}u_2u_1u_2\subset U$.
			\qedhere
			\end{itemize}
	\end{proof}
	We can now prove the main theorem of this section.
	\begin{thm} $H_{G_6}=U$.
		\label{thm66}
		\end{thm}
		\begin{proof}
			Since $1\in U$, it is enough to prove that $U$ is a left-sided ideal of $H_{G_6}$. For this purpose, one may check that $sU$ and $tU$ are subsets of $U$. However, by the definition of $U$, we only have to prove that $sU\subset U$. By the definition of $U$ and by remark \ref{r66}, we must prove that for every $k\in\{0,\dots,3\}$, $z^ksu_2u_1\subset U$. However, this holds for every $k\in\{1,2\}$, by lemma \ref{sts66}(iii).
			
			 For $k=0$ we have: $su_2u_1\subset s(R+Rt+Rt^2)(R+Rs)\subset u_2u_1u_2+u_1tu_1+Rst^2s$. By the definition of $U$ and lemma \ref{sts66}(i), it will be sufficient to prove that $st^2s\in U$. We have:
			$$\small{\begin{array}[t]{lcl}st^2s&=&(st)^3t^{-1}s^{-1}t^{-1}s^{-1}ts\\
			&\in& zu_2s^{-1}t^{-1}(R+Rs)ts\\
			&\in&zu_2u_1u_2+zu_2s^{-1}t^{-1}(st)^3t^{-1}s^{-1}t^{-1}\\
			&\in& U+u_2(z^2u_1u_2u_1)u_2.			
			\end{array}}$$
			The result follows from lemma \ref{sts66}(iii) and remark \ref{r66}.
			
			It remains to prove the case where $k=3$. We have: $z^3su_2u_1\subset s(R+Rt^{-1}+Rt^{-2})(R+Rs^{-1})\subset u_2u_1u_2+u_1t^{-1}u_1+Rst^{-2}s^{-1}$. By the definition of $U$ and lemma \ref{sts66}(ii), we need to prove that $st^{-2}s^{-1}\in U$. We have:
			$$\small{\begin{array}[t]{lcl}
			z^3st^{-2}s^{-1}&\in&z^3(R+Rs^{-1})t^{-2}s^{-1}\\
			&\in&z^3u_2u_1u_2+Rz^3(ts)^{-3}tstst^{-1}s^{-1}\\
			&\in&U+z^2u_2st(R+Rs^{-1})t^{-1}s^{-1}\\
			&\in&U+z^2u_2u_1u_2+z^2u_2st(ts)^{-3}tst\\
			&\in&U+u_2(zu_1u_2u_1)u_2.
			\end{array}}$$
			The result follows again from lemma \ref{sts66}(iii) and remark \ref{r66}.
		\end{proof}
	\begin{cor}
		The BMR freeness conjecture holds for the generic Hecke algebra $H_{G_6}$.
	\end{cor}
	\begin{proof}
		By theorem \ref{thm66} we have that $H_{G_6}=U$. By definition, $U=\sum\limits_{k=0}^3z^ku_2u_1u_2$. We expand $u_1$ as $R+Rs$ and we have that $U=\sum\limits_{k=0}^3(z^ku_2+z^ku_2su_2)$. Therefore, $U$ is generated as $u_2$-module by 16 elements. Since $u_2$ is generated as $R$-module by 3 elements, we have that $H_{G_6}$ is generated as
		$R$-module by $|G_6|=48$ elements and the result follows from proposition \ref{BMR PROP}.
	\end{proof}
\subsection{The case of $G_7$}

Let $R=\ZZ[u_{s,i}^{\pm},u_{t,j}^{\pm},u_{u,l}^{\pm}]_{\substack{1\leq i\leq 2 \\1\leq j,l\leq 3}}$. We also let $$H_{G_{7}}=\langle s,t,u\;|\; stu=tus=ust, \;\prod\limits_{i=1}^{2}(s-u_{s,i})=\prod\limits_{j=1}^{3}(t-u_{t,j})=\prod\limits_{l=1}^{3}(u-u_{u,l})=0\rangle$$ be the generic Hecke algebra associated to $G_{7}$. Let $u_1$ be the subalgebra of $H_{G_{7}}$ generated by $s$, $u_2$ the subalgebra of $H_{G_{7}}$ generated by $t$ and $u_3$ the subalgebra of $H_{G_{7}}$ generated by $u$. We recall that $z:=stu=tus=ust$  generates the center of the associated complex braid group and that $|Z(G_7)|=12$.
We set $U=\sum\limits_{k=0}^{11}(z^ku_3u_2+z^ktu^{-1}u_2).$
By the definition of $U$, we have the following remark.
\begin{rem}
	$Uu_2 \subset U$.
	\label{r77}
\end{rem}
To make it easier for the reader to follow the calculations, we will underline the elements that  belong to $U$ by definition.  Moreover, we will use directly  remark \ref{r77}; this means that every time we have a power of $t$ at the end of an element, we may ignore it. In order to remind that to the reader, we put a parenthesis around the part of the element we consider.

Our goal is to prove that $H_{G_{7}}=U$ (theorem \ref{thm77}). In order to do so, we first need  to prove some preliminary results.

\begin{lem}
	\mbox{}
	\vspace*{-\parsep}
	\vspace*{-\baselineskip}\\
	\begin{itemize}[leftmargin=0.8cm]
		\item[(i)]For every $k\in\{0,\dots,10\}$, $z^ku_1\subset U$.
		\item [(ii)] For every $k\in\{0,\dots,9\}$, $z^ku_2u\subset U$.
		\item[(iii)]For every $k\in\{1,\dots,10\}$, $z^ku_1u_3\subset U$.
	\label{l77}
\end{itemize}
\end{lem}
\begin{proof}
	\mbox{}
	\vspace*{-\parsep}
	\vspace*{-\baselineskip}\\
	\begin{itemize}[leftmargin=0.8cm]
		\item[(i)] 
		
		$z^ku_1=z^k(R+Rs)\subset \underline{z^ku_3}+z^k(stu)u^{-1}t^{-1}\subset U+\underline{z^{k+1}u_3u_2}$.
		\item[(ii)] $\begin{array}[t]{lcl}
			z^ku_2u&=&z^k(R+Rt+Rt^2)u\\
			&\subset& \underline{z^ku_3}+z^k(tus)s^{-1}+z^kt(tus)s^{-1}\\
			&\subset& U+z^{k+1}u_1+z^{k+1}t(R+Rs)\\
			&\stackrel{(i)}{\subset}&U+\underline{z^{k+1}u_2}+z^{k+1}t(stu)u^{-1}t^{-1}\\&\subset& U+\underline{z^{k+2}tu^{-1}u_2}.
			\end{array}$
			\item[(iii)] $z^ku_1u_3=z^k(R+Rs^{-1})u_3\subset \underline{z^ku_3}+z^k(s^{-1}u^{-1}t^{-1})tu_3\subset U+z^{k-1}tu_3$.
			
			The result follows from the definition of $U$ and (i), if we expand $u_3$ as $R+Ru^{-1}+Ru$.
		\qedhere
			\end{itemize}
\end{proof}

From now on, we will double-underline the elements described in lemma \ref{l77} and we will use directly the fact that these elements are inside $U$. 
\begin{prop}$u_3U\subset U$.
	\label{pr77}
\end{prop}
\begin{proof}
	Since $u_3=R+Ru+Ru^2$, it is enough to prove that $uU\subset U$. By the definition of $U$ and remark \ref{r77} it will be sufficient to prove that for every $k\in\{0,\dots,11\}$, $z^kutu^{-1}\in U$. We distinguish the following cases:
	\begin{itemize}[leftmargin=*]
		\item \underline{$k\in\{0,\dots,7\}$}:
		
		$\small{\begin{array}[t]{lcl}
			z^kutu^{-1}&\in&z^kut(R+Ru+Ru^2)\\
			&\in&\underline{z^ku_3}+z^ku(tus)s^{-1}+z^ku(tus)s^{-1}u\\
			&\in&U+Rz^{k+1}u(R+Rs)+Rz^{k+1}u(R+Rs)u\\
			&\in&U+\underline{z^{k+1}u_3}+Rz^{k+1}(ust)t^{-1}+Rz^{k+1}(ust)t^{-1}u\\
			&\in&U+\underline{z^{k+2}u_2}+\underline{\underline{z^{k+2}u_2u}}.
			\end{array}}$
			\item \underline{$k\in\{8,\dots,11\}$}:
			
			$\small{\begin{array}[t]{lcl}
			z^kutu^{-1}&\in&z^ku(R+Rt^{-1}+Rt^{-2})u^{-1}\\
			&\in&\underline{z^ku_3}+Rz^ku(t^{-1}s^{-1}u^{-1})usu^{-1}+Rz^kut^{-2}(u^{-1}t^{-1}s^{-1})st\\
			&\in&U+z^{k-1}u_3(R+Rs^{-1})u^{-1}+Rz^{k-1}ut^{-2}(R+Rs^{-1})t\\
			&\in&U+\underline{z^{k-1}u_3}+Rz^{k-1}u_3(s^{-1}u^{-1}t^{-1})t+\underline{z^{k-1}u_3u_2}+
			Rz^{k-1}ut^{-1}(t^{-1}s^{-1}u^{-1})ut\\
				&\in&U+\underline{z^{k-2}u_3u_2}+Rz^{k-2}(R+Ru^{-1}+Ru^{-2})t^{-1}ut\\
				&\in&U+\underline{\underline{(z^{k-2}u_2u)t}}+Rz^{k-2}(u^{-1}t^{-1}s^{-1})sut+Rz^{k-2}u^{-1}(u^{-1}t^{-1}s^{-1})sut\\
				&\in&U+\underline{\underline{(z^{k-3}u_1u_3)t}}+Rz^{k-3}u^{-1}s(R+Ru^{-1}+Ru^{-2})t\\
				&\in&U+Rz^{k-3}u^{-1}(stu)u^{-1}+Rz^{k-3}u^{-1}(R+Rs^{-1})u^{-1}t+
				Rz^{k-3}u^{-1}(R+Rs^{-1})u^{-2}t\\
				&\in&U+\underline{z^{k-2}u_3}+\underline{z^{k-3}u_3u_2}+Rz^{k-3}u^{-1}(s^{-1}u^{-1}t^{-1})t^2+Rz^{k-3}u^{-1}(s^{-1}u^{-1}t^{-1})tu^{-1}t\\
				&\in&U+\underline{z^{k-4}u_3u_2}+Rz^{k-4}u^{-1}(R+Rt^{-1}+Rt^{-2})u^{-1}t\\
				&\in&U+\underline{z^{k-4}u_3u_2}+Rz^{k-4}(u^{-1}t^{-1}s^{-1})su^{-1}t+
				Rz^{k-4}(u^{-1}t^{-1}s^{-1})st^{-1}u^{-1}t\\
				&\in&U+\underline{\underline{(z^{k-5}u_1u_3)t}}+Rz^{k-5}s(t^{-1}s^{-1}u^{-1})usu^{-1}t\\
				&\in&U+Rz^{k-6}su(R+Rs^{-1})u^{-1}t\\
				&\in&U+\underline{\underline{(z^{k-6}u_1)t}}+Rz^{k-6}su(s^{-1}u^{-1}t^{-1})t^2\\
					&\in&U+\underline{\underline{(z^{k-7}u_1u_3)t}}.
		\end{array}}$\\
		\qedhere
	\end{itemize}
\end{proof}
We can now prove the main theorem of this section.
\begin{thm} $H_{G_7}=U$.
	\label{thm77}
\end{thm}

\begin{proof}
	Since $1\in U$, it will be sufficient to prove that $U$ is a left-sided ideal of $H_{G_7}$. For this purpose, one may check that  $sU$, $tU$ and $uU$ are subsets of $U$. However, by proposition \ref{pr77} we
	only have to prove that $tU$ and $sU$ are subsets of $U$. We recall that $z=stu$, therefore $s=zu^{-1}t^{-1}$ and $s^{-1}=z^{-1}tu$. We notice that $U=
	\sum\limits_{k=0}^{10}z^k(u_3u_2+tu^{-1}u_2)+z^{11}(u_3u_2+tu^{-1}u_2).$
	Hence, \\
	$\small{\begin{array}[t]{lcl}sU&\subset&\sum\limits_{k=0}^{10}
	z^ks(u_3u_2+tu^{-1}u_2)+
	z^{11}s(u_3u_2+tu^{-1}u_2)\\
	&\subset& \sum\limits_{k=0}^{10}z^{k+1}u^{-1}t^{-1}(u_3u_2+tu^{-1}u_2)+z^{11}(R+Rs^{-1})(u_3u_2+
	tu^{-1}u_2)\\
	&\subset& \sum\limits_{k=0}^{10}u^{-1}t^{-1}(\underline{z^{k+1}u_3u_2}+\underline{z^{k+1}tu^{-1}u_2})+
	\underline{z^{11}u_3u_2}+\underline{z^{11}tu^{-1}u_2}+z^{11}s^{-1}u_3u_2+
	z^{11}s^{-1}tu^{-1}u_2\\ 
	&\subset&u_3u_2U+z^{10}tu_3u_2+
	z^{10}tutu^{-1}u_2\\ 
	&\subset&u_3u_2U+
	t\underline{z^{10}u_3u_2}+tu(\underline{z^{10}tu^{-1}u_2})\\
	&\subset&u_3u_2u_3U.
	\end{array}}$\\\\
	By proposition \ref{pr77} we have that $u_3U\subset U$. 
	 If we also suppose that $u_2U\subset U$ then, obviously, we have  $tU\subset U$ but we  also have  $sU\subset U$ (since $sU\subset u_3u_2u_3U$). Hence, in order to prove that $U=H_{G_7}$ we restrict ourselves to proving that $u_2U \subset U$.
	 
	 By definition, $u_2=R+Rt^{-1}+Rt^{-2}$, therefore it will be sufficient to prove that $t^{-1}U\subset U$.
	By the definition of $U$ and remark \ref{r77}, this is the same as proving that, for every $k\in\{0,\dots,11\}$, $z^kt^{-1}u_3\subset U$. 
	For $k\in\{2,\dots,11\}$ the result is obvious, since $z^kt^{-1}u_3=z^k(t^{-1}s^{-1}u^{-1})usu_3\subset u_3(\underline{\underline{z^{k-1}u_1u_3}})\subset u_3U\stackrel{\ref{pr77}}{\subset}U$. For $k\in\{0,1\}$, we have:
	$$\small{\begin{array}{lcl}
	z^kt^{-1}u_3&=&z^kt^{-1}(R+Ru+Ru^2)\\
	&\subset& \underline{z^ku_2}+\underline{\underline{z^ku_2u}}+z^k(R+Rt+Rt^2)u^2\\
	&\subset&U+\underline{z^ku_3}+Rz^kt(R+Ru+Ru^{-1})+Rz^kt(tus)s^{-1}u\\
	&\subset&U+\underline{z^ku_2}+\underline{\underline{z^ku_2u}}+\underline{z^ktu^{-1}u_2}+Rz^{k+1}(tus)s^{-1}u^{-1}s^{-1}u\\
	&\subset&U+Rz^{k+2}(R+Rs)u^{-1}(R+Rs)u\\
	&\subset&U+\underline{\underline{z^{k+2}u_1}}+u_3\underline{\underline{z^{k+2}u_1u_3}}+Rz^{k+2}su^{-1}su\\
		&\subset&U+u_3U+Rz^{k+2}s(R+Ru+Ru^2)su\\
		&\subset&U+u_3U+\underline{\underline{z^{k+2}u_1u_3}}+Rz^{k+2}s(ust)t^{-1}u+Rz^{k+2}su(ust)t^{-1}u\\
		&\subset&U+u_3U+Rz^{k+3}(stu)u^{-1}t^{-2}u+Rz^{k+3}su(R+Rt+Rt^2)u\\
		&\subset&U+u_3U+u_3\underline{\underline{z^{k+4}u_2u}}+\underline{\underline{z^{k+3}u_1u_3}}+Rz^{k+3}su(tus)s^{-1}+Rz^{k+3}sut(tus)s^{-1}\\
		&\subset&U+u_3U+Rz^{k+4}su(R+Rs)+Rz^{k+4}sut(R+Rs)\\
		&\subset&U+u_3U+\underline{\underline{z^{k+4}u_1u_3}}+Rz^{k+4}s(ust)t^{-1}+
		\underline{\underline{(z^{k+4}u_1u_3)t}}+Rz^{k+4}(stu)u^{-1}t^{-1}uts\\
			&\subset&U+u_3U+\underline{\underline{(z^{k+5}u_1)t}}+z^{k+5}u_3(R+Rt+Rt^2)uts\\
			&\subset&U+u_3U+z^{k+5}u_3t(stu)u^{-1}t^{-1}+z^{k+5}u_3(tus)s^{-1}ts+
			z^{k+5}u_3t(tus)s^{-1}ts\\
				&\subset&U+u_3U+u_3\underline{z^{k+5}tu^{-1}u_2}+z^{k+5}u_3(R+Rs)ts+
				z^{k+6}u_3t(R+Rs)ts\\
				&\subset&U+u_3U+z^{k+5}u_3t(stu)u^{-1}t^{-1}+z^{k+5}u_3(stu)u^{-1}s+
				z^{k+6}u_3t^2(R+Rs^{-1})+\\&&+z^{k+6}u_3t(stu)u^{-1}s\\
				&\subset&U+u_3U+u_3\underline{z^{k+6}tu^{-1}u_2}+u_3\underline{z^{k+6}u_1}+\underline{z^{k+6}u_3u_2}+z^{k+6}u_3t^2(s^{-1}u^{-1}t^{-1})tu+\\&&+
				z^{k+6}u_3t(R+Ru+Ru^2)s\\
					&\subset&U+u_3U+u_3\underline{\underline{z^{k+5}u_2u}}+z^{k+6}u_3t(stu)u^{-1}t^{-1}+z^{k+6}u_3(tus)+
					z^{k+6}u_3(tus)s^{-1}(ust)t^{-1}\\
						&\subset&U+u_3U+u_3\underline{z^{k+7}tu^{-1}u_2}+\underline{z^{k+7}u_3}+u_3\underline{\underline{(z^{k+8}u_1)t^{-1}}}\\
						&\subset&U+u_3U.
	\end{array}}$$
The result follows from  proposition \ref{pr77}.	
\end{proof}
\begin{cor}
	The BMR freeness conjecture holds for the generic Hecke algebra $H_{G_7}$.
\end{cor}
\begin{proof}
	By theorem \ref{thm77} we have that $H_{G_7}=U$. The result follows from proposition \ref{BMR PROP}, since by definition $U$ is generated as a right $u_2$-module by 48 elements and, hence, as $R$-module by $|G_7|=144$ elements (recall that $u_2$ is generated as $R$-module by 3 elements).
\end{proof}

\section{The Octahedral family}
\indent

In this family we encounter the exceptional groups $G_8$, $G_9$, $G_{10}$, $G_{11}$, $G_{12}$, $G_{13}$, $G_{14}$ and $G_{15}$. In Chapter 2 we proved that the BMR freeness conjecture holds for $G_8$ (see corollary \ref{G8}). Moreover, we know the validity of the conjecture for the group $G_{12}$ (see theorem \ref{2case}). In this section we prove the conjecture for the rest of the exceptional groups in this family using a case-by-case analysis. We also reprove the validity of the BMR freeness conjecture for the exceptional group $G_{12}$. As in the tetrahedral case, we use directly the relations (\ref{ooo}) and (\ref{oooo}).

\subsection{The case of $G_9$}
\indent

Let $R=\ZZ[u_{s,i}^{\pm},u_{t,j}^{\pm}]_{\substack{1\leq i\leq 2 \\1\leq j\leq 4}}$ and let $H_{G_{9}}=\langle s,t\;|\; ststst=tststs, \prod\limits_{i=1}^{2}(s-u_{s,i})=\prod\limits_{j=1}^{4}(t-u_{t,j})=0\rangle$ be the generic Hecke algebra associated to $G_{9}$. Let $u_1$ be the subalgebra of $H_{G_{9}}$ generated by $s$ and $u_2$ be the subalgebra of $H_{G_{9}}$ generated by $t$.
We recall that $z:=(st)^3=(ts)^3$  generates the center of the associated complex braid group and that $|Z(G_9)|=8$. We set $U=\sum\limits_{k=0}^7(z^ku_2u_1u_2+z^ku_2st^{-2}s)$. By the definition of $U$, we have the following remark.
\begin{rem} $u_2U\subset U$.
	\label{rem9}
\end{rem}
From now on, we will underline the elements that by definition belong to $U$.  Moreover, we will use directly the remark \ref{rem9}; this means that every time we have a power of $t$ at the beginning of an element, we may ignore it. In order to remind that to the reader, we put a parenthesis around the part of the element we consider.

Our goal is to prove that $H_{G_{9}}=U$ (theorem \ref{thm9}). The next proposition provides the necessary conditions for this to be true.

\begin{prop}
	If $z^ku_1u_2u_1\subset U$ and $z^kst^{-2}st\in U$ for every $k\in\{0,\dots,7\}$, then $H_{G_9}=U$.
	\label{pr9}
\end{prop}
\begin{proof}
	Since $1\in U$, it is enough to prove that $U$ is a right-sided ideal of $H_{G_9}$. For this purpose, one may check that $Us$ and $Ut$ are subsets of $U$. By the definition of $U$ we have that  $Us\subset \sum\limits_{k=0}^7u_2(z^ku_1u_2u_1)$ and $Ut\subset \sum\limits_{k=0}^7(\underline{z^ku_2u_1u_2}+z^ku_2st^{-2}st)\subset U+\sum\limits_{k=0}^7u_2(z^kst^{-2}st).$ The result follows from hypothesis and remark \ref{rem9}.
\end{proof}
As a first step, we prove the conditions of the above proposition for a smaller range of the values of $k$.
\begin{lem}
	\mbox{}
	\vspace*{-\parsep}
	\vspace*{-\baselineskip}\\
	\begin{itemize}[leftmargin=0.8cm]	
		\item[(i)] For every $k\in\{0,\dots,6\}$, $z^ku_1tu_1\subset U$.
			\item[(ii)] For every $k\in\{1,\dots,7\}$, $z^ku_1t^{-1}u_1\subset U$.
				\item[(iii)] For every $k\in\{1,\dots,6\}$, $z^ku_1u_2u_1\subset U$.
				
			\item[(iv)]For every $k\in\{0,\dots,5\}$, $z^ku_1t^2u_1t\subset U+z^{k+2}u_2u_1u_2u_1$. Therefore, for every $k\in\{0,\dots,4\}$, 
			$z^ku_1t^2st\subset U$.
			\item[(v)]For every $k\in\{1,\dots,5\}$, $z^ku_1u_2u_1t\subset U+z^{k+2}u_2u_1u_2u_1$. Therefore, for every $k\in\{0,\dots,4\}$, $z^ku_1u_2u_1t\subset U$.
			\qedhere
	\end{itemize}
	\label{lem9}
\end{lem}
\begin{proof}
	\mbox{}
	\vspace*{-\parsep}
	\vspace*{-\baselineskip}\\
	\begin{itemize}[leftmargin=0.8cm]
		\item [(i)] $z^ku_1tu_1=z^k(R+Rs)t(R+Rs)\subset \underline{z^ku_2u_1u_2}+Rz^ksts\subset U+Rz^k(st)^3t^{-1}s^{-1}t^{-1}\subset U+\underline{z^{k+1}u_2u_1u_2}$.
		\item [(ii)] $z^ku_1t^{-1}u_1=z^k(R+Rs^{-1})t^{-1}(R+Rs^{-1})\subset \underline{z^ku_2u_1u_2}+Rz^ks^{-1}t^{-1}s^{-1}\subset U+Rz^k(ts)^{-3}tst\subset U+\underline{z^{k-1}u_2u_1u_2}$.
		\item[(iii)]Since $u_2=R+Rt+Rt^{-1}+Rt^{-2}$, by the definition of $U$ and by (i) and (ii) we only have to prove that, for every $k\in\{1,\dots,6\}$,  $z^ku_1t^{-2}u_1 \subset U$. Indeed, $z^ku_1t^{-2}u_1=z^k(R+Rs)t^{-2}(R+Rs)\subset \underline{z^ku_2u_1u_2}+\underline{z^ku_2st^{-2}s}$.
		\item [(iv)] 
		$\small{\begin{array}[t]{lcl}
		 z^ku_1t^2u_1t&=&z^ku_1t^2(R+Rs)t\\
		 &\subset& \underline{z^ku_1u_2}+z^ku_1t^2st\\
		 &\subset& U+z^k(R+Rs)t^2st\\
		 &\subset& U+ \underline{z^ku_2u_1u_2}+Rz^kst(ts)^3s^{-1}t^{-1}s^{-1}\\
		 &\subset& U+Rz^{k+1}st(R+Rs)t^{-1}s^{-1}\\
		 &\subset& U+\underline{z^{k+1}u_2}+Rz^{k+1}(st)^3t^{-1}s^{-1}t^{-2}s^{-1}\\
		 &\subset& U+z^{k+2}u_2u_1u_2u_1.
		 \end{array}}$
		 
		 Since $z^{k+2}u_2u_1u_2u_1=u_2(z^{k+2}u_1u_2u_1)$ and since, for every $k\in\{0,\dots, 4\}$, we have $k+2\in \{2,\dots,6\}$, we can use (iii) and we have that for every $k\in\{0,\dots,4\}$, $z^ku_1t^2st\subset U$.

		 \item [(v)] 
		 $\small{\begin{array}[t]{lcl}
		 z^ku_1u_2u_1t&=&z^ku_1(R+Rt+Rt^{-1}+Rt^2)u_1t\\
		 &\subset& \underline{z^ku_1u_2}+
		 z^k(R+Rs)t(R+Rs)t+z^k(R+Rs^{-1})t^{-1}(R+Rs^{-1})t+z^ku_1t^2u_1t\\
		 &\stackrel{(iv)}{\subset}& U+\underline{z^ku_2u_1u_2}+Rz^kstst+Rz^ks^{-1}t^{-1}s^{-1}t+z^{k+2}u_2u_1u_2u_1\\
		 &\subset& U+Rz^k(st)^3t^{-1}s^{-1}+Rz^k(ts)^{-3}tst^2+z^{k+2}u_2u_1u_2u_1\\
		 &\subset& U+\underline{(z^{k+1}+z^{k-1})u_2u_1u_2}+z^{k+2}u_2u_1u_2u_1\\
		 &\subset& U+z^{k+2}u_2u_1u_2u_1.
		 \end{array}}$
		 
		 Using the same arguments as in (iv), we can use (iii) and we have that for every $k\in\{0,\dots,4\}$, $z^ku_1t^2st\subset U$.
		  \qedhere
	\end{itemize}
\end{proof}
To make it easier for the reader to follow the calculations, we will double-underline the elements described in lemma \ref{lem9} and  we will use directly the fact that these elements are inside $U$. The next proposition proves the first condition of \ref{pr9}.
\begin{prop}For every $k\in\{0,\dots,7\}$, $z^ku_1u_2u_1\subset U$.
	\label{prr9}
\end{prop}
\begin{proof}
	By lemma \ref{lem9}(iii), we need to prove the cases where $k\in\{0,7\}$. We have:
	\begin{itemize}[leftmargin=*]
		\item \underline{$k=0$}:
		
		$\small{\begin{array}[t]{lcl}u_1u_2u_1&=&(R+Rs)(R+Rt+Rt^{-2}+Rt^3)(R+Rs)\\
		&\subset& \underline{u_2u_1u_2}+\underline{\underline{u_1tu_1}}+\underline{Rst^{-2}s}+Rst^3s\\
		&\subset&U+Rt^{-1}(ts)^3s^{-1}t^{-1}s^{-1}t^2s\\
		&\subset&U+zu_2s^{-1}t^{-1}(R+Rs)t^2s\\
		&\subset&U+\underline{\underline{u_2(zu_1tu_1)}}+zu_2(R+Rs)t^{-1}st^2s\\
				&\subset&U+\underline{\underline{u_2(zu_1u_2u_1)}}+zu_2s(R+Rt+Rt^2+Rt^3)st^2s\\
				&\subset&U+\underline{\underline{u_2(zu_1u_2u_1)}}+zu_2(st)^3t^{-1}s^{-1}ts+zu_2st(ts)^3s^{-1}t^{-1}s^{-1}ts+zu_2st^2(ts)^3s^{-1}t^{-1}s^{-1}ts\\
				&\subset&U+\underline{\underline{u_2(z^2u_1tu_1)}}+z^2u_2sts^{-1}t^{-1}(R+Rs)ts+z^2u_2st^2s^{-1}t^{-1}(R+Rs)ts\\
				&\subset&U+\underline{z^2u_2u_1u_2}+z^2u_2st(R+Rs)t^{-1}sts+z^2u_2st^2(R+Rs)t^{-1}sts\\
				&\subset&U+\underline{\underline{u_2(z^2u_1tu_1})}+z^2u_2stst^{-1}sts+z^2u_2(st)^3t^{-1}+
				z^2u_2st^2st^{-1}sts\\
					&\subset&U+z^2u_2stst^{-1}sts+\underline{z^3u_2}+
					z^2u_2st^2st^{-1}sts
				\end{array}}$
				
				It remains to prove that $A:=z^2u_2stst^{-1}sts+
				z^2u_2st^2st^{-1}sts$ is a subset of $U$. For this purpose, we expand $t^{-1}$ as a linear combination of 1, $t$, $t^2$ and $t^3$ and we have:

				$\small{\begin{array}{lcl}
				A&\subset&z^2u_2sts(R+Rt+Rt^2+Rt^3)sts+
				z^2u_2st^2s(R+Rt+Rt^2+Rt^3)sts\\

				&\subset&z^2u_2sts^2ts+z^2u_2(st)^3s+z^2u_2(ts)^3s^{-2}(st)^3t^{-1}+
				z^2u_2(ts)^3s^{-1}t(ts)^3s^{-1}t^{-1}+z^2u_2st^2s^2ts+\\&&+Rz^2st(ts)^3+z^2u_2st(ts)^3s^{-1}t^{-1}s^{-1}(ts)^3s^{-1}t^{-1}+z^2u_2st^2st^2(ts)^3s^{-1}t^{-1}\\ 
				&\subset&z^2u_2st(R+Rs)ts+\underline{z^3u_2u_1}+\underline{z^4u_2u_1u_2}+
				z^4u_2(R+Rs)t(R+Rs)t^{-1}+z^2u_2st^2(R+Rs)ts+\\&&+\underline{z^3u_2u_1u_2}+z^4u_2sts^{-1}t^{-1}s^{-2}t^{-1}+z^3u_2st^2st^2(R+Rs)t^{-1}\\ 
				&\subset&U+\underline{\underline{u_2(z^2u_1u_2u_1)}}+z^2u_2(st)^3t^{-1}+\underline{z^4u_2u_1u_2}+z^4u_2stst^{-1}+\underline{\underline{u_2(z^2u_1u_2u_1)}}+
				z^2u_2st^2sts+\\&&+z^4u_2sts^{-1}t^{-1}(R+Rs^{-1})t^{-1}+\underline{\underline{u_2(z^3u_1t^2st)}}+z^3u_2st^2st^2st^{-1}\\ 
				&\subset&U+\underline{z^3u_2}+z^4u_2(ts)^3s^{-1}t^{-2}+z^2u_2st(ts)^3s^{-1}t^{-1}+z^4u_2st(R+Rs)t^{-2}+z^4u_2st^2(st)^{-3}s+\\&&+z^3u_2st(ts)^3s^{-1}t^{-1}s^{-1}tst^{-1}\\ 
				
				&\subset&U+\underline{z^5u_2u_1u_2}+z^3u_2st(R+Rs)t^{-1}+\underline{z^4u_2u_1u_2}+z^4u_2(st)^3t^{-1}s^{-1}t^{-3}+
				\underline{\underline{u_2(z^3u_1u_2u_1)}}+\\&&+z^4u_2st(R+Rs)t^{-1}s^{-1}tst^{-1}\\
				&\subset&U+\underline{(z^3+z^4+z^5)u_2u_1u_2}+z^3u_2(ts)^3s^{-1}t^{-2}+
				z^4u_2(ts)^3s^{-1}t^{-2}(R+Rs)tst^{-1}\\
				&\subset&U+\underline{z^4u_2u_1u_2}+z^5u_2(st)^{-3}sts^2t^{-1}+z^5u_2s^{-1}t^{-3}(ts)^3s^{-1}t^{-2}\\
				&\subset&U+z^4u_2st(R+Rs)t^{-1}+z^6u_2s^{-1}(R+Rt^{-1}+Rt^{-2}+Rt)s^{-1}t^{-2}\\
				&\subset&U+\underline{z^4u_2u_1}+z^4u_2(ts)^3s^{-1}t^{-2}+
				\underline{z^6u_2u_1u_2}+z^6u_2(st)^{-3}st^{-1}+z^6u_2(st)^{-3}sts^2(ts)^{-3}tst^{-1}+\\&&+z^6u_2(R+Rs)t(R+Rs)t^{-2}\\ 
				&\subset&U+\underline{(z^5+z^6)u_2u_1u_2}+z^4u_2st(R+Rs)tst^{-1}+
				z^6u_2stst^{-2}\\
				&\subset&U+z^4u_2st^2st^{-1}+z^4u_2(ts)^3+z^6u_2(ts)^3s^{-1}t^{-3}\\
				&\subset&U+z^4u_2s(R+Rt+Rt^{-1}+Rt^{-2})st^{-1}+\underline{z^5u_2}+\underline{z^7u_2u_1u_2}\\
				&\subset&U+\underline{z^4u_2u_1u_2}+z^4u_2(ts)^3s^{-1}t^{-2}+
				z^4u_2(R+Rs^{-1})t^{-1}(R+Rs^{-1})t^{-1}+\\&&+z^4u_2(R+Rs^{-1})t^{-2}(R+Rs^{-1})t^{-1}\\
				&\subset&U+\underline{(z^4+z^5)u_2u_1u_2}+z^4u_2s^{-1}t^{-1}s^{-1}t^{-1}+
				z^4u_2s^{-1}t^{-2}s^{-1}t^{-1}\\
				&\subset&U+z^4u_2(st)^{-3}s+z^4u_2s^{-1}t^{-1}(ts)^{-3}sts\\
				&\subset&U+\underline{z^3u_2u_1}+z^3u_2s^{-1}t^{-1}(R+Rs^{-1})ts\\
				&\subset&U+\underline{z^3u_2}+z^3u_2(st)^{-3}st^2s\\
				&\subset& U+\underline{\underline{u_2(z^2u_1u_2u_1)}}.
						\end{array}}$
		\item \underline{$k=7$}:
		
			$\small{\begin{array}[t]{lcl}z^7u_1u_2u_1&=&z^7(R+Rs)(R+Rt^{-1}+Rt^{-2}+Rt^{-3})(R+Rs)\\
		&\subset& \underline{z^7u_2u_1u_2}+\underline{\underline{z^7u_1t^{-1}u_1}}+\underline{Rz^7st^{-2}s}+Rz^7st^{-3}s\\
		&\subset&U+Rz^7(R+Rs^{-1})t^{-3}(R+Rs^{-1})\\
		&\subset&U+\underline{z^7u_2u_1u_2}+Rz^7s^{-1}t^{-3}s^{-1}\\
		&\subset&U+Rz^7(ts)^{-3}tstst^{-2}s^{-1}\\
		&\subset&U+Rz^6tst(R+Rs^{-1})t^{-2}s^{-1}\\
		&\subset&U+\underline{\underline{t(z^6u_1u_2u_1)}}+Rz^6t(R+Rs^{-1})ts^{-1}t^{-2}s^{-1}\\
		&\subset&U+\underline{\underline{t^2(z^6u_1u_2u_1)}}++Rz^6ts^{-1}(R+Rt^{-1}+Rt^{-2}+Rt^{-3})s^{-1}t^{-2}s^{-1}\\

		&\subset&U+\underline{\underline{t(z^6u_1u_2u_1)}}+Rz^6t^2(st)^{-3}st^{-1}s^{-1}+Rz^6ts^{-1}t^{-1}(st)^{-3}stst^{-1}s^{-1}+\\&&+Rz^6ts^{-1}t^{-2}(st)^{-3}stst^{-1}s^{-1}\\
		&\subset&U+\underline{\underline{t^2(z^5u_1u_2u_1)}}+Rz^5ts^{-1}t^{-1}(R+Rs^{-1})t(R+Rs^{-1})t^{-1}s^{-1}+\\&&+
		Rz^5ts^{-1}t^{-2}st(R+Rs^{-1})t^{-1}s^{-1}\\

		&\subset&U+\underline{\underline{t(z^5u_1u_2u_1)}}+Rz^5t(sts)^{-1}t(sts)^{-1}+Rz^5ts^{-1}t^{-2}(R+Rs^{-1})ts^{-1}t^{-1}s^{-1}\\
		&\subset&U+Rz^5t(ts)^{-3}tst^3(st)^{-3}st+Rz^5t(ts)^{-3}t+Rz^5ts^{-1}t^{-2}s^{-1}t^2(st)^{-3}st\\
			&\subset&U+\underline{\underline{t^2(z^3u_1u_2u_1t)}}+\underline{z^4u_2}+
			Rz^4ts^{-1}t^{-2}s^{-1}t^2st.
			\end{array}}$
		
		It remains to prove that the element $z^4ts^{-1}t^{-2}s^{-1}t^2st$ is inside $U$. For this purpose, we expand $t^2$ as a linear combination of 1, $t$, $t^{-1}$ and $t^{-2}$ and we have:
		
		$\small{\begin{array}{lcl}
		z^4ts^{-1}t^{-2}s^{-1}t^2st&\in&				z^4ts^{-1}t^{-2}s^{-1}(R+Rt+Rt^{-1}+Rt^{-2})st\\
			&\in&U+\underline{z^4u_2u_1u_2}+Rz^4ts^{-1}t^{-2}(R+Rs)tst+
			Rz^4ts^{-1}t^{-2}s^{-1}t^{-1}(R+Rs^{-1})t+\\&&+Rz^4ts^{-1}t^{-2}s^{-1}t^{-2}(R+Rs^{-1})t\\
			&\in&U+\underline{\underline{t(z^4u_1u_2u_1t)}}+Rz^4ts^{-1}t^{-3}(ts)^3s^{-1}+\underline{\underline{t(z^4u_1u_2u_1)}}+
			Rz^4ts^{-1}t^{-1}(st)^{-3}st^2+\\&&+Rz^4ts^{-1}t^{-1}(st)^{-3}sts+
			Rz^4ts^{-1}t^{-1}(st)^{-3}stst^{-1}s^{-1}t\\
			&\in&U+\underline{\underline{t(z^5u_1u_2u_1)}}+Rz^3ts^{-1}t^{-1}(R+Rs^{-1})t^2+
			Rz^3ts^{-1}t^{-1}(R+Rs^{-1})ts+\\&&+Rz^3ts^{-1}t^{-1}st(R+Rs^{-1})t^{-1}s^{-1}t\\
			&\in&U+\underline{z^3u_2u_1u_2}+Rz^3t^2(st)^{-3}st^3+Rz^3t^2(ts)^{-3}st^2s+\\&&+Rz^3ts^{-1}t^{-1}(R+Rs^{-1})ts^{-1}t^{-1}s^{-1}t\\
			&\in&U+\underline{z^2u_2u_1u_2}+\underline{z^2u_2st^2s}+Rz^3ts^{-2}t^{-1}s^{-1}t+Rz^3t^2(st)^{-3}st^3(st)^{-3}st^2\\
			&\in&U+Rz^3t(R+Rs^{-1})t^{-1}s^{-1}t+Rzt^2s(R+Rt+Rt^2+Rt^{-1})st^2\\
			&\in&U+\underline{(z+z^3)u_2u_1u_2}+Rz^3t^2(st)^{-3}st^2+Rzt(ts)^3s^{-1}t+Rzt(ts)^3s^{-1}t^{-1}s^{-1}tst^2+\\&&+Rzt^2(R+Rs^{-1})t^{-1}(R+Rs^{-1})t^2\\
			&\in&U+\underline{(z+z^2)u_2u_1u_2}+Rz^2ts^{-1}t^{-1}(R+Rs)tst^2+
			Rzt^2s^{-1}t^{-1}s^{-1}t^2\\
				&\in&U+\underline{z^2u_2}+Rz^2ts^{-1}t^{-2}(ts)^3s^{-1}t+Rzt^3(st)^{-3}st^3\\
				&\in&U+\underline{\underline{t(z^3u_1u_2u_1t)}}+\underline{u_2u_1u_2}.
		\end{array}}$\\
		\qedhere
	\end{itemize}
\end{proof}
\begin{cor} $Uu_1\subset U$
	\label{cor9}
\end{cor}
\begin{proof}
	By the definition of $U$ we have that
	$Uu_1\subset \sum\limits_{k=0}^7(z^ku_2u_1u_2u_1+z^ku_2st^{-2}u_1)$.
	By proposition \ref{prr9}, we only have to prove that for every $k\in\{0,\dots,7\}$, $z^ku_2st^{-2}u_1\subset U$, which follows from the definition of $U$ if we expand $u_1$ as $R+Rs$.
\end{proof}
For the rest of this section, we will use directly corollary \ref{cor9}; this means that every time we have a power of $s$ at the end of an element, we may ignore it, as we did for the powers of $t$ in the beginning of the elements. In order to remind that to the reader, we put again a parenthesis around the part of the element we consider.
\begin{thm}$H_{G_9}=U$.
	\label{thm9}
\end{thm}
\begin{proof}                       
	By propositions \ref{pr9} and \ref{prr9}, we only have to prove that  $z^kst^{-2}st\in U$, for every $k\in\{0,\dots,7\}$,
However, by lemma \ref{lem9}(v) we have that  $z^kst^{-2}st\in U+u_2(z^{k+2}u_1u_2u_1)$, for every $k\in\{1,\dots, 5\}$. Therefore,  by proposition \ref{prr9}, we only the cases where $k\in\{0\}\cup\{6,7\}$.
\begin{itemize}[leftmargin=*]
	\item \underline{$k=0$}: 
	
	$\small{\begin{array}[t]{lcl}
	st^{-2}st&\in&s(R+Rt+Rt^2+Rt^3)st\\
	&\in& \underline{u_1u_2}+R(st)^3t^{-1}s^{-1}+Rst(ts)^3s^{-1}t^{-1}s^{-1}+Rst^2(ts)^3s^{-1}t^{-1}s^{-1}\\
	&\in&U+\underline{zu_2u_1}+zst(R+Rs)t^{-1}u_1+zst^2(R+Rs)t^{-1}u_1\\
	&\in&U+\underline{zu_1}+z(st)^3t^{-1}s^{-1}t^{-2}u_1+\underline{\underline{zu_1u_2u_1}}+zst^2s(R+Rt+Rt^2+Rt^3)u_1\\
	&\in&U+\underline{(z^2u_2u_1u_2)u_1}+\underline{\underline{zu_1u_2u_1}}+\underline{\underline{(zu_1t^2u_1t)u_1}}+zst(ts)^3s^{-1}t^{-1}s^{-1}tu_1+
	zst(ts)^3s^{-1}t^{-1}s^{-1}t^2u_1\\
	&\in&U+z^2st(R+Rs)t^{-1}s^{-1}tu_1+z^2st(R+Rs)t^{-1}s^{-1}t^2u_1\\
	&\in&U+\underline{z^2u_2u_1}+z^2(st)^3t^{-1}s^{-1}t^{-2}s^{-1}tu_1+z^2(st)^3t^{-1}s^{-1}t^{-2}s^{-1}t^2u_1\\
	&\in&U+\underline{\underline{t^{-1}(z^3u_1u_2u_1t)u_1}}+z^3t^{-1}s^{-1}t^{-2}s^{-1}t^2u_1
	
\end{array}}$

It remains to prove that the element $z^3t^{-1}s^{-1}t^{-2}s^{-1}t^2$ is inside $U$. For this purpose, we expand $t^{-2}$ as a linear combination of 1, $t^{-1}$, $t$ and $t^2$ and we have:

$\small{\begin{array}{lcl}
	z^3t^{-1}s^{-1}t^{-2}s^{-1}t^2	&\in&z^3t^{-1}s^{-1}(R+Rt^{-1}+Rt+Rt^2)s^{-1}t^2\\
	&\in&\underline{z^3u_2u_1u_2}+Rz^3(st)^{-3}st^3+
	Rz^3t^{-1}(R+Rs)t(R+Rs)t^2+\\&&+
	Rz^3t^{-1}(R+Rs)t^2(R+Rs)t^2\\
	&\in&U+\underline{\underline{u_2\big((z^2+z^3)u_1u_2\big)}}+	Rz^3t^{-1}stst^2+Rz^3t^{-1}st^2st^2\\
	&\in&U+Rz^3t^{-2}(ts)^3s^{-1}t+Rz^3t^{-1}st(ts)^3s^{-1}t^{-1}s^{-1}t\\
	&\in&U+\underline{\underline{t^{-2}(z^4u_1u_2)}}+rz^4t^{-1}st(R+Rs)t^{-1}s^{-1}t\\
	&\in&U+\underline{z^4u_1}+Rz^4t^{-2}(ts)^3s^{-1}t^{-2}s^{-1}t\\
	&\in&U+z^5u_2(u_1u_2u_1t).
	\end{array}}$

However, by lemma 	\ref{lem9}(v) we have that $z^5u_2(u_1u_2u_1t)\subset u_2U+
\underline{(z^7u_2u_1u_2)u_1}$. The result follows from remark \ref{rem9}.

	\item \underline{$k\in\{6,7\}$}: 
	
	$\small{\begin{array}[t]{lcl}
	z^kst^{-2}st&\in&z^k(R+Rs^{-1})t^{-2}(R+Rs^{-1})t\\
	&\in&\underline{z^ku_2u_1u_2}+Rz^{k}s^{-1}t^{-2}s^{-1}t\\
	&\in&U+Rz^kt(st)^{-3}stst^{-1}s^{-1}t\\
	&\in&U+Rz^{k-1}tst(R+Rs^{-1})t^{-1}s^{-1}t\\
	&\in&U+\underline{z^{k-1}u_2}+
	Rz^{k-1}tst^2(st)^{-3}st^2\\
	&\in&U+Rz^{k-2}ts(R+Rt+Rt^{-1}+Rt^{-2})st^2\\
	&\in&U+\underline{z^{k-2}u_2u_1u_2}+Rz^{k-2}(ts)^3s^{-1}t+Rz^{k-2}t(R+Rs^{-1})t^{-1}(R+Rs^{-1})t^2+\\&&+Rz^{k-2}t(R+Rs^{-1})t^{-2}(R+Rs^{-1})t^2\\
	&\in&U+\underline{(z^{k-1}+z^{k-2})u_2u_1u_2}+Rz^{k-2}ts^{-1}t^{-1}s^{-1}t^2+Rz^{k-2}ts^{-1}t^{-2}s^{-1}t^2\\
	&\in&U+Rz^{k-2}t^2(st)^{-3}st^3+Rz^{k-2}t^2(st)^{-3}stst^{-1}s^{-1}t^2\\
	&\in&U+\underline{z^{k-3}u_2u_1u_2}+z^{k-3}u_2st(R+Rs^{-1})t^{-1}s^{-1}t^2\\
	&\in&U+\underline{z^{k-3}u_2}+z^{k-3}u_2st^2(st)^{-3}st^3\\
	&\in&U+z^{k-4}u_2st^2s(R+Rt+Rt^2+Rt^{-1})\\
	&\in&U+\underline{\underline{u_2(z^{k-4}u_1u_2u_1)}}+\underline{\underline{z^{k-4}u_2(u_1t^2u_1t)}}+z^{k-4}u_2st(ts)^3s^{-1}t^{-1}s^{-1}t+\\&&+z^{k-4}u_2st^2(R+Rs^{-1})t^{-1}\\
	&\in&U+z^{k-3}u_2st(R+Rs)t^{-1}s^{-1}t+\underline{z^{k-4}u_2u_1u_2}+z^{k-4}u_2st^3(st)^{-3}sts\\
	&\in&U+\underline{z^{k-3}u_2}+z^{k-3}u_2(ts)^3s^{-1}t^{-2}s^{-1}t+\underline{\underline{u_2(z^{k-5}u_1u_2u_1t)s}}\\
		&\in&U+z^{k-2}u_2(u_1u_2u_1t).
	\end{array}}$\\

Again, by lemma \ref{lem9}(v) we have that $z^{k-2}u_2(u_1u_2u_1t)\subset u_2U+
\underline{(z^ku_2u_1u_2)u_1}$. The result follows from remark \ref{rem9}.
	\qedhere
\end{itemize}
\end{proof}

\begin{cor}
	The BMR freeness conjecture holds for the generic Hecke algebra $H_{G_9}$.
\end{cor}
\begin{proof}
	By theorem \ref{thm9} we have that $H_{G_9}=U=\sum\limits_{k=0}^7z^k(u_2+u_2su_2+u_2st^{-2})$. The result then follows from proposition \ref{BMR PROP}, since $H_{G_9}$ is generated as left $u_2$-module by 48 elements and, hence, as $R$-module by $|G_9|=192$ elements (recall that $u_2$ is generated as $R$-module by 4 elements.)
\end{proof}

\subsection{The case of $G_{10}$}
\indent

Let $R=\ZZ[u_{s,i}^{\pm},u_{t,j}^{\pm}]_{\substack{1\leq i\leq 3 \\1\leq j\leq 4}}$ and let $H_{G_{10}}=\langle s,t\;|\; stst=tsts,\prod\limits_{i=1}^{3}(s-u_{s,i})=\prod\limits_{j=1}^{4}(t-u_{t,i})=0\rangle$ be the generic Hecke algebra associated to $G_{10}$. Let $u_1$ be the subalgebra of $H_{G_{10}}$ generated by $s$ and $u_2$ the subalgebra of $H_{G_{10}}$ generated by $t$. We recall that $z:=(st)^2=(ts)^2$  generates the center of the associated complex braid group and that $|Z(G_{10})|=12$. We set $U=\sum\limits_{k=0}^{11}(z^ku_2u_1+z^ku_2st^{-1}+z^ku_2s^{-1}t+z^ku_2s^{-1}ts^{-1}).$

From now on, we will underline the elements that belong to $U$  by definition. Our goal is to prove that $H_{G_{10}}=U$ (theorem \ref{thm10}). 
 Since $1\in U$, it is enough to prove that $U$ is a right-sided ideal of $H_{G_{10}}$ or, equivalently, that $Us$ and $Ut$ are subsets of $U$. 
For this purpose, we first need to prove some preliminary results. 

In the following lemmas we prove that some subsets of $z^ku_2u_1u_2u_1$, where $k$ belongs in a smaller range of $\{0,\dots,11\}$, are also subsets of $U$.
\begin{lem} \mbox{}
	\vspace*{-\parsep}
	\vspace*{-\baselineskip}\\
	\begin{itemize}[leftmargin=0.8cm]
		\item [(i)] For every $k\in\{1,\dots, 10\}$, $z^ku_2u_1u_2\subset U$.
		\item[(ii)]For every $k\in\{1,\dots, 11\}$, $z^ku_2st^{-1}s\subset U$.
		\item[(iii)]For every $k\in\{0,\dots,9\}$, $z^ku_2st^2s\subset U$.
		\item [(iv)]For every $k\in\{1,\dots,9\}$, $z^ku_2su_2s\subset U$.
	\end{itemize}
	\label{lem10}
	\end{lem}
	\begin{proof}
		\mbox{}
		\vspace*{-\parsep}
		\vspace*{-\baselineskip}\\
		\begin{itemize}[leftmargin=0.8cm]
			\item[(i)] 
			$\small{\begin{array}[t]{lcl}
			z^ku_2u_1u_2&=&z^ku_2(R+Rs+Rs^{-1})u_2\\
			&=&\underline{z^ku_2}+z^ku_2s(R+Rt+Rt^{-1}+Rt^2)+z^ku_2s^{-1}(R+Rt+Rt^{-1}+Rt^{-2})\\
			&\subset& U+\underline{z^ku_2u_1}+Rz^ku_2(st)^2t^{-1}s^{-1}+\underline{z^ku_2st^{-1}}+
			z^ku_2(st)^2t^{-1}s^{-1}t+\underline{z^ku_2s^{-1}t}+\\&&+z^ku_2(ts)^{-2}ts+z^ku_2(ts)^{-2}tst^{-1}\\
			&\subset& U+\underline{z^{k+1}u_2u_1}+\underline{z^{k+1}u_2s^{-1}t}+
			\underline{z^{k-1}u_2u_1}+\underline{z^{k-1}u_2st^{-1}}.
			\end{array}}$
			\item[(ii)] 
			$\small{\begin{array}[t]{lcl}
			z^ku_2st^{-1}s&=&z^ku_2st^{-1}(R+Rs^{-1}+Rs^{-2})\\
			&\subset&\underline{z^ku_2st^{-1}}+z^ku_2s^2(ts)^{-2}t+
			z^ku_2(R+Rs^{-1}+Rs^{-2})t^{-1}s^{-2}\\
			&\subset&U+z^{k-1}u_2(R+Rs+Rs^{-1})t+\underline{z^{k}u_2u_1}+z^ku_2(ts)^{-2}ts^{-1}+z^ku_2s^{-1}(ts)^{-2}ts^{-1}\\
			&\subset& U+\underline{z^{k-1}u_2u_1}+z^{k-1}u_2(st)^2t^{-1}s^{-1}+\underline{z^{k-1}
			u_2s^{-1}t}+\underline{z^{k-1}u_2s^{-1}ts^{-1}}\\
		&\subset&U+\underline{z^ku_2u_1}.
			\end{array}}$
			
			\item[(iii)]$\small{\begin{array}[t]{lcl}
			z^ku_2st^2s&=&z^ku_2st(ts)^2s^{-1}t^{-1}\\
			&\subset&U+z^{k+1}u_2st(R+Rs+Rs^{2})t^{-1}\\
			&\subset&U+\underline{z^{k+1}u_2u_1}+z^{k+1}u_2(st)^2t^{-2}+z^{k+1}u_2(st)^2t^{-1}st^{-1}\\
			&\subset&U+\underline{z^{k+2}u_2}+\underline{z^{k+2}u_2st^{-1}}.
			\end{array}}$
			\item [(iv)] $\small{\begin{array}[t]{lcl}
		z^ku_2su_2s&=&z^ku_2s(R+Rt+Rt^2+Rt^{-1})\\
		&\subset& \underline{z^ku_2u_1}+
			z^ku_2(st)^2t^{-1}s^{-1}+z^ku_2st^2s+z^ku_2st^{-1}s\\
			&\subset& U+\underline{z^{k+1}u_2u_1}+z^ku_2st^2s+z^ku_2st^{-1}s\\
			&\subset& U+z^ku_2st^2s+z^ku_2st^{-1}s.
			\end{array}}$
			
			The result follows then from (ii) and (iii).
			\qedhere
		\end{itemize}
	\end{proof}
	\begin{lem}
		\mbox{}
		\vspace*{-\parsep}
		\vspace*{-\baselineskip}\\
		\begin{itemize}[leftmargin=0.8cm]
			\item[(i)] For every $k\in\{0,\dots,10\}$, $z^ku_2u_1tu_1\subset U$.
			
				\item[(ii)] For every $k\in\{2,\dots,11\}$, $z^ku_2s^{-1}u_2s^{-1}\subset U$.
				\item[(iii)] For every $k\in\{2,\dots,10\}$, $z^ku_2u_1u_2s^{-1}\subset U$.
				\label{lem210}
		\end{itemize}
	\end{lem}
	\begin{proof}
		\mbox{}
		\vspace*{-\parsep}
		\vspace*{-\baselineskip}\\
\begin{itemize}[leftmargin=0.8cm]
	\item[(i)] 
	$\small{\begin{array}[t]{lcl}
	z^ku_2u_1tu_1&=&z^ku_2(R+Rs+Rs^{-1})tu_1\\
	&\subset&\underline{z^ku_2u_1}+z^ku_2(st)^2t^{-1}s^{-1}u_1+z^ku_2s^{-1}t(R+Rs+Rs^{-1})\\
	&\subset&U+\underline{z^{k+1}u_2u_1}+\underline{z^ku_2s^{-1}t}+z^ku_2s^{-2}(st)^2t^{-1}+\underline{z^ku_2s^{-1}ts^{-1}}\\
	&\subset&U+z^{k+1}u_2(R+Rs+Rs^{-1})t^{-1}\\
	&\subset&U+\underline{z^{k+1}u_2}+\underline{z^{k+1}u_2st^{-1}}+z^{k+1}u_2(ts)^{-2}ts\\
	&\subset&U+\underline{z^ku_2u_1}.
	\end{array}}$

	\item[(ii)]
	$\small{\begin{array}[t]{lcl}
	z^ku_2s^{-1}u_2s^{-1}&=&z^ku_2s^{-1}(R+Rt+Rt^{-1}+Rt^{-2})s^{-1}\\
	&\subset&\underline{z^ku_2u_1}+\underline{z^ks^{-1}ts^{-1}}+z^ku_2(ts)^{-2}t+z^ku_2(ts)^{-2}tst^{-1}s^{-1}\\
	&\subset&U+\underline{z^{k-1}u_2}+z^{k-1}u_2s(st)^{-2}st\\
	&\subset&U+z^{k-2}u_2s^2t\\
	&\subset&U+z^{k-2}u_2(R+Rs+Rs^{-1})t\\
	&\subset&U+\underline{z^{k-2}u_2}+z^{k-2}u_2(st)^2t^{-1}s^{-1}+\underline{z^{k-2}u_2s^{-1}t}\\
	&\subset& U+\underline{z^{k-1}u_2u_1}.
	\end{array}}$
	\item[(iii)]
	$\small{\begin{array}[t]{lcl}
	z^ku_2u_1u_2s^{-1}&=&z^ku_2(R+Rs+Rs^{-1})u_2s^{-1}\\
	&\subset&\underline{z^ku_2u_1}+z^ku_2su_2s^{-1}+z^ku_2s^{-1}u_2s^{-1}\\
	&\stackrel{(ii)}{\subset}&U+z^ku_2s(R+Rt+Rt^2+Rt^{-1})s^{-1}\\
	&\subset&U+\underline{z^ku_2}+z^ku_2u_1tu_1+z^ku_2(st)^2t^{-1}s^{-1}ts^{-1}+z^ku_2s(st)^{-2}st\\
	&\stackrel{(i)}{\subset}&U+\underline{z^{k+1}u_2s^{-1}ts^{-1}}+z^{k-1}u_2u_1u_2\
	\end{array}}$
	
	The result follows from lemma \ref{lem10}(i).
	\qedhere
	
\end{itemize}
	\end{proof}
To make it easier for the reader to follow the calculations, we will double-underline the elements as described in lemmas \ref{lem10} and \ref{lem210} and  we will use directly the fact that these elements are inside $U$. In the following lemma we prove that some subsets of $z^ku_2u_1u_2u_1u_2$, where $k$ belongs in a smaller range of $\{0,\dots,11\}$ are also subsets of $U$.

	\begin{lem}
		 \mbox{}
		 \vspace*{-\parsep}
		 \vspace*{-\baselineskip}\\
		 \begin{itemize}[leftmargin=0.8cm]
		 	\item [(i)] For every $k\in\{1,\dots,8\}$, $z^ku_2u_1tu_1t\subset U$.
		 	\item [(ii)] For every $k\in\{1,\dots,6\}$, $z^ku_2su_2st^2\subset U$.
		 	\item[(iii)] For every $k\in\{1,\dots,5\}$, $z^ku_2u_1tu_1t^2\subset U$.
		 	\item[(iv)] For every $k\in\{1,\dots,3\}$, $z^ku_2su_2su_2\subset U$.
		 	\item[(v)]For every $k\in\{3,\dots,5\}$, $z^ku_2s^{-1}u_2s^{-1}u_2\subset U$.
		 \end{itemize}
		 \label{stst10}
	\end{lem}
	\begin{proof}
		\mbox{}
		\vspace*{-\parsep}
		\vspace*{-\baselineskip}\\
		\begin{itemize}[leftmargin=0.8cm]
			\item[(i)]
			$\small{\begin{array}[t]{lcl}
			z^ku_2u_1tu_1t&=&z^ku_2u_1t(R+Rs+Rs^2)t\\
			&\subset&\underline{\underline{z^ku_2u_1u_2}}+z^ku_2u_1(st)^2+z^ku_2(R+Rs+Rs^2)ts^2t\\
				&\subset&U+\underline{z^{k+1}u_2u_1}+
				\underline{\underline{z^ku_2u_1u_2}}+z^ku_2(st)^2t^{-1}st+z^ku_2s(st)^2t^{-1}st\\
				&\subset&U+\underline{\underline{z^{k+1}u_2u_1u_2}}+z^{k+1}u_2s(R+Rt+Rt^2+Rt^3)st\\
				&\subset&U+\underline{\underline{z^{k+1}u_2u_1u_2}}+z^{k+1}u_2(st)^2+
				z^{k+1}u_2st(ts)^2s^{-1}+z^{k+1}u_2st^2(ts)^2s^{-1}\\
				&\subset&U+\underline{z^{k+2}u_2}+z^{k+2}u_2(st)^2t^{-1}s^{-2}+z^{k+2}u_2(st)^2t^{-1}s^{-1}ts^{-1}\\
						&\subset&U+\underline{z^{k+3}u_2u_1}+\underline{z^{k+3}u_2s^{-1}ts^{-1}}.
			
			\end{array}}$
			\item[(ii)]
			$\small{\begin{array}[t]{lcl}
			z^ku_2su_2st^2&=&z^ku_2s(R+Rt+Rt^2+Rt^3)st^2\\
			&\subset&\underline{\underline{z^ku_2u_1u_2}}+z^ku_2(st)^2t+z^ku_2st(ts)^2s^{-1}t+z^ku_2st^2(ts)^2s^{-1}t\\
				&\subset&U+\underline{z^{k+1}u_2}+z^{k+1}u_2u_1tu_1t+z^{k+1}u_2(st)^2t^{-1}s^{-1}ts^{-1}t\\
				&\subset&U+(z^{k+1}+z^{k+2})u_2u_1tu_1t.
			\end{array}}$
			
			The result follows from (i).
			\item[(iii)]
			$\small{\begin{array}[t]{lcl}
			z^ku_2u_1tu_1t^2&=&z^ku_2u_1t(R+Rs+Rs^2)t^2\\
			&\subset&\underline{\underline{z^ku_2u_1u_2}}+z^ku_2u_1(ts)^2s^{-1}t+
			z^ku_2(R+Rs+Rs^2)ts^2t^2\\
			&\subset&U+\underline{\underline{z^{k+1}u_2u_1u_2}}+
			\underline{\underline{z^ku_2u_1u_2}}+z^ku_2(st)^2t^{-1}st^2+z^ku_2s(st)^2t^{-1}st^2\\
			&\subset&U+\underline{\underline{z^{k+1}u_2u_1u_2}}+z^{k+1}u_2su_2st^2.
			\end{array}}$
			
			The result follows from (ii).
			\item[(iv)]
			$\small{\begin{array}[t]{lcl}
			z^ku_2su_2su_2&=&z^ku_2s(R+Rt+Rt^2+Rt^3)su_2\\
			&\subset&\underline{\underline{z^ku_2u_1u_2}}+z^ku_2(st)^2u_2+z^ku_2(st)^2t^{-1}s^{-2}(st)^2u_2+
			z^ku_2st^3s(R+Rt+Rt^2+Rt^3)\\
			&\subset&U+\underline{z^{k+1}u_2}+\underline{\underline{z^{k+2}u_2u_1u_2}}+
			\underline{\underline{z^ku_2su_2s}}+z^ku_2st^2(ts)^2s^{-1}+z^ku_2su_2st^2+\\&&+
			z^ku_2st^2(ts)^2s^{-1}t^2\\
			&\stackrel{(ii)}{\subset}&U+z^{k+1}u_2(st)^2t^{-1}s^{-1}ts^{-1}+z^{k+1}u_2(st)^2t^{-1}s^{-1}ts^{-1}t^2\\
			&\subset&U+\underline{z^{k+2}u_2s^{-1}ts^{-1}}+z^{k+2}u_2u_1tu_1t^2.
			\end{array}}$
			
			The result follows from (iii).
			\item[(v)] $z^ku_2s^{-1}u_2s^{-1}u_2=
			z^ku_2(ts)^{-1}tsu_2(ts)^{-1}tsu_2\subset z^{k-2}u_2su_2su_2$.
			The result follows from (iv).
			\qedhere
		\end{itemize}
		\end{proof}
		To make it easier for the reader to follow the calculations, we will also double-underline the elements as described in lemmas \ref{lem10} and \ref{lem210} and  we will use directly the fact that these elements are inside $U$.
		The following lemma helps us to prove that $Uu_1\subset U$ (see proposition \ref{Us10}).
		\begin{lem} For every $k\in\{8,9\}$, $z^ks^{-1}u_2s^{-1}u_2s^{-1}\subset U$.
			\label{ll10}
		\end{lem}
		\begin{proof} We expand $\bold{u_2}$ as $R+Rt^{-1}+Rt^{-2}+Rt^{-3}$ and we have:
				$$\small{\begin{array}[t]{lcl}
			z^ks^{-1}\bold{u_2}s^{-1}u_2s^{-1}
			&\subset&\underline{\underline{z^ku_2u_1u_2s^{-1}}}+z^k(ts)^{-2}u_2s^{-1}+
			z^k(ts)^{-2}tst^{-1}s^{-1}u_2s^{-1}+
			z^k(ts)^{-2}tst^{-2}s^{-1}u_2s^{-1}\\
			&\subset&U+\underline{z^{k-1}u_2u_1}+z^{k-1}u_2s^2(ts)^{-2}u_2s^{-1}+
			z^{k-1}u_2st^{-1}(st)^{-2}su_2s^{-1}\\
			&\subset&U+\underline{\underline{z^{k-2}u_2u_1u_2s^{-1}}}+z^{k-2}u_2st^{-1}s(R+Rt+Rt^{-1}+Rt^{-2})s^{-1}\\
			&\subset&U+\underline{z^{k-2}u_2st^{-1}}+z^{k-2}u_2st^{-1}(st)^2t^{-1}s^{-2}+z^{k-2}u_2st^{-1}s(st)^{-2}st+\\&&+
			z^{k-2}u_2st^{-1}(R+Rs^{-1}+Rs^{-2})t^{-2}s^{-1}\\
			&\subset&U+z^{k-1}u_2st^{-2}(R+Rs+Rs^{-1})+z^{k-3}u_2st^{-1}s^2t+\underline{\underline{z^{k-2}u_2u_1u_2s^{-1}}}+\\&&+
			z^{k-2}u_2s(st)^{-2}st^{-1}s^{-1}+z^{k-2}u_2s(st)^{-2}sts^{-1}t^{-2}s^{-1}\\
			&\subset&U+\underline{z^{k-1}u_2u_1}+\underline{\underline{z^{k-1}u_2su_2s}}+\underline{\underline{z^{k-1}u_2u_1u_2s^{-1}}}+
			z^{k-3}u_2st^{-1}(R+Rs+Rs^{-1})t+\\&&+\underline{\underline{z^{k-3}u_2u_1u_2s^{-1}}}+
			z^{k-3}u_2s^2t(ts)^{-2}tst^{-1}s^{-1}\\
			&\subset&U+\underline{z^{k-3}u_2u_1}+z^{k-3}u_2st^{-1}(st)^2t^{-1}s^{-1}+
			z^{k-3}u_2s(st)^{-2}st^2+\\&&+
			z^{k-4}u_2s^2t^2(R+Rs^{-1}+Rs^{-2})t^{-1}s^{-1}\\
			&\subset&U+\underline{\underline{z^{k-2}u_2u_1u_2s^{-1}}}+\underline{\underline{z^{k-4}u_2u_1u_2}}+\underline{\underline{z^{k-4}u_2u_1u_2s^{-1}}}+
			z^{k-4}u_2s^2t^2(ts)^{-2}t+\\&&+z^{k-4}u_2s^2t^2s^{-1}(ts)^{-2}t\\
			&\subset&U+\underline{\underline{z^{k-5}u_2u_1u_2}}+z^{k-5}u_2s^2t^3(st)^{-2}st^2\\
			&\subset&U+z^{k-6}u_2(R+Rs+Rs^{-1})t^3st^2\\
			&\subset&U+\underline{\underline{z^{k-6}u_2u_1u_2}}+\underline{\underline{z^{k-6}u_2su_2st^2}}+z^{k-6}u_2(ts)^{-2}tst^4st^2\\
			&\subset&U+\underline{\underline{z^{k-7}u_2su_2st^2}}.
			\end{array}}$$
		\end{proof}
		\begin{prop} $Uu_1\subset U$.
			\label{Us10}
		\end{prop}
		\begin{proof}
			Since $u_1=R+Rs+Rs^2$, it is enough to prove that $Us\subset U$. By the definition of $U$ we only have to prove that for every $k\in\{0,\dots,11\}$, $z^ku_2st^{-1}s$ and $z^ku_2s^{-1}ts$ are subsets of $U$. 
			 However, by lemma \ref{lem10}(ii) we have that for every $k\in\{1,\dots,11\}$, $z^ku_2st^{-1}s\subset U$. Hence, it will be sufficient to prove that $u_2st^{-1}s\subset U$. We have:
			 $$\small{\begin{array}{lcl}
			 u_2st^{-1}s&=&u_2s(R+Rt+Rt^2+Rt^3)s\\
			 &\subset&\underline{u_2u_1}+u_2(ts)^2+\underline{\underline{u_2st^2s}}+u_2st^2(ts)^2s^{-1}t^{-1}\\
			 &\subset&U+\underline{zu_2}+zu_2st^2(R+Rs+Rs^2)t^{-1}\\
			 &\subset&U+\underline{\underline{u_2u_1u_2}}+zu_2(ts)^2s^{-2}(st)^2t^{-2}+zu_2st(ts)^2s^{-1}t^{-1}st^{-1}\\
			 &\subset&U+\underline{\underline{z^3u_2u_1u_2}}+z^2u_2st(R+Rs+Rs^2)t^{-1}st^{-1}\\
			 &\subset&U+\underline{\underline{z^2u_2u_1u_2}}+z^2u_2(st)^2t^{-2}st^{-1}+
			 z^2u_2(ts)^2st^{-1}st^{-1}\\
			 &\subset&U+\underline{z^3u_2st^{-1}}+\underline{\underline{z^3u_2su_2su_2}}.
			 \end{array}}$$
				It remains to prove that for every $k\in\{0,\dots,11\}$, $z^ku_2s^{-1}ts\subset U$. For $k\not=11$, the result is obvious since 
				$z^ku_2s^{-1}ts\subset z^ku_2(R+Rs+Rs^2)ts\subset \underline{z^ku_2u_1}+z^ku_2(st)^2t^{-1}+z^ku_2s(st)^2t^{-1}\subset U+\underline{z^{k+1}u_2}+\underline{z^{k+1}u_2st^{-1}}.$ Therefore, we only have to prove that $z^{11}u_2s^{-1}ts\subset U.$ 
					$$\small{\begin{array}{lcl}
						z^{11}u_2s^{-1}ts&\subset&z^{11}u_2s^{-1}t(R+Rs^{-1}+Rs^{-2})\\
						&\subset&\underline{z^{11}u_2s^{-1}t}+\underline{z^{11}u_2s^{-1}ts^{-1}}+z^{11}u_2s^{-1}(R+Rt^{-1}+Rt^{-2}+Rt^{-3})s^{-2}\\
						&\subset&U+\underline{z^{11}u_2u_1}+z^{11}u_2(ts)^{-2}ts^{-1}+
						z^{11}u_2(ts)^{-2}tst^{-1}s^{-2}+	z^{11}u_2(ts)^{-2}tst^{-2}s^{-2}\\
						&\subset&U+\underline{z^{10}u_2u_1}+z^{10}u_2s(st)^{-2}sts^{-1}+
						z^{10}u_2(R+Rs^{-1}+Rs^{-2})t^{-2}s^{-2}\\
						&\subset&U+\underline{\underline{z^9u_2u_1u_2s^{-1}}}+\underline{z^{10}u_2u_1}+z^{10}u_2(ts)^{-2}ts(st)^{-2}sts^{-1}+\\&&+
						z^{10}u_2s^{-1}(ts)^{-2}tst^{-1}s^{-2}\\
						&\subset&U+\underline{\underline{z^8u_2u_1u_2s^{-1}}}+z^9u_2s^{-1}t(R+Rs^{-1}+Rs^{-2})t^{-1}s^{-2}\\
					&\subset&U+\underline{z^9u_2u_1}+z^9u_2s^{-1}t(ts)^{-2}ts^{-1}+z^9u_2s^{-1}ts^{-1}(ts)^{-2}ts^{-1}\\
					&\subset&U+\underline{\underline{z^8s^{-1}u_2s^{-1}}}+z^8u_2s^{-1}u_2s^{-1}u_2s^{-1}
						\end{array}}$$
					The result follows from lemma \ref{ll10}.
							\end{proof}
							For the rest of this section, we will use directly proposition  \ref{Us10}; this means that every time we have a power of $s$ at the end of an element, we may ignore it. In order to remind that to the reader, we put  a parenthesis around the part of the element we consider.
							We can now prove the main theorem of this section.
							
\begin{thm} $H_{G_{10}}=U$.
	\label{thm10}
\end{thm}
\begin{proof}
	As we explain in the beginning of this section, since $1\in U$, it will be sufficient to prove that $U$ is a right-sided ideal of $H_{G_{10}}$. For this purpose one may check that $Us$ and $Ut$ are subsets of $U$. By proposition \ref{Us10} it is enough to prove that $Ut\subset U$. Since $t\in R+Rt^{-1}+Rt^{-2}+Rt^{-3}$, we only have to prove that $Ut^{-1}\subset U$. By the definition of $U$ we have:
	$$Ut^{-1}\subset\sum\limits_{k=0}^{11}(z^ku_2u_1t^{-1}+z^ku_2st^{-2}+\underline{z^ku_2s^{-1}}+z^ku_2s^{-1}ts^{-1}t^{-1}).$$
	As a result, we restrict ourselves to proving that for every $k\in\{0,\dots,11\}$, $z^ku_2u_1t^{-1}$, $z^ku_2st^{-2}$ and $z^ku_2s^{-1}ts^{-1}t^{-1}$ are subsets of $U$. We distinguish the following cases:
	\begin{itemize}[leftmargin=0.8cm]
		\item[C1.] \underline{The case of $z^ku_2u_1t^{-1}$}: 
		\begin{itemize}[leftmargin=*]
			\item \underline{$k\not=0$}: \\
			$z^ku_2u_1t^{-1}=z^ku_2(R+Rs+Rs^{-1})t^{-1}\subset
		\underline{z^ku_2}+\underline{z^ku_2st^{-1}}+z^ku_2(ts)^{-2}ts\subset
	U+\underline{z^{k-1}u_2u_1}.$
		
		\item \underline {$k=0$}:\\
		$\small{\begin{array}[t]{lcl}
		u_2u_1t^{-1}&=&u_2(R+Rs+Rs^2)t^{-1}\\
		&\subset& 
		\underline{u_2}+\underline{u_2st^{-1}}+u_2s^2t^{-1}\\
		&\subset&U+u_2s^2(R+Rt+Rt^2+Rt^3)\\
		&\subset&U+\underline{u_2u_1}+u_2s(st)^2t^{-1}s^{-1}+u_2s(st)^2t^{-1}s^{-1}t+u_2s(st)^2t^{-1}s^{-1}t^2\\
		&\subset&U+\underline{(zu_2st^{-1})s^{-1}}+zu_2st^{-1}(R+Rs+Rs^2)t+
		zu_2st^{-1}(R+Rs+Rs^2)t^2\\
		&\subset&U+\underline{zu_2s}+zu_2st^{-1}(st)^2t^{-1}s^{-1}+zu_2st^{-1}s(st)^2t^{-1}s^{-1}+\underline{\underline{zu_2u_1u_2}}+\underline{\underline{zu_2su_2st^2}}+\\&&+zu_2s(R+Rt+Rt^2+Rt^3)s^2t^2\\
		&\subset&U+\underline{\underline{(z^2u_2u_1u_2)s^{-1}}}+\underline{\underline{(z^2u_2su_2su_2)s^{-1}}}+\underline{\underline{zu_2u_1u_2}}+zu_2(st)^2t^{-1}st^2+\\&&+
		zu_2(st)^2t^{-1}s^{-1}ts^2t^2+zu_2(st)^2t^{-1}s^{-1}t^2s^2t^2\\
		&\subset&U+\underline{\underline{z^2u_2u_1u_2}}+z^2u_2(R+Rs+Rs^2)ts^2t^2+
		z^2u_2(R+Rs+Rs^2)t^2s^2t^2\\
			&\subset&U+\underline{\underline{z^2u_2u_1u_2}}+z^2u_2(st)^2t^{-1}st^2+
			z^2u_2s(st)^2t^{-1}st^2+z^2u_2(st)^2t^{-1}s^{-1}ts^2t^2+\\&&+
			z^2u_2s(st)^2t^{-1}s^{-1}ts^2t^2\\
			&\subset&U+\underline{\underline{z^3u_2u_1u_2}}+\underline{\underline{z^3u_2su_2st^2}}+z^3u_2s^{-1}t(R+Rs+Rs^{-1})t^2+\\&&+z^3u_2st^{-1}(R+Rs+Rs^2)ts^2t^2\\
			&\subset&U+\underline{\underline{z^3u_2u_1u_2}}+z^3u_2s^{-1}(ts)^2s^{-1}t+\underline{\underline{z^3u_2s^{-1}u_2s^{-1}u_2}}+z^3u_2st^{-1}(st)^2t^{-1}st^2+\\&&+z^3u_2st^{-1}s(st)^2t^{-1}st^2\\
			&\subset&U+\underline{\underline{z^4u_2u_1u_2}}+\underline{\underline{z^4u_2su_2st^2}}+z^4u_2st^{-1}st^{-2}(ts)^2s^{-1}t\\
			&\subset&U+z^5u_2st^{-1}s(R+Rt+Rt^{-1}+Rt^2)s^{-1}t\\
			&\subset&U+\underline{z^5u_2u_1}+z^5u_2st^{-1}(st)^2t^{-1}s^{-2}t+z^5u_2st^{-1}s(st)^{-2}st^2+z^5u_2st^{-1}st^2s^{-1}t\\
			&\subset&U+z^6st^{-2}(R+Rs+Rs^{-1})t+z^4u_2s(st)^{-2}sts^3t^2+
			z^5u_2st^{-1}st^2(R+Rs+Rs^2)t\\
			&\subset&U+\underline{\underline{z^6u_2u_1u_2}}+z^6u_2st^{-2}(st)^2t^{-1}s^{-1}+z^6u_2st^{-1}(ts)^{-2}st^2+\underline{\underline{z^3u_2u_1tu_1t^2}}+\\&&+
			z^5u_2st^{-1}st^3+z^5u_2st^{-1}(st)^2t^{-1}s^{-1}(ts)^2s^{-1}+\\&&+z^5u_2s(R+Rt+Rt^2+Rt^3)st^2s^2t\\
				&\subset&U+\underline{\underline{(z^7u_2u_1u_2)u_1}}+\underline{\underline{z^5u_2su_2st^2}}+z^5u_2st^{-1}s(R+Rt+Rt^{-1}+Rt^{2})+\\&&
			+z^5u_2s^2t^2s^2t+z^5u_2(st)^2ts^2t+z^5u_2st(ts)^2s^{-1}ts^2t+z^5u_2st^2(ts)^2s^{-1}ts^2t\\
			&\subset&U+\underline{(z^5u_2st^{-1})s}+z^5u_2st^{-1}(st)^2t^{-1}s^{-1}+z^5u_2st^{-1}st^{-1}+\\&&+z^5u_2st^{-1}(st)^2t^{-1}s^{-1}t+z^5u_2s^2t^2(R+Rs+Rs^{-1})t+\underline{\underline{z^6u_2u_1u_2}}+\\&&+
			z^6u_2(st)^2t^{-1}s^{-2}ts^2t+z^6u_2st^2(R+Rs+Rs^2)ts^2t\\
				&\subset&U+\underline{\underline{(z^6u_2u_1u_2)s^{-1}}}+
				z^5u_2st^{-1}(R+Rs^{-1}+Rs^{-2})t^{-1}+
				z^5u_2st^{-2}s^{-1}t+\\&&+
				\underline{\underline{z^5u_2u_1u_2}}+z^5u_2s^2t(ts)^2s^{-1}+z^5u_2(R+Rs+Rs^{-1})t^2s^{-1}t+
				\\&&+z^7u_2s^{-2}(ts)^2s^{-1}t^{-1}(st)^2t^{-1}s^{-1}+
				z^6u_2st^3s^2t+z^6u_2st(ts)^2st+\\&&+
				z^6u_2(st)^2t^{-1}s^{-1}ts(st)^2t^{-1}(st)^2t^{-1}s^{-1}\\
				&\subset&U+\underline{\underline{z^5u_2u_1u_2}}+z^5u_2s(st)^{-2}s+
				z^5u_2s(st)^{-2}st(ts)^{-2}ts+\\&&+
				z^5u_2s(R+Rt+Rt^{-1}+Rt^2)s^{-1}t+
				\underline{\underline{(z^6u_2u_1u_2)s^{-1}}}+
				\underline{z^6u_2s^{-1}t}+\\&&+
				z^5u_2st^2s^{-1}t+\underline{\underline{z^5u_2s^{-1}u_2s^{-1}u_2}}+
					\underline{\underline{(z^9u_2u_1u_2)s^{-1}}}+\\&&+
					z^6u_2st^3(R+Rs+Rs^{-1})t+z^7u_2(st)^2+
					z^9u_2s^{-1}(ts)^2s^{-1}t^{-3}s^{-1}\\
					&\subset&U+\underline{z^4u_2u_1}+	\underline{\underline{(z^3u_2u_1u_2)}}+\underline{z^5u_2}+
					z^5u_2(st)^2t^{-1}s^{-2}t+z^5u_2s(ts)^{-2}st^2+\\&&+
					
					z^5u_2st^2s^{-1}t+
					\underline{\underline{(z^6u_2u_1u_2)}}+z^6u_2st^2(ts)^2s^{-1}+
					z^6u_2st^3(ts)^{-2}tst^2+\underline{z^8u_2}+\\&&+\underline{\underline{(z^{10}u_2u_1u_2)s^{-1}}}\\
					&\subset&U+\underline{\underline{(z^4+z^6)u_2u_1u_2}}+z^5u_2(st)^2t^{-1}s^{-1}ts^{-1}t+\underline{\underline{(z^7u_2u_1u_2)s^{-1}}}+
					\underline{\underline{z^5u_2su_2st^2}}\\
					&\subset&U+z^6u_2s^{-1}t(R+Rs+Rs^2)t\\
					&\subset&U+\underline{z^6u_2u_1}+z^6u_2s^{-1}(ts)^2s^{-1}+
					z^6u_2s^{-1}(ts)^2s^{-1}t^{-1}(st)^2t^{-1}s^{-1}\\
					&\subset&U+\underline{z^7u_2u_1}+\underline{\underline{(z^8u_2u_1u_2)s^{-1}}}.
				\end{array}}$
			\end{itemize}
		\item [C2.] \underline{The case of $z^ku_2st^{-2}$}:
		
		For $k\not=11$, we expand $t^{-2}$ as a linear combination of 1, $t$, $t^{-1}$ and $t^{2}$ and  we have that 
		$z^ku_2st^{-2}\subset \underline{z^ku_2u_1}+z^ku_2(st)^2t^{-1}s^{-1}+\underline{z^ku_2st^{-1}}+\underline{\underline{(z^ku_2st^2s)s^{-1}}}\subset U+\underline{z^{k+1}u_2u_1}\subset U$.
			It remains to prove that $z^{11}u_2st^{-2}\subset U$. We have:
			
		$\small{\begin{array}{lcl}
		z^{11}u_2st^{-2}&\subset&z^{11}u_2(R+Rs^{-1}+Rs^{-2})t^{-2}\\
		&\subset&\underline{z^{11}u_2}+z^{11}u_2(ts)^{-2}tst^{-1}+z^{11}u_2s^{-1}(ts)^{-2}tst^{-1}\\
		&\subset&U+\underline{z^{10}u_2st^{-1}}+z^{10}u_2s^{-1}t(R+Rs^{-1}+Rs^{-2})t^{-1}\\
		&\subset&U+\underline{z^{10}u_2u_1}+z^{10}u_2s^{-1}t(ts)^{-2}ts+z^{10}u_2s^{-1}ts^{-1}(ts)^{-2}ts\\
		&\subset&U+\underline{\underline{(z^9u_2u_1u_2)s}}+(z^9u_2s^{-1}u_2s^{-1}u_2s^{-1})s^2.
		\end{array}}$
	
	The result follows from lemma \ref{ll10}.
		
		\item [C3.] \underline{The case of $z^ku_2s^{-1}ts^{-1}t^{-1}$}:
		
		For $k\not \in\{0,1\}$, we have that $z^ku_2s^{-1}ts^{-1}t^{-1}=z^ku_2s^{-1}t(ts)^{-2}ts=\underline{\underline{(z^{k-1}u_2u_1u_2)s}}\subset U.$ It remains to prove the case where $k\in\{0,1\}$. We have:
		
		$\small{\begin{array}{lcl}
		z^ku_2s^{-1}ts^{-1}t^{-1}&\subset&z^ku_2s^{-1}t(R+Rs+Rs^2)t^{-1}\\
		&\subset&\underline{z^ku_2u_1}+z^ku_2s^{-1}(ts)^2s^{-1}t^{-2}+z^ku_2(R+Rs+Rs^2)ts^2t^{-1}\\
		&\subset&U+\underline{\underline{z^{k+1}u_2u_1u_2}}+z^ku_2u_1t^{-1}+z^ku_2(st)^2t^{-1}st^{-1}+z^ku_2s(st)^2t^{-1}st^{-1}\\
		&\stackrel{C1}{\subset}&U+\underline{z^{k+1}u_2st^{-1}}+\underline{\underline{z^{k+1}u_2su_2su_2}}.
		\end{array}}$
		
		\qedhere
	\end{itemize}
\end{proof}
\begin{cor}
	The BMR freeness conjecture holds for the generic Hecke algebra $H_{G_{10}}$.
\end{cor}
\begin{proof}
	By theorem \ref{thm10} we have that $H_{G_{10}}=U$. The result follows from proposition \ref{BMR PROP}, since by definition $U$ is generated as left $u_2$-module by 28 elements and, hence, as $R$-module by $|G_{10}|=288$ elements (recall that $u_2$ is generated as $R$-module by 4 elements).
\end{proof}
\subsection{The case of $G_{11}$}
\indent

Let $R=\ZZ[u_{s,i}^{\pm},u_{t,j}^{\pm},u_{u,l}^{\pm}]$, where $1\leq i\leq 2$, $1\leq j\leq 3$  and $1\leq l\leq 4$. We also let $$H_{G_{11}}=\langle s,t,u\;|\; stu=tus=ust,\prod\limits_{i=1}^{2}(s-u_{s,i})=\prod\limits_{j=1}^{3}(t-u_{t,j})=\prod\limits_{l=1}^{4}(u-u_{u,l})=0\rangle$$ be the generic Hecke algebra associated to $G_{11}$. Let $u_1$ be the subalgebra of $H_{G_{11}}$ generated by $s$, $u_2$ the subalgebra of $H_{G_{11}}$ generated by $t$ and $u_3$ the subalgebra of $H_{G_{11}}$ generated by $u$. We recall that $z:=stu=tus=ust$  generates the center of the associated complex braid group and that $|Z(G_{11})|=24$.
We set $U=\sum\limits_{k=0}^{23}(z^ku_3u_2+z^ku_3tu^{-1}u_2).$
By the definition of $U$, we have the following remark.

\begin{rem}
	$Uu_2 \subset U$.
	\label{rem11}
\end{rem}
 From now on, we will underline the elements that by definition belong to $U$.  Moreover, we will use directly the remark \ref{rem11}; this means that every time we have a power of $t$ at the end of an element, we may ignore it. 
 To remind that to the reader, we put a parenthesis around the part of the element we consider.

Our goal is to prove that $H_{G_{11}}=U$ (theorem \ref{thm11}). Since $1\in U$, it will be sufficient to prove that $U$ is a left-sided ideal of $H_{G_{11}}$. For this purpose, one may check that $sU$, $tU$ and $uU$ are subsets of $U$. The following proposition states that it is enough to prove $tU\subset U$.
\begin{prop}
If $tU\subset U$ then $H_{G_{11}}=U$.
\label{Ut11}
\end{prop}
\begin{proof}
As we explained above, we have to prove that $sU$, $tU$ and $uU$ are subsets of $U$. However, by the definition of $U$ we have $uU\subset U$ and, hence, by hypothesis we
only have to prove that $sU\subset U$. We recall that $z=stu$, therefore $s=zu^{-1}t^{-1}$ and $s^{-1}=z^{-1}tu$. We notice that $$U=
\sum\limits_{k=0}^{22}z^k(u_3u_2+u_3tu^{-1}u_2)+z^{23}(u_3u_2+u_3tu^{-1}u_2).$$
Hence, we have:\\ \\
$\small{\begin{array}[t]{lcl}sU&\subset&
\sum\limits_{k=0}^{22}z^ks(u_3u_2+u_3tu^{-1}u_2)+z^{23}s(u_3u_2+u_3tu^{-1}u_2)\\
&\subset& \sum\limits_{k=0}^{22}z^{k+1}u^{-1}t^{-1}(u_3u_2+u_3tu^{-1}u_2)+z^{23}(R+Rs^{-1})(u_3u_2
+u_3tu^{-1}u_2)\\
&\subset& \sum\limits_{k=0}^{22}u^{-1}t^{-1}(\underline{z^{k+1}u_3u_2}+\underline{z^{k+1}u_3tu^{-1}u_2})+
\underline{z^{23}u_3u_2}+\underline{z^{23}u_3tu^{-1}u_2}+z^{23}s^{-1}u_3u_2+
z^{23}s^{-1}u_3tu^{-1}u_2\\
&\subset&u_3u_2U+z^{22}tu_3u_2+z^{22}tu_3tu^{-1}u_2\\
&\subset&u_3u_2U+
t(\underline{z^{22}u_3u_2}+\underline{z^{22}u_3t^{-1}u_2})\\
&\subset&u_3u_2U.
\end{array}}$\\\\
By hypothesis, $tU\subset U$ and, hence, $u_2U\subset U$, since $u_2=R+Rt+Rt^2$. Moreover, $u_3U\subset U$, by the definition of $U$. Therefore, $sU\subset U$.
\end{proof}

\begin{cor}
	If $z^ktu_3$ and $z^ktu_3tu^{-1}$ are subsets of $U$ for every $k\in\{0,\dots,23\}$, then $H_{G_{11}}=U$.
	\label{corr11}
\end{cor}
\begin{proof}
	
The result follows directly from the definition of $U$, proposition \ref{Ut11} and remark \ref{rem11}.
\end{proof}

As a first step we will prove the conditions of corollary \ref{corr11} for some shorter range of the values of $k$, as we can see in proposition  \ref{tu11} and corollary \ref{tuts111}.

\begin{prop}
	\mbox{}
	\vspace*{-\parsep}
	\vspace*{-\baselineskip}\\
	\begin{itemize}[leftmargin=0.8cm]
		\item[(i)] For every $k\in\{0,\dots,21\}$, $z^ku_3tu\subset U$.
		\item[(ii)] For every $k\in\{0,\dots,19\}$, $z^ku_3tu^2\subset U$.
		\item[(iii)] For every $k\in\{0,\dots,19\}$, $z^ku_3tu_3\subset U$.
		\item[(iv)] For every $k\in\{2,\dots,21\}$, $z^ku_3u_2u_3\subset U$.
	\end{itemize}
	\label{tu11}
\end{prop}
\begin{proof}
	Since $u_3=R+Ru+Ru^{-1}+Ru^2$, (iii) follows from (i) and (ii) and the definition of $U$. Moreover, (iv) follows directly from (iii), since: \\
	$\small{\begin{array}{lcl}
	z^ku_3u_2u_3&=&z^ku_3(R+Rt+Rt^{-1})u_3\\
	&\subset& \underline{z^ku_3}+z^ku_3tu_3+z^ku_3(t^{-1}s^{-1}u^{-1})usu_3\\
	&\subset& U+z^ku_3tu_3+
	z^{k-1}u_3(R+Rs^{-1})u_3\\
	&\subset& U+z^ku_3tu_3+\underline{z^{k-1}u_3}+z^{k-1}u_3(s^{-1}u^{-1}t^{-1})tu_3\\
	&\subset& U+(z^k+z^{k-2})tu_3.\end{array}}$\\
	Therefore, it is enough to prove (i) and (ii). 
For every $k\in\{0,\dots,21\}$ we have $z^ku_3tu=z^ku_3(tus)s^{-1}\subset z^{k+1}u_3(R+Rs)\subset \underline{z^{k+1}u_3}+z^{k+1}u_3(ust)t^{-1}\subset U+\underline{z^{k+2}u_3t^{-1}}\subset U$
 and, hence, we prove (i).
 
For (ii), we notice that, for every $k\in\{0,\dots,19\}$, $z^ku_3tu^2=z^ku_3(tus)s^{-1}u\subset z^{k+1}u_3(R+Rs)u\subset \underline{z^{k+1}u_3}+z^{k+1}u_3su\subset U+z^{k+1}u_3(ust)t^{-1}u\subset U+z^{k+2}u_3t^{-1}u$. However, if we expand $t^{-1}$ as a linear combination of 1, $t$ and $t^2$ we have that $z^{k+2}u_3t^{-1}u\subset
\underline{z^{k+2}u_3}+z^{k+2}u_3tu+z^{k+2}u_3t(tus)s^{-1}\stackrel{(i)}{\subset}U+z^{k+3}u_3t(R+Rs)\subset U+\underline{z^{k+3}u_3t}+z^{k+3}u_3t(stu)u^{-1}t^{-1}\subset
 U+\underline{z^{k+4}u_3tu^{-1}u_2}$.

\end{proof}

\begin{lem}
	\mbox{}
	\vspace*{-\parsep}
	\vspace*{-\baselineskip}\\
	\begin{itemize}[leftmargin=0.8cm]
		\item[(i)] For every $k\in\{1,\dots,22\}$, $z^ku_3u_2u_1\subset U$.
		\item[(ii)] For every $k\in\{0,\dots,18\}$, $z^ku_3tu_3u_1\subset U$.
		\item[(iii)] For every $k\in\{3,\dots,22\}$, $z^ku_3tu_1u_3\subset U$.
	\end{itemize}
	\label{tsu11}
\end{lem}
\begin{proof}
	\mbox{}
	\vspace*{-\parsep}
	\vspace*{-\baselineskip}\\
	\begin{itemize}[leftmargin=0.8cm]
		\item[(i)]
		$\begin{array}[t]{lcl}
		z^ku_3u_2u_1&=&z^ku_3u_2(R+Rs)\\
		&\subset& \underline{z^ku_3u_2}+z^ku_3(R+Rt^{-1}+Rt)s\\
		&\subset& z^ku_3(ust)t^{-1}+z^ku_3t^{-1}s+z^ku_3t(stu)u^{-1}t^{-1}\\
		& \subset& U+\underline{z^{k+1}u_3u_2}+z^ku_3t^{-1}(R+Rs^{-1})+\underline{z^{k+1}u_3tu^{-1}u_2}\\
		&\subset& U+\underline{z^{k}u_3u_2}+z^ku_3(u^{-1}t^{-1}s^{-1})\\
		&\subset& U+\underline{z^{k-1}u_3}.
		\end{array}$
		\item[(ii)]
		$z^ku_3tu_3u_1=z^ku_3tu_3(R+Rs)\subset z^ku_3tu_3+z^ku_3tu_3(ust)t^{-1}
		\subset \big(z^ku_3tu_3+z^{k+1}u_3tu_3\big)u_2$. The result follows from proposition \ref{tu11}(iii) .
		\item[(iii)]	
		$z^ku_3tu_1u_3=z^ku_3t(R+Rs^{-1})u_3=z^ku_3tu_3+z^ku_3t^2(t^{-1}s^{-1}u^{-1})u_3\stackrel{\ref{tu11}(iii)}{\subset}U+z^{k-1}u_3u_2u_3$.
		The result follows directly from proposition \ref{tu11}(iv).
		\qedhere
	\end{itemize}
\end{proof}
In order to make it easier for the reader to follow the calculations, from now on we will double-underline the elements as described in proposition \ref{tu11} and in lemma \ref{tsu11} and we will use directly the fact that these elements are inside $U$.

\begin{prop}
	\mbox{}
	\vspace*{-\parsep}
	\vspace*{-\baselineskip}\\
	\begin{itemize}[leftmargin=0.8cm]
	\item[(i)] For every $k\in\{2,\dots, 23\}$, $z^kt^2u^{-1}\in U$.	
	\item[(ii)]For every $k\in\{0,\dots,21\}$, $z^ktutu^{-1} \in U$.
		\item[(iii)]For every $k\in\{0,\dots,15\}$, $z^ktu^2tu^{-1} \in U$.
		\item[(iv)]For every $k\in\{6,\dots,23\}$, $z^ktu^{-1}tu^{-1} \in U+z^ku_3tu_3$.
		Therefore, for every $k\in\{6, \dots, 19\}$, $z^ktu^{-1}tu^{-1}\in U$.
		\item[(v)]For every $k\in\{0,\dots, 5\}$, $z^ktu^3tu^{-1} \in U$.
		\item[(vi)]For every $k\in\{16,\dots,23\}$, $z^ktu^{-2}tu^{-1} \in U+(z^k+z^{k-1}+z^{k-2})u_3tu_3u_2$.
		Therefore, for every $k\in\{16, \dots, 19\}$, $z^ktu^{-2}tu^{-1}\in U$.
		
	\end{itemize}
	\label{tuts11}
\end{prop}
\begin{proof}
	\mbox{}
	\vspace*{-\parsep}
	\vspace*{-\baselineskip}\\
	\begin{itemize}[leftmargin=0.8cm]
		\item[(i)]$\small{\begin{array}[t]{lcl}
		z^kt^2u^{-1}&\in& z^k(R+Rt+Rt^{-1})u^{-1}\\
		&\in& \underline{z^ku_3}+\underline{z^ku_3tu^{-1}}+z^ku_3(u^{-1}t^{-1}s^{-1})su^{-1}\\
	&\in& U+z^{k-1}u_3(R+Rs^{-1})u^{-1}\\
	&\in& U+\underline{z^{k-1}u_3}+z^{k-1}u_3(s^{-1}u^{-1}t^{-1})t\\
	&\in& U+\underline{z^{k-2}u_3u_2}.
	\end{array}}$
		
		\item[(ii)] 
		$z^ktutu^{-1}=z^k(tus)s^{-1}tu^{-1}\in z^{k+1}(R+Rs)tu^{-1}\subset
		\underline{z^{k+1}u_3tu^{-1}}+z^{k+1}u_3(stu)u^{-2}$. Hence, 
		$z^ktutu^{-1}\subset U+\underline{z^{k+2}u_3}\subset U$.
		\item[(iii)] 
		$\small{\begin{array}[t]{lcl}
		z^ktu^2tu^{-1}&=&z^k(tus)s^{-1}utu^{-1}\\
		&\in& z^{k+1}(R+Rs)utu^{-1}\\
		&\in& \underline{z^{k+1}u_3tu^{-1}}+z^{k+1}u_3(ust)t^{-1}utu^{-1}\\
		
		&\in& U+z^{k+2}u_3(R+Rt+Rt^2)utu^{-1}\\
		&\in& U+\underline{z^{k+2}u_3tu^{-1}}+z^{k+2}u_3tutu^{-1}+z^{k+2}u_3t(tus)s^{-1}tu^{-1}\\
		
		&\stackrel{(ii)}{\in}&U+z^{k+3}u_3t(R+Rs)tu^{-1}\\
		&\in& U+\underline{\underline{z^{k+3}u_3u_2u_3}}+z^{k+3}u_3t(stu)u^{-2}\\
		&\in& U+\underline{\underline{z^{k+4}u_3tu_3}}.
		\end{array}}$
		
		\item[(iv)] 
		$\begin{array}[t]{lcl}
		z^ktu^{-1}tu^{-1}&\in& z^ktu^{-1}(R+Rt^{-1}+Rt^{-2})u^{-1}\\
		&\in& z^ku_3tu_3+z^kt(u^{-1}t^{-1}s^{-1})su^{-1}+z^kt(u^{-1}t^{-1}s^{-1})st^{-1}u^{-1}\\
		&\in& z^ku_3tu_3+z^{k-1}t(R+Rs^{-1})u^{-1}+z^{k-1}tsu(u^{-1}t^{-1}s^{-1})su^{-1}\\
		&\in& z^ku_3tu_3+\underline{z^{k-1}u_3tu^{-1}}+z^{k-1}u_3t^2(t^{-1}s^{-1}u^{-1})+
		z^{k-2}tsu(R+Rs^{-1})u^{-1}\\
		&\in& U+z^ku_3tu_3+\underline{z^{k-2}u_3u_2}+z^{k-2}u_3t(stu)u^{-1}t^{-1}+
		z^{k-2}u_3tsu(s^{-1}u^{-1}t^{-1})t\\
		&\in&U+ z^ku_3tu_3+\underline{z^{k-1}u_3tu^{-1}u_2}+\underline{\underline{(z^{k-3}u_3tu_1u_3)u_2}}\\
		&\in& U+z^ku_3tu_3.
		\end{array}$

The result follows from proposition \ref{tu11}(iii).
		
		\item[(v)] 
		$\small{\begin{array}[t]{lcl}
		z^ktu^3tu^{-1}&=&z^k(tus)s^{-1}u^2tu^{-1}\\
		&\in&z^{k+1}(R+Rs)u^2tu^{-1}\\
		&\in&\underline{z^{k+1}u_3tu^{-1}}+z^{k+1}u_3(ust)t^{-1}u^2tu^{-1}\\
		&\in&U+z^{k+2}u_3(R+Rt+Rt^2)u^2tu^{-1}\\
		&\in&U+\underline{z^{k+2}u_3tu^{-1}}+z^{k+2}u_3tu^2tu^{-1}+z^{k+2}u_3t(tus)s^{-1}utu^{-1}\\
		&\stackrel{(iii)}{\in}&U+z^{k+3}u_3t(R+Rs)utu^{-1}\\
		&\in&U+z^{k+3}u_3tutu^{-1}+z^{k+3}u_3tu^{-1}(ust)t^{-2}(tus)s^{-1}tu^{-1}\\
		&\stackrel{(ii)}{\in}&U+z^{k+5}u_3tu^{-1}t^{-2}(R+Rs)tu^{-1}\\
		&\in&U+z^{k+5}u_3t(u^{-1}t^{-1}s^{-1})su^{-1}+z^{k+5}u_3tu^{-1}t^{-2}(stu)u^{-2}\\ 
		&\in&U+\underline{\underline{z^{k+4}u_3tu_1u_3}}+
		z^{k+6}u_3tu^{-1}t^{-2}u^{-2}.
		\end{array}}$
	
	We expand $t^{-2}$ as a linear combination of 1, $t^{-1}$ and $t$ and we have:\\
	$\small{\begin{array}{lcl}
		z^{k+6}u_3tu^{-1}t^{-2}u^{-2}	&\in&\underline{\underline{z^{k+6}u_3tu_3}}+z^{k+6}u_3t(u^{-1}t^{-1}s^{-1})su^{-2}+
		z^{k+6}u_3tu^{-1}t(R+Ru+Ru^{-1}+Ru^{2})\\ 
		&\in&U+\underline{\underline{z^{k+5}u_3tu_1u_3}}+\underline{z^{k+6}u_3tu^{-1}u_2}+
		z^{k+6}u_3tu^{-1}(tus)s^{-1}+z^{k+6}u_3tu^{-1}tu^{-1}+\\&&+
		z^{k+6}u_3tu^{-1}(tus)s^{-1}(ust)t^{-1}s^{-1} \\ 
		&\stackrel{(iv)}{\in}&U+\underline{\underline{z^{k+7}u_3tu_3u_1}}+z^{k+8}u_3tu^{-1}(R+Rs)t^{-1}s^{-1}\\
		&\in& U+z^{k+8}u_3t(u^{-1}t^{-1}s^{-1})+z^{k+8}u_3tu^{-2}(ust)t^{-2}s^{-1}\\
		&\in&U+\underline{z^{k+7}u_3u_2}+z^{k+9}u_3tu^{-2}(R+Rt^{-1}+Rt)s^{-1}\\
		&\in&U+\underline{\underline{z^{k+9}u_3tu_3u_1}}+z^{k+9}u_3tu^{-1}(u^{-1}t^{-1}s^{-1})+z^{k+9}u_3tu^{-2}t(R+Rs)\\
		&\in&U+\underline{z^{k+8}u_3tu^{-1}}+\underline{\underline{(z^{k+9}u_3tu_3)u_2}}
		+z^{k+9}u_3tu^{-2}t(stu)u^{-1}t^{-1}\\
		&\in& U+(z^{k+10}u_3tu^{-2}tu^{-1})u_2.
		\end{array}}$
		
		 The result follows from (i), (ii), (iii) and (iv), if we expand $u^{-2}$ as a linear combination of 1, $u$, $u^2$ and $u^{-1}$.

	\item[(vi)] 
	$\small{\begin{array}[t]{lcl}
	z^ktu^{-2}tu^{-1}&\in& z^ktu^{-2}(R+Rt^{-1}+Rt^{-2})u^{-1}\\
	&\in& z^ku_3tu_3u_2+z^ku_3tu^{-1}(u^{-1}t^{-1}s^{-1})su^{-1}+z^ku_3tu^{-1}(u^{-1}t^{-1}s^{-1})st^{-1}u^{-1}\\
	&\in&z^ku_3tu_3u_2+z^{k-1}u_3tu^{-1}(R+Rs^{-1})u^{-1}+z^{k-1}u_3tu^{-1}s(t^{-1}s^{-1}u^{-1})usu^{-1}\\
	&\in&(z^k+z^{k-1})u_3tu_3u_2+z^{k-1}u_3tu^{-1}(s^{-1}u^{-1}t^{-1})t+
	z^{k-2}u_3tu^{-1}su(R+Rs^{-1})u^{-1}\\
	&\in&(z^k+z^{k-1})u_3tu_3u_2+\underline{z^{k-2}u_3tu^{-1}u_2}+z^{k-2}u_3tu^{-1}(stu)u^{-1}t^{-1}+\\&&+z^{k-2}u_3tu^{-1}su(s^{-1}u^{-1}t^{-1})t\\
		&\in&U+(z^k+z^{k-1})u_3tu_3u_2+\underline{z^{k-1}u_3}+z^{k-3}u_3tu^{-1}(R+Rs^{-1})ut\\
			&\in&U+(z^k+z^{k-1})u_3tu_3u_2+\underline{z^{k-3}u_3u_2}+z^{k-3}u_3tu^{-1}(s^{-1}u^{-1}t^{-1})tu^2t\\
		&\in&U+(z^k+z^{k-1})u_3tu_3u_2+z^{k-4}u_3tu^{-1}t(R+Ru+Ru^{-1}+Ru^{-2})t\\
			&\in&U+(z^k+z^{k-1})u_3tu_3u_2+\underline{z^{k-4}u_3tu^{-1}u_2}+
			z^{k-4}u_3tu^{-1}(tus)s^{-1}t+
			\\&&+(z^{k-4}u_3
			tu^{-1}tu^{-1})t+z^{k-4}u_3tu^{-1}(R+Rt^{-1}+Rt^{-2})u^{-2}t\\
			&\stackrel{(iv)}{\in}&U+(z^k+z^{k-1})u_3tu_3u_2+z^{k-3}u_3tu^{-1}(R+Rs)t+
			\underline{\underline{(z^{k-4}u_3tu_3)u_2}}+\\&&+
			z^{k-4}u_3t(u^{-1}t^{-1}s^{-1})su^{-2}t+z^{k-4}u_3t(u^{-1}t^{-1}s^{-1})st^{-1}u^{-2}t\\
			&\in&U+(z^k+z^{k-1})u_3tu_3u_2+\underline{z^{k-3}u_3tu^{-1}u_2}+z^{k-3}u_3tu^{-2}(ust)+\\&&+
			z^{k-5}u_3t(R+Rs^{-1})u^{-2}t+z^{k-5}u_3ts(t^{-1}s^{-1}u^{-1})usu^{-2}t\\
			&\in&U+(z^k+z^{k-1}+z^{k-2})u_3tu_3u_2+\underline{\underline{(z^{k-5}u_3tu_3)u_2}}+
			z^{k-5}u_3t(s^{-1}u^{-1}t^{-1})tu^{-1}t+\\&&+z^{k-6}u_3t(R+Rs^{-1})usu^{-2}t\\
			&\in&U+(z^k+z^{k-1}+z^{k-2})u_3tu_3u_2+\underline{\underline{z^{k-6}u_3u_2u_3}}+z^{k-6}u_3(tus)u^{-2}t+\\&&+
			z^{k-6}u_3ts^{-1}u(R+Rs^{-1})u^{-2}t\\
			&\in&U+(z^k+z^{k-1}+z^{k-2})u_3tu_3u_2+\underline{z^{k-5}u_3u_2}+
			z^{k-6}u_3t(s^{-1}u^{-1}t^{-1})t^2+\\&&+
			z^{k-6}u_3t^2(t^{-1}s^{-1}u^{-1})u^2(s^{-1}u^{-1}t^{-1})tu^{-1}t\\
			&\in&U+(z^k+z^{k-1}+z^{k-2})u_3tu_3u_2+\underline{z^{k-7}u_3u_2}+
			z^{k-8}u_3(R+Rt+Rt^{-1})u^2tu^{-1}t\\
			&\in&U+(z^k+z^{k-1}+z^{k-2})u_3tu_3u_2+\underline{z^{k-8}u_3tu^{-1}u_2}+
			(z^{k-8}u_3tu^2tu^{-1})t+\\&&+
			z^{k-8}u_3(u^{-1}t^{-1}s^{-1})su^2tu^{-1}t\\
			&\stackrel{(iii)}{\in}&U+(z^k+z^{k-1}+z^{k-2})u_3tu_3u_2+z^{k-9}u_3(R+Rs^{-1})u^2tu^{-1}t\\
			&\in&U+(z^k+z^{k-1}+z^{k-2})u_3tu_3u_2+\underline{z^{k-9}u_3tu^{-1}u_2}+
		z^{k-9}u_3(s^{-1}u^{-1}t^{-1})tu^3tu^{-1}t\\
		&\in&U+(z^k+z^{k-1}+z^{k-2})u_3tu_3u_2+(z^{k-10}u_3tu^3tu^{-1})t.
	
	\end{array}}$
	
	However, if we expand $u^3$ as a linear combination of 1, $u$, $u^2$ and $u^{-1}$, we can use (i), (ii), (iii) and (iv) and we have that, for every $k\in\{16,\dots 23\}$, $z^{k-10}u_3tu^3tu^{-1}\subset U$.
	Therefore, for every $k\in\{16,\dots, 23\}$, $z^ktu^{-2}tu^{-1}\in U+(z^k+z^{k-1}+z^{k-2})u_3tu_3u_2$. 
	Moreover, by proposition \ref{tu11}(iii),  we have that  $\big((z^k+z^{k-1}+z^{k-2})u_3tu_3\big)u_2\subset U$, for every $k\in\{16,\dots, 19\}$ and, hence, $z^ktu^{-2}tu^{-1}\in U$ for every $k\in\{16,\dots, 19\}$.
	\qedhere
\end{itemize}
\end{proof}
	\begin{cor}
		\mbox{}
		\vspace*{-\parsep}
		\vspace*{-\baselineskip}\\
		\begin{itemize}[leftmargin=0.8cm]
			\item[(i)] For every $k\in\{2, \dots, 19\}$, $z^ktu_3tu^{-1}\in U$.
				\item[(ii)] For every $k\in\{1, \dots, 23\}$, $z^ku_2u^mtu^{-1}\subset U+z^ku_3tu^mtu^{-1}+z^{k-2}u_3tu^{m+1}tu^{-1}$, where $m\in \mathbb{Z}$. Therefore, for every $k\in\{4,\dots, 19\}$, $z^ku_3u_2u_3tu^{-1}\subset U$.
		\end{itemize}
		\label{tuts111}
	\end{cor}
	
	\begin{proof}
	We first prove (i). We use different definitions of $u_3$ and we have:
		\begin{itemize}[leftmargin=*]
			\item For $k\in\{2,\dots,5\}$, we write $u_3=R+Ru+Ru^2+Ru^3$. The result then follows from proposition \ref{tuts11} (i), (ii), (iii) and (v).
				\item For $k\in\{6,\dots,15\}$, we write $u_3=R+Ru+Ru^2+Ru^{-1}$. The result then follows from proposition \ref{tuts11} (i), (ii), (iii) and (iv).
				\item For $k\in\{16,\dots,19\}$, we write $u_3=R+Ru+Ru^{-2}+Ru^{-1}$. The result then follows from proposition \ref{tuts11} (i), (ii), (iv) and (vi).
		\end{itemize}
		 For the first part of (ii) we have:
		 \\
		$\small{\begin{array}[t]{lcl}
	 z^ku_2u^mtu^{-1}&=&z^k(R+Rt+Rt^{-1})u^mtu^{-1}\\
	 &\subset& \underline{z^ku_3tu^{-1}}+
		z^ku_3tu^mtu^{-1}+z^ku_3(u^{-1}t^{-1}s^{-1})su^mtu^{-1}\\
		&\subset& U+z^{k}u_3tu^mtu^{-1}+z^{k-1}u_3(R+Rs^{-1})u^mtu^{-1}\\
		&\subset& U+z^ku_3tu^mtu^{-1}
		+\underline{z^{k-1}u_3tu^{-1}}+z^{k-1}u_3(s^{-1}u^{-1}t^{-1})tu^{m+1}tu^{-1}\\
		&\subset& U+z^ku_3tu^mtu^{-1}+z^{k-2}u_3tu^{m+1}tu^{-1}.
		\end{array}}$\\
		Therefore, for every $k\in\{4,\dots,19\}$ we have that
		$z^ku_3u_2u_3tu^{-1}\subset z^ku_3(R+Rt+Rt^{-1})u_3tu^{-1}\subset u_3U+
		u_3t^{\pm}u_3tu^{-1}\subset u_3U+u_3(z^k+z^{k-2})tu_3tu^{-1}\stackrel{(i)}{\subset}u_3U$. The result follows from the definition of $U$.
		\end{proof}
	We now prove a lemma that leads us to the main theorem of this section (theorem \ref{thm11}).
	
	\begin{lem}
		\mbox{}
		\vspace*{-\parsep}
		\vspace*{-\baselineskip}\\
		\begin{itemize}[leftmargin=0.8cm]
			\item[(i)] For every $k\in\{0,\dots, 15\}$, $z^ku_3t^2u^2\subset U$.
			\item[(ii)] For every $k\in\{6,\dots,16\}$, $z^ktu^{-1}tu^{-1}tu^{-1}\in U+z^{k-6}u_3^{\times}t^2u^3t
		.$ 
			\item [(iii)]For every $k\in\{12,\dots,16\}$, 
			$z^ktu^{-1}tu^{2}tu^{-1}\in U+z^{k-12}u_3^{\times}t^2u^3t.$
			\item [(iv)]For every $k\in\{1,\dots,15\}$, 
			$z^ktu^{-1}tu^{2}tu^{-1}\in U+(z^{k+6}+z^{k+7}+z^{k+8})u_3tu_3u_2+z^{k+8}u_3^{\times}tu^{-3}tu^{-1}t.$
		\end{itemize}
		\label{ttuttu}
	\end{lem}
	\begin{proof}
		\mbox{}
		\vspace*{-\parsep}
		\vspace*{-\baselineskip}\\
		\begin{itemize}[leftmargin=0.8cm]
		\item[(i)]	
		$\small{\begin{array}[t]{lcl}
		z^ku_3t^2u^2&=&z^ku_3t(tus)s^{-1}u\\
		&\subset& z^{k+1}u_3t(R+Rs)u\\
		&\subset& \underline{\underline{z^{k+1}u_3tu}}+z^{k+1}u_3tu^{-1}(ust)t^{-1}u\\
		&\subset&U+z^{k+2}u_3tu^{-1}(R+Rt+Rt^2)u\\
		&\subset&U+\underline{z^{k+2}u_3u_2}+z^{k+2}u_3tu^{-1}(tus)s^{-1}+
		z^{k+2}u_3tu^{-1}t(tus)s^{-1}\\
		&\subset&\underline{\underline{z^{k+3}u_3tu_3u_1}}+z^{k+3}u_3tu^{-1}t(R+Rs)\\
		&\subset&U+\underline{z^{k+3}u_3tu^{-1}u_2}+z^{k+3}u_3tu^{-1}t(stu)u^{-1}t^{-1}\\
		&\in& U+(z^{k+4}u_3tu_3tu^{-1})t.
		\end{array}}$
	
			The result follows from corollary \ref{tuts111}(i).
			\item[(ii)] 
			$\small{\begin{array}[t]{lcl}
			z^ktu^{-1}tu^{-1}tu^{-1}&\in&
			z^{k}tu^{-1}(R+Rt^{-1}+R^{\times}t^{-2})u^{-1}tu^{-1}\\
			&\in&z^{k}u_3tu_3tu^{-1}+z^{k}u_3t(u^{-1}t^{-1}s^{-1})su^{-1}tu^{-1}+
			z^ku_3^{\times}t(u^{-1}t^{-1}s^{-1})st^{-1}u^{-1}tu^{-1}\\
			&\stackrel{\ref{tuts111}(i)}{\in}&U+z^{k-1}u_3t(R+Rs^{-1})u^{-1}tu^{-1}+z^{k-1}u_3^{\times}ts(t^{-1}s^{-1}u^{-1})usu^{-1}tu^{-1}\\ 
			&\in&U+z^{k-1}u_3tu_3tu^{-1}+z^{k-1}u_3t(s^{-1}u^{-1}t^{-1})t^2u^{-1}+\\&&+
			z^{k-2}u_3^{\times}tsu(R+R^{\times}s^{-1})u^{-1}tu^{-1}\\ 
				&\stackrel{\ref{tuts111}(i)}{\in}&
		U+\underline{\underline{z^{k-2}u_3u_2u_3}}+z^{k-2}u_3t(stu)u^{-2}+z^{k-2}u_3^{\times}tsu(s^{-1}u^{-1}t^{-1})t^2u^{-1}\\ 
			&\in&U+\underline{\underline{z^{k-1}u_3tu_3}}+z^{k-3}u_3^{\times}t(R+R^{\times}s^{-1})ut^2u^{-1}\\ 	
		&\in&U+z^{k-3}u_3(tus)s^{-1}t^2u^{-1}+
		z^{k-3}u_3^{\times}t^2(t^{-1}s^{-1}u^{-1})u^2t^2u^{-1}\\	
		&\in&U+z^{k-2}u_3(R+Rs)t^2u^{-1}+z^{k-4}u_3^{\times}t^2u^2(R+Rt+R^{\times}t^{-1})u^{-1}\\ 
			&\in&U+\underline{\underline{z^{k-2}u_3u_2u_3}}
			+z^{k-2}u_3(ust)tu^{-1}+\underline{\underline{z^{k-4}u_3u_2u_3}}+z^{k-4}u_3u_2u^2tu^{-1}+\\&&+z^{k-4}u_3^{\times}t^2u^2t^{-1}(u^{-1}t^{-1}s^{-1})st	\\ 
		&\stackrel{\ref{tuts111}(ii)}{\in}&U+\underline{z^{k-1}u_3tu^{-1}}+
		z^{k-6}u_3tu^3tu^{-1}+
		z^{k-5}u_3^{\times}t^2u^2t^{-1}(R+R^{\times}s^{-1})t\\
		&\stackrel{\ref{tuts11}(v)}{\in}&U+z^{k-5}u_3t^2u^2+z^{k-5}u_3^{\times}t^2u^3(u^{-1}t^{-1}s^{-1})t\\
		&\stackrel{(i)}{\subset}&U+z^{k-6}u_3^{\times}t^2u^3t.
			
			\end{array}}$
		\item[(iii)] 
		$\small{\begin{array}[t]{lcl}
		z^ktu^{-1}tu^{2}tu^{-1}&\in&z^ktu^{-1}t(R+Ru+Ru^{-1}+R^{\times}u^{-2})tu^{-1}\\
		&\in& z^ku_3tu^{-1}t^2u^{-1}+z^ku_3tu^{-1}(tus)s^{-1}tu^{-1}+z^ku_3tu^{-1}tu^{-1}tu^{-1}+\\&&+
		z^ku_3^{\times}tu^{-1}tu^{-2}tu^{-1}
		\\ 
		&\stackrel{(ii)}{\in}&z^ku_3tu^{-1}(R+Rt+Rt^{-1})u^{-1}+z^{k+1}u_3tu^{-1}(R+Rs)tu^{-1}+\\&&+\underline{\underline{(z^{k-6}u_3u_2u_3)t}} +z^{k}u_3^{\times}tu^{-1}(R+Rt^{-1}+R^{\times}t^{-2})u^{-2}tu^{-1}\\ 
		&\in& U+\underline{\underline{z^ku_3tu_3}}+(z^k+z^{k+1})u_3tu_3tu^{-1}+z^ku_3t(u^{-1}t^{-1}s^{-1})su^{-1}+\\&&+z^ku_3tu^{-1}(stu)u^{-2}
		+z^{k}u_3t(u^{-1}t^{-1}s^{-1})su^{-2}tu^{-1}+\\&&+
		z^{k}u_3^{\times}t(u^{-1}t^{-1}s^{-1})st^{-1}u^{-2}tu^{-1}\\ 
		&\stackrel{\ref{tuts111}(i)}{\in}&U+\underline{\underline{z^{k-1}u_3tu_1u_3}}+\underline{\underline{z^{k+1}u_3tu_3}}
		+z^{k-1}u_3t(R+Rs^{-1})u^{-2}tu^{-1}+\\&&+z^{k-1}u_3^{\times}tsu(u^{-1}t^{-1}s^{-1})su^{-2}tu^{-1}\\ 
		&\in&U+z^{k-1}u_3tu_3tu^{-1}+z^{k-1}u_3t(s^{-1}u^{-1}t^{-1})tu^{-1}tu^{-1}+\\&&+
		z^{k-2}u_3^{\times}tsu(R+R^{\times}s^{-1})u^{-2}tu^{-1}\\ 
		&\stackrel{\ref{tuts111}(i)}{\in}&U+z^{k-2}u_3u_2u_3tu^{-1}+
		z^{k-2}u_3tsu^{-1}tu^{-1}+
		z^{k-2}u_3^{\times}tsu(s^{-1}u^{-1}t^{-1})tu^{-1}tu^{-1}\\  
		&\stackrel{\ref{tuts111}(ii)}{\in}&U+z^{k-2}u_3t(R+Rs^{-1})u^{-1}tu^{-1}+z^{k-3}u_3^{\times}t(R+Rs^{-1})utu^{-1}tu^{-1}\\ 
		&\in&U+z^{k-2}u_3tu_3tu^{-1}+z^{k-2}u_3t(s^{-1}u^{-1}t^{-1})t^2u^{-1}+\\&&+
		z^{k-3}u_3(tus)s^{-1}tu^{-1}tu^{-1}+z^{k-3}u_3^{\times}t^2(t^{-1}s^{-1}u^{-1})u^2tu^{-1}tu^{-1}\\ 
			&\stackrel{\ref{tuts111}(i)}{\in}&U+\underline{\underline{z^{k-3}u_3u_2u_3}}+
			z^{k-2}u_3(R+Rs)tu^{-1}tu^{-1}+\\&&+
			z^{k-4}u_3^{\times}(R+Rt+Rt^{-1})u^2tu^{-1}tu^{-1}\\  \
		&\in&U+(z^{k-2}+z^{k-4})u_3u_2u_3tu^{-1}+z^{k-2}u_3(stu)u^{-2}tu^{-1}+z^{k-4}u_3tu^2tu^{-1}tu^{-1}+\\&&+z^{k-4}u_3^{\times}(u^{-1}t^{-1}s^{-1})su^2tu^{-1}tu^{-1}\\ 
			&\stackrel{\ref{tuts111}(ii)}{\in}&U+\underline{z^{k-1}u_3tu^{-1}}+
			z^{k-4}u_3tu^2tu^{-1}tu^{-1}+
			z^{k-5}u_3^{\times}(R+R^{\times}s^{-1})u^2tu^{-1}tu^{-1}\\ 
		&\in&U+z^{k-4}u_3tu^2tu^{-1}tu^{-1}+z^{k-5}u_3tu_3tu^{-1}+
		z^{k-5}u_3^{\times}(s^{-1}u^{-1}t^{-1})tu^3tu^{-1}tu^{-1}\\ 
	&\stackrel{\ref{tuts111}(i)}{\in}&U+z^{k-4}u_3tu^2tu^{-1}tu^{-1}+
	z^{k-6}u_3^{\times}tu^3tu^{-1}tu^{-1}.
	\end{array}}$
	
	It remains to prove that $B:=z^{k-4}u_3tu^2tu^{-1}tu^{-1}+
	z^{k-6}u_3^{\times}tu^3tu^{-1}tu^{-1}$ is a subset of $U+z^{k-12}u_3^{\times}t^2u^3t.$ We have:
	
	$\small{\begin{array}{lcl}
		B&=&z^{k-4}u_3tu^2tu^{-1}tu^{-1}+
		z^{k-6}u_3^{\times}tu^3tu^{-1}tu^{-1}\\
		&\subset&z^{k-4}u_3tu^2tu^{-1}tu^{-1}+
		z^{k-6}u_3^{\times}t(R+Ru+Ru^2+R^{\times}u^{-1})tu^{-1}tu^{-1}\\
		&\subset&U+(z^{k-4}+z^{k-6})u_3tu^2tu^{-1}tu^{-1}+z^{k-6}u_3u_2u_3tu^{-1}+ 
		z^{k-6}u_3(tus)s^{-1}tu^{-1}tu^{-1}+\\&&+z^{k-6}u_3^{\times}tu^{-1}tu^{-1}tu^{-1}\\ 	
	
			&\stackrel{\ref{tuts111}(ii)}{\subset}&U+(z^{k-4}+z^{k-6})u_3(tus)s^{-1}utu^{-1}tu^{-1}+
			z^{k-5}u_3(R+Rs)tu^{-1}tu^{-1}+\\&&+z^{k-6}u_3^{\times}tu^{-1}tu^{-1}tu^{-1}\\ 
		
	 &\stackrel{(ii)}{\subset}&U+(z^{k-3}+z^{k-5})u_3(R+Rs)utu^{-1}tu^{-1}+z^{k-5}u_3tu_3tu^{-1}+z^{k-5}u_3(ust)u^{-1}tu^{-1}+\\&&+z^{k-12}u_3^{\times}t^2u^3t\\
	 	&\stackrel{\ref{tuts111}(i)}{\subset}&U+(z^{k-3}+z^{k-5})u_3tu_3tu^{-1}+
	 	(z^{k-3}+z^{k-5})u_3(ust)t^{-2}(tus)s^{-1}tu^{-1}tu^{-1}+
	 \underline{z^{k-4}u_3tu^{-1}}+\\&&+
	 z^{k-12}u_3^{\times}t^2u^3t\\
	 &\stackrel{\ref{tuts111}(i)}{\subset}&U+(z^{k-1}+z^{k-3})u_3t^{-2}(R+Rs)tu^{-1}tu^{-1}+
	 z^{k-12}u_3^{\times}t^2u^3t\\
	 &\subset&U+(z^{k-1}+z^{k-3})u_3u_2u_3tu^{-1}+
	 (z^{k-1}+z^{k-3})u_3t^{-2}(stu)u^{-2}tu^{-1}+z^{k-12}u_3^{\times}t^2u^3t\\
	 &\subset&U+(z^{k-1}+z^{k-3}+z+z^{k-2})u_3u_2u_3tu^{-1}+z^{k-12}u_3^{\times}t^2u^3t\\
	 &\stackrel{\ref{tuts111}(ii)}{\in}&U+z^{k-12}u_3^{\times}t^2u^3t.
	 \end{array}}$
			
	\item[(iv)]	
	$\small{\begin{array}[t]{lcl}
	z^ktu^{-1}tu^2tu^{-1}&=&z^ktu^{-1}(tus)s^{-1}utu^{-1}\\
	&\in& z^{k+1}tu^{-1}(R+R^{\times}s)utu^{-1}\\
	&\in&\underline{\underline{z^{k+1}u_3u_2u_3}}+z^{k+1}u_3^{\times}tu^{-2}(ust)t^{-2}(tus)s^{-1}tu^{-1}\\
	&\in&U+z^{k+3}u_3^{\times}tu^{-2}t^{-2}(R+R^{\times}s)tu^{-1}\\
	&\in&U+z^{k+3}u_3tu^{-2}t^{-1}(u^{-1}t^{-1}s^{-1})st+z^{k+3}u_3^{\times}tu^{-2}t^{-2}(stu)u^{-2}\\
	&\in&U+
	z^{k+2}u_3tu^{-2}t^{-1}(R+Rs^{-1})t+z^{k+4}u_3^{\times}tu^{-2}(R+Rt^{-1}+R^{\times}t)u^{-2}\\
	&\in&U+\underline{\underline{z^{k+2}u_3tu_3}}+z^{k+2}u_3tu^{-1}(u^{-1}t^{-1}s^{-1})
	t+\underline{\underline{z^{k+4}u_3tu_3}}+\\&&+
	z^{k+4}u_3tu^{-1}(u^{-1}t^{-1}s^{-1})su^{-2}+z^{k+4}u_3^{\times}tu^{-2}tu^{-2}\\
	&\in&U+\underline{z^{k+1}u_3tu^{-1}u_2}+z^{k+3}u_3tu^{-1}(R+Rs^{-1})u^{-2}+\\&&+
	z^{k+4}u_3^{\times}tu^{-2}t(R+Ru+Ru^{-1}+R^{\times}u^{2})\\
	&\in&U+\underline{\underline{z^{k+3}u_3tu_3}}+
	z^{k+3}u_3tu^{-1}(s^{-1}u^{-1}t^{-1})tu^{-1}+\underline{\underline{(z^{k+4}u_3tu_3)t}}+\\&&+
	z^{k+4}u_3tu^{-2}(tus)s^{-1}+z^{k+4}u_3tu_3tu^{-1}+z^{k+4}u_3^{\times}tu^{-2}(tus)s^{-1}u\\
	&\stackrel{\ref{tuts111}(i)}{\in}&U+z^{k+2}u_3tu_3tu^{-1}+z^{k+5}u_3tu^{-2}(R+Rs)+z^{k+5}u_3^{\times}tu^{-2}(R+Rs)u\\
	&\stackrel{\ref{tuts111}(i)}{\in}&U+\underline{\underline{z^{k+5}u_3tu_3}}+z^{k+5}u_3tu^{-2}(stu)u^{-1}t^{-1}+\underline{z^{k+5}u_3tu^{-1}}+\\&&+
	z^{k+5}u_3^{\times}tu^{-3}(ust)t^{-1}u\\
	&\in&U+z^{k+6}u_3tu_3u_2+z^{k+6}u_3^{\times}tu^{-3}(R+Rt+R^{\times}t^2)u\\
	&\in&U+z^{k+6}u_3tu_3u_2+z^{k+6}u_3tu^{-3}(tus)s^{-1}+
	z^{k+6}u_3^{\times}tu^{-3}t(tus)s^{-1}\\
		&\in&U+z^{k+6}u_3tu_3u_2+z^{k+7}u_3tu^{-3}(R+Rs)+
		z^{k+7}u_3^{\times}tu^{-3}t(R+R^{\times}s)\\
			&\in&U+(z^{k+6}+z^{k+7})u_3tu_3u_2+z^{k+7}u_3tu^{-4}(ust)t^{-1}
			+z^{k+7}u_3^{\times}tu^{-3}t(stu)u^{-1}t\\
			&\in&U+(z^{k+6}+z^{k+7}+z^{k+8})u_3tu_3u_2+z^{k+8}u_3^{\times}tu^{-3}tu^{-1}t.

	\end{array}}$\\
\qedhere
			
		\end{itemize}
	\end{proof}

		\begin{thm}	
			\mbox{}
			\vspace*{-\parsep}
			\vspace*{-\baselineskip}\\
			\begin{itemize}[leftmargin=0.8cm]
			\item[(i)]For every $k\in\{0,\dots,23\}$, $z^ktu_3\subset U$.
			\item[(ii)]For every $k\in\{0,\dots,23\}$, $z^ktu_3tu^{-1}\subset U$.
			\item[(iii)]$H_{G_{11}}=U$.
			\label{thm11}
		\end{itemize}
		\end{thm}
		
	\begin{proof}
		\mbox{}
		\vspace*{-\parsep}
		\vspace*{-\baselineskip}\\
		
	\begin{itemize}[leftmargin=0.8cm]
		\item [(i)] By proposition \ref{tu11} (iii), we have to prove that $z^ktu_3\subset U$, for every $k\in\{20,\dots,23\}$. We use different definitions of $u_3$ and we have:
		\begin{itemize}[leftmargin=*]
			\item $\underline{k\in\{20,21\}}$: $z^ktu_3=z^kt(R+Ru+Ru^{-1}+Ru^{-2})\subset \underline{z^ku_3t}+\underline{\underline{z^ku_3tu}}+\underline{z^ku_3tu^{-1}}+z^ku_3tu^{-2}$. Therefore, $z^ktu_3\subset U+ \mathbf{z^ku_3tu^{-2}}$.
			\item $\underline{k\in\{22,23\}}$: $z^ktu_3=z^kt(R+Ru^{-1}+Ru^{-2}+Ru^{-3})\subset \underline{z^ku_3t}+{\underline{z^ku_3tu^{-1}}}+z^ku_3tu^{-2}+z^ku_3tu^{-3}$. Therefore, $z^ktu_3\subset U+ \mathbf{z^ku_3tu^{-2}+z^ku_3tu^{-3}}$.
		\end{itemize}
		As a result, it will be sufficient to prove that, for every $k\in\{20,\dots,23\}$, $z^ku_3tu^{-2}$ is a subset of $U$,  and that, for every $k\in\{22,23\}$,  $z^ku_3tu^{-3}$ is also a subset of $U$.
		We have:
		\\ \\
		$\begin{array}{lcl}
		z^ku_3tu^{-2}&=&z^ku_3t^2(t^{-1}s^{-1}u^{-1})usu^{-2}\\
		&\subset& z^{k-1}u_3t^2u(R+Rs^{-1})u^{-2}\\
		&\subset&z^{k-1}u_3t^2(u^{-1}t^{-1}s^{-1})st+z^{k-1}u_3t^2u(s^{-1}u^{-1}t^{-1})tu^{-1}\\
		&\subset& U+z^{k-2}u_3t^2(R+Rs^{-1})t+z^{k-2}u_3u_2utu^{-1}\\
		&\stackrel{\ref{tuts111}(ii)}{\subset}&U+\underline{\underline{z^{k-2}u_3u_2u_3}}+z^{k-2}u_3t^3(t^{-1}s^{-1}u^{-1})ut+z^{k-2}u_3tutu^{-1}+z^{k-4}u_3tu^2tu^{-1}.
			\end{array}$
	\\ \\
		Therefore, by proposition  \ref{tuts11}(ii) and by corollary \ref{tuts111}(i) we have
		\begin{equation}
		z^ku_3tu^{-2}\subset U,
		\label{tuu111}
		\end{equation}
		for every $k\in\{20,\dots,23\}$.
		Moreover, since $z^ktu_3 \subset U+z^ku_3tu^{-2}$, $k\in\{20,21\}$,  we use proposition \ref{tu11}(iii) and we have that
		\begin{equation}
		z^ku_3tu_3\subset U,
		\label{tuu1111}
		\end{equation}
		for every $k\in\{0,\dots,21\}$.
		
		We now prove that $z^ku_3tu^{-3}\subset U$, for every $k\in\{22,23\}$.
	We have:\\ \\
	$\begin{array}[t]{lcl}
	z^ku_3tu^{-3}&\subset&z^ku_3(R+Rt^{-1}+Rt^{-2})u^{-3}\\
	&\subset&\underline{z^ku_3}+z^ku_3(u^{-1}t^{-1}s^{-1})su^{-3}+
	z^ku_3t^{-1}u(u^{-1}t^{-1}s^{-1})su^{-3}\\ 
	&\subset&U+z^{k-1}u_3(R+Rs^{-1})u^{-3}+z^{k-1}u_3t^{-1}u(R+Rs^{-1})u^{-3}\\ 
	&\subset&U+\underline{z^{k-1}u_3}+z^{k-1}u_3(s^{-1}u^{-1}t^{-1})tu^{-2}+ z^{k-1}u_3(u^{-1}t^{-1}s^{-1})su^{-2}+\\&&+
	z^{k-1}u_3(u^{-1}t^{-1}s^{-1})su(s^{-1}u^{-1}t^{-1})tu^{-2}\\ 
	&\subset& U+z^{k-2}u_3tu^{-2}+z^{k-2}u_3(R+Rs^{-1})u^{-2}+
	z^{k-3}u_3(R+Rs^{-1})utu^{-2}\\
	&\stackrel{(\ref{tuu111})}{\subset}&U+
		\underline{z^{k-2}u_3}+z^{k-2}u_3(s^{-1}u^{-1}t^{-1})tu^{-1}+
		z^{k-3}u_3tu_3+z^{k-3}u_3(s^{-1}u^{-1}t^{-1})tu^2tu^{-2}\\ 
	&\stackrel{(\ref{tuu1111})}{\subset}& U+\underline{z^{k-3}u_3tu^{-1}}+
	z^{k-4}u_3tu^2tu^{-2}.
	\end{array}$

Therefore, it will be sufficient to prove that $z^{k-4}u_3tu^2tu^{-2}$ is a subset of $U$. For this purpose, we expand $u^2$ as a linear combination of 1, $u$, $u^{-1}$ and $u^{-2}$ and we have:

$\small{\begin{array}{lcl}
	z^{k-4}u_3tu^2tu^{-2}
	&\subset& U+\underline{\underline{z^{k-4}u_3u_2u_3}}+
	z^{k-4}u_3tutu^{-2}+z^{k-4}u_3tu^{-1}tu^{-2}+
	z^{k-4}u_3tu^{-2}tu^{-2}
	\end{array}}$

However, $	z^{k-4}u_3tutu^{-2}=z^{k-4}u_3(tus)s^{-1}tu^{-2}=z^{k-3}u_3s^{-1}tu^{-2}$. If we expand $s^{-1}$ as a linear combination of 1 and $s$ we have that
$z^{k-3}u_3s^{-1}tu^{-2}\subset  z^{k-3}u_3tu^{-2}+z^{k-3}u_3(stu)u^{-3}=z^{k-3}u_3tu^{-2}+\underline{z^{k-2}u_3}$.
Therefore, by relation (\ref{tuu1111}) we have that $z^{k-3}u_3s^{-1}tu^{-2}\subset U$ and, hence, $z^{k-4}u_3tutu^{-2}\subset U$. 

It remains to prove that $C:=z^{k-4}u_3tu^{-1}tu^{-2}+
z^{k-4}u_3tu^{-2}tu^{-2}$ is a subset of $U$. We have:

$\small{\begin{array}{lcl}
C&=&z^{k-4}u_3tu^{-1}tu^{-2}+
z^{k-4}u_3tu^{-2}tu^{-2}\\
&\subset&z^{k-4}u_3t^2(t^{-1}s^{-1}u^{-1})usu^{-1}tu^{-2}+
	z^{k-4}u_3tu^{-2}(R+Rt^{-1}+Rt^{-2})u^{-2}\\
	
	&\subset&z^{k-5}u_3t^2u(R+Rs^{-1})u^{-1}tu^{-2}+
	\underline{\underline{z^{k-4}u_3tu_3}}+
	z^{k-4}u_3tu^{-1}(u^{-1}t^{-1}s^{-1})su^{-2}+\\&&+z^{k-4}u_3tu^{-2}
	t^{-1}(t^{-1}s^{-1}u^{-1})usu^{-2}\\

	&\subset& U+\underline{\underline{z^{k-5}u_3u_2u_3}}+z^{k-5}u_3t^2u(s^{-1}u^{-1}t^{-1})t^2u^{-2}+z^{k-5}u_3tu^{-1}(R+Rs^{-1})u^{-2}+\\&&+z^{k-5}u_3tu^{-2}t^{-1}u(R+Rs^{-1})u^{-2}\\ 
	&\subset& U+z^{k-6}u_3t^2ut^2u^{-2}+\underline{\underline{z^{k-5}u_3u_2u_3}}+
	z^{k-5}u_3tu^{-1}(s^{-1}u^{-1}t^{-1})tu^{-1}+\\&&+
	z^{k-5}u_3tu^{-1}(u^{-1}t^{-1}s^{-1})su^{-1}+
	z^{k-5}u_3tu^{-1}(u^{-1}t^{-1}s^{-1})su(s^{-1}u^{-1}t^{-1})tu^{-1}\\ 
	&\subset& U+z^{k-6}u_3t^2u(R+Rt+Rt^{-1})u^{-2}
	+z^{k-6}u_3tu_3tu^{-1}+z^{k-6}u_3tu^{-1}(R+Rs^{-1})u^{-1}+\\&&+
	z^{k-7}u_3tu^{-1}(R+Rs^{-1})utu^{-1}\\ 
	&\stackrel{\ref{tuts111}(i)}{\subset}&U+
	\underline{\underline{z^{k-6}u_3u_2u_3}}+z^{k-6}u_3t(tus)s^{-1}tu^{-2}+
	z^{k-6}u_3t^2u^2(u^{-1}t^{-1}s^{-1})su^{-2}+\\&&+

	\underline{\underline{(z^{k-6}+z^{k-7})u_3u_2u_3}}+z^{k-6}u_3tu^{-1}(s^{-1}u^{-1}t^{-1})t
	+z^{k-7}u_3tu^{-1}(s^{-1}u^{-1}t^{-1})tu^2tu^{-1}\\ 
	&\subset&U+z^{k-5}u_3t(R+Rs)tu^{-2}+z^{k-7}u_3t^2u^2(R+Rs^{-1})u^{-2}+
	\underline{z^{k-7}u_3tu^{-1}u_2}+\\&&+
	z^{k-8}u_3tu^{-1}tu^2tu^{-1}\\ &\stackrel{\ref{ttuttu}(iii)}{\subset}&U+\underline{\underline{z^{k-5}u_3u_2u_3}}+z^{k-5}u_3t(stu)u^{-3}+
	\underline{z^{k-7}u_3u_2}+z^{k-7}u_3t^2u^2(s^{-1}u^{-1}t^{-1})tu^{-1}+\\&&+
	\underline{\underline{(z^{k-20}u_3u_2u_3)t}}
	 \\
	 &\subset&U+\underline{\underline{z^{k-4}tu_3}}+
	z^{k-8}u_3u_2u_3tu^{-1}.
	\end{array}}$
	
	The result follows from corollary \ref{tuts111}(ii).
\item[(ii)] By corollary \ref{tuts111}(i), we restrict ourselves to proving the cases where $k\in\{0,1\}\cup\{20,\dots,23\}$.
We distinguish the following cases:
\begin{itemize}[leftmargin=*]
	\item [C1.] \underline{$k\in\{20,21\}$}: We expand $u_3$ as $R+Ru+Ru^{-1}+Ru^{-2}$ and by proposition \ref{tuts11}(i), (ii), (iv) and (v) we have that
	$z^ktu_3tu^{-1}\subset U+(z^k+z^{k-1}+z^{k-2})u_3tu_3$. The result follows from (i).
	\item [C2.]\underline{$k\in\{22,23\}$}:  We expand $u_3$ as $R+Ru^{-1}+Ru^{-2}+Ru^{-3}$ and by proposition \ref{tuts11}(i), (iv) and (v) we have that $z^ktu_3tu^{-1}\subset U+(z^k+z^{k-1}+z^{k-2})u_3tu_3+z^ktu^{-3}tu^{-1}\stackrel{(i)}{\subset}U+z^ktu^{-3}tu^{-1}$.  However, since $k-8\in\{1,\dots,15\}$, we can apply lemma \ref{ttuttu}(iv) and we have that
	$z^ktu^{-3}tu^{-1}\in
	u_3^{\times}Ut^{-1}+(z^k+z^{k-1}+z^{k-2})u_3tu_3u_2+z^{k-8}u_3^{\times}tu^{-1}tu^2tu^{-1}t^{-1}$.
	However, by (i) and by lemma \ref{ttuttu}(iii) we also have that
	$z^ktu^{-3}tu^{-1}\in u_3^{\times}Uu_2+\underline{\underline{(z^{k-20}u_3u_2u_3)t}}$.
	The result then follows from the definition of $U$ and remark \ref{rem11}.
	\item [C3.]\underline{$k\in\{0,1\}$}: We expand $u_3$ as $R+Ru+Ru^2+Ru^{3}$ and by proposition \ref{tuts11}(ii), (iii) and (v) we have that
	 $z^ktu_3tu^{-1}\subset U+z^kt^2u^{-1}\subset U+z^kt^2(R+Ru+Ru^2+Ru^3)\subset U+\underline{z^ku_2}+z^ku_3t(tus)s^{-1}+z^ku_3t^2u^2+z^ku_3t^2u^3 \stackrel{\ref{ttuttu}(i)}{\subset} U+\underline{\underline{z^{k+1}u_3u_2u_1}}+z^ku_2t^2u^3$. 
	However, by lemma \ref{ttuttu}(iii), we have that $z^kt^2u_3\in u_3^{\times}Ut^{-1}+z^{k+12}u_3^{\times}tu^{-1}tu^2tu^{-1}t^{-1}$. Since $k+12\in\{1,\dots,15\}$, we can apply lemma \ref{ttuttu}(iv) we have that
	 $z^{k+12}tu^{-1}tu^2tu^{-1}t^{-1}\in Ut^{-1}+( z^{k+18}+z^{k+19}+z^{k+20})u_3tu_3u_2+z^{k+20}u_3tu_3tu^{-1}$.
	 However, by (i) and by case C1, we have that $z^{k+12}tu^{-1}tu^2tu^{-1}t^{-1}\in Uu_2$.
	Therefore, $z^ktu_3tu^{-1}\subset u_3Ut^{-1}+U$. The result follows from the definition of $U$ and remark \ref{rem11}.

\end{itemize}

		\item[(iii)]
		The result follows immediately from (i) and (ii) (see corollary \ref{corr11}(iii)).
		\qedhere	
\end{itemize}	
	\end{proof}	
\begin{cor}
	The BMR freeness conjecture holds for the generic Hecke algebra $H_{G_{11}}$.
\end{cor}
\begin{proof}
	By theorem \ref{thm11}(iii) we have that $H_{G_{11}}=U$. The result follows from proposition \ref{BMR PROP}, since by definition $U$ is generated as left $u_3$ module by 144 elements and, hence, as $R$-module by $|G_{11}|=576$ elements (recall that $u_3$ is generated as $R$-module by 4 elements).
\end{proof}

\subsection{The case of $G_{12}$}
\indent

Let $R=\ZZ[u_{s,i}^{\pm},u_{t,j}^{\pm},u_{u,l}^{\pm}]_{\substack{1\leq i,j,l\leq 2 }}$ and let $$ H_{G_{12}}=\langle s,t,u\;|\; stus=tust=ustu,\prod\limits_{i=1}^{2}(s-u_{s,i})=\prod\limits_{j=1}^{2}(t-u_{t,j})=\prod\limits_{l=1}^{2}(u-u_{u,l})=0\rangle.$$ 
 Let also $\bar R=\ZZ[u_{s,i}^{\pm},]_{\substack{1\leq i\leq 2 }}$. Under the specialization $\grf: R\twoheadrightarrow \bar R$, defined by $u_{t,j}\mapsto u_{s,i}$, $u_{u,l}\mapsto u_{s,i}$ the algebra $\bar H_{G_{12}}:=H_{G_{12}}\otimes_{\grf}\bar R$ is the generic Hecke algebra associated to $G_{12}$. 
 
 Let $u_1$ be the subalgebra of $H_{G_{12}}$ generated by $s$, $u_2$ the subalgebra of $H_{G_{12}}$ generated by $t$ and $u_3$ the subalgebra of $H_{G_{12}}$ generated by $u$. We recall that $z:=(stu)^4$ generates the center of the  complex braid group associated to $G_{12}$ and that $|Z(G_{12})|=2$.  We set $U=\sum\limits_{k=0}^1(z^ku_1u_3u_1u_2+z^ku_1u_2u_3u_2+z^ku_2u_3u_1u_2+z^ku_2u_1u_3u_2+z^ku_3u_2u_1u_2).$
By the definition of $U$ we have the following remark:
\begin{rem} $Uu_2\subset U$.
	\label{r12}
\end{rem}
To make it easier for the reader to follow the next calculations, we will underline the elements that by definition belong to $U$. Moreover, we will use directly remark \ref{r12}; this means that every time we have a power of $t$ at the end of an element we may ignore it. In order to remind that to the reader, we put a parenthesis around the part of the element we consider.

\begin{prop}$u_1U\subset U$.
	\label{su12}
\end{prop}
\begin{proof}
	Since $u_1=R+Rs$ we have to prove that $sU\subset U$. By the definition of $U$ and by remark \ref{r12}, it is enough to prove that for every $k\in\{0,1\}$, $z^ksu_2u_3u_1$, $z^ksu_2u_1u_3$ and $z^ksu_3u_2u_1$ are subsets of $U$.
	We have:
	\begin{itemize}[leftmargin=*]
	\item $\small{\begin{array}[t]{lcl}
	z^ksu_2u_3u_1&=&z^ks(R+Rt)u_3u_1\\
	&\subset&\underline{z^ku_1u_3u_1}+z^kst(R+Ru)u_1\\
	&\subset&U+z^ku_1u_2u_1+z^kstu(R+Rs)\\
	&\subset&U+z^ku_1u_2u_1+\underline{z^ku_1u_2u_3}+Rz^kstus\\
	&\subset&U+z^ku_1u_2u_1+Rz^ktust\\
		&\subset&U+z^ku_1u_2u_1+\underline{u_2u_3u_1u_2}\\
		&\subset&U+z^ku_1u_2u_1.
		\end{array}}$
		
		\item $\small{\begin{array}[t]{lcl}
		z^ksu_2u_1u_3&=&z^ks(R+Rt^{-1})u_1u_3\\
		&\subset&\underline{z^ku_1u_3}+z^kst^{-1}(R+Rs^{-1})u_3\\
		&\subset&U+\underline{z^ku_1u_2u_3}+z^kst^{-1}s^{-1}(R+Ru^{-1})\\
		&\subset&U+z^ku_1u_2u_1+Rz^ks(t^{-1}s^{-1}u^{-1}t^{-1})t\\
			&\subset&U+z^ku_1u_2u_1+Rz^ku^{-1}t^{-1}s^{-1}t\\
			&\subset&U+z^ku_1u_2u_1+\underline{z^ku_3u_2u_1u_2}\\
				&\subset&U+z^ku_1u_2u_1.
		\end{array}}$
		\item $\small{\begin{array}[t]{lcl}
		z^ksu_3u_2u_1&\subset&z^k(R+Rs^{-1})u_3u_2u_1\\
		&\subset&\underline{z^ku_3u_2u_1}+z^ks^{-1}(R+Ru^{-1})u_2u_1\\
		&\subset&U+z^ku_1u_2u_1+z^ks^{-1}u^{-1}(R+Rt^{-1})u_1\\
		&\subset&U+z^ku_1u_2u_1+\underline{z^ku_1u_3u_1}+z^ks^{-1}u^{-1}t^{-1}(R+Rs^{-1})\\
		&\subset&U+z^ku_1u_2u_1+\underline{z^ku_1u_3u_2}+Rz^k(s^{-1}u^{-1}t^{-1}s^{-1})\\
		&\subset&U+z^ku_1u_2u_1+Rz^kt^{-1}s^{-1}u^{-1}t^{-1}\\
		&\subset&U+z^ku_1u_2u_1+\underline{z^ku_2u_1u_3u_2}\\
		&\subset&U+z^ku_1u_2u_1.
		\end{array}}$
	\end{itemize}
	Therefore, we have to prove that for every $k\in\{0,1\}$, $z^ku_1u_2u_1\subset U$.
	We distinguish the following cases:
	\begin{itemize}[leftmargin=*]
		\item $\underline{k=0}$: 
		$\small{\begin{array}[t]{lcl}
		u_1u_2u_1&=&(R+Rs)(R+Rt)(R+Rs)\\
		&\subset&\underline{u_2u_1u_2}+Rsts\\
		&\subset&U+
		(stu)^4u^{-1}(t^{-1}s^{-1}u^{-1}t^{-1})s^{-1}(u^{-1}t^{-1}s^{-1}u^{-1})s\\
		&\subset&U+zu^{-2}t^{-1}s^{-1}u^{-1}s^{-2}u^{-1}t^{-1}\\
			&\subset&U+zu^{-2}t^{-1}s^{-1}u^{-1}(R+Rs^{-1})u^{-1}t^{-1}\\
				&\subset&U+zu^{-2}t^{-1}s^{-1}u^{-2}t^{-1}+u^{-2}t^{-1}s^{-1}zu^{-1}s^{-1}u^{-1}t^{-1}\\
				&\subset&U+zu^{-2}t^{-1}s^{-1}(R+Ru^{-1})t^{-1}+u^{-2}t^{-1}s^{-1}(stu)^4u^{-1}s^{-1}u^{-1}t^{-1}\\
				&\subset&U+\underline{zu_3u_2u_1u_2}+zu^{-1}(u^{-1}t^{-1}s^{-1}u^{-1})t^{-1}+
				u^{-1}(stus)(tust)s^{-1}u^{-1}t^{-1}\\
				&\subset&U+z(u^{-1}t^{-1}s^{-1}u^{-1})t^{-2}+(stus)\\
					&\subset&U+zt^{-1}s^{-1}u^{-1})t^{-3}+tust\\
					&\subset&U+\underline{zu_2u_1u_3u_2}+\underline{u_2u_3u_1u_2}.
		\end{array}}$
			\item $\underline{k=1}$:
			$\small{\begin{array}[t]{lcl}
			zu_1u_2u_1&=&(R+Rs^{-1})(R+Rt^{-1})(R+Rs^{-1})\\
			&\subset&\underline{zu_2u_1u_2}+zRs^{-1}t^{-1}s^{-1}\\
			&\subset&U+Rs^{-1}t^{-1}s^{-1}z\\
			&\subset&U+Rs^{-1}t^{-1}s^{-1}(stu)^4\\
			&\subset&U+Rs^{-1}(ustu)s(tust)u\\
			&\subset&U+Rtus^2ustu^2\\
			&\subset&U+Rtu(R+Rs)ustu^2\\
			&\subset&U+Rtu^2st(R+Ru)+Rtus(ustu)u\\
			&\subset&U+\underline{u_2u_3u_1u_2}+Rtu^2(stu)^4(stu)^{-3}+Rtustustu\\

			&\subset&U+Rztu(t^{-1}s^{-1}u^{-1}t^{-1})s^{-1}u^{-1}t^{-1}s^{-1}+Rtu(stu)^4(stu)^{-2}\\

			&\subset&U+Rzs^{-1}u^{-1}(s^{-1}u^{-1}t^{-1}s^{-1})+	Rz(s^{-1}u^{-1}t^{-1}s^{-1})\\

			&\subset&U+Rz(s^{-1}u^{-1}t^{-1}s^{-1})u^{-1}t^{-1}+R	zt^{-1}s^{-1}u^{-1}t^{-1}\\

			&\subset&U+Rzu^{-1}t^{-1}s^{-1}u^{-2}t^{-1}+\underline{zu_2u_1u_3u_2}\\
			&\subset&U+zu^{-1}t^{-1}s^{-1}(R+Ru^{-1})t^{-1}\\
			&\subset&U+\underline{zu_3u_2u_1u_2}+z(u^{-1}t^{-1}s^{-1}u^{-1})t^{-1}\\
			&\subset&U+zt^{-1}s^{-1}u^{-1}t^{-2}\\
			&\subset&U+\underline{zu_2u_1u_3u_2}.
			\end{array}}$\\
			\qedhere
	\end{itemize}
	
\end{proof}

Our goal is to prove that $H_{G_{12}}=U$ (theorem \ref{thm12}). The following proposition provides a criterium for this to be true.

\begin{prop} If $u_3U\subset U$, then $H_{G_{12}}=U$.
	\label{prop12}
	\end{prop}
	\begin{proof}
		Since  $1\in U$, it is enough to prove that $U$ is a left-sided ideal of $H_{G_{11}}$. For this purpose, one may check that $sU$, $tU$ and $uU$ are subsets of $U$. By hypothesis and proposition \ref{prop12}, it is enough to prove that $tU\subset U$. By the definition of $U$ and remark \ref{r12} we have to prove that for every $k\in\{0,1\}$ $z^ktu_1u_3u_1$, $z^ktu_1u_2u_3$ and $z^ktu_3u_2u_1$ are subsets of $U$. We have:
		\begin{itemize}[leftmargin=*]
			\item $\small{\begin{array}[t]{lcl}
			z^ktu_1u_3u_1&\subset&z^k(R+Rt^{-1})u_1u_3u_1\\
			&\subset&\underline{z^ku_1u_3u_1}+z^kt^{-1}(R+Rs^{-1})(R+Ru)(R+Rs)\\
			&\subset&U+\underline{z^ku_2u_3u_1}+\underline{z^ku_2u_1u_3}+Rz^kt^{-1}s^{-1}us\\
			&\subset&U+Rz^kt^{-1}s^{-1}(ustu)u^{-1}t^{-1}\\
			&\subset&U+Rz^kusu^{-1}t^{-1}\\
			&\subset&U+u_3\underline{z^ku_1u_3u_2}\\
			&\subset& U+u_3U.
			\end{array}}$
			\item $\small{\begin{array}[t]{lcl}
			z^ktu_1u_2u_3&\subset&z^k(R+Rt^{-1})u_1u_2u_3\\
			&\subset&\underline{z^ku_1u_2u_3}+z^kt^{-1}(R+Rs)(R+Rt)(R+Ru)\\
			&\subset&U+\underline{z^ku_2u_1u_3}+\underline{z^ku_2u_1u_2}+Rz^kt^{-1}stu\\
			&\subset&U+Rz^kt^{-1}(stus)s^{-1}\\
			&\subset&U+Rz^kusts^{-1}\\
			&\subset&U+u_3u_1\underline{z^ku_2u_1}\\
			&\subset&U+u_3u_1U.
			\end{array}}$
			\item $\small{\begin{array}[t]{lcl}
			z^ktu_3u_2u_1&=&z^kt(R+Ru^{-1})(R+Rt^{-1})(R+Rs^{-1})\\
			&\subset&\underline{z^ku_2u_3u_1u_2}+Rz^ktu^{-1}t^{-1}s^{-1}\\
			&\subset&U+Rz^kt(u^{-1}t^{-1}s^{-1}u^{-1})u\\
			&\subset&U+Rz^ks^{-1}u^{-1}t^{-1}u\\
			&\subset&U+u_1u_3\underline{z^ku_2u_3}\\
			&\subset&U+u_1u_3U.
		\end{array}}$
	
	The result follows from proposition \ref{su12} and remark \ref{r12}.
			\qedhere
		\end{itemize}
	\end{proof}
	We can now prove the main theorem of this section.
	\begin{thm}$H_{G_{12}}=U$.
		\label{thm12}
	\end{thm}
	\begin{proof}
		By proposition \ref{prop12} it is enough to prove that $u_3U\subset U$. Since $u_3=R+Ru$, it will be sufficient to check that $uU\subset U$. By the definition of $U$ and remark \ref{r12}, we only have to prove that for every $k\in\{0,1\}$, $z^kuu_1u_3u_1$, $z^kuu_1u_2u_3$, $z^kuu_2u_3u_1$ and $z^kuu_2u_1u_3$ are subsets of $U$. 
		\begin{itemize}[leftmargin=0.8cm]
			\item[C1.] We will prove that $z^kuu_1u_3u_1\subset U$,  $k\in\{0,1\}$.
			\begin{itemize}[leftmargin=*]
				\item [(i)]$\underline{k=0}$: 
				
	$\small{\begin{array}[t]{lcl}
				uu_1u_3u_1&=&u(R+Rs)(R+Ru)(R+Rs)\\
				&\subset&\underline{u_3u_1}+Rusu+Rusus\\
				&\subset&U+Rzuz^{-1}su+Rzusuz^{-1}s\\
				&\subset&U+Rzu(stu)^{-4}su+Rzusu(stu)^{-4}s\\
				&\subset&U+Rzt^{-1}(s^{-1}u^{-1}t^{-1}s^{-1})u^{-1}(t^{-1}s^{-1}u^{-1}t^{-1}u)+\\&&+
				Rzus(t^{-1}s^{-1}u^{-1}t^{-1})s^{-1}(u^{-1}t^{-1}s^{-1}u^{-1})t^{-1}\\
				&\subset&U+Rzt^{-2}s^{-1}u^{-1}t^{-1}u^{-2}t^{-1}s^{-1}+Rzt^{-1}s^{-2}t^{-1}s^{-1}u^{-1}t^{-2}\\
				&\subset&U+Rz(R+Rt^{-1})s^{-1}u^{-1}t^{-1}u^{-2}t^{-1}s^{-1}+
				Rzt^{-1}(R+Rs^{-1})t^{-1}s^{-1}u^{-1}t^{-2}\\
				&\subset&U+Rzs^{-1}u^{-1}t^{-1}(R+Ru^{-1})t^{-1}s^{-1}+
				Rt^{-1}s^{-1}u^{-1}t^{-1}u^{-2}t^{-1}s^{-1}z+
				\underline{zu_2u_1u_3u_2}+\\&&+Rzt^{-1}s^{-1}(t^{-1}s^{-1}u^{-1}t^{-1})t^{-1}\\
				&\subset&U+u_1\underline{zu_3u_2u_1}+Rs^{-1}u^{-1}t^{-1}u^{-1}t^{-1}s^{-1}z+Rt^{-1}s^{-1}u^{-1}t^{-1}u^{-2}t^{-1}s^{-1}(stu)^4+\\&&+
				Rt^{-1}s^{-2}u^{-1}t^{-1}s^{-1}zt^{-1}\\ 

				&\subset&U+u_1U+Rs^{-1}u^{-1}t^{-1}u^{-1}t^{-1}s^{-1}(stu)^4+
				R(t^{-1}s^{-1}u^{-1}t^{-1})u^{-1}(stus)(tust)u+\\&&+
				Rt^{-1}s^{-2}u^{-1}t^{-1}s^{-1}(stu)^4t^{-1}\\
				
				&\subset&U+u_1U+Rs^{-1}u^{-1}t^{-1}(stus)t(ustu)+Rs^{-1}t(ustu)+
				Rt^{-1}s^{-1}(tust)(ustu)t^{-1}\\
				
				&\subset&U+u_1U+Rt^3ust+u_1t^2ust+R(ustu)s\\
				&\subset&U+u_1U+u_1\underline{u_2u_3u_1u_2}+Rstus^2\\
				&\subset&U+u_1U+u_1\underline{u_2u_3u_1}\subset U+u_1U.
				\end{array}}$
			
			The result follows from proposition \ref{su12}.
				\item [(ii)]\underline{$k=1$}:
				
				$\small{\begin{array}[t]{lcl}
					zuu_1u_3u_1&\subset& z(R+Ru^{-1})u_1u_3u_1\\
					&\subset&\underline{zu_1u_3u_1}+z
					u^{-1}(R+Rs^{-1})(R+Ru^{-1})(R+Rs^{-1})\\
					&\subset&U+\underline{u_3u_1}+Rzu^{-1}s^{-1}u^{-1}+Rzu^{-1}s^{-1}u^{-1}s^{-1}\\
					&\subset&U+Ru^{-1}s^{-1}zu^{-1}+Ru^{-1}s^{-1}zu^{-1}s^{-1}\\
					
					&\subset&U+Ru^{-1}s^{-1}(stu)^4u^{-1}+Ru^{-1}s^{-1}(stu)^4u^{-1}s^{-1}\\
					
					&\subset& U+Ru^{-1}(tust)u(stus)t+Ru^{-1}(tust)us(tust)s^{-1}\\
					
					&\subset&U+Rstu^3stut+R(stus)stu\\
					
					&\subset&U+u_1t(R+Ru)stut+Rust(ustu)\\

					&\subset&U+u_1t(stu)^4(stu)^{-3}t+u_1(tust)ut+Rust^2ust\\

				&\subset&U+zu_1t(u^{-1}t^{-1}s^{-1}u^{-1})t^{-1}s^{-1}u^{-1}t^{-1}s^{-1}t+u_1ustu^2+rus(R+Rt)ust\\
				
				&\subset&U+zu_1u^{-1}t^{-2}s^{-1}u^{-1}t^{-1}s^{-1}t+u_1ust(R+Ru)+
			(	uu_1u_3u_1)u_2+R(ustu)st\\ 
				
				&\stackrel{(i)}{\subset}&U+zu_1u^{-1}(R+Rt^{-1})s^{-1}u^{-1}t^{-1}s^{-1}t+u_1\underline{u_3u_1u_2}+
				u_1(ustu)+Uu_2+\\&&+Rstus^2t\\
				
				&\stackrel{\ref{r12}}{\subset}&U+u_1U+zu_1u^{-1}(s^{-1}u^{-1}t^{-1}s^{-1})t+u_1(stu)^{-2}zt+
			u_1tust+u_1\underline{u_2u_3u_1u_2}\\
			
			&\subset& U+u_1U+zu_1(u^{-1}t^{-1}s^{-1}u^{-1})+u_1(stu)^{-2}(stu)^4t+
			u_1\underline{u_2u_3u_1u_2}\\
			
			&\subset&U+u_1U+zu_1t^{-1}s^{-1}u^{-1}t^{-1}+u_1t(ustu)t\\
			
			&\subset&U+u_1U+u_1\underline{u_2u_1u_3u_2}+u_1t^2ust^2\\
			
			&\subset&U+u_1U+u_1\underline{u_2u_3u_1u_2}\\
				&\subset&U+u_1U.
					\end{array}}$
				
				The result follows again from proposition \ref{su12}.
			\end{itemize}
			\item [C2.] We will prove that $z^kuu_1u_2u_3\subset U$,  $k\in\{0,1\}$. We notice that  $z^kuu_1u_2u_3=z^ku(R+Rs)(R+Rt)(R+Ru)\subset \underline{z^ku_3u_1u_2}+z^ku_3u_2u_3+z^kuu_1u_3u_1+Rz^k(ustu)\stackrel{C1}{\subset} U+z^ku_3u_2u_2+Rz^ktust\subset U+z^ku_3u_2u_2+\underline{z^ku_2u_3u_1u_2}\subset U+z^ku_3u_2u_2.$
			Therefore, we must prove that, for every $k\in\{0,1\}$, $z^ku_3u_2u_2\subset U$. We distinguish the following cases:
			\begin{itemize}[leftmargin=*]
				\item [(i)] $\underline{k=0}:$
				$\small{\begin{array}[t]{lcl}
				u_3u_2u_3&=&(R+Ru)(R+Rt)(R+Ru)\\
				&\subset&\underline{u_2u_3u_2}+Rutu\\
				&\subset&U+Rus^{-1}stu\\
				&\subset&U+Ru(R+Rs)stu\\
				&\subset&U+R(ustu)+Rus(stu)^4(stu)^{-3}\\
				&\subset&U+Rtust+Rzus(u^{-1}t^{-1}s^{-1}u^{-1})(t^{-1}s^{-1}u^{-1}t^{-1})s^{-1}\\
				&\subset&U+\underline{u_2u_3u_1u_2}+Rzt^{-1}s^{-2}u^{-1}t^{-1}s^{-2}\\
				&\subset&U+Rzt^{-1}s^{-2}u^{-1}t^{-1}(R+Rs^{-1})\\
				&\subset&U+\underline{zu_2u_1u_3u_2}+Rt^{-1}s^{-1}u^{-1}t^{-1}s^{-1}z\\
				&\subset&U+Rt^{-1}s^{-2}u^{-1}t^{-1}s^{-1}(stu)^4\\
				&\subset&U+Rt^{-1}s^{-1}(tust)(ustu)\\
				&\subset&U+R(ustu)st\\
				&\subset&U+Rstus^2t\\
				&\subset&U+u_1\underline{u_2u_3u_1u_2}\\
					&\subset&U+u_1U.
				\end{array}}$
			
			The result follows from proposition \ref{su12}.
				\item[(ii)]$\underline{k=1}:$
				$\small{\begin{array}[t]{lcl}
				zu_3u_2u_3&=&z(R+Ru^{-1})(R+Rt^{-1})(R+Ru^{-1})\\
				&\subset&\underline{zu_2u_3u_2}+Rzu^{-1}t^{-1}u^{-1}\\
				&\subset&U+R(stu)^4u^{-1}t^{-1}u^{-1}\\
				&\subset&U+Rs(tust)u(stus)u^{-1}\\
				&\subset&U+Rs^2tusu^2st\\
				&\subset&U+u_1tus(R+Ru)st\\
				&\subset&U+u_1\underline{u_2u_3u_1u_2}+u_1(stu)^4(stu)^{-3}sust\\
				&\subset&U+u_1U+zu_1u^{-1}t^{-1}s^{-1}u^{-1}(t^{-1}s^{-1}u^{-1}t^{-1})ust\\
				&\subset&U+u_1U+zu_1u^{-1}t^{-1}s^{-1}u^{-2}\\
				&\subset&U+u_1U+zu_1u^{-1}t^{-1}s^{-1}(R+Ru^{-1})\\
				&\subset&U+u_1U+u_1\underline{zu_3u_2u_1}+zu_1(u^{-1}t^{-1}s^{-1}u^{-1})\\
				&\subset&U+u_1U+zu_1t^{-1}s^{-1}u^{-1}t^{-1}\\
				&\subset&U+u_1U+u_1\underline{zu_2u_1u_3u_2}\\
			&\subset&	U+u_1U.
				\end{array}}$
			
			he result follows again from proposition \ref{su12}.
			\end{itemize}
			\item[C3.] For every $k\in\{0,1\}$ we have:
			$z^kuu_2u_3u_1=z^ku(R+Rt)(R+Ru)(R+Rs)\subset \underline{z^ku_3u_2u_1}+Rz^kutu+Rz^kutus\subset U+z^kuu_1u_2u_3+Rz^ku(tust)t^{-1}$.
			Therefore, by C2 we have 	$z^kuu_2u_3u_1\subset U+Rz^k(ustu)st^{-1}\subset U+Rz^kstus^2t^{-1}\subset U+u_1\underline{z^ku_2u_3u_1u_2}\subset U+u_1U$. The result follows from \ref{su12}.
			\item[C4.]For every $k\in\{0,1\}$ we have:
			$z^kuu_2u_1u_3=z^ku(R+Rt^{-1})(R+Rs^{-1})(R+Ru^{-1})\subset U+\underline{z^ku_3u_2u_1}+z^kuu_1u_2u_3+z^kuu_1u_3u_1+Rz^kut^{-1}s^{-1}u^{-1}$. Therefore, by C1 and C2 we have  that $z^kuu_2u_1u_3\subset U+Rz^ku(t^{-1}s^{-1}u^{-1}t^{-1})t\subset U+
			Rz^kt^{-1}s^{-1}u^{-1}t\subset U+\underline{z^ku_2u_1u_3u_2}.$
			\qedhere
		\end{itemize}
	\end{proof}
		\begin{cor}
			The BMR freeness conjecture holds for the generic Hecke algebra $\bar H_{G_{12}}$.
		\end{cor}
		\begin{proof}
			By theorem \ref{thm12} we have that $H_{G_{12}}=U=\sum\limits_{k=0}^1z^k(u_2+su_2+uu_2+suu_2+usu_2+tuu_2+tsu_2+susu_2+stuu_2+tusu_2+tsuu_2+utsu_2)$ and, hence, $H_{G_{12}}$ is generated as right $u_2$-module by 24 elements and, hence, as $R$-module by $|G_{12}|=48$ elements (recall that $u_2$ is generated as $R$-module by 2 element). Therefore, $\bar H_{G_{12}}$ is generated as $\bar R$-module by $|G_{12}|=48$ elements, since the action of $\bar R$ factors through $R$. The result follows from proposition \ref{BMR PROP}.
\end{proof}

\subsection{The case of $G_{13}$}
\indent

Let $R=\ZZ[u_{s,i}^{\pm},u_{t,j}^{\pm},u_{u_{u,l}}]_{\substack{1\leq i,j,l\leq 2 }}$ and let  $$H_{G_{13}}=\langle s,t,u\;|\; ustu=tust,stust=ustus, \prod\limits_{i=1}^{2}(s-u_{s,i})=\prod\limits_{j=1}^{2}(t-u_{t,l})=\prod\limits_{l=1}^{2}(u-u_{u_k})=0\rangle$$ be the generic Hecke algebra associated to $G_{13}$. Let $u_1$ be the subalgebra of $H_{G_{13}}$ generated by $s$, $u_2$ the subalgebra of $H_{G_{13}}$ generated by $t$ and $u_3$ be the subalgebra of $H_{G_{13}}$ generated by $u$. We recall that $z:=(stu)^3=(tus)^3=(ust)^3$ generates the center of the associated complex braid group and that $|Z(G_{13})|=4$. We set
 $U=\sum\limits_{k=0}^{3}z^k(u_1u_2u_1u_2+u_1u_2u_3u_2+u_2u_3u_1u_2+u_2u_1u_3u_2+u_3u_2u_1u_2)$. 
By the definition of $U$, we have the following remark:

\begin{rem}
	$Uu_2 \subset U$.
	\label{rem13}
\end{rem}
From now on, we will underline the elements that by definition belong to $U$. Moreover, we will use directly the remark \ref{rem13}; this means that every time we have a power of  $t$ at the end of an element, we may ignore it. In order to remind that to the reader, we put a parenthesis around the part of the element we consider.

Our goal is to prove that $H_{G_{13}}=U$ (theorem \ref{thm13}). Since $1\in U$, it is enough to prove that $U$ is a left-sided ideal of $H_{G_{13}}$. For this purpose, one may check that $sU$, $tU$ and $uU$ are subsets of $U$. We set $$\small{\begin{array}{lcl}
v_1&=&1\\
v_2&=&u\\
v_3&=&s\\
v_4&=&ts\\
v_5&=&su\\
v_6&=&us\\
v_7&=&tu\\
v_8&=&tsu\\
v_9&=&tus\\
v_{10}&=&sts\\
v_{11}&=&stu\\
v_{12}&=&uts
\end{array}}$$
By the definition of $U$, one may notice that $U=\sum\limits_{k=0}^3\sum\limits_{i=1}^{12}(Rz^kv_i+Rz^kv_it)$. Hence, by remark \ref{rem13} we only have to prove that for every $k\in\{0,\dots,3\}$, $sz^kv_i$, $tz^kv_i$ and $uz^kv_i$ are elements inside $U$, $i=1,\dots,12$. 
 As a first step, we prove this argument for $i=8,\dots, 12$ and for a smaller range of the values of $k$, as we can see in the following proposition.

\begin{prop}
	\mbox{}
	\vspace*{-\parsep}
	\vspace*{-\baselineskip}\\
	\begin{itemize}[leftmargin=0.8cm]
	\item[(i)] For every $k\in\{1,2,3\}$, $z^ksv_{8}\in U.$	
		\item[(ii)] For every $k\in\{0,1,2\}$, $z^ksv_{9}\in U.$	
			\item[(iii)] For every $k\in\{0,1,2\}$, $z^kuv_{10}\in U.$
		
		\item[(iv)]For every $k\in\{0,1,2,3\}$, $z^ktv_{10}\in U.$
		
		\item[(v)]For every $k\in\{0,1,2,3\}$, $z^ktv_{11}\in U.$

	\item[(vi)] For every $k\in\{1,2,3\}$, $z^ksv_{12}\in U.$
	\item[(vii)] For every $k\in\{0,1,2,3\}$, $z^ktv_{12}\in U.$

	\end{itemize}
	\label{easy}
	\end{prop}
	\begin{proof}
		
	\mbox{}
	\vspace*{-\parsep}
	\vspace*{-\baselineskip}\\
	\begin{itemize}	[leftmargin=0.8cm]
		
	\item[(i)]
	$\small{\begin{array}[t]{lcl}z^ksv_{8}&=&z^kstsu\\
	&=&z^kz^k(R+Rs^{-1})(R+Rt^{-1})(R+Rs^{-1})(R+Ru^{-1})\\
	&\in&
	\underline{z^ku_2u_1u_3}+\underline{z^ku_1u_2u_3}+\underline{z^ku_1u_2u_1}+
	Rz^ks^{-1}t^{-1}s^{-1}u^{-1}\\
	&\in& U+Rz^ks^{-1}(ust)^{-3}ustust\\
	&\in& U+Rz^{k-1}s^{-1}(ustus)t\\
	&\in& U+Rz^{k-1}tust^2\\
	&\in& U+\underline{(z^{k-1}u_2u_3u_1)t^2}.
	\end{array}}$
	
	\item[(ii)]
	$z^ksv_{9}=z^kstus=z^k(stu)^3(u^{-1}t^{-1}s^{-1}u^{-1})t^{-1}=z^{k+1}t^{-1}s^{-1}u^{-1}t^{-2}\in \underline{(z^{k+1}u_2u_1u_3)t^{-2}}$.
	\item[(iii)] 
	$z^kuv_{10}=z^kusts=z^k(ust)^3t^{-1}(s^{-1}u^{-1}t^{-1}s^{-1}u^{-1})s=
	z^{k+1}t^{-2}s^{-1}u^{-1}t^{-1}\in \underline{(z^{k+1}u_2u_1u_3)t}.$
	
	\item[(iv)]
	\begin{itemize}[leftmargin=*]
		\item[$\bullet$]\underline{$k\in\{0,1\}$}:
	$\small{\begin{array}[t]{lcl}
	z^ktv_{10}&=&z^ktsts\\
	&=&z^kt(stu)^3(u^{-1}t^{-1}s^{-1}u^{-1})t^{-1}s^{-1}u^{-1}s\\
	&=&z^{k+1}s^{-1}u^{-1}t^{-2}s^{-1}u^{-1}s\\
	&\in&z^{k+1}s^{-1}u^{-1}(R+Rt^{-1})s^{-1}u^{-1}s\\
	&\in&Rz^{k+1}s^{-1}u^{-1}s^{-1}(R+Ru)s+Rz^{k+1}(s^{-1}u^{-1}t^{-1}s^{-1}u^{-1})s\\
	&\in&\underline{z^{k+1}u_1u_3}+Rz^{k+1}s^{-1}u^{-1}(R+Rs)us+Rz^{k+1}t^{-1}s^{-1}u^{-1}t^{-1}\\
	&\in&U+\underline{z^{k+1}u_2}+Rz^{k+1}s^{-1}(R+Ru)sus+\underline{(z^{k+1}u_2u_1u_3)t}\\
	&\in&U+\underline{z^{k+1}u_3u_1}+Rz^{k+1}s^{-1}us(ust)^3(t^{-1}s^{-1}u^{-1}t^{-1}s^{-1})u^{-1}t^{-1}\\
	&\in&U+Rz^{k+2}s^{-1}t^{-1}s^{-1}u^{-2}t^{-1}\\
	&\in&U+Rz^{k+2}s^{-1}t^{-1}s^{-1}(R+Ru^{-1})t^{-1}\\
	&\in&U+\underline{(z^{k+2}u_1u_2u_1)t^{-1}}+Rz^{k+1}s^{-1}t^{-1}s^{-1}u^{-1}(ust)^3t^{-1}\\
	&\in&U+Rz^{k+1}s^{-1}(ustus)\\
	&\in& U+Rz^{k+1}tust\\
	&\in& U+\underline{(z^{k+2}u_2u_3u_1)t}.
	\end{array}}$	
	 \item[$\bullet$] \underline{$k\in\{2,3\}$}:
	 $\small{\begin{array}[t]{lcl}
	 z^ktv_{10}&=&z^ktsts\\
	 &\in&z^k(R+Rt^{-1})(R+Rs^{-1})(R+Rt^{-1})(R+Rs^{-1})\\
	 &\in& \underline{(z^ku_1u_2u_1)u_2}+Rz^kt^{-1}s^{-1}t^{-1}s^{-1}\\
	 &\in& U+Rz^{k-1}t^{-1}s^{-1}t^{-1}s^{-1}(stu)^3\\
	 &\in&U+Rz^{k-1}t^{-1}s^{-1}(ustus)tu\\
	 &\in&U+Rz^{k-1}ust^2u\\
	 &\in&U+Rz^{k-1}us(R+Rt)u\\
	  &\in&U+Rz^{k-1}usu+Rz^{k-1}(ustu)\\
	  &\in&U+Rz^{k-1}(R+Ru^{-1})(R+Rs^{-1})(R+Ru^{-1})+Rz^{k-1}tust\\
	  &\in&U+\underline{z^{k-1}u_3u_1}+\underline{z^{k-1}u_1u_3}+Rz^{k-1}u^{-1}s^{-1}u^{-1}+\underline{(z^{k-1}u_2u_3u_1)t}\\
	  	&\in&U+Rz^{k-2}u^{-1}s^{-1}u^{-1}(ust)^3\\
	  	&\in&U+Rz^{k-2}u^{-1}(tust)ust\\
	  	&\in&U+Rz^{k-2}stu^2st\\
	  	&\in&U+Rz^{k-2}st(R+Ru)st\\
	  	&\in& U+\underline{(z^{k-2}u_1u_2u_1)t}+(Rz^{k-2}sv_9)t
	  	\stackrel{(ii)}{\subset}U.
	 \end{array}}$
\end{itemize}

\item[(v)] 
\begin{itemize}[leftmargin=*]
	
\item[$\bullet$]\underline{$k\in\{0,1\}$}: 
$\small{\begin{array}[t]{lcl}
z^ktv_{11}&=&z^ktstu\\
&=&z^kt(stu)^3(u^{-1}t^{-1}s^{-1}u^{-1})t^{-1}s^{-1}\\
&=&z^{k+1}s^{-1}u^{-1}t^{-2}s^{-1}\\
&\in&z^{k+1}s^{-1}u^{-1}(R+Rt^{-1})s^{-1}\\
 &\in& Rz^{k+1}s^{-1}u^{-1}s^{-1}+Rz^{k+1}s^{-1}u^{-1}t^{-1}s^{-1}\\
 &\in& Rz^{k+1}(R+Rs)(R+Ru)(R+Rs)+Rz^ks^{-1}u^{-1}t^{-1}s^{-1}(stu)^3\\
 &\in& \underline{z^{k+1}u_1u_3}+\underline{z^{k+1}u_3u_1}+Rz^{k+1}sus+Rz^kt(ustu)\\
 &\in& U+Rz^{k+1}s(ust)^3(t^{-1}s^{-1}u^{-1}t^{-1}s^{-1})u^{-1}t^{-1}+Rz^kt^2ust\\
 &\in&U+Rz^{k+2}u^{-1}t^{-1}s^{-1}u^{-2}+\underline{(z^ku_2u_3u_1)t}\\
 &\in&U+Rz^{k+2}u^{-1}t^{-1}s^{-1}(R+Ru^{-1})\\
 &\in&U+\underline{z^{k+2}u_3u_2u_1}+Rz^{k+2}(u^{-1}t^{-1}s^{-1}u^{-1})\\
 &\in& U+Rz^{k+2}t^{-1}s^{-1}u^{-1}t^{-1}\\
 &\in& U+\underline{(z^{k+2}u_2u_1u_3)t^{-1}}.
\end{array}}$

\item[$\bullet$]\underline{$k\in\{2,3\}$}:
 $\small{\begin{array}[t]{lcl}
z^ktv_{11}&=&z^ktstu\\
&\in& z^k(R+Rt^{-1})(R+Rs^{-1})(R+Rt^{-1})(R+Ru^{-1})\\
&\in&\underline{z^ku_1u_2u_3}+\underline{z^ku_2u_1u_3}+\underline{(z^ku_2u_1)u_2}+
Rz^kt^{-1}s^{-1}t^{-1}u^{-1}\\
&\in& U+Rz^{k-1}t^{-1}s^{-1}(stu)^3t^{-1}u^{-1}\\
&\in&U+Rz^{k-1}ust(ustu)t^{-1}u^{-1}\\
&\in&U+Rz^{k-1}ust^2usu^{-1}\\
&\in&U+Rz^{k-1}us(Rt+R)usu^{-1}\\
&\in&U+Rz^{k-1}(ust)^3(t^{-1}s^{-1}u^{-1}t^{-1})u^{-1}+Rz^{k-1}usu(R+Rs^{-1})u^{-1}\\
&\in&U+Rz^ku^{-1}t^{-1}s^{-1}u^{-2}+\underline{z^{k-1}u_3u_1}+Rz^{k-1}us(R+Ru^{-1})s^{-1}u^{-1}\\
&\in&U+Rz^ku^{-1}t^{-1}s^{-1}(R+Ru^{-1})+\underline{z^{k-1}u_2}+Rz^{k-1}u(R+Rs^{-1})u^{-1}s^{-1}u^{-1}\\
&\in& U+\underline{z^ku_3u_2u_1}+Rz^k(u^{-1}t^{-1}s^{-1}u^{-1})+\underline{z^{k-1}u_1u_3}+\\&&+
Rz^{k-2}us^{-1}(stu)^3u^{-1}s^{-1}u^{-1}\\
&\in& U+Rz^kt^{-1}s^{-1}u^{-1}t^{-1}+Rz^{k-2}utu(stust)s^{-1}u^{-1}\\
&\in&U+\underline{(z^ku_2u_1u_3)t^{-1}}+Rz^{k-2}utu^2st\\
&\in&U+Rz^{k-2}ut(R+Ru)st\\
&\in& U+\underline{(z^{k-2}u_3u_2u_1)t}+Rz^{k-2}u(tust)\\
&\in&U+Rz^{k-2}u^2stu\\
&\in&U+Rz^{k-2}(R+Ru)stu\\
&\in& U+\underline{z^{k-2}u_1u_2u_3}+Rz^{k-2}(ustu)\\
&\in& U+Rz^{k-2}tust\\
&\in& U+\underline{(z^{k-2}u_2u_3u_1)t}.
\end{array}}$
\end{itemize}

		\item[(vi)] 
		$\small{\begin{array}[t]{lcl}z^ksv_{12}&=&z^ksuts\\
			&\in& z^k(R+Rs^{-1})(R+Ru^{-1})(R+Rt^{-1})(R+Rs^{-1})\\
		&\in&
		\underline{z^ku_3u_2u_1}+\underline{z^ku_1u_2u_1}+\underline{(z^ku_2u_1u_3)u_2}+
		Rz^ks^{-1}u^{-1}s^{-1}+Rz^ks^{-1}u^{-1}t^{-1}s^{-1}\\
		&\in& U+ Rz^{k-1}s^{-1}(stu)^3u^{-1}s^{-1}+Rz^{k-1}s^{-1}u^{-1}t^{-1}s^{-1}(stu)^3\\
		
		&\in& U+Rz^{k-1}tu(stust)s^{-1}+Rz^{k-1}tustu\\
		&\in& U+Rz^{k-1}tu^2stu+Rz^{k-1}tustu\\
		&\in&U+z^{k-1}tu_3stu\\
		&\in& U+z^{k-1}t(R+Ru)stu\\
		&\in& U+ Rz^{k-1}tv_{11}+Rz^{k-1}t(ustu)\\
		&\stackrel{(v)}{\in}&U+Rz^{k-1}t^2ust\\&\in& U+\underline{(z^{k-1}u_2u_3u_1)t}.
		\end{array}}$
		\item[(vii)] 
		\begin{itemize}[leftmargin=*]
			
		\item[$\bullet$]\underline{$k\in\{0,1\}$}:
		$\small{\begin{array}[t]{lcl}
		z^ktv_{12}&=&z^ktuts\\
		&=&z^ks^{-1}(stu)^3(u^{-1}t^{-1}s^{-1}u^{-1})t^{-1}s^{-1}ts\\
		&=&z^{k+1}s^{-1}t^{-1}s^{-1}u^{-1}t^{-2}s^{-1}ts\\
		&\in&z^{k+1}s^{-1}t^{-1}s^{-1}u^{-1}(R+Rt^{-1})s^{-1}ts\\
		&\in&Rz^{k+1}s^{-1}t^{-1}s^{-1}u^{-1}(R+Rs)ts+Rz^{k+1}s^{-1}(ust)^{-3}(ustu)ts\\
		&\in&Rz^{k+1}s^{-1}(ust)^{-3}(ustus)ts+Rz^{k+1}s^{-1}t^{-1}s^{-1}(R+Ru)sts+Rz^ks^{-1}tust^2s\\
		&\in&Rz^ktust^2s+\underline{z^{k+1}u_2}+Rz^{k+1}s^{-1}t^{-1}s^{-1}usts+Rz^k(R+Rs)tus(R+Rt)s\\
		&\in&U+Rz^ktus(R+Rt)s+Rz^{k+1}s^{-1}t^{-1}s^{-1}(ust)^3t^{-1}(s^{-1}u^{-1}t^{-1}s^{-1}u^{-1})s+\\&&+
		\underline{z^ku_2u_3u_1}+Rz^k(tust)s+Rz^kstus^2+Rz^k(stust)s\\
		&\in& U+
		\underline{z^ku_2u_3u_1}+Rz^k(tust)s+Rz^{k+2}s^{-1}t^{-1}s^{-1}t^{-2}s^{-1}u^{-1}t^{-1}+Rz^k(ustus)+\\&&+Rz^kstu(R+Rs)+Rz^kustus^2\\
		&\in&U+Rz^k(ustus)+Rz^{k+2}s^{-1}t^{-1}s^{-1}(R+Rt^{-1})s^{-1}u^{-1}t^{-1}+
		(Rz^ksv_9)t+\\&&+
		\underline{z^ku_1u_2u_3}+Rz^ksv_9+Rz^kustu(R+Rs)\\
		&\stackrel{(ii)}{\in}&U+(Rz^ksv_9)t+Rz^{k+2}s^{-1}t^{-1}s^{-2}u^{-1}t^{-1}+
		Rz^{k+2}s^{-1}t^{-1}s^{-1}t^{-1}s^{-1}u^{-1}t^{-1}+\\&&+
		Rz^k(ustu)+Rz^k(ustus).
		\end{array}}$
	
	However, by (ii) we have that $z^ksv_9\in U$. Moreover, $z^k(ustu)=z^ktust\in\underline{(z^ku_2u_3u_1)t}$ and $z^k(ustus)=z^ksv_9t\stackrel{(ii)}{\in}U$. Therefore, it remains to prove that $D:=Rz^{k+2}s^{-1}t^{-1}s^{-2}u^{-1}t^{-1}+
	Rz^{k+2}s^{-1}t^{-1}s^{-1}t^{-1}s^{-1}u^{-1}t^{-1}$ is a subset of $U$. We have:
	\\ \\
	$\small{\begin{array}{lcl}
		D&=&Rz^{k+2}s^{-1}t^{-1}s^{-2}u^{-1}t^{-1}+
		Rz^{k+2}s^{-1}t^{-1}s^{-1}t^{-1}s^{-1}u^{-1}t^{-1}\\
	&\subset&Rz^{k+2}s^{-1}t^{-1}(R+Rs^{-1})u^{-1}t^{-1}+
	Rz^{k+2}s^{-1}t^{-1}s^{-1}(ust)^{-3}(ustus)\\
	&\subset&U+\underline{(z^{k+2}u_1u_2u_3)t^{-1}}+Rz^{k+2}s^{-1}(ust)^{-3}(ustus)+Rz^{k+1}s^{-1}ust+
	\underline{(z^ku_2u_3u_1)t}\\
	&\in&U+Rz^{k+1}tust+Rz^{k+1}s^{-1}(R+Ru^{-1})(R+Rs^{-1})t\\
	&\in&U+\underline{(z^{k+1}u_2u_3u_1)t}+\underline{(z^{k+1}u_1u_3)t}+
	Rz^{k+1}s^{-1}u^{-1}s^{-1}t\\
	&\in&U+Rz^{k}s^{-1}(stu)^3u^{-1}s^{-1}t\\
	&\in&U+Rz^ktu(stust)s^{-1}t\\
	&\in&U+Rz^ktu^2stut\\
	&\in&U+Rz^kt(R+Ru)stut\\
	&\in&U+(Rz^ktv_{11})t+Rz^{k}(tust)ut\\
	&\stackrel{(v)}{\in}&U+Rz^kustu^2t\\
	&\in&U+Rz^kust(R+Ru)t\\
	&\in& U+\underline{(z^ku_3u_1)t^2}+Rz^k(ustu)t\\&\in& U+\underline{(z^ku_2u_3u_1)t}.
	\end{array}}$
		
\			\item[$\bullet$]\underline{$k\in\{2,3\}$}: 
			$\small{\begin{array}[t]{lcl}
				z^ktv_{12}&=&z^ktuts\\
				&\in& z^k(R+Rt^{-1})(R+Ru^{-1})(R+Rt^{-1})(R+Rs^{-1})\\
				&\in&
				\underline{z^ku_3u_2u_1}+\underline{(z^ku_2u_3u_1)u_2}+Rz^kt^{-1}u^{-1}t^{-1}s^{-1}\\
				&\in& U+Rz^{k-1}t^{-1}u^{-1}t^{-1}s^{-1}(stu)^3\\
				&\in&U+Rz^{k-1}t^{-1}s(tust)u\\
				&\in& U+Rz^{k-1}t^{-1}sustu^2\\
				& \in& U+Rz^{k-1}t^{-1}sust(Ru+R)\\
				&\in& U+Rz^{k-1}t^{-1}s(ustu)+Rz^{k-1}t^{-1}(R+Rs^{-1})ust\\
				&\in&U+Rz^{k-1}t^{-1}stust+\underline{(z^{k-1}u_2u_3u_1)t}+Rz^{k-1}t^{-1}s^{-1}(R+Ru^{-1})st\\
				&\in&U+Rz^{k-1}t^{-1}(stu)^3(u^{-1}t^{-1}s^{-1}u^{-1})+\underline{z^{k-1}u_2}+\\&&+
				Rz^{k-1}t^{-1}s^{-1}u^{-1}(R+Rs^{-1})t\\
				&\in&U+Rz^kt^{-2}s^{-1}u^{-1}t^{-1}+\underline{(z^{k-1}u_2u_1u_3)t}+
				Rz^{k-1}(ust)^{-3}u(stust)s^{-1}t\\
				&\in&U+\underline{(z^ku_2u_1u_3)t^{-1}}+Rz^{k-2}u^2stut\\
				&\in&U+Rz^{k-2}(R+Ru)stut\\
				&\in&U+\underline{(z^{k-2}u_1u_2u_3)t}+Rz^{k-2}(ustu)\\
				&\in& U+Rz^{k-2}tust\\
				&\in& U+\underline{(z^{k-2}u_2u_3u_1)t}.
				
				\end{array}}$

\end{itemize}	
		
\end{itemize}
		\qedhere

		\end{proof}

\begin{cor} $u_2U\subset U$.
	\label{cor13}
	\end{cor}	
	
	\begin{proof}
		Since $u_2=R+Rt$, it is enough to prove that $tU\subset U$. However, by the definition of $U$ and remark \ref{rem13}(i), this is the same as proving that for every $k\in\{0,\dots,3\}$, $z^ktv_i\in U$, which follows directly from the definition of $U$ and proposition \ref{easy}(v), (vi), (vii).
	\end{proof}
By remark \ref{rem13} (ii), we have that for every $k\in\{0,\dots,3\}$, $z^ku_iu_ju_k\subset U$, for some combinations of $i,j,l \in\{1,2,3\}$. We can generalize this argument for every $i,j,l \in\{1,2,3\}$.
\begin{prop}
	For every $k\in\{0,\dots,3\}$ and for every combination of $i,j,l\in\{1,2,3\}$, $z^ku_iu_ju_l\subset U$.
	\label{prr13}
\end{prop}
\begin{proof}
By the definition of $U$  we only have to prove that, for every $k\in\{0,\dots,3\}$, $z^ku_1u_3u_1$, $z^ku_3u_1u_3$ and $z^ku_3u_2u_3$ are subsets of $U$. We distinguish the following 3 cases:
	\begin{itemize}[leftmargin=0.8cm]
		\item[C1.] 
		\begin{itemize}[leftmargin=*]
			\item [$\bullet$]\underline{$k\in\{0,1,2\}$}:
			$\small{\begin{array}[t]{lcl}z^ku_1u_3u_1&=&z^k(R+Rs)(R+Ru)(R+Rs)\\
			 &\subset& z^kv_5+z^kv_6+\underline{z^ku_1}+Rz^ksus\\
			 &\subset& U+Rz^ks(ust)^3(t^{-1}s^{-1}u^{-1}t^{-1}s^{-1})u^{-1}\\
			 &\subset& U+Rz^{k+1}u^{-1}t^{-1}s^{-1}u^{-2}\\
			 &\subset& U+z^{k+1}u^{-1}t^{-1}s^{-1}(R+Ru^{-1})\\
			 &\subset& U+\underline{z^{k+1}u_3u_2u_1}+Rz^{k+1}(u^{-1}t^{-1}s^{-1}u^{-1})\\
			 &\subset& U+z^{k+1}t^{-1}s^{-1}u^{-1}t^{-1}\\
			 &\subset& U+\underline{(z^{k+1}u_2u_1u_3)t^{-1}}.
			 \end{array}}$
			 \item[$\bullet$]\underline{$k=3$}:
			 $\small{\begin{array}[t]{lcl}z^3u_1u_3u_1&=&z^k(R+Rs^{-1})(R+Ru{-1})s\\
			 &\subset& \underline{z^3u_3u_1}+\underline{z^3u_1u_3}+\underline{z^3u_1}+Rz^3s^{-1}u^{-1}s\\
			 &\subset& U+Rz^{2}s^{-1}(stu)^3u^{-1}s\\
			 &\subset& U+Rz^{2}tu(stust)s\\
			 &\subset& U+Rz^{2}tu^2stu\\
			 &\subset& U+Rz^{2}t(R+Ru)stu\\
			 &\subset& U+
			 t(\underline{Rz^{2}u_1u_2u_3})+Rz^{2}(ustu)\\
			 &\subset& U+u_2U+Rz^{2}v_9t.
			 \end{array}}$
			
			The result follows from corollary \ref{cor13} and the definition of $U$.
		\end{itemize}
		\item [C2.] 
		\begin{itemize}[leftmargin=*]
			\item[$\bullet$]\underline{$k\in\{1,2,3\}$}:
				$\small{\begin{array}[t]{lcl}z^ku_3u_1u_3&=&z^k(R+Ru^{-1})(R+Rs^{-1})(R+Ru^{-1})\\&\subset& \underline{z^ku_1u_3}+\underline{z^ku_3u_1}+Rz^ku^{-1}s^{-1}u^{-1}\\
			&\subset& U+Rz^{k-1}u^{-1}s^{-1}u^{-1}(ust)^3\\&\subset& U+Rz^{k-1}u^{-1}(tust)ust\\
			&\subset& U+Rz^{k-1}stu^2st\\
			&\subset& U+Rz^{k-1}st(R+Ru)st\\
			&\subset& U+\underline{(z^{k-1}u_1u_2u_1)t}+(Rz^{k-1}sv_9)t.\end{array}}$.
		
	The result follows from proposition \ref{easy}(ii).
			\item[$\bullet$] \underline{$k=0$}:
				$\small{\begin{array}[t]{lcl}
			 	u_3u_1u_3&=&(R+Ru)(R+Rs)(R+Ru)\\
			 	&\subset& \underline{u_1u_3}+\underline{u_3u_1}+Rusu\\
			 	&\subset& U+Rustt^{-1}u\\
			 	&\subset& U+Rust(R+Rt)u\\
			 	&\subset& U+R(ustu)+R(ust)^3t^{-1}s^{-1}(u^{-1}t^{-1}s^{-1}u^{-1})tu\\
			 &\subset&U+Rtust+Rzt^{-1}s^{-1}t^{-1}s^{-1}\\
			 &\subset& U+\underline{(u_2u_3u_1)t}+u_2(\underline{zu_1u_2u_1})\\
			 &\subset& U+u_2U.\end{array}}$
			
			The result follows from corollary \ref{cor13}.
			 
		\end{itemize}
		\item[C3.] 
		\begin{itemize}[leftmargin=*]
			\item[$\bullet$] \underline{$k\in\{0,1\}$}:
			$\small{\begin{array}[t]{lcl}
			z^ku_3u_1u_3&=&z^k(R+Ru)(R+Rt)(R+Ru)\\
			&\subset& \underline{z^ku_2u_3}+\underline{z^ku_3u_2}+Rz^kutu\\
			&\subset&U+Rz^kus^{-1}stu\\
			&\subset&U+Rz^ku(R+Rs)stu\\
			&\subset&U+Rz^k(ustu)+Rz^kus(stu)^3(u^{-1}t^{-1}s^{-1}u^{-1})t^{-1}s^{-1}\\
			&\subset&U+Rz^ktust+Rz^{k+1}ust^{-1}s^{-1}u^{-1}t^{-2}s^{-1}\\
			&\subset&U+\underline{(z^ku_2u_3u_1)t}+Rz^{k+1}ust^{-1}s^{-1}u^{-1}(R+Rt^{-1})s^{-1}\\
			&\subset&U+Rz^{k+1}us(R+Rt)s^{-1}u^{-1}s^{-1}+Rz^{k+1}us(t^{-1}s^{-1}u^{-1}t^{-1}s^{-1})\\
			&\subset&U+\underline{z^{k+1}u_1}+Rz^{k+1}usts^{-1}(R+Ru)s^{-1}+Rz^{k+1}t^{-1}s^{-1}u^{-1}\\
			&\subset&U+z^{k+1}ustu_1+Rz^{k+1}ust(R+Rs)us^{-1}+\underline{z^{k+1}u_2u_1u_3}\\
			&\subset&U+z^{k+1}ustu_1+Rz^{k+1}(ustu)s^{-1}+\\&&+
			Rz^{k+1}(ust)^3t^{-1}(s^{-1}u^{-1}t^{-1}s^{-1}u^{-1})sus^{-1}\\
				&\subset&U+z^{k+1}u_2ustu_1+Rz^{k+2}t^{-1}(t^{-1}s^{-1}u^{-1}t^{-1})us^{-1}\\
				&\subset&U+z^{k+1}u_2ust(R+Rs)+Rz^{k+2}t^{-1}u^{-1}t^{-1}s^{-2}\\
				&\subset&U+u_2\underline{(z^{k+1}u_3u_1)t}+u_2z^{k+1}uv_{10}+
				u_2\underline{z^{k+2}u_3u_2u_1}\\
&\stackrel{\ref{easy}(iv)}{\subset}&U+u_2U.\end{array}}$

The result follows from proposition \ref{cor13}(i).

			\item[$\bullet$]  \underline{$k\in\{2,3\}$}:
			$\small{\begin{array}[t]{lcl}
				z^ku_3u_1u_3&=&z^k(R+Ru^{-1})(R+Rt^{-1})(R+Ru^{-1})\\
				&\subset& \underline{z^ku_2u_3}+\underline{z^ku_3u_2}+Rz^ku^{-1}t^{-1}u^{-1}\\
				&\subset&U+Rz^ku^{-1}t^{-1}s^{-1}su^{-1}\\
				&\subset&U+Rz^ku^{-1}t^{-1}s^{-1}(R+Rs^{-1})u^{-1}\\
				&\subset&U+Rz^k(u^{-1}t^{-1}s^{-1}u^{-1})+Rz^k(stu)^{-3}st(ustu)s^{-1}u^{-1}\\
				&\subset&U+Rz^kt^{-1}s^{-1}u^{-1}t^{-1}+Rz^{k-1}st^2usts^{-1}u^{-1}\\
				&\subset&U+\underline{(z^{k}u_2u_1u_3)t^{-1}}+Rz^{k-1}s(R+Rt)usts^{-1}u^{-1}\\
				&\subset&U+Rz^{k-1}sus(R+Rt^{-1})s^{-1}u^{-1}+Rz^{k-1}(stust)s^{-1}u^{-1}\\
				&\subset&U+\underline{z^{k-1}u_1}+Rz^{k-1}s(R+Ru^{-1})st^{-1}s^{-1}u^{-1}+Rz^{k-1}ust\\
				&\subset&U+z^{k-1}u_1t^{-1}s^{-1}u^{-1}+Rz^{k-1}su^{-1}(R+Rs^{-1})t^{-1}s^{-1}u^{-1}+\\&&+
				\underline{(z^{k-1}u_3u_1)t}\\
				&\subset&U+z^{k-1}u_1t^{-1}s^{-1}u^{-1}+Rz^{k-1}s(u^{-1}t^{-1}s^{-1}u^{-1})+\\&&+Rz^{k-1}su^{-1}s^{-1}(ust)^{-3}(ustus)t\\
					&\subset&U+z^{k-1}u_1t^{-1}s^{-1}u^{-1}u_2+Rz^{k-2}su^{-1}(tust)t\\
					&\subset&U+z^{k-1}(R+Rs)t^{-1}s^{-1}u^{-1}u_2+Rz^{k-2}s^2tut\\
					&\subset&U+\underline{z^{k-1}u_2u_1u_3}+Rz^{k-1}s(R+Rt)(R+Rs)(R+Ru)u_2+\\&&+\underline{(z^{k-2}u_1u_2u_3)t}\\
					&\subset&U+\underline{(z^{k-1}u_3)u_2}+\underline{(z^{k-1}u_1u_2u_3)u_2}+\underline{(z^{k-1}u_1u_2u_1)u_2}+
					+Rz^{k-1}stsuu_2\\
					&\subset&U+(Rz^{k-1}sv_8)u_2.
			\end{array}}$
	
	The result follows from proposition \ref{easy}(i).
	\qedhere
	\end{itemize}
\end{itemize}
\end{proof}

\begin{prop} If $u_3U\subset U$, then $H_{G_{13}}=U$.
	\label{prop13}
\end{prop}
\begin{proof}
	As we explained in the beginning of this section, in order to prove that $H_{G_{13}}=U$, it will be sufficient to prove that   $sU$, $tU$ and $uU$ are subsets of $U$. By corollary \ref{cor13} and hypothesis, we only have to prove that $sU\subset U$. By the definition of $U$ and remark \ref{rem13}, we have to prove that for every $k\in\{0,\dots,3\}$, $z^ksv_i\in U$, $k=1,\dots,12$. However, for every $i\in\{1,\dots,7\}\cup\{10,11\}$, we have that $z^ksv_i\in z^ku_1u_ju_l$, where $j,l\in\{1,2,3\}$ and not necessarily distinct. Therefore, by proposition \ref{prr13} we only have to prove that $z^ksv_i\in U$, for $i=8,9,12$.
For this purpose, for every $k\in\{0,\dots,3\}$ we have to check the following cases:
	\begin{itemize}[leftmargin=*]
		\item For the element $z^ksv_8$ we only have to prove that $z^0sv_8\in U$, since by proposition \ref{easy}(i) we have the rest of the cases. We have:
		$sv_8=stsu=sts(ust)^3(t^{-1}s^{-1}u^{-1}t^{-1}s^{-1})u^{-1}t^{-1}s^{-1}=
		zst(u^{-1}t^{-1}s^{-1})u^{-1})u^{-1}t^{-1}s^{-1}
	=zu^{-1}t^{-1}u^{-1}t^{-1}s^{-1}
		\in u_3u_2\underline{zu_3u_2u_1}$.
		Therefore, 	$sv_8
		\in u_3u_2U$.
	The result follows from corollary \ref{cor13} and from hypothesis.
		\item For the element $z^ksv_9$ we use proposition \ref{easy}(ii) and we only have to prove that $z^3sv_9\in U$.
			$z^3sv_9=z^3stus\in z^3(R+Rs^{-1})(R+Rt^{-1})(R+Ru^{-1})(R+Rs^{-1})\subset \underline{z^3u_2u_3u_1}+\underline{z^3u_1u_2u_1}+z^ku_1u_2u_1+Rz^3s^{-1}t^{-1}u^{-1}s^{-1}\stackrel{\ref{prr13}}{\subset}U+Rz^2s^{-1}(stu)^3t^{-1}u^{-1}s^{-1}\subset U+Rz^2tu(stustut^{-1}u^{-1}s^{-1}).$
			By hypothesis and by corollary \ref{cor13}, it is enough to prove that 	$z^2stustut^{-1}u^{-1}s^{-1}\in U$. 
			
	$\small{\begin{array}{lcl}
	z^2st(ustu)t^{-1}u^{-1}s^{-1}&=&z^2st^2usu^{-1}s^{-1}\\&\in& z^2s(R+Rt)usu^{-1}s^{-1}\\
	&\in&Rz^2su(R+Rs^{-1})u^{-1}s^{-1}+
	Rz^2(stu)^3u^{-1}(t^{-1}s^{-1}u^{-1}t^{-1})u^{-1}s^{-1}\\
	&\in& U+\underline{z^2u_2}+Rz^2s(R+Ru^{-1})s^{-1}u^{-1}s^{-1}+Rz^3u^{-2}t^{-1}s^{-1}u^{-2}s^{-1}\\
	&\in&U+\underline{z^2u_3u_1}+Rz^2(R+Rs^{-1})u^{-1}s^{-1}u^{-1}s^{-1}+u_3u_2(z^3u_1u_3u_1)\\
	&\stackrel{\ref{prr13}}{\in}&U+u_3(z^2u_1u_3u_1)+Rzs^{-1}u^{-1}s^{-1}(stu)^3u^{-1}s^{-1}+u_3u_2U\\
	&\stackrel{\ref{prr13}}{\in}&U+u_3u_2U+Rzs^{-1}u^{-1}(tust)usts^{-1}\\
	&\in&U+u_3u_2U+Rztu^2sts^{-1}\\
	&\in&U+u_3u_2U+u_2u_3(\underline{zu_1u_2u_1})\\
	&\subset& U+u_3u_2u_3U.
	\end{array}}$
		
	The result follows from hypothesis and from corollary \ref{cor13}.
		
	\item For the element $z^ksv_{12}$  we only have to prove that $z^0sv_{12}\in U$, since the rest of the cases have been proven in proposition \ref{easy}(vi).
	
	$\small{\begin{array}[t]{lcl}
		sv_{12}&=&suts\\
		&=&(stu)^3u^{-1}t^{-1}s^{-1}u^{-1}(t^{-1}s^{-1}u^{-1}t^{-1})uts\\
		&=&z u^{-1}t^{-1}s^{-1}u^{-2}t^{-1}s^{-1}ts\\
		&\in&zu^{-1}t^{-1}s^{-1}(R+Ru^{-1})t^{-1}s^{-1}ts\\
		&\in&Rzu^{-1}t^{-1}s^{-1}t^{-1}s^{-1}ts+Rz(stu)^{-3}u^{-1}(ustu)ts\\
		&\in&Rzu^{-1}t^{-1}s^{-1}t^{-1}(R+Rs)ts+Ru^{-1}tust^2s\\
		&\in&\underline{(zu_3)u_2}+Rzu^{-1}t^{-1}s^{-1}(R+Rt)sts+u_3u_2u_3\underline{(u_1u_2u_1)}\\
		&\in&U+\underline{zu_3u_1}+Rzu^{-1}t^{-1}(R+Rs)tsts+u_3u_2u_3U\\
		&\in&U+u_2u_3U+u_3\underline{(zu_1u_2u_1)}+zu_3u_2(stu)^3u^{-1}t^{-1}(s^{-1}u^{-1}t^{-1}s^{-1}u^{-1})sts\\
			&\in&U+u_3u_2u_3U+z^2u_3u_2u^{-1}t^{-2}s^{-1}u^{-1}s\\
			&\in&U+u_3u_2u_3U+z^2u_3u_2u_3u_2\underline{(z^2u_1u_3u_1)}\\
			&\subset& U+u_3u_2u_3u_2U.
		\end{array}}$
		
			The result follows again from hypothesis and from corollary \ref{cor13}.
			\qedhere
	\end{itemize}
\end{proof}
We can now prove the main theorem of this section.
\begin{thm} $H_{G_{13}}=U$.
	\label{thm13}
	\end{thm}
\begin{proof}
	By proposition \ref{prop13} it will be sufficient to prove that $u_3U\subset U$. By the definition of $U$ and remark \ref{rem13}(i) we have to prove that for every $k\in\{0,\dots,3\}$, $z^ku_3v_i\subset U$, $k=1,\dots,12$. However, for every $i\in\{1,\dots,7\}\cup\{12\}$, $z^ku_3v_i\subset z^ku_3u_ju_l$, where $j,l\in\{1,2,3\}$ and not necessarily distinct. Therefore, by proposition \ref{prr13} we restrict ourselves to proving that $z^ku_3v_i\subset U$, for $i=8,9,10,11$.
	For every $k\in\{0,\dots,3\}$ we have:
	\begin{itemize}[leftmargin=*]
		
	\item$\small{\begin{array}[t]{lcl}
		z^ku_3v_8&\subset&z^k(R+Ru^{-1})tsu\\
		&\subset& \underline{z^ku_2u_1u_3}+Rz^ku^{-1}(R+Rt^{-1})(R+Rs^{-1})(R+Ru^{-1})\\
		&\subset& U+z^ku_3u_1u_3+z^ku_3u_2u_3+\underline{z^ku_3u_2u_1}+Rz^k(u^{-1}t^{-1}s^{-1}u^{-1})\\
	&\stackrel{\ref{prr13}}{\subset}&U+Rz^kt^{-1}s^{-1}u^{-1}t^{-1}\\&\subset& U+\underline{(z^ku_2u_1u_3)t^{-1}}.
	\end{array}}$
	
	\item $\small{\begin{array}[t]{lcl}
		z^ku_3v_9&\subset& z^k(R+Ru)tus\\
		&\subset& \underline{z^ku_2u_3u_1}+Rz^ku(tust)t^{-1}\\
		&\subset& U+z^ku^2stut^{-1}\\
		&\subset& U+z^k(R+Ru)stut^{-1}\\
		&\subset& U+\underline{(z^ku_1u_2u_3)t^{-1}}+Rz^k(ustu)t^{-1}\\
		&\subset& U+Rz^ktus\\
		&\subset& U+\underline{z^ku_2u_3u_1}.\end{array}}$
	
	\item $z^ku_3v_{10}\subset z^k(R+Ru)v_{10}\subset \underline{Rz^kv_{10}}+Rz^kuv_{10}.$ Therefore, by proposition \ref{easy}(iii), we only have to prove that $z^3uv_{10}\in U$. However, $z^3v_{10}=
	z^3t^{-1}(tust)s=z^3t^{-1}(ustus)=z^3t^{-1}stust\in u_2(z^3stus)u_2$. Hence, by remark \ref{rem13} and corollary \ref{cor13}, we need to prove that $z^3stus\in U$.\\ \\
	$\small{\begin{array}{lcl}
	z^3stus&=&z^3(R+Rs^{-1})(R+Rt^{-1})(R+Ru^{-1})(R+Rs^{-1})\\
	&\in&\underline{z^3u_1u_2u_3}+\underline{z^3u_1u_2u_1}+\underline{z^3u_2u_3u_1}+
	z^3u_2u_3u_1+Rs^{-1}t^{-1}u^{-1}s^{-1}\\
	&\stackrel{\ref{prr13}}{\in}& U+Rz^2s^{-1}(stu)^3t^{-1}u^{-1}s^{-1}\\
	&\in&U+Rz^2tust(ustu)t^{-1}u^{-1}s^{-1}\\
	&\in&U+Rz^2tust^2usu^{-1}s^{-1}\\
	&\in&U+Rz^2tus(Rt+R)usu^{-1}s^{-1}\\
	&\in&U+Rz^2t(ust)^3(t^{-1}s^{-1}u^{-1}t^{-1})u^{-1}s^{-1}+Rz^2tusu(R+Rs^{-1})u^{-1}s^{-1}\\
	&\in&U+Rz^3tu^{-1}t^{-1}s^{-1}u^{-2}s^{-1}+\underline{z^2u_2u_3}+Rz^2tus(R+Ru^{-1})s^{-1}u^{-1}s^{-1}\\
	&\in&U+Rz^3tu^{-1}t^{-1}s^{-1}(R+Ru^{-1})s^{-1}+\underline{z^2u_2u_1}+
	Rz^2tu(R+Rs^{-1})u^{-1}s^{-1}u^{-1}s^{-1}\\
	&\in&U+u_2\underline{z^3u_3u_2u_1}+Rz^2tu^{-1}t^{-1}s^{-1}u^{-1}(ust)^3s^{-1}+
	u_2(z^2u_1u_3u_1)+\\&&+Rztus^{-1}u^{-1}(ust)^3s^{-1}u^{-1}s^{-1}\\
	&\stackrel{\ref{prr13}}{\in}& U+u_2U+Rz^3t(stust)s^{-1}+Rztutu(stust)s^{-1}u^{-1}s^{-1}\\

	&\stackrel{\ref{cor13}}{\in}& U+Rz^3t(ustu)+Rztutu^2sts^{-1}\\
	&\in&U+Rz^3t^2ust+Rztut(R+Ru)sts^{-1}\\

	&\in&U+u_2\underline{(z^3u_3u_1)t}+Rztu(R+Rt^{-1})sts^{-1}+Rztu(tust)s^{-1}\\
	&\in&U+Rztust(R+Rs)+Rztut^{-1}(R+Rs^{-1})ts^{-1}+
	Rztu^2stus^{-1}\\
	&\in&U+u_2U+\underline{(zu_2u_3u_1)t}+zu_2uv_{10}+\underline{zu_2u_3u_1}+Rztut^{-1}s^{-1}(R+Rt^{-1})s^{-1}+\\&&+
	Rzt(R+Ru)stus^{-1}\\
	&\stackrel{\ref{easy}(iii)}{\in}& U+u_2U+u_2\underline{zu_3u_2u_1}+Rzt(R+Ru^{-1})t^{-1}s^{-1}t^{-1}s^{-1}+
	Rztstu(R+Rs)+\\&&+
	Rztustu(R+Rs)\\
	&\stackrel{\ref{cor13}}{\in}& U+\underline{zu_1u_2u_1}+Rzt(stu)^{-3}st(ustu)t^{-1}s^{-1}+Rztv_{11}+u_2zsv_9+
	
	Rzt(ustu)+\\&&+Rzt(ustu)s\\
	&\stackrel{\ref{easy}}{\in}&U+Rst^2u+Rzt^2ust+Rzt^2usts\\
	&\in&U+\underline{u_1u_2}+\underline{(zu_2u_3u_1)t}+zu_2uv_{10}\\
	&\stackrel{\ref{easy}(iii)}{\subset}&U+u_2U.
	\end{array}}$

The result follows from corollary \ref{cor13}.
	
	\item $z^ku_3v_{11}\subset z^k(R+Ru)stu\subset \underline{z^ku_1u_2u_3}+Rz^k(ustu)
	\subset U+Rz^ktust\subset U+\underline{(z^ku_2u_3u_1)t}\subset U.$
	\qedhere
\end{itemize}
\end{proof}

\begin{cor}
	The BMR freeness conjecture holds for the generic Hecke algebra $H_{G_{13}}$.
\end{cor}
\begin{proof}
	By theorem \ref{thm13} we have that $H_{G_{13}}=U=\sum\limits_{k=0}^3\sum\limits_{i=1}^{12}(Rz^kv_i+Rz^kv_it)$. The result follows from proposition \ref{BMR PROP}, since by definition $H_{G_{13}}$ is generated as $R$-module by $|G_{13}|=96$ elements.
\end{proof}

\subsection{The case of $G_{14}$}

Let $R=\ZZ[u_{s,i}^{\pm},u_{t,j}^{\pm}]_{\substack{1\leq i\leq 2 \\1\leq j\leq 4}}$ and let $H_{G_{14}}=\langle s,t\;|\; stststst=tstststs,\prod\limits_{i=1}^{2}(s-u_{s,i})=\prod\limits_{j=1}^{3}(t-u_{t,j})=0\rangle$ be the generic Hecke algebra associated to $G_{14}$. Let $u_1$ be the subalgebra of $H_{G_{14}}$ generated by $s$ and $u_2$ the subalgebra of $H_{G_{14}}$ generated by $t$. We recall that that $z:=(st)^4=(ts)^4$  generates the center of the associated complex braid group and that $|Z(G_{14})|=6$.
We set $U=\sum\limits_{k=0}^{5}(z^ku_1u_2u_1u_2+z^ku_1tst^{-1}su_2).$
By the definition of $U$, we have the following remark.

\begin{rem}
$Uu_2 \subset U$.
	\label{rem14}
\end{rem}
 From now on, we will underline the elements that by definition belong to $U$.  Moreover, we will use directly remark \ref{rem14} and the fact that $u_1U\subset U$; this means that every time we have a power of $s$ in the beginning of an element or a power of $t$ at the end of it, we may ignore it. In order to remind that to the reader, we put a parenthesis around the part of the element we consider.

Our goal is to prove that $H_{G_{14}}=U$ (theorem \ref{thm14}). Since $1\in U$, it is enough to prove that $U$ is a left-sided ideal of $H_{G_{14}}$. For this purpose, one may check that $sU$ and $tU$ are subsets of $U$. However, by the definition of $U$ and remark \ref{rem14}, we only have to prove that for every $k\in\{0,\dots,5\}$,
$z^ktu_1u_2u_1$ and $z^ktu_1tst^{-1}s$ are subsets of $U$. In the following proposition we first prove this statement for a smaller range of the values of $k$. 

\begin{prop}
	\mbox{}
	\vspace*{-\parsep}
	\vspace*{-\baselineskip}\\
	\begin{itemize}[leftmargin=0.8cm]
	\item[(i)]For every $k\in\{0,\dots,4\}$, $z^ktu_1u_2u_1\subset U$.
	\item[(ii)]For every $k\in\{0,\dots,4\}$, $z^ku_2u_1tu_1\subset U+z^{k+1}tu_1u_2u_1u_2$. Therefore, for every $k\in\{0,\dots,3\}$, $z^ku_2u_1tu_1 \subset U$.
	\item[(iii)]For every $k\in\{1,\dots,5\}$, $z^ku_2u_1t^{-1}u_1\subset U$.
	\item[(iv)]For every $k\in\{1,\dots,4\}$, $z^ku_2u_1u_2u_1\subset U+z^{k+1}tu_1u_2u_1u_2$. Therefore, for evey $k\in\{1,\dots,3\}$, $z^ku_2u_1u_2u_1\subset U$.
	
\end{itemize}
	\label{pr14}
\end{prop}

\begin{proof}
\mbox{}
\vspace*{-\parsep}
\vspace*{-\baselineskip}\\
\begin{itemize}	[leftmargin=0.8cm]
\item [(i)]	
$z^ktu_1u_2u_1=z^kt(R+Rs)(R+Rt+Rt^{-1})(R+Rs)\subset \underline{z^ku_2u_1u_2}+\underline{Rz^ktst^{-1}s}+Rz^ktsts\subset U+Rz^k(ts)^4s^{-1}t^{-1}s^{-1}t^{-1}\subset U+\underline{z^{k+1}u_1u_2u_1u_2}\subset U.$
	
\item[(ii)] We notice that $z^ku_2u_1tu_1=z^k(R+Rt+Rt^2)u_1tu_1\subset
\underline{z^ku_1u_2u_1}+z^ktu_1u_2u_1+
z^kt^2u_1tu_1$. However, $z^ktu_1u_2u_1\subset U$, by (i). Therefore, we have to prove that $z^kt^2u_1tu_1$ is a subset of $U$. Indeed, $z^kt^2u_1tu_1=
z^kt^2(R+Rs)t(R+Rs)\subset U+\underline{z^ku_2u_1u_2}+Rz^kt^2sts\subset
Rz^kt(ts)^4s^{-1}t^{-1}s^{-1}t^{-1}\subset U+
(z^{k+1}tu_1u_2u_1)u_2$. The result follows from (i).

\item[(iii)]We expand $u_2$ as $R+Rt+Rt^{-1}$ and we have that
$z^ku_2u_1t^{-1}u_1\subset \underline{z^ku_1u_2u_1}+z^ktu_1t^{-1}u_1+
	z^kt^{-1}u_1t^{-1}u_1\subset U+z^kt(R+Rs)t^{-1}(R+Rs)+
	z^kt^{-1}(R+Rs^{-1})t^{-1}(R+Rs^{-1})\subset
 U+\underline{z^ku_2u_1u_2}+
	\underline{Rz^ktst^{-1}s}+
	Rz^kt^{-1}s^{-1}t^{-1}s^{-1}\subset U+
	Rz^k(st)^{-4}stst\subset U+\underline{z^{k-1}u_1u_2u_1u_2}.$
	
\item[(iv)] The result follows from the definition of U and from (ii) and (iii), since $z^ku_2u_1u_2u_1=z^ku_2u_1(R+Rt+Rt^{-1})u_1$.
\qedhere
\end{itemize}	
\end{proof}
To make it easier for the reader to follow the calculations, from now on we will double-underline the elements as described in the above proposition (proposition \ref{pr14}) and we will use directly the fact that these elements are inside $U$.
We can now prove the following lemmas that lead us to the main theorem of this section.
\begin{lem}
 For every $k\in\{3,4\}$, $z^ktsu_2sts\subset
U+z^{k+1}u_1tu_1u_2u_1u_2$.
		
	\label{lm14}
\end{lem}
\begin{proof}We have:
	\\
	$\small{\begin{array}[t]{lcl}
	z^ktsu_2sts&=&z^kts(R+Rt+Rt^{-1})sts\\
	
	&\subset& \underline{\underline{z^ktu_1u_2u_1}}+Rz^k(ts)^4s^{-1}t^{-1}+z^ktst^{-1}sts\\
	&\subset&U+\underline{z^{k+1}u_1u_2}+z^kt(R+Rs^{-1})t^{-1}(R+Rs^{-1})t(R+Rs^{-1})\\
	&\subset&U+\underline{z^ku_1u_2u_1}+\underline{\underline{(z^ktu_1u_2u_1)t}}+
	Rz^kts^{-1}t^{-1}s^{-1}ts^{-1}\\
	&\subset&U+Rz^kts^{-1}t^{-1}s^{-1}(R+Rt^{-1}+Rt^{-2})s^{-1}\\
	&\subset&U+\underline{\underline{z^ktu_1u_2u_1}}+Rz^kt^2(st)^{-4}st+
	Rz^kt^2(st)^{-4                                                 }stst^{-1}s^{-1}\\
	&\subset&U+\underline{z^{k-1}u_2u_1u_2}+Rz^{k-1}t^2st(R+Rs^{-1})t^{-1}s^{-1}\\
	&\subset&U+\underline{z^{k-1}u_2}+Rz^{k-1}t^2st^2(st)^{-4}stst\\
	&\subset&U+Rz^{k-2}(R+Rt+Rt^{-1})st^2stst\\
	&\subset&U+\underline{\underline{z^{k-2}s(u_2u_1tu_1)t}}+Rz^{k-2}tst(ts)^4s^{-1}t^{-1}s^{-1}+
	Rz^{k-2}t^{-1}s(R+Rt+Rt^{-1})stst\\
	&\subset&U+Rz^{k-1}tst(R+Rs)t^{-1}s^{-1}+\underline{\underline{(z^{k-2}u_2u_1tu_1)t}}+Rz^{k-2}t^{-2}(ts)^4s^{-1}+\\&&+Rz^{k-2}t^{-1}(R+Rs^{-1})t^{-1}(R+Rs^{-1})tst\\
	&\subset&U+\underline{(z^{k-1}+z^{k-2})u_2u_1u_2}+Rz^{k-1}(ts)^4s^{-1}t^{-1}s^{-1}t^{-2}s^{-1}+
	\underline{\underline{(z^{k-2}u_2u_1tu_1)t}}+\\&&+Rz^{k-2}t^{-1}s^{-1}t^{-1}s^{-1}tst\\
	&\subset&U+z^ku_1u_2u_1u_2u_1+
	Rz^{k-2}(st)^{-4}stst^2st\\
	&\stackrel{\ref{pr14}(iv)}{\subset}&U+z^{k+1}u_1tu_1u_2u_1u_2+\underline{\underline{z^{k-3}s(tu_1u_2u_1)t}}.
	\end{array}}$
		
\end{proof}
\begin{lem}
	$z^kt^{-2}s^{-1}t^{-2}s^{-1}\in U$.
	\label{lem114}
\end{lem}
\begin{proof}We have:
	\\
	$\small{\begin{array}[t]{lcl}
	z^5t^{-2}s^{-1}t^{-2}s^{-1}&=&z^5t^{-1}(st)^{-4}ststst^{-1}s^{-1}\\
	&\in&z^4t^{-1}stst(R+Rs^{-1})t^{-1}s^{-1}\\
	&\in&\underline{z^4u_2u_1u_2}+Rz^4t^{-1}sts(R+Rt^{-1}+Rt^{-2})s^{-1}t^{-1}s^{-1}\\
	&\in&U+\underline{z^4u_2}+Rz^4t^{-1}st(R+Rs^{-1})t^{-1}s^{-1}t^{-1}s^{-1}+Rz^4t^{-1}stst^{-1}(st)^{-4}stst\\
	&\in&U+\underline{z^4u_2u_1}+Rz^4t^{-1}st^2(st)^{-4}st+
	Rz^3t^{-1}st(R+Rs^{-1})t^{-1}(R+Rs^{-1})tst\\
	&\in&U+\underline{\underline{(z^3u_2u_1u_2u_1)u_2}}+
	Rz^3t^{-1}sts^{-1}t^{-1}s^{-1}tst\\
	&\in&U+
		Rz^3t^{-1}(R+Rs^{-1})ts^{-1}t^{-1}s^{-1}tst\\

	&\in&U+\underline{\underline{u_1(z^3u_2u_1u_2u_1)u_2}}+Rz^3t^{-1}s^{-1}ts^{-1}t^{-1}s^{-1}tst\\
	
	&\in&U+
	Rz^3t^{-1}s^{-1}t^2(st)^{-4}stst^2st\\
	&\in&U+Rz^2t^{-1}s^{-1}(R+Rt+Rt^{-1})stst^2st\\
	&\in&U+\underline{z^2u_1u_2u_1u_2}+Rz^2t^{-1}s^{-2}(st)^4t^{-1}s^{-1}tst+
	z^2t^{-1}s^{-1}t^{-1}(R+Rs^{-1})tst^2st\\
	&\in&U+Rz^3t^{-1}(R+Rs^{-1})t^{-1}s^{-1}tst+\underline{z^2u_2u_1u_2}+
	Rz^2(st)^{-4}stst^2st^2st\\
		&\in&U+\underline{\underline{(z^3u_2u_1tu_1)t}}+Rz^3(st)^{-4}stst^2st+
		Rzstst^2s(R+Rt+Rt^{-1})st\\
			&\in&U+\underline{\underline{s\big((z+z^2)tu_1u_2u_1\big)t}}+Rzstst(ts)^4s^{-1}t^{-1}s^{-1}t^{-1}+
			Rzsts(R+Rt+Rt^{-1})st^{-1}st\\
	
		&\in&U+Rz^2stst(R+Rs)t^{-1}s^{-1}t^{-1}+\underline{\underline{s(ztu_1u_2u_1)t}}+Rz(st)^4t^{-1}s^{-1}t^{-2}st+\\&&+Rzstst^{-1}st^{-1}st\\
	&\in&U+\underline{z^2u_1}+Rz^2(st)^4t^{-1}s^{-1}t^{-2}s^{-1}t^{-1}+
		\underline{\underline{(z^2u_2u_1u_2u_1)t}}+\\&&+
		Rzst(R+Rs^{-1})t^{-1}(R+Rs^{-1})t^{-1}st\\
		&\in&U+	\underline{\underline{(z^3u_2u_1u_2u_1)t^{-1}}}+\underline{zu_1u_2u_1u_2}+
		\underline{\underline{s(z^3tu_1u_2u_1)t}}+Rzsts^{-1}t^{-1}s^{-1}t^{-1}st\\
	
	\end{array}}$
\\
However, $zsts^{-1}t^{-1}s^{-1}t^{-1}st\in	zsts^{-1}t^{-1}s^{-1}t^{-1}(R+Rs^{-1})t\subset
\underline{\underline{s(ztu_1u_2u_1)}}+Rzst^2(st)^{-4}st^2\subset
U+\underline{u_1u_2u_1u_2}$.
\end{proof}
\begin{lem}
	For every $k\in\{0,\dots,5\}$, $z^ktu_1u_2u_1\subset U$.
	\label{cor14}
\end{lem}
\begin{proof}
	By proposition \ref{pr14} (i), we only have to prove that $z^5tu_1u_2u_1\subset U$. We have:
	$z^5tu_1u_2u_1=z^5t(R+Rs^{-1})(R+Rt^{-1}+Rt^{-2})u_1\subset \underline{z^5u_2u_1}+
	\underline{\underline{z^5u_2u_1t^{-1}u_1}}+z^5ts^{-1}t^{-2}u_1\subset U+z^5ts^{-1}t^{-2}(R+Rs^{-1})\subset U+\underline{z^5u_2u_1u_2}+Rz^5ts^{-1}t^{-2}s^{-1}.$
	We expand $t$ as a linear combination of 1, $t^{-1}$ and $t^{-2}$ and by the definition of $U$ and lemma \ref{lem114}, we only have to prove that $z^5t^{-1}s^{-1}t^{-2}s^{-1}\in U$. Indeed, we have:
	$z^5t^{-1}s^{-1}t^{-2}s^{-1}=z^5(st)^{-4}
	ststst^{-1}s^{-1}\in z^4stst(R+Rs^{-1})t^{-1}s^{-1}
	\subset \underline{z^4u_1u_2}+z^4stst^2(st)^{-4}stst\subset U+s(z^3tsu_2sts)t$.
	We use lemma \ref{lm14} and we have that $z^3tsu_2sts\subset U+z^4u_1tu_1u_2u_1u_2$. The result follows from proposition \ref{pr14}(i). 
	
\end{proof}
\begin{thm}
	$H_{G_{14}}=U$.
	\label{thm14}
\end{thm}
\begin{proof}
	As we explained in the beginning of this section, it is enough to prove that, for every $k\in\{0,\dots,5\}$, $z^ktu_1u_2u_1$ and $z^ktu_1tst^{-1}s$ are subsets of $U$. The first part is exactly what we proved in lemma \ref{cor14}. It remains to prove the second one. 
Since $z^ktu_1tst^{-1}s=z^kt(R+Rs)tst^{-1}s$, we must prove that for every $k\in\{0,\dots,5\}$, the elements $z^kt^2st^{-1}s$ and $z^ktstst^{-1}s$ are inside $U$. 
We distinguish the following cases:
\begin{itemize}[leftmargin=*]
	\item \underline{The element $z^kt^2st^{-1}s$}:
	
	By proposition \ref{pr14} (iii), we only have to prove the case where $k=0$. We have:
	$$\small{\begin{array}{lcl}
	t^2st^{-1}s&\in&t^2s(R+Rt+Rt^2)s\\
	&\in& \underline{u_2u_1}+\underline{\underline{u_2u_1tu_1}}+t(ts)^4s^{-1}t^{-1}s^{-1}t^{-1}s^{-1}ts\\
	&\in&U+zts^{-1}t^{-1}s^{-1}t^{-1}(R+Rs)ts\\
	&\in&U+\underline{zu_2u_1u_2}+Rzts^{-1}t^{-1}s^{-1}(R+Rt+Rt^2)sts\\
	&\in&U+\underline{zu_2}+Rzts^{-1}t^{-1}(R+Rs)tsts+Rzts^{-1}t^{-1}(R+Rs)t^2sts\\
	&\in& U+\underline{zu_2u_1}+Rzts^{-1}t^{-2}(ts)^4s^{-1}t^{-1}+
	Rzts^{-2}(st)^4t^{-1}s^{-1}t^{-1}+\\&&+
	Rzt(R+Rs)t^{-1}st^2sts\\
	&\in&U+\underline{\underline{(z^2tu_1u_2u_1)t^{-1}}}+\underline{\underline{s(zu_2u_1tu_1)}}+Rzts(R+Rt+Rt^2)st^2sts\\
	&\in&U+Rzts^2t^2sts+Rz(ts)^4s^{-1}t^{-1}s^{-1}tsts+Rztst^2st(ts)^4s^{-1}t^{-1}s^{-1}t^{-1}\\
	&\in& U+Rzt(R+Rs)t^2sts+Rz^2s^{-1}t^{-1}(R+Rs)tsts+Rz^2tst^2st(R+Rs)t^{-1}s^{-1}t^{-1}\\
	&\in& U+\underline{\underline{zu_2u_1tu_1}}+Rztst(ts)^4s^{-1}t^{-1}s^{-1}t^{-1}+
	\underline{z^2u_2u_1}+
	Rz^2s^{-1}t^{-2}(ts)^4s^{-1}t^{-1}+\\&&+\underline{z^2u_2u_1u_2}+
	Rz^2tst(ts)^4s^{-1}t^{-1}s^{-1}t^{-2}s^{-1}\\
	&\in&U+Rz^2tst(R+Rs)t^{-1}s^{-1}t^{-1}+\underline{z^3u_1u_2u_1u_2}+
	z^3tsts^{-1}t^{-1}s^{-1}(R+Rt^{-1}+Rt)s^{-1}\\
	&\in& U+\underline{z^2u_1}+Rz^2(ts)^4s^{-1}t^{-1}s^{-1}t^{-2}s^{-1}t^{-1}+
	Rz^3tsts^{-1}t^{-1}s^{-2}+
	
	Rz^3tst^2(st)^{-4}st+\\&&+
	Rz^3tst(R+Rs)t^{-1}s^{-1}ts^{-1}\\

	&\in&U+\underline{\underline{s^{-1}(z^3u_2u_1u_2u_1)t^{-1}}}+Rz^3tst(R+Rs)t^{-1}s^{-2}+\underline{\underline{(z^2tu_1u_2u_1)t}}+
	\underline{\underline{z^3u_2u_1}}+\\&&+Rz^3tstst^{-1}(R+Rs)ts^{-1}\\
	&\in&U+\underline{z^3u_2u_1}+Rz^3(ts)^4s^{-1}t^{-1}s^{-1}t^{-2}s^{-2}+
	\underline{z^3u_2u_1u_2}+Rz^3(ts)^4s^{-1}t^{-1}s^{-1}t^{-2}sts^{-1}\\
	&\in&U+s^{-1}(z^4u_2u_1u_2u_1)+Rz^4s^{-1}t^{-1}s^{-1}(R+Rt^{-1}+Rt)sts^{-1}\\
	&\stackrel{\ref{pr14}(iv)}{\in}&U+s^{-1}z^5(tu_1u_2u_1)u_2+\underline{z^4u_1}+
	Rz^4(ts)^{-4}tsts^2ts^{-1}+
	Rz^4s^{-1}t^{-1}(R+Rs)tsts^{-1}\\
	&\stackrel{\ref{cor14}}{\in}&U+Rz^3tst(R+Rs)ts^{-1}+\underline{z^4u_2u_1}+
	Rz^4s^{-1}t^{-2}(ts)^4s^{-1}t^{-1}s^{-2}\\
	&\in&U+\underline{\underline{z^3tu_1u_2u_1}}+Rz^3(ts)^4s^{-1}t^{-1}s^{-2}+
	s^{-1}(z^4u_2u_1u_2u_1)\\
		&\stackrel{\ref{cor14}}{\in}&U+\underline{z^4u_1u_2u_1}+s^{-1}(z^5tu_1u_2u_1)u_2\\
		&\stackrel{\ref{cor14}}{\in}&U.

		\end{array}}$$

	\item \underline{The element $z^ktstst^{-1}s$}:
	
	\underline{For $k\in\{4,5\}$}:
	$\begin{array}[t]{lcl}
	z^ktstst^{-1}s&\in& z^ktst(R+Rs^{-1})t^{-1}s\\
	&\subset&
	\underline{z^ku_2u_1}+z^ktsts^{-1}t^{-1}(R+Rs^{-1})\\
	&\subset& U+(z^ktsts^{-1})t^{-1}+
	z^ktst^2(st)^{-4}stst\\
		&\stackrel{\ref{cor14}}{\subset}&U+(z^{k-1}tst^2sts)t\\
		&\stackrel{\ref{lm14}}{\subset}&
	U+u_1(z^ktu_1u_2u_1)u_2\\&\stackrel{\ref{cor14}}{\subset}&U.
	\end{array}$
	
\underline{	For $k\in\{0,\dots,3\}$}: $z^ktstst^{-1}s=z^k(ts)^4s^{-1}t^{-1}s^{-1}t^{-2}s\in s^{-1}(z^{k+1}u_2u_1u_2u_1)$. However, by 
\ref{pr14}(iv) we have that $z^{k+1}u_2u_1u_2u_1\subset U+(z^{k+2}tu_1u_2
u_1)u_2\stackrel{\ref{cor14}}{\subset}U.$
\qedhere	
\end{itemize}

\end{proof}

\begin{cor}
	The BMR freeness conjecture holds for the generic Hecke algebra $H_{G_{14}}$.
\end{cor}
\begin{proof}
	By theorem \ref{thm14} we have that $H_{G_{14}}=U=\sum\limits_{k=0}^{5}(z^ku_1u_2+z^ku_1tsu_2+z^ku_1t^{-1}su_2+z^ku_1tst^{-1}su_2)$. The result follows from proposition \ref{BMR PROP}, since $H_{G_{14}}$ is generated as a left $u_1$-module by 72 elements and, hence, as
	$R$-module by $|G_{14}|=144$ elements (recall that $u_1$ is generated as $R$ module by 2 elements).
\end{proof}

\subsection{The case of $G_{15}$}
\indent

Let $R=\ZZ[u_{s,i}^{\pm},u_{t,j}^{\pm},u_{u,l}^{\pm}]_{\substack{1\leq i,j\leq 2 \\1\leq l\leq 3}}$ and let $$H_{G_{15}}=\langle s,t,u\;|\; stu=tus,ustut=stutu,\prod\limits_{i=1}^{2}(s-u_{s,i})=\prod\limits_{j=1}^{2}(t-u_{t,j})=\prod\limits_{l=1}^{3}(u-u_{u,l})=0\rangle$$ be the generic Hecke algebra associated to $G_{15}$. Let $u_1$ be the subalgebra of $H_{G_{15}}$ generated by $s$, $u_2$ the subalgebra of $H_{G_{15}}$ generated by $t$ and $u_3$ the subalgebra of $H_{G_{15}}$ generated by $u$. We recall that $z:=stutu=ustut=tustu=tutus=utust$  generates the center of the associated complex braid group and that $|Z(G_{15})|=12$.
We set $U=\sum\limits_{k=0}^{11}z^ku_3u_2u_1u_2u_1$. 
By the definition of $U$ we have the following remark.

\begin{rem}
	$Uu_1 \subset U$.
	\label{rem15}
\end{rem}
From now on, we will underline the elements that belong to $U$ by definition.
Our goal is to prove that $H_{G_{15}}=U$ (theorem \ref{thm 15}). Since $1\in U$, it will be sufficient to prove that $U$ is a left-sided ideal of $H_{G_{15}}$. For this purpose, one must check that $sU$, $tU$ and $uU$ are subsets of $U$. The following proposition states that it is enough to prove $tU\subset U$.
\begin{prop}
If $tU\subset U$, then $H_{G_{15}}=U$.
\label{prrrr1}
\end{prop}
\begin{proof}
As we explained above, it is enough to prove that $sU$, $tU$ and $uU$ are subsets of $U$. However, by hypothesis and by the definition of $U$, we
can restrict ourselves to proving that $sU\subset U$. We recall that $z=stutu$. Therefore, $s=zu^{-1}t^{-1}u^{-1}t^{-1}$ and $s^{-1}=z^{-1}tutu$. 
We notice that $$U=
\sum\limits_{k=0}^{10}z^ku_3u_2u_1u_2u_1+z^{11}u_3u_2u_1u_2u_1.$$
Hence, $\begin{array}[t]{lcl}sU&\subset& \sum\limits_{k=0}^{10}z^ksu_3u_2u_1u_2u_1+
z^{11}su_3u_2u_1u_2u_1\\
&\subset& \sum\limits_{k=0}^{10}z^{k+1}u^{-1}t^{-1}u^{-1}t^{-1}u_3u_2u_1u_2u_1+z^{11}(R+Rs^{-1})u_3u_2u_1u_2u_1\\
&\subset& u^{-1}t^{-1}u^{-1}t^{-1}\sum\limits_{k=0}^{10}z^{k+1}u_3u_2u_1u_2u_1+z^{11}u_3u_2u_1u_2u_1+(z^{11}s^{-1})u_3u_2u_1u_2u_1\\
&\subset&u_3u_2u_3u_2U+
z^{10}tutz^ku_3u_2u_1u_2u_1\\
&\subset&U+u_3u_2u_3u_2U+
tut(z^{10}z^ku_3u_2u_1u_2u_1)\\&\subset& U+u_3u_2u_3u_2U.
\end{array}$\\\\
Since $tU\subset U$ we have $u_2U\subset U$ (recall that $u_2=R+Rt$) and, by the definition of $U$, we also have $u_3U\subset U$. The result then is obvious.
 \end{proof}

A first step to prove our main theorem is analogous to \ref{rem15} (see proposition \ref{prrr15}). For this purpose, we first prove some preliminary results.

\begin{lem}
	\mbox{}
	\vspace*{-\parsep}
	\vspace*{-\baselineskip}\\
	\begin{itemize}[leftmargin=0.8cm]
		\item[(i)] For every $k\in\{0,...,11\}$, $z^ku_3u_1u^{-1}u_2\subset U$.
			\item[(ii)] For every $k\in\{3,...,11\}$, $z^ku_3u_1u^{-2}t^{-1}\subset U$.
				\item[(iii)] For every $k\in\{3,...,11\}$, $z^ku_3u_1u_3t^{-1}\subset U$.
	\end{itemize}
	\label{sutt}
\end{lem}
\begin{proof}
Since $u_3=R+Ru^{-1}+Ru^{-2}$, (iii) follows from (i) and (ii) and from the definition of $U$. For (i) we have: $z^ku_3u_1u^{-1}u_2=z^ku_3(R+Rs^{-1})u^{-1}u_2\subset \underline{z^ku_3u_2}+z^ku_3(s^{-1}u^{-1}t^{-1})u_2\subset U+z^ku_3(u^{-1}t^{-1}s^{-1})u_2\subset U+\underline{z^ku_3u_2u_1u_2}\subset U.$
It remains to prove (ii). We have:
\\\\$\small{\begin{array}{lcl}
z^ku_3u_1u^{-2}t^{-1}&=&z^ku_3(R+Rs^{-1})u^{-2}t^{-1}\\
&\subset&\underline{z^ku_3t^{-1}}+z^ku_3(s^{-1}u^{-1}t^{-1}u^{-1}t^{-1})tutu^{-1}t^{-1}\\
	&\subset&U+z^{k-1}u_3tu(R+Rt^{-1})u^{-1}t^{-1}\\
	&\subset&U+\underline{z^{k-1}u_3}+z^{k-1}u_3tu^2(u^{-1}t^{-1}u^{-1}t^{-1}s^{-1})s\\
	&\subset&U+z^{k-2}u_3t(R+Ru+Ru^{-1})s\\
	&\subset&U+\underline{z^{k-2}u_3ts}+z^{k-2}u_3(tustu)u^{-1}t^{-1}+z^{k-2}u_3(R+Rt^{-1})u^{-1}s\\
	&\subset&U+\underline{z^{k-1}u_3t^{-1}}+\underline{z^{k-2}u_3s}+z^{k-2}u_3(u^{-1}t^{-1}u^{-1}t^{-1}s^{-1})st\\
	&\subset&U+\underline{z^{k-3}u_3st}.
\end{array}}$

\qedhere
\end{proof}

\begin{lem}
\mbox{}
\vspace*{-\parsep}
\vspace*{-\baselineskip}\\
\begin{itemize}[leftmargin=0.8cm]
	\item[(i)] For every $k\in\{0,...,10\}$, $z^ku_3u_2uu_2\subset U$.
	\item[(ii)] For every $k\in\{1,...,11\}$, $z^ku_3u_2u^{-1}u_2\subset U$.
	\item[(iii)] For every $k\in\{1,...,10\}$, $z^ku_3u_2u_3u_2\subset U$.
\end{itemize}
\label{tutt}
\end{lem}
\begin{proof}
	Since $u_3=R+Ru^{-1}+Ru^{-2}$, (iii) follows from (i) and (ii) and from the definition of $U$. For (i) we have: $z^ku_3u_2uu_2=z^ku_3(R+Rt)uu_2\subset \underline{z^ku_3u_2}+z
	^ku_3(utust)t^{-1}s^{-1}u_2\subset U+\underline{z^{k+1}u_3t^{-1}s^{-1}u_2}\subset U.$
	For (ii), we use similar kind of calculations: $z^ku_3u_2u^{-1}u_2=z^ku_3(R+Rt^{-1})u^{-1}u_2\subset \underline{z^ku_3u_2}+z
	^ku_3(u^{-1}t^{-1}u^{-1}t^{-1}s^{-1})su_2\subset U+\underline{z^{k-1}u_3su_2}.$
\end{proof}
\begin{lem}
	\mbox{}
	\vspace*{-\parsep}
	\vspace*{-\baselineskip}\\
	\begin{itemize}[leftmargin=0.8cm]
		\item[(i)] For every $k\in\{0,...,10\}$, $z^ku_3u_1tu\subset U$.
		\item[(ii)] For every $k\in\{0,...,9\}$, $z^ku_3u_1tu^2\subset U$.
		\item[(iii)] For every $k\in\{0,...,9\}$, $z^ku_3u_1tu_3\subset U$.
	\end{itemize}
	\label{stuu}
\end{lem}
\begin{proof}
	Since $u_3=R+Ru+Ru^{2}$, (iii) follows from (i) and (ii) and from the definition of $U$. For (i) we have: $z^ku_3u_1tu=z^ku_3(R+Rs)tu\subset z^ku_3tu+z^ku_3(stutu)u^{-1}t^{-1}\stackrel{\ref{tutt}(i)}{\subset}U+\underline{z^{k+1}u_3t^{-1}}\subset U.$ Similarly, for (ii): $z^ku_3u_1tu^2=z^ku_3(R+Rs)tu^2\subset z^ku_3tu^2+z^ku_3(stutu)u^{-1}t^{-1}u\stackrel{\ref{tutt}(iii)}{\subset}U+z^{k+1}u_3t^{-1}u\stackrel{\ref{tutt}(i)}{\subset}U.$
\end{proof}

\begin{lem}
	\mbox{}
	\vspace*{-\parsep}
	\vspace*{-\baselineskip}\\
	\begin{itemize}
		\item[(i)] For every $k\in\{0,...,10\}$, $z^ku_3u_2uu_1\subset U$.
		\item[(ii)] For every $k\in\{1,...,11\}$, $z^ku_3u_2u^{-1}u_1\subset U$.
		\item[(iii)] For every $k\in\{1,...,10\}$, $z^ku_3u_2u_3u_1\subset U$.
	\end{itemize}
	\label{tuss}
\end{lem}
\begin{proof}
	Since $u_3=R+Ru+Ru^{-1}$, (iii) follows from (i) and (ii) and from the definition of $U$. For (i) we have: $z^ku_3u_2uu_1=z^ku_3(R+Rt)u(R+Rs)
	\subset \underline{z^ku_3u_1}+z^ktu+z^k(tus)\stackrel{\ref{tutt}(i)}{\subset}U+z^kstu\stackrel{\ref{stuu}(i)}{\subset}U.$ Similarly, for (ii): $z^ku_3u_2u^{-1}u_1=z^ku_3(R+Rt^{-1})u^{-1}u_1\subset \underline{z^ku_3u_1}+z^ku_3(u^{-1}t^{-1}u^{-1}t^{-1}s^{-1})stu_1\subset U+\underline{z^{k-1}u_3stu_1}.$
	
\end{proof}

\begin{lem}
	\mbox{}
	\vspace*{-\parsep}
	\vspace*{-\baselineskip}\\
	\begin{itemize}[leftmargin=0.8cm]
		\item[(i)] For every $k\in\{0,...,8\}$, $z^ku_3u_1uu_1\subset U$.
		\item[(ii)] For every $k\in\{0,...,11\}$, $z^ku_3u_1u^{-1}u_1\subset U$.
		\item[(iii)] For every $k\in\{0,...,8\}$, $z^ku_3u_1u_3u_1\subset U$.
	\end{itemize}
	\label{suu}
\end{lem}

\begin{proof} 
	By remark \ref{rem15} we can ignore the $u_1$ in the end.
	Moreover, since $u_3=R+Ru+Ru^{-1}$, (iii) follows from (i) and (ii). However, $z^ku_3u_1u^{-1}\subset z^ku_3u_1u^{-1}u_2$ and, hence, (ii) follows from lemma \ref{sutt} (i). Therefore, it will be sufficient to prove (i). We have:\\ 
	$\small{\begin{array}{lcl}
	z^ku_3u_1u&=&z^ku_3(R+Rs)u\\
	&\subset& \underline{z^ku_3}+z^ku_3(ustut)t^{-1}u^{-1}t^{-1}u\\
	&\subset& U+z^{k+1}u_3(R+Rt)u^{-1}(R+Rt)u\\
	&\subset& U+z^{k+1}u_3u_2u+\underline{z^{k+1}u_3}+\underline{z^{k+1}u_3t}+z^{k+1}u_3tu^{-1}tu\\
	&\stackrel{\ref{tutt}(i)}{\subset}&U+z^{k+1}u_3t(R+Ru+Ru^2)tu\\
	&\subset&U+z^{k+1}u_3u_2u+z^{k+1}u_3(tutus)s^{-1}+z^{k+1}u_3(utust)t^{-1}s^{-1}(utust)t^{-1}s^{-1}\\
	&\stackrel{\ref{tutt}(i)}{\subset}&U+\underline{z^{k+2}u_3s^{-1}}+\underline{z^{k+3}u_3t^{-1}s^{-1}t^{-1}s^{-1}}.
	\end{array}}$
	
	\qedhere
\end{proof}
To make it easier for the reader to follow the calculations, from now on we will double-underline the elements described in the above lemmas (lemmas \ref{sutt} - \ref{suu}) and we will use directly the fact that these elements are inside $U$.

\begin{prop}
	\mbox{}
	\vspace*{-\parsep}
	\vspace*{-\baselineskip}\\
	\begin{itemize}[leftmargin=0.8cm]
		\item[(i)] For every $k\in\{2,\dots,11\}$, $z^ks^{-1}u^{-2}\in U+z^{k-3}u_3^{\times}stst$.
		\item[(ii)] For every $k\in\{0,\dots,11\}$, $z^ku_3u_1u_2u_1u_2\subset U$.
		\item[(iii)]For every $k\in\{0,\dots,11\}$, $z^ku_3u_1u_3\subset U$.
	\end{itemize}
	\label{ststt}
\end{prop}
\begin{proof}\mbox{}
	\vspace*{-\parsep}
	\vspace*{-\baselineskip}\\
	\begin{itemize}[leftmargin=0.8cm]
		\item[(i)]$\small{\begin{array}[t]{lcl}
		z^ks^{-1}u^{-2}&=&z^k(s^{-1}u^{-1}t^{-1}u^{-1}t^{-1})tutu^{-1}\\
		&\in& z^{k-1}(R+R^{\times}t^{-1})u(R+R^{\times}t^{-1})u^{-1}\\
		&\in& \underline{z^{k-1}u_3}+\underline{z^{k-1}u_3t^{-1}}+\underline{\underline{z^{k-1}u_3u_2u^{-1}u_2}}+R^{\times}z^{k-1}t^{-1}u^2(u^{-1}t^{-1}u^{-1}t^{-1}s^{-1})st\\
		&\in& U+R^{\times}z^{k-2}t^{-1}(R+Ru+R^{\times}u^{-1})st\\
		&\in& U+\underline{z^{k-2}u_3t^{-1}st}+Rz^{k-2}t^{-1}(ustut)t^{-1}u^{-1}+R^{\times}z^{k-2}u(u^{-1}t^{-1}u^{-1}t^{-1}s^{-1})stst\\
		&\in& U+\underline{\underline{z^{k-1}u_3u_2u^{-1}u_2}}+z^{k-3}u_3^{\times}stst\\ &\subset&
		U+z^{k-3}u_3^{\times}stst.
		\end{array}}$
	\item[(ii)] For $k\in\{0,...,5\}$, we have:
	$z^ku_3u_1u_2u_1u_2=z^ku_3(R+Rs)(R+Rt)(R+Rs)(R+Rt)\subset  \underline{z^ku_3u_1u_2u_1}+\underline{z^ku_3u_2u_1u_2}+z^ku_3stst
			\stackrel{(i)}{\subset}U+z^{k+3}u_3^{\times}s^{-1}u^{-2}\stackrel{\ref{suu}(iii)}{\subset}U.$ It remains to prove the case where $k\in\{6,\dots,11\}$. We have:
			
			$\small{\begin{array}[t]{lcl}	z^ku_3u_1u_2u_1u_2&=&z^ku_3(R+Rs^{-1})(R+Rt^{-1})(R+Rs^{-1})(R+Rt^{-1})\\
				&\subset& \underline{z^ku_3u_1u_2u_1}+\underline{z^ku_3u_2u_1u_2}+z^ku_3s^{-1}t^{-1}s^{-1}t^{-1}\\
			&\subset&z^ku_3s^{-1}(t^{-1}s^{-1}u^{-1}t^{-1}u^{-1})utut^{-1}\\
			&\subset& z^{k-1}u_3s^{-1}(R+Ru^{-1}+Ru^{-2})tut^{-1}\\
			&\subset& z^{k-1}u_3(R+Rs)tut^{-1}+z^{k-1}u_3s^{-1}u^{-1}(R+Rt^{-1})ut^{-1}+\\&&+
			z^{k-1}u_3s^{-1}u^{-2}(R+Rt^{-1})ut^{-1}\\
			&\subset&\underline{\underline{z^{k-1}u_3u_2uu_2}}+z^{k-1}u_3(ustut)t^{-2}+
			\underline{z^{k-1}u_3s^{-1}t^{-1}}+z^{k-1}u_3(s^{-1}u^{-1}t^{-1}u^{-1}t^{-1})tu^2t^{-1}+\\&&+\underline{\underline{z^{k-1}u_3u_1u_3t^{-1}}}+
			z^{k-1}u_3s^{-1}u^{-2}t^{-1}(R+Ru^{-1}+Ru^{-2})t^{-1}\\
			&\subset&U+\underline{z^kt^{-2}}+\underline{\underline{z^{k-2}u_3u^2u_3u_2}}+z^{k-1}u_3s^{-1}u^{-2}t^{-2}+\\&&+z^{k-1}u_3(s^{-1}u^{-1}t^{-1})t(u^{-1}t^{-1}u^{-1}t^{-1}s^{-1})s+\\&&+z^{k-1}u_3(s^{-1}u^{-1}t^{-1}u^{-1}t^{-1})tut(u^{-1}t^{-1}u^{-1}t^{-1}s^{-1})stu^{-1}t^{-1}\\
			&\subset& U+z^{k-1}u_3s^{-1}u^{-2}(R+Rt^{-1})+\underline{z^{k-2}u_3t^{-1}s^{-1}ts}+\\&&+
			z^{k-3}u_3tu(R+Rt^{-1})stu^{-1}t^{-1}\\
			&\subset&U+z^{k-1}u_3s^{-1}u^{-2}+\underline{\underline{z^{k-1}u_3u_1u_3t^{-1}}}+z^{k-3}u_3(utust)u^{-1}t^{-1}+\\&&+z^{k-3}u_3tut^{-1}(R+Rs^{-1})tu^{-1}t^{-1}\\
			&\subset& U+z^{k-1}u_3s^{-1}u^{-2}
			+\underline{(z^{k-2}+z^{k-3})u_3u_2}+z^{k-3}u_3tut^{-1}s^{-1}tu^{-1}t^{-1}.
			\end{array}}$
		
		However, we notice that $z^{k-3}u_3tut^{-1}s^{-1}tu^{-1}t^{-1}$ is a subset of $U$. Indeed, we have:
		
		$\small{\begin{array}{lcl}
		z^{k-3}u_3tut^{-1}s^{-1}tu^{-1}t^{-1}	&\subset&z^{k-3}u_3tut^{-1}s^{-1}(R+Rt^{-1})u^{-1}t^{-1}\\
		&\subset& z^{k-3}tu^2(u^{-1}t^{-1}s^{-1}u^{-1}t^{-1})+
			\\&&+z^{k-3}u_3tu^2(u^{-1}t^{-1}s^{-1}u^{-1}t^{-1})tu^2(u^{-1}t^{-1}u^{-1}t^{-1}s^{-1})s\\
			&\subset& \underline{\underline{z^{k-4}u_3u_2u_3u_2}}+z^{k-5}u_3
			t(R+Ru+Ru^{-1})tu^2s\\
			&\subset&U+\underline{\underline{z^{k-5}u_3u_2u_3u_1}}+
			z^{k-5}u_3(tutus)s^{-1}us+\\&&+z^{k-5}u_3tu^{-1}t(R+Ru+Ru^{-1})s\\
			&\subset&U+\underline{\underline{z^{k-4}u_3u_1uu_1}}+z^{k-5}u_3tu^{-1}(R+Rt^{-1})(R+Rs^{-1})+\\&&+z^{k-5}u_3tu^{-2}(utust)t^{-1}+z^{k-5}u_3tu^{-1}(R+Rt^{-1})u^{-1}s\\
			&\subset&U+\underline{\underline{(z^{k-5}+z^{k-4})u_3u_2u_3u_2}}+\underline{\underline{z^{k-5}u_3u_2u_3u_1}}
		+z^{k-5}u_3tu^{-1}t^{-1}s^{-1}+\\&&+
			z^{k-5}u_3t(u^{-1}t^{-1}u^{-1}t^{-1}s^{-1})sts\\

			&\subset&U+z^{k-5}u_3t^2u(u^{-1}t^{-1}u^{-1}t^{-1}s^{-1})+\underline{z^{k-6}u_3tsts}\\
			&\subset&U+\underline{\underline{z^{k-6}u_3u_2uu_2}}.
			\end{array}}$
			
			Hence, \begin{equation}
			z^ku_3u_1u_2u_1u_2\subset U+z^{k-1}u_3s^{-1}u^{-2},\; k\in\{6,\dots,11\}.
			\label{pp}
			\end{equation} 
		For $k\in\{6,...,9\}$ we rewrite (\ref{pp}) and we have $z^ku_3u_1u_2u_1u_2\subset U+z^{k-1}u_3u_1u_3u_1$. Therefore, by lemma \ref{suu}(iii) we have $z^ku_3u_1u_2u_1u_2\subset U$. For $k\in\{10,11\}$ we use (i) and (\ref{pp}) becomes $z^ku_3u_1u_2u_1u_2\subset U+z^{k-4}u_3stst$. However, since $k-4\in\{6,7\}$, we can apply (\ref{pp})  and we have that $z^{k-4}u_3stst\subset U+z^{k-5}u_3s^{-1}u^{-2}$. The result follows from lemma \ref{suu}(iii).

		\item [(iii)] By lemma \ref{suu} (iii), it is enough to prove that for $k\in\{9,10,11\}$, $z^ku_1u_3\subset U$. We expand $u_1$ as $R+Rs^{-1}$ and $u_3$ as $R+Ru^{-1}+Ru^{-2}$ and we have 
		$z^ku_1u_3\subset \underline{z^ku_3}+z^ku_1u^{-1}+z^ku_1u^{-2}$. Hence, by lemma  \ref{suu}(ii) we only have to prove that
		 $z^ku_3s^{-1}u^{-2}\subset U$ . However, by (i) we have $z^ku_3s^{-1}u^{-2}\subset U+z^{k-3}u_3stst$ and the result  follows directly from (ii).
		\qedhere
	\end{itemize}
	
\end{proof}
\begin{lem}
		\mbox{}
		\vspace*{-\parsep}
		\vspace*{-\baselineskip}\\
		\begin{itemize}[leftmargin=0.8cm]
			\item[(i)] For every $k\in\{3,\dots,8\}$, $z^ku_3tu^{-1}u_1u\subset U$.
				\item[(ii)] For every $k\in\{3,4\}$, $z^ku_3tu^{-1}u_1u^2\subset U$.
					\item[(iii)] For every $k\in\{5,\dots,8\}$, $z^ku_3tu^{-1}u_1u^{-1}\subset U$.
						\item[(iv)] For every $k\in\{3,\dots,8\}$, $z^ku_3tu^{-1}u_1u_3\subset U$.
		\end{itemize}
	\label{stuuuu}
\end{lem}
\begin{proof}
	\mbox{}
	\vspace*{-\parsep}
	\vspace*{-\baselineskip}\\
	\begin{itemize}[leftmargin=0.8cm]
		\item[(i)]$\small{\begin{array}[t]{lcl}
		z^ku_3tu^{-1}u_1u&\subset&z^{k}u_3t(R+Ru+Ru^2)(R+Rs)u\\
		&\subset& U+\underline{\underline{z^{k}u_3u_1tu_3}}+z^{k}u_3tsu+z^{k}u_3(tus)u+z^{k}u_3tu^2su\\ 
		&\subset&U+z^{k}u_3(R+Rt^{-1})(R+Rs^{-1})u+z^{k}u_3stu+z^{k}u_3tu(ustut)t^{-1}u^{-1}t^{-1}u\\ 
		&\subset&U+\underline{\underline{z^{k}u_3u_1u_3u_1}}+\underline{\underline{z^{k}u_3u_2u_3u_1}}+
		z^{k}u_3(t^{-1}s^{-1}u^{-1}t^{-1}u^{-1})utu^2+\\&&+\underline{\underline{z^{k}u_3u_1tu_3}}+
		z^{k+1}u_3tu(R+Rt)u^{-1}t^{-1}u\\ 		&\subset&U+\underline{\underline{z^{k-1}u_3u_2u_3u_1}}+\underline{z^{k+1}u_3}+
		z^{k+1}u_3tut(R+Ru+Ru^2)t^{-1}u\\ 
		&\subset&U+
		
		\underline{\underline{z^{k+1}u_3u_2u_3u_2}}+
		z^{k+1}u_3(tutus)s^{-1}t^{-1}u+z^{k+1}u_3(tutus)s^{-1}ut^{-1}u\\
		
		&\subset&U+z^{k+2}u_3s^{-1}(R+Rt)u+z^{k+2}u_3s^{-1}u(R+Rt)u\\ 
		
		&\subset&U+z^{k+2}u_3u_1u_3+\underline{\underline{z^{k+2}u_3u_1tu}}+z^{k+2}u_3s^{-1}(utust)t^{-1}s^{-1}\\ 
		&\stackrel{\ref{ststt}(iii)}{\subset}&U+
	\underline{z^{k+3}u_3u_1u_2u_1}.
	\end{array}}$

	\item[(ii)]$\small{\begin{array}[t]{lcl}
		z^ku_3tu^{-1}u_1u^2&\subset&z^{k}u_3t(R+Ru+Ru^2)(R+Rs)u^2\\
		&\subset& U+\underline{\underline{z^{k}u_3u_1tu_3}}+z^{k}u_3tsu^2+z^{k}u_3(tus)u^2+z^{k}u_3tu^2su^2\\ 
		&\subset&U+z^{k}u_3(R+Rt^{-1})(R+Rs^{-1})u^2+z^{k}u_3stu^2+z^{k}u_3tu(ustut)t^{-1}u^{-1}t^{-1}u^2\\ 
		&\subset&U+\underline{\underline{z^{k}u_3u_1u_3u_1}}+\underline{\underline{z^{k}u_3u_2u_3u_1}}+
		z^{k}u_3(t^{-1}s^{-1}u^{-1}t^{-1}u^{-1})utu^3+\\&&+\underline{\underline{z^{k}u_3u_1tu_3}}+
		z^{k+1}u_3tu(R+Rt)u^{-1}t^{-1}u^2\\ 		&\subset&U+\underline{\underline{z^{k-1}u_3u_2u_3u_1}}+\underline{z^{k+1}u_3}+
		z^{k+1}u_3tut(R+Ru+Ru^2)t^{-1}u^2\\ 
		&\subset&U+
		
		\underline{\underline{z^{k+1}u_3u_2u_3u_2}}+
		z^{k+1}u_3(tutus)s^{-1}t^{-1}u^2+z^{k+1}u_3(tutus)s^{-1}ut^{-1}u^2\\
		
		&\subset&U+z^{k+2}u_3s^{-1}(R+Rt)u^2+z^{k+2}u_3s^{-1}u(R+Rt)u^2\\ 
		
		&\subset&U+z^{k+2}u_3u_1u_3+\underline{\underline{z^{k+2}u_3u_1tu_3}}+z^{k+2}u_3s^{-1}(utust)t^{-1}s^{-1}u\\ 
		&\stackrel{\ref{ststt}(iii)}{\subset}&U+
		z^{k+3}u_3(R+Rs)t^{-1}s^{-1}u\\
		&\subset&z^{k+3}(u^{-1}t^{-1}s^{-1}u^{-1}t^{-1})tu^2+z^{k+3}u_3s(R+Rt)(R+Rs)u\\
		&\subset&U+\underline{\underline{z^{k+2}u_3u_1tu_3}}+\underline{\underline{z^{k+3}u_3u_1u_3u_1}}+\underline{\underline{z^{k+3}u_3u_1tu_3}}+z^{k+3}u_3stsu\\
		&\subset&U+z^{k+3}u_3(ustut)t^{-1}u^{-1}su\\
			&\subset&U+z^{k+4}u_3(R+Rt)u^{-1}su\\
			&\subset&U+z^{k+4}u_3u_1u_3+z^{k+4}tu^{-1}u_1u.
				\end{array}}$
			
The result follows from proposition \ref{ststt}(iii) and from (i).

\item[(iii)]$\small{\begin{array}[t]{lcl}
	z^ku_3tu^{-1}u_1u^{-1}&\subset&z^ku_3(R+Rt^{-1})u^{-1}(R+Rs^{-1})u^{-1}\\
	&\subset&\underline{\underline{z^ku_3u_1u_3u_1}}+\underline{\underline{z^ku_3u_2u_3u_2}}+z^ku_3t^{-1}u^{-1}s^{-1}u^{-1}\\
	&\subset&U+z^ku_3(u^{-1}t^{-1}u^{-1}t^{-1}s^{-1})sts^{-1}u^{-1}\\
	&\subset&U+z^{k-1}u_3s(R+Rt^{-1})s^{-1}u^{-1}\\
		&\subset&\underline{\underline{z^{k-1}u_3u_1u_3u_1}}+z^{k-1}u_3s(t^{-1}s^{-1}u^{-1}t^{-1}u^{-1})ut\\
		&\subset&U+z^{k-2}u_3su(R+Rt^{-1})\\
		&\subset&\underline{\underline{z^{k-2}u_3u_1u_3u_1}}+\underline{\underline{z^{k-2}u_3u_1u_3t^{-1}}}.
	
\end{array}}$
\item[(iv)] For $k\in\{3,4\}$ we have $z^ku_3tu^{-1}u_1u_3\subset z^ku_3tu^{-1}u_1(R+Ru+Ru^2) \stackrel{(i),(ii)}{\subset}U+ \underline{\underline{z^ku_3u_2u^{-1}u_1}}\subset U$.
Similarly, for $k\in\{5,\dots,8\}$ we have $z^ku_3tu^{-1}u_1u_3\subset z^ku_3tu^{-1}u_1(R+Ru+Ru^{-1}) \stackrel{(i),(iii)}{\subset}U+ \underline{\underline{z^ku_3u_2u^{-1}u_1}}$.
\qedhere
	\end{itemize}
\end{proof}

\begin{prop}$Uu_2\subset U$.
	\label{prrr15}
\end{prop}

\begin{proof}
	By the definition of $U$ and the fact that $u_2=R+Rt$, we have to prove that $z^ku_3u_2u_1u_2u_1t\subset U$, for every $k\in\{0,...,11\}$. If we expand $u_1$ as $R+Rs$ and $u_2$ as $R+Rt$ we notice that $z^ku_3u_2u_1u_2u_1t\subset z^ku_3u_2u_1u_2u_1+z^ku_3u_1u_2u_1u_2+z^ku_3tstst$. Therefore, by the definition of $U$ and by proposition \ref{ststt}(ii), 
	we only have to prove  $z^ku_3tstst\subset U$, for every $k\in\{0,\dots,11\}$. We distinguish the following cases:
	\begin{itemize}[leftmargin=*]
		\item \underline{$k\in\{0,\dots,5\}$}:
		
		$\small{\begin{array}[t]{lcl}
		z^ku_3tstst&=&z^ku_3tst(stutu)u^{-1}t^{-1}u^{-1}\\
		&\subset&z^{k+1}u_3tst(R+Ru+Ru^2)t^{-1}u^{-1}\\
		&\subset&z^{k+1}u_3tsu^{-1}+z^{k+1}u_3t(stutu)u^{-1}t^{-2}u^{-1}+z^{k+1}u_3tstu^2
		t^{-1}u^{-1}\\
		&\subset& z^{k+1}u_3(R+Rt^{-1})(R+Rs^{-1})u^{-1}+z^{k+2}u_3tu^{-1}(R+Rt^{-1})u^{-1}+
		z^{k+1}u_3tstu^2(R+Rt)u^{-1}\\
		&\subset&\underline{\underline{z^{k+1}u_3u_1u_3u_1}}+\underline{\underline{(z^{k+1}+z^{k+2})u_3u_2u_3u_1}}+
		
		z^{k+1}u_3(t^{-1}s^{-1}u^{-1}t^{-1}u^{-1})ut+\\&&+
		z^{k+2}u_3t(u^{-1}t^{-1}u^{-1}t^{-1}s^{-1})st+z^{k+1}u_3t(stutu)u^{-1}t^{-1}+
		z^{k+1}u_3t(stutu)u^{-1}t^{-1}utu^{-1}\\
		&\subset&U+\underline{(z^k+z^{k+1})u_3u_2u_1u_2}+\underline{\underline{z^{k+2}u_3u_2u^{-1}u_2}}+z^{k+2}u_3tu^{-1}(R+Rt)utu^{-1}\\
		&\subset& U+\underline{\underline{z^{k+2}u_3u_2u^{-1}u_2}}+z^{k+2}u_3tu^{-1}(tutus)s^{-1}u^{-2}\\
		&\subset&U+z^{k+3}u_3tu^{-1}u_1u_3.
		\end{array}}$
	
	The result follows from lemma \ref{stuuuu}(iii).
		\item 	\underline{$k\in\{6,\dots,11\}$}:
		
		$\small{\begin{array}[t]{lcl}
		z^ku_3v_8t&=&z^ku_3tstst\\
		&\subset&z^ku_3(R+Rt^{-1})(R+Rs^{-1})(R+Rt^{-1})(R+Rs^{-1})(R+Rt^{-1})\\
		&\subset&z^ku_3u_1u_2u_1u_2+z^kt^{-1}s^{-1}t^{-1}s^{-1}\\
		&\stackrel{\ref{ststt}(ii)}{\subset}& U+z^ku_3(t^{-1}s^{-1}u^{-1}t^{-1}u^{-1})utu^2tu(u^{-1}t^{-1}u^{-1}t^{-1}s^{-1})t^{-1}\\
		&\subset&U+z^{k-2}u_3t(R+Ru+Ru^{-1})tut^{-1}\\
		&\subset&U+\underline{\underline{z^{k-2}u_3t^2ut^{-1}}}+z^{k-2}u_3(tutus)s^{-1}t^{-1}+
		z^{k-2}u_3tu^{-1}(R+Rt^{-1})ut^{-1}\\
		&\subset&U+\underline{z^{k-1}u_3s^{-1}t^{-1}+z^{k-2}u_3}+z^{k-2}u_3(R+Rt^{-1})u^{-1}t^{-1}ut^{-1}\\
		&\subset&U+\underline{\underline{z^{k-2}u_3t^{-1}ut^{-1}}}+z^{k-2}u_3(u^{-1}t^{-1}u^{-1}t^{-1}s^{-1})su^{-1}t^{-1}\\
		&\subset&U+\underline{\underline{z^{k-3}u_3su^{-1}t^{-1}}}.
		\end{array}}$
	
	\qedhere
	\end{itemize}
\end{proof}
We can now prove the main theorem of this section.
\begin{thm} $H_{G_{15}}=U$.
	
	\label{thm 15}
\end{thm}
\begin{proof}
	By proposition \ref{prrrr1} it is enough to prove that $tU\subset U$. By remark \ref{rem15} and proposition \ref{prrr15}, we only have to prove that $z^ktu_3\subset U$. By lemma \ref{tutt} (iii), we only have to check the cases where $k\in\{0,11\}$. We have:
	\begin{itemize}
		\item \underline{$k=0$}:
		$\begin{array}[t]{lcl}
		tu_3&=&t(R+Ru+Ru^2)\\
		&\subset&\underline{t}+\underline{\underline{tu}}+tu^2\\
		&\subset& U+s^{-1}(stutu)u^{-1}t^{-1}u\\
		&\subset& U+zs^{-1}u^{-1}(R+Rt)u\\
		&\subset& U+\underline{zs}+zs^{-1}u^{-2}(utust)t^{-1}s^{-1}\\
		&\subset&U+zu_1u_3t^{-1}s^{-1}\\
		&\stackrel{\ref{ststt}(iii)}{\subset}&U+Uu_2u_1.
		\end{array}$
	
		\item \underline{$k=11$}:
		$\begin{array}[t]{lcl}
		z^{11}tu_3&\subset&z^{11}(R+Rt^{-1})(R+Ru^{-1}+Ru^{-2})\\
		&\subset&\underline{z^{11}u_3}+\underline{\underline{z^{11}t^{-1}u^{-1}}}+z^{11}t^{-1}u^{-2}\\
		&\subset& U+z^{11}u(u^{-1}t^{-1}u^{-1}t^{-1}s^{-1})stu^{-1}\\
		&\subset&U+z^{10}u_3s(R+Rt^{-1})u^{-1}\\
		&\subset&U+\underline{\underline{z^{10}u_3su^{-1}}}+z^{10}u_3su(u^{-1}t^{-1}u^{-1}t^{-1}s^{-1})st\\
		&\subset&U+z^9u_3sust\\
		&\stackrel{\ref{ststt}(iii)}{\subset}&U+Uu_1u_2.
		\end{array}$
	
	\end{itemize}
	The result follows from remark \ref{rem15} and proposition \ref{prrr15}.
\end{proof}

\begin{cor}
	The BMR freeness conjecture holds for the generic Hecke algebra $H_{G_{15}}$.
\end{cor}
\begin{proof}
	By theorem \ref{thm 15} we have that $H_{G_{15}}=U=\sum\limits_{k=0}^{11}z^k(u_3+u_3s+u_3t+u_3ts+u_3st+u_3tst+u_3sts+u_3tsts)$. The result follows from proposition \ref{BMR PROP}, since $H_{G_{15}}$ is generated as left $u_3$-module by 96 elements and, hence, as
	$R$-module by $|G_{15}|=288$ elements (recall that $u_3$ is generated as $R$-module by 3 elements).
\end{proof}

\begin{appendices}
	\chapter{The BMR and ER presentations}
\indent

Let $W$ be an exceptional group of rank 2 and $B$ the complex braid group associated to $W$. We know that for every complex reflection group we have  a Coxeter-like presentation (see remark \ref{coxeterlike}(ii)) and for the associated complex braid group we have an Artin-like presentation (see theorem \ref{Presentt}). For the exceptional groups of rank 2 we
 call these presentations \emph{the BMR presentations}, due to M. Brou\'e, G. Malle and R. Rouquier. 
 
 In 2006 P. Etingof and E. Rains gave different presentations of $W$ and $B$, based on the BMR presentations associated to the maximal groups $G_7$, $G_{11}$ and $G_{19}$ (see \textsection 6.1 of \cite{ERrank2}). We call these presentations \emph{the ER presentations}. In table \ref{t3} we give the ER and BMR presentation for every exceptional group of rank 2 and also the isomorphism between them.
 
  We will now prove that we have  indeed such an isomorphism. An argument in \cite{ERrank2} is that the ER presentation of $W$ has some defining relations of the form $g^p=1$ and that the ER presentation of $B$ is obtained by removing such relations. As we can see in the following pages, when we prove that the ER presentation of $W$ is isomorphic to its BMR presentation we don't use the relations $g^p=1$ in order to prove the correspondence between the other kind of relations in ER presentation and the braid relations in BMR presentation. Hence, we restrict ourselves to proving that the ER presentation and BMR presentation of $W$ are isomorphic, since the case of $B$ can be proven in the same way.
\subsection*{The groups $G_4$, $G_8$ and $G_{16}$}

\indent

Let $W$ be the group $G_4$, $G_8$ or $G_{16}$. For these groups the BMR presentation is of the form $$\langle s,t\;|\;s^k=t^k=1, sts=tst\rangle,$$ where $k=3,4,5$, for each group respectively. We also know that the center of the group in each case is generated by the element $z=(st)^3$ (see Tables in Appendix 1 in \cite{broue2000}). 
According to \S6.1 in \cite{ERrank2}, the ER presentation of these groups is of the form
$$\langle a,b,c\;|\; a^2=b^{-3}=\text{central}, c^k=1, abc=1\rangle,$$ where $k=3,4,5$, for each group respectively. Let $\tilde W$ be a group having this ER presentation. We will prove that $\tilde W\simeq W$.
Let $\grf_1:W\rightarrow \tilde W$ and $\grf_2: \tilde W \rightarrow W$, defined by $
\grf_1(s)=c$, $\grf_1(t)=c^{-1}b$ and $\grf_2(a)=(sts)^{-1}$, $\grf_2(b) =st$ and $\grf_2(c)=s$. 

 We now prove that $\grf_1$ is a well-defined group homomorphism.
Since $abc=1$ and $a^{-2}=b^3$ we have: $a(abc)bc=1\an bcbc=a^{-2}\an cbc=b^{-1}a^{-2} \an cbc=b^2\an bc=c^{-1}b^2$, meaning that $\grf_1(s)\grf_1(t)\grf_1(s)=\grf_1(t)\grf_1(s)\grf_1(t)$. Moreover, $\grf_1(s)^k=c^k=1$ and by the relations $abc=1$ and $a^{-2}=b^3$ we also have $(c^{-1}b)^k=(ab^2)^k=(aa^{-2}b^{-1})^k=(a^{-1}b^{-1})^k=a^{-1}(b^{-1}a^{-1})^ka=a^{-1}c^ka=1$, meaning that $\grf_1(t)^k=1$. 
Similarly, we can prove that $\grf_2$ is a well-defined group homomorphism, since $\grf_2(a)^2=(sts)^{-1}(sts)^{-1}=(tst)^{-1}(sts)^{-1}=z^{-1}=\grf_2(b)^3$, $\grf_2(c)^k=s^k=1$ and $\grf_2(a)\grf_2(b)\grf_2(c)=1$.

We also notice that $\grf_1(\grf_2(a))=c^{-1}b^{-1}=a$, $\grf_1(\grf_2(b))=b$ and $\grf_1(\grf_2(c))=c$. Moreover, $\grf_2(\grf_1(s))=s$ and $\grf_2(\grf_1(t))=t$, meaning that $\grf_1\circ\grf_2=\text{id}_{\tilde W}$ and $\grf_2\circ\grf_1=\text{id}_{W}$ and, hence, $W\simeq \tilde W$.

\subsection*{The groups $G_5$, $G_{10}$ and $G_{18}$}

\indent

Let $W$ be the group $G_5$, $G_{10}$ or $G_{18}$. For these groups the BMR presentation is of the form $$\langle s,t\;|\;s^3=t^k=1, stst=tsts\rangle,$$ where $k=3,4,5$, for each group respectively. We also know that the center of the group in each case is generated by the element $z=(st)^2$ (see Tables in Appendix 1 in \cite{broue2000}). 
According to \S6.1 in \cite{ERrank2}, the ER presentation of these groups is of the form
$$\langle a,b,c\;|\; a^2=\text{central}, b^3=c^k=1, abc=1\rangle,$$ where $k=3,4,5$, for each group respectively. Let $\tilde W$ be a group having this ER presentation. We will prove that $\tilde W\simeq W$.
Let $\grf_1:W\rightarrow \tilde W$ and $\grf_2: \tilde W \rightarrow W$, defined by $
\grf_1(s)=b$, $\grf_1(t)=c$ and $\grf_2(a)=(st)^{-1}$, $\grf_2(b) =s$ and $\grf_2(c)=t$. 

We now prove that $\grf_1$ is a well-defined group homomorphism.
Since $abc=cab=1$ and $a^{2}:=\grz$ is central we have: $cbcb=ca^{-1}(abc)b=ca^{-1}b=\grz^{-1}cab=\grz^{-1}=a^{-2}=a^{-1}(abc)a^{-1}(abc)=bcbc$, meaning that $\grf_1(t)\grf_1(s)\grf_1(t)\grf_1(s)=\grf_1(s)\grf_1(t)\grf_1(s)\grf_1(t)$.
Moreover, $\grf_1(s)^3=b^3=1$ and  $\grf_1(t)^k=c^k=1$. Similarly, we can prove that $\grf_2$ is a well-defined  group homomorphism, since $\grf_2(a)^2=(st)^{-2}=z^{-1}$, $\grf_2(b)^3=s^3=1$, $\grf_2(c)^k=t^k=1$ and $\grf_2(a)\grf_2(b)\grf_2(c)=1$.

We also notice that $\grf_1(\grf_2(a))=c^{-1}b^{-1}=a$, $\grf_1(\grf_2(b))=b$ and $\grf_1(\grf_2(c))=c$. Moreover, $\grf_2(\grf_1(s))=s$ and $\grf_2(\grf_1(t))=t$, meaning that $\grf_1\circ\grf_2=\text{id}_{\tilde W}$ and $\grf_2\circ\grf_1=\text{id}_{W}$ and, hence, $W\simeq \tilde W$.

\subsection*{The groups $G_6$, $G_{9}$ and $G_{17}$}

\indent

Let $W$ be the group $G_6$, $G_{9}$ or $G_{17}$. For these groups the BMR presentation is of the form $$\langle s,t\;|\;s^2=t^k=1, (st)^3=(ts)^3\rangle,$$ where $k=3,4,5$, for each group respectively. We also know that the center of the group in each case is generated by the element $z=(st)^3$ (see Tables in Appendix 1 in \cite{broue2000}). 
According to \S6.1 in \cite{ERrank2}, the ER presentation of these groups is of the form
$$\langle a,b,c\;|\; a^2=c^k=1, b^3=\text{central}, abc=1\rangle,$$ where $k=3,4,5$, for each group respectively. Let $\tilde W$ be a group having this ER presentation. We will prove that $\tilde W\simeq W$.
Let $\grf_1:W\rightarrow \tilde W$ and $\grf_2: \tilde W \rightarrow W$, defined by $
\grf_1(s)=a$, $\grf_1(t)=c$ and $\grf_2(a)=s$, $\grf_2(b) =(ts)^{-1}$ and $\grf_2(c)=t$. 

We now prove that $\grf_1$ is a well-defined group homomorphism.
Since $abc=cab=1$ and $b^{3}:=\grz$ is central we have: $cacaca=(cab)b^{-1}(cab)b^{-1}(cab)b^{-1}=b^{-3}=\grz^{-1}=\grz^{-1}abc=a\grz^{-1}bc=ab^{-2}c=a(cab)b^{-1}(cab)b^{-1}c=acacac$, meaning that $\big(\grf_1(t)\grf_1(s)\big)^3=\big(\grf_1(s)\grf_1(t)\big)^3$.
Moreover, $\grf_1(s)^2=s^2=1$ and  $\grf_1(t)^k=c^k=1$. Similarly, we can prove that $\grf_2$ is a well-defined group homomorphism, since $\grf_2(a)^2=s^2=1$, $\grf_2(b)^3=(ts)^{-3}=z^{-1}$, $\grf_2(c)^k=t^k=1$ and $\grf_2(a)\grf_2(b)\grf_2(c)=1$.

We also notice that $\grf_1(\grf_2(a))=a$, $\grf_1(\grf_2(b))=a^{-1}c^{-1}=b$ and $\grf_1(\grf_2(c))=c$. Moreover, $\grf_2(\grf_1(s))=s$ and $\grf_2(\grf_1(t))=t$, meaning that $\grf_1\circ\grf_2=\text{id}_{\tilde W}$ and $\grf_2\circ\grf_1=\text{id}_{W}$ and, hence, $W\simeq \tilde W$.

\subsection*{The groups $G_7$, $G_{11}$ and $G_{19}$}

\indent

This is the case of the maximal groups. We notice that in this case the BMR presentation and the ER presentation coincide. 
\subsection*{The group $G_{12}$}

\indent

Let $W$ be the group $G_{12}$, whose BMR presentation is of the form $$\langle s,t,u\;|\;s^2=t^2=u^2=1, stus=tust=ustu\rangle.$$ We also know that the center of the group  is generated by the element $z=(stu)^4$ (see Tables in Appendix 1 in \cite{broue2000}). 
According to \S6.1 in \cite{ERrank2}, the ER presentation of $W$ is of the form
$$\langle a,b,c\;|\; a^2=1, b^3=c^{-4}=\text{central}, abc=1\rangle.$$ Let $\tilde W$ be a group having this ER presentation. We will prove that $\tilde W\simeq W$.
Let $\grf_1:W\rightarrow \tilde W$ and $\grf_2: \tilde W \rightarrow W$, defined by $
\grf_1(s)=a$, $\grf_1(t)=a^{-1}c^2b$, $\grf_1(u)=(cb)^{-1}$ and $\grf_2(a)=s$, $\grf_2(b) =(stus)^{-1}$ and $\grf_2(c)=stu$. 

We now prove that $\grf_1$ is a well-defined group homomorphism.
Since $cab=1$ we have $ca=b^{-1}$ meaning that $\grf_1(s)\grf_1(t)\grf_1(u)\grf_1(s)=\grf_1(u)\grf_1(s)\grf_1(t)\grf_1(u)$. Moreover, since $bca=1$ and $b^{-3}=c^{4}$ we have  $a^{-1}c^3b=a^{-1}c^{-1}b^{-2}=(bca)^{-1}b^{-1}=b^{-1}$, meaning that we also have  $\grf_1(t)\grf_1(u)\grf_1(s)\grf_1(t)= \grf_1(u)\grf_1(s)\grf_1(t)\grf_1(u)$.
Moreover, $\grf_1(s)^2=a^2=1$ and since $a^2=1$ and $bca=cab=1$ we have $\grf_1(u)^2=b^{-1}c^{-1}b^{-1}c^{-1}=b^{-1}a(a^{-1}c^{-1}b^{-1})c^{-1}=b^{-1}a^{-1}c^{-1}=1$. Since $a^{-1}=bc$, $c^3=c^{-1}b^{-3}$, $a^{-1}=a$ and $abc=bca=1$ we also have 
$\grf_1(t)^2=a^{-1}c^2ba^{-1}c^2b=bc^3b^2c^3b=bc^{-1}b^{-1}c^{-1}b^{-2}=b(c^{-1}b^{-1}a^{-1})ac^{-1}b^{-2}=bac^{-1}b^{-2}=b(a^{-1}c^{-1}b^{-1})b^{-1}=1$. Similarly, we can prove that $\grf_2$ is a well-defined group homomorphism, since $\grf_2(a)^2=s^2=1$, $\grf_2(b)^3=\big(stus(stus)(stus)\big)^{-1}=\big(stus(tust)(ustu)\big)^{-1}=
(stu)^{-4}=\grf_2(c)^{-4}=z^{-1}$ and $\grf_2(a)\grf_2(b)\grf_2(c)=1$.

We also notice that $\grf_1(\grf_2(a))=a$, $\grf_1(\grf_2(b))=a^{-1}c^{-1}=b$ and $\grf_1(\grf_2(c))=c$. Moreover, $\grf_2(\grf_1(s))=s$, $\grf_2(\grf_1(t))=tustu(s^{-1}u^{-1}t^{-1}s^{-1})=tustu(u^{-1}t^{-1}s^{-1}u^{-1})=t$ and $\grf_2(\grf_1(u))=(stus)u^{-1}t^{-1}s^{-1}=(ustu)u^{-1}t^{-1}s^{-1}=u$,  meaning that $\grf_1\circ\grf_2=\text{id}_{\tilde W}$ and $\grf_2\circ\grf_1=\text{id}_{W}$ and, hence, $W\simeq \tilde W$.

\subsection*{The group $G_{13}$}

\indent

Let $W$ be the group $G_{13}$, whose BMR presentation is of the form $$\left\langle
\begin{array}{l|cl}
&s^2=t^2=u^2=1 &\\s,t,u&stust=ustus&\\
&tust=ustu&
\end{array}\right \rangle.$$ We also know that the center of the group  is generated by the element $z=(stu)^3=(tus)^3=(ust)^3$ (see Tables in Appendix 1 in \cite{broue2000}). 
According to \S6.1 in \cite{ERrank2}, the ER presentation of $W$ is of the form
$$\left	\langle \begin{array}{l|cl} &a^2=d^2=1&\\a,b,c,d &b^3=dc^{-2}=\text{central}&\\ &abc=1&
\end{array}\right\rangle.$$ Let $\tilde W$ be a group having this ER presentation. We will prove that $\tilde W\simeq W$.
Let $\grf_1:W\rightarrow \tilde W$ and $\grf_2: \tilde W \rightarrow W$, defined by $
\grf_1(s)=d$, $\grf_1(t)=a$, $\grf_1(u)=b(da)^{-1}$ and $\grf_2(a)=t$, $\grf_2(b) =ust$, $\grf_2(c)=(tust)^{-1}$ and $\grf_2(d)=s$. 

We now prove that $\grf_1$ is a well-defined group homomorphism.
Since $abc=1$ and $dc^{-1}=b^3c$ we have $dab=d(abc)c^{-1}=dc^{-1}=b^3c=b^2a^{-1}(abc)=b^2a^{-1}$ meaning that $\grf_1(s)\grf_1(t)\grf_1(u)\grf_1(s)\grf_1(t)=\grf_1(u)\grf_1(s)\grf_1(t)\grf_1(u)\grf_1(s)$. Moreover, since $cab=1$ and $dc^{-2}=b^3:=\grz$ we have  $b^2a^{-1}d^{-1}=b^3(b^{-1}a^{-1}c^{-1})cd^{-1}=dc^{-1}d^{-1}=\grz cd^{-1}=c\grz d^{-1}=cc^{-2}=c^{-1}=ab$, meaning that we also have  $\grf_1(u)\grf_1(s)\grf_1(t)\grf_1(u)=\grf_1(t)\grf_1(u)\grf_1(s)\grf_1(t)$.
Moreover, $\grf_1(s)^2=d^2=1$, $\grf_1(t)^2=a^2=1$ and since $b^3=dc^{-2}:=\grz\an d^{-1}\grz=c^{-2}$ and $abc=cab=1$ we also have  $ba^{-1}d^{-1}ba^{-1}d^{-1}=b^2(b^{-1}a^{-1}c^{-1})cd^{-1}b^2(b^{-1}a^{-1}c^{-1})cd^{-1}=b^{-1}b^3cd^{-1}b^{-1}b^3cd^{-1}=b^{-1}\grz cd^{-1}b^{-1}\grz c d^{-1}$. However, we have $b^{-1}\grz cd^{-1}b^{-1}\grz c d^{-1}=
b^{-1}c(d^{-1}\grz) b^{-1}c(d^{-1}\grz)=b^{-1}c^{-1}b^{-1}c^{-1}=b^{-1}(c^{-1}b^{-1}a^{-1})ac^{-1}=b^{-1}a^{-1}c^{-1}=1$, meaning that $\grf_1(u)^2=1$. Similarly, we can prove that $\grf_2$ is a well-defined group homomorphism: $\grf_2(f)^2=d^2=1,$ $\grf_2(d)\grf_2(c)^{-2}=(stust)tust=ustustust=\grf_2(b)^3=z$ and $\grf_2(a)\grf_2(b)\grf_2(c)=1$.

We also notice that $\grf_1(\grf_2(a))=a$, $\grf_1(\grf_2(b))=b$, $\grf_1(\grf_2(c))=b^{-1}a^{-1}=c$ and $\grf_1(\grf_2(d))=d$. Moreover, $\grf_2(\grf_1(s))=s$, $\grf_2(\grf_1(t))=t$ and $\grf_2(\grf_1(u))=u$,  meaning that $\grf_1\circ\grf_2=\text{id}_{\tilde W}$ and $\grf_2\circ\grf_1=\text{id}_{W}$ and, hence, $W\simeq \tilde W$.

\subsection*{The group $G_{14}$}

\indent

Let $W$ be the group $G_{14}$, whose BMR presentation is of the form $$\langle s,t\;|\;s^2=t^3=1, (st)^4=(ts)^4\rangle.$$ We also know that the center of the group  is generated by the element $z=(st)^4$ (see Tables in Appendix 1 in \cite{broue2000}). 
According to \S6.1 in \cite{ERrank2}, the ER presentation of $W$ is of the form
$$\langle a,b,c\;|\; a^2=b^3=1, c^{4}=\text{central}, abc=1\rangle.$$ Let $\tilde W$ be a group having this ER presentation. We will prove that $\tilde W\simeq W$.
Let $\grf_1:W\rightarrow \tilde W$ and $\grf_2: \tilde W \rightarrow W$, defined by $
\grf_1(s)=a$, $\grf_1(t)=b$ and $\grf_2(a)=s$, $\grf_2(b) =t$ and $\grf_2(c)=(st)^{-1}$. 

We now prove that $\grf_1$ is a well-defined group homomorphism.
Since $abc=1$ and $c^{4}:=\grz$ we have $(ab)^4=c^{-4}=\grz^{-1}=a^{-1}\grz^{-1} a=a^{-1}(abc)c^{-4}a=bc^{-3}a=b(abc)c^{-1}(abc)c^{-1}(abc)c^{-1}a=(ba)^4$ 
meaning that $\big(\grf_1(s)\grf_1(t)\big)^4=\big(\grf_1(s)\grf_1(t)\big)^4$. We also have  
$\grf_1(s)^2=a^2=1$ and $\grf_2(b)^3=t^3=1$.  Similarly, we can prove that $\grf_2$ is a well-defined group homomorphism, since $\grf_2(a)^2=s^2=1$, $\grf_2(b)^3=t^3=1$, $\grf_2(c)^4=(st)^{-4}=z^{-1}$ and $\grf_2(a)\grf_2(b)\grf_2(c)=1$.

We also notice that $\grf_1(\grf_2(a))=a$, $\grf_1(\grf_2(b))=a^{-1}c^{-1}=b$ and $\grf_1(\grf_2(c))=b^{-1}a^{-1}=c$. Moreover, $\grf_2(\grf_1(s))=s$ and  $\grf_2(\grf_1(t))=t$,  meaning that $\grf_1\circ\grf_2=\text{id}_{\tilde W}$ and $\grf_2\circ\grf_1=\text{id}_{W}$ and, hence, $W\simeq \tilde W$.

\subsection*{The group $G_{15}$}

\indent

Let $W$ be the group $G_{15}$, whose BMR presentation is of the form $$\left\langle
\begin{array}{l|cl}
&s^2=t^2=u^3=1 &\\s,t,u&ustut=stutu&\\
&tus=stu&
\end{array}\right \rangle.$$ We also know that the center of the group  is generated by the element $z=stutu$ (see Tables in Appendix 1 in \cite{broue2000}). 
According to \S6.1 in \cite{ERrank2}, the ER presentation of $W$ is of the form
$$\left	\langle \begin{array}{l|cl} &a^2=b^3=d^2=1&\\a,b,c,d &dc^{-2}=\text{central}&\\ &abc=1&
\end{array}\right\rangle.$$ Let $\tilde W$ be a group having this ER presentation. We will prove that $\tilde W\simeq W$.
Let $\grf_1:W\rightarrow \tilde W$ and $\grf_2: \tilde W \rightarrow W$, defined by $
\grf_1(s)=d$, $\grf_1(t)=a$, $\grf_1(u)=b$ and $\grf_2(a)=t$, $\grf_2(b) =u$, $\grf_2(c)=(tu)^{-1}$ and $\grf_2(d)=s$. 

We now prove that $\grf_1$ is a well-defined group homomorphism.
Since $abc=(bca)=1$ and $dc^{-2}:=\grz$ we have $bdaba=bd(abc)c^{-1}a=b(dc^{-2})c^{-1}ca=
b\grz ca=(bca)\grz=dc^{-2}=d(abc)c^{-1}(abc)c^{-1}=dabab$,
 meaning that $\grf_1(u)\grf_1(s)\grf_1(t)\grf_1(u)\grf_1(t)=\grf_1(s)\grf_1(t)\grf_1(u)\grf_1(t)\grf_1(u)$. Moreover, since $abc=1$ and $d=\grz c^2$ we have $abd=(abc)c^{-1}d=c^{-1}\grz c^2=(\grz c^2)c^{-1}=dc^{-1}=d(abc)c^{-1}=dab$, meaning that we also have  $\grf_1(t)\grf_1(u)\grf_1(s)=\grf_1(s)\grf_1(t)\grf_1(u)$.
Moreover, $\grf_1(s)^2=d^2=1$, $\grf_1(t)^2=a^2=1$ and  $\grf_1(u)^3=b^3=1$. Similarly, we can prove that $\grf_2$ is a well-defined group homomorphism: $\grf_2(a)^2=t^2=1,$ $\grf_2(b)^3=u^3=1$, $\grf_2(d)^2=s^2=1$, 
$\grf_2(d)\grf_2(c)^{-2}=stutu=z$ and $\grf_2(a)\grf_2(b)\grf_2(c)=1$.

We also notice that $\grf_1(\grf_2(a))=a$, $\grf_1(\grf_2(b))=b$, $\grf_1(\grf_2(c))=b^{-1}a^{-1}=c$ and $\grf_1(\grf_2(d))=d$. Moreover, $\grf_2(\grf_1(s))=s$, $\grf_2(\grf_1(t))=t$ and $\grf_2(\grf_1(u))=u$,  meaning that $\grf_1\circ\grf_2=\text{id}_{\tilde W}$ and $\grf_2\circ\grf_1=\text{id}_{W}$ and, hence, $W\simeq \tilde W$.

\subsection*{The group $G_{20}$}

\indent

Let $W$ be the group $G_{14}$, whose BMR presentation is of the form $$\langle s,t\;|\;s^3=t^3=1, ststs=tstst\rangle.$$ We also know that the center of the group  is generated by the element $z=(st)^5$ (see Tables in Appendix 1 in \cite{broue2000}). 
According to \S6.1 in \cite{ERrank2}, the ER presentation of $W$ is of the form
$$\left\langle \begin{array}{l|cl} &a^2=\text{central}&\\
a,b,c&b^3=1,a^4=c^{-5}&\\ &abc=1&
\end{array}\right\rangle .$$ Let $\tilde W$ be a group having this ER presentation. We will prove that $\tilde W\simeq W$.
Let $\grf_1:W\rightarrow \tilde W$ and $\grf_2: \tilde W \rightarrow W$, defined by $
\grf_1(s)=b$, $\grf_2(a)=(ststs)^{-1}$, $\grf_2(b) =s$ and $\grf_2(c)=tsts$. 

We now prove that $\grf_1$ is a well-defined group homomorphism.
Since $abc=1$, $a^2:=\grz$ and $\grz^{-2}c^{-5}=1$ we have $bc^{-1}a^{-1}bc^{-1}a^{-1}b=bc^{-1}a^{-2}(abc)c^{-2}a^{-2}(abc)c^{-1}=bc^{-1}\grz^{-1} c^{-2}\grz^{-1} c^{-1}=b(\grz^{-2} c^{-4})=bc=a^{-1}(abc)=\grz^{-2} c^{-5} a^{-1}= c^{-1}\grz^{-1} c^{-2} \grz^{-1} c^{-2} a^{-1}$. However, since $\grz^{-1}=a^{-2}$, $c^{-1}\grz^{-1} c^{-2} \grz^{-1} c^{-2} a^{-1}= c^{-1}a^{-2}c^{-2}a^{-2}c^{-2}a^{-1}=c^{-1}a^{-2}(abc)c^{-2}a^{-2}(abc)c^{-2}a^{-1}=(c^{-1}a^{-1}b)^2c^{-1}a^{-1}$,  
meaning that $\grf_1(s)\grf_1(t)\grf_1(s)\grf_1(t)\grf_1(s)=\grf_1(t)\grf_1(s)\grf_1(t)\grf_1(s)\grf_1(t)$. We also have  
$\grf_1(s)^3=b^3=1$ and since $abc=bca=1$ and $b^2=b^{-1}$ we also have $\grf_2(t)^3=c^{-1}(a^{-1}c^{-1}b^{-1})b(a^{-1}c^{-1}b^{-1})ba^{-1}=c^{-1}b^2a^{-1}=c^{-1}b^{-1}a^{-1}=1$. Similarly, we can prove that $\grf_2$ is a well-defined group homomorphism, since $\grf_2(a)^2=(ststs)^{-1}(ststs)^{-1}=(st)^{-5}=z^{-1}$, $\grf_2(b)^3=s^3=1$, $\grf_1(c)^{-5}=(ts)^{10}=(st)^{10}=z^{2}=\grf_2(a)^4$ and $\grf_2(a)\grf_2(b)\grf_2(c)=1$.

Since $abc=1$, $a^2:=\grz$ and $\grz^{-2}c^{-5}=1$ we have that  $\grf_1(\grf_2(a))=(bc^{-1}a^{-1}bc^{-1}a^{-1}b)^{-1}=(bc^{-1}a^{-2}(abc)c^{-2}a^{-2}(abc)c^{-1})^{-1}=(bc^{-1}\grz^{-1} c^{-2}\grz^{-1} c^{-1}=b(\grz^{-2} c^{-4}))^{-1}=(bc)^{-1}=a$. Moreover, $\grf_1(\grf_2(b))=b$ and since $a^2:=\grz$ and $\grz^{-2}c^{-5}=1$ we have $\grf_1(\grf_2(c))=c^{-1}a^{-2}(abc)c^{-2}a^{-2}(abc)c^{-1}=c^{-1}\grz^{-1}c^{-2}\grz^{-1}c^{-1}=\grz^{-2}c^{-4}=c$. Furthermore, $\grf_2(\grf_1(s))=s$ and  $\grf_2(\grf_1(t))=(ts)^{-2}(ststs)=(ts)^{-2}(tstst)=t$,  meaning that $\grf_1\circ\grf_2=\text{id}_{\tilde W}$ and $\grf_2\circ\grf_1=\text{id}_{W}$ and, hence, $W\simeq \tilde W$.

\subsection*{The group $G_{21}$}

\indent

Let $W$ be the group $G_{14}$, whose BMR presentation is of the form $$\langle s,t\;|\;s^3=t^3=1, (st)^5=(ts)^5\rangle.$$ We also know that the center of the group  is generated by the element $z=(st)^5$ (see Tables in Appendix 1 in \cite{broue2000}). 
According to \S6.1 in \cite{ERrank2}, the ER presentation of $W$ is of the form
$$\langle a,b,c\;|\; a^2=b^3=1, c^{5}=\text{central}, abc=1\rangle.$$ Let $\tilde W$ be a group having this ER presentation. We will prove that $\tilde W\simeq W$.
Let $\grf_1:W\rightarrow \tilde W$ and $\grf_2: \tilde W \rightarrow W$, defined by $
\grf_1(s)=a$, $\grf_1(t)=b$ and $\grf_2(a)=s$, $\grf_2(b) =t$ and $\grf_2(c)=(st)^{-1}$. 

We now prove that $\grf_1$ is a well-defined group homomorphism.
Since $abc=1$ and $c^{5}:=\grz$ we have $(ab)^5=c^{-5}=\grz^{-1}=a^{-1}\grz^{-1} a=a^{-1}(abc)c^{-5}a=bc^{-4}a=b(abc)c^{-1}(abc)c^{-1}(abc)c^{-1}(abc)c^{-1}a=(ba)^4$
meaning that $\big(\grf_1(s)\grf_1(t)\big)^5=\big(\grf_1(s)\grf_1(t)\big)^5$. We also have  
$\grf_1(s)^2=a^2=1$ and $\grf_2(b)^3=t^3=1$. Similarly, we can prove that $\grf_2$ is a well-defined group homomorphism, since $\grf_2(a)^2=s^2=1$, $\grf_2(b)^3=t^3=1$, $\grf_2(c)^5=(st)^{-5}=z^{-1}$ and $\grf_2(a)\grf_2(b)\grf_2(c)=1$.

We also notice that $\grf_1(\grf_2(a))=a$, $\grf_1(\grf_2(b))=a^{-1}c^{-1}=b$ and $\grf_1(\grf_2(c))=b^{-1}a^{-1}=c$. Moreover, $\grf_2(\grf_1(s))=s$ and  $\grf_2(\grf_1(t))=t$,  meaning that $\grf_1\circ\grf_2=\text{id}_{\tilde W}$ and $\grf_2\circ\grf_1=\text{id}_{W}$ and, hence, $W\simeq \tilde W$.

\subsection*{The group $G_{22}$}

\indent

Let $W$ be the group $G_{12}$, whose BMR presentation is of the form $$\langle s,t,u\;|\;s^2=t^2=u^2=1, stust=tustu=ustus\rangle.$$ We also know that the center of the group  is generated by the element $z=(stu)^5$ (see Tables in Appendix 1 in \cite{broue2000}). 
According to \S6.1 in \cite{ERrank2}, the ER presentation of $W$ is of the form
$$\langle a,b,c\;|\; a^2=1, b^3=\text{central}, c^5=b^{-6}, abc=1\rangle.$$ Let $\tilde W$ be a group having this ER presentation. We will prove that $\tilde W\simeq W$.
Let $\grf_1:W\rightarrow \tilde W$ and $\grf_2: \tilde W \rightarrow W$, defined by $
\grf_1(s)=a$, $\grf_1(t)=(b^{-1}cb^2)^{-1}$, $\grf_1(u)=(cb)^{-1}$ and $\grf_2(a)=s$, $\grf_2(b) =tustu$ and $\grf_2(c)=(stustu)^{-1}$. 

We now prove that $\grf_1$ is a well-defined group homomorphism.
Since $abc=1$, $b^3:=\grz$ and $c^{-5}=b^6$ we have:
$ab^{-2}c^{-2}ab^{-2}c^{-1}b=a\grz^{-1}bc^{-2}a\grz^{-1}bc^{-1}b=\grz^{-2}(abc)c^{-3}(abc)c^{-2}b=\grz^{-2}c^{-5}b=(\grz^{-2}b^6)b=b$. Moreover, $b^{-2}c^{-2}ab^{-2}c^{-2}=\grz^{-1}bc^{-2}a\grz^{-1}bc^{-2}=\grz^{-2}bc^{-2}(abc)c^{-3}=\grz^{-2}bc^{-5}=\grz^{-2}b^7=b$ and, analogously,  $b^{-1}c^{-1}ab^{-2}c^{-2}a=
b^{-1}c^{-1}a\grz^{-1}bc^{-2}a=\grz^{-1}b^{-1}c^{-1}(abc)c^{-3}(abc)c^{-1}b^{-1}=\grz^{-1}b^{-1}c^{-5}b^{-1}=\grz^{-1}b^4=b$. Hence, 
 $\grf_1(s)\grf_1(t)\grf_1(u)\grf_1(s)\grf_1(t)=\grf_1(t)\grf_1(u)\grf_1(s)\grf_1(t)\grf_1(u)=\grf_1(u)\grf_1(s)\grf_1(t)\grf_1(u)\grf_1(s)$. 
 Moreover, $\grf_1(s)^2=a^2=1$ and since $abc=bca=cab=1$ and $\gra=\gra^{-1}$ we also have $\grf_2(t)^2=b{-2}c^{-1}b^{-1}c^{-1}b=b^{-2}(c^{-1}b^{-1}a^{-1})a^{-1}c^{-1}b=b^{-2}(a^{-1})c^{-1}b^{-1})b^2=1$ and $\grf_1(u)^2=b^{-1}c^{-1}b^{-1}c^{-1}=b^{-1}(c^{-1}b^{-1}a^{-1})a^{-1}c^{-1}=(b^{-1}a^{-1}c^{-1})=1$. 
 
Similarly, we can prove that $\grf_2$ is a well-defined group homomorphism, since $\grf_2(a)^2=s^2=1$, $\grf_2(b)^3=(tustu)(tustu)(tustu)=(stust)(ustus)(tustu)=(stu)^5=z$,  $\grf_2(c)^{-5}=(stu)^{10}=z^2=\grf_2(b)^3$ and $\grf_2(a)\grf_2(b)\grf_2(c)=1$.

We also notice that $\grf_1(\grf_2(a))=a$ and since  $abc=1$, $b^3:=\grz$ and $c^{-5}=b^6$ we have: $\grf_1(\grf_2(b))=\grf_1(tustu)=b^{-2}c^{-2}ab^{-2}c^{-2}=\grz^{-1}bc^{-2}a\grz^{-1}bc^{-2}=\grz^{-2}bc^{-2}(abc)c^{-3}=\grz^{-2}bc^{-5}=\grz^{-2}b^7=b$. Moreover,  $\grf_1(\grf_2(c))=(ab^{-2}c^{-2}ab^{-2}c^{-2})^{-1}=(a\grz^{-1}bc^{-2}a\grz^{-1}c^{-2}=\grz^{-2}(abc)c^{-3}(abc)c^{-3})^{-1}=\grz^{2}c^{6}=c^{-5}c^6=c$. 
We also have that  $\grf_2(\grf_1(s))=s$, $\grf_2(\grf_1(t))=(tustu)^{-1}(tustu)^{-1}(stu)^2(tustu)=(ustus)^{-1}(tustu)^{-1}(tustu)u(stust)=t$ and $\grf_2(\grf_1(u))=(tustu)^{-1}(stust)u=(tustu)^{-1}(tustu)u=u$,  meaning that $\grf_1\circ\grf_2=\text{id}_{\tilde W}$ and $\grf_2\circ\grf_1=\text{id}_{W}$ and, hence, $W\simeq \tilde W$.
\newpage
~
\chapter{Tables}
\begin{table}[h]
	\begin{center}
		\small
		\caption{
			\bf{BMR and ER presentations for the complex braid groups associated to the exceptional groups of rank 2}}
		\label{t2}
		
		\scalebox{0.65}
		{\begin{tabular}{|c|c|c|c|c|c|}
				\hline
				Group & BMR presentation & ER presentation & $\grf_1$: BMR $\leadsto$ ER& $\grf_2$: ER $\leadsto$ BMR \\
				\hline
				$\begin{array}[t]{lcl} G_4\\ \\
				G_8\\ \\
				G_{16} 
				\end{array}$&
				$\begin{array}[t]{lcl} \\\\	\langle s,t\;|\; sts=tst\rangle\\
				\end{array}$& 
				$\begin{array}[t]{lcl} \\\\	\langle \gra,\grb,\grg\;|\; \gra^2=\grb^{-3}=\text{central},\gra\grb\grg=1\rangle \end{array}$&
				$\begin{array}[t]{lcl} \\
				
				s&\mapsto &\grg\\t&\mapsto& \grg^{-1}\grb
				\end{array}$ &
				$\begin{array}[t]{lcl}\\
				\gra&\mapsto& (sts)^{-1}\\\grb&\mapsto& st\\\grg&\mapsto&s
				\end{array}$\\
				
				\hline \hline
				$\begin{array}[t]{lcl} G_5\\ \\
				G_{10}\\ \\
				G_{18}
				\end{array}$&
				$\begin{array}[t]{lcl} \\\\	\langle s,t\;|\; stst=tsts\rangle\\
				\end{array}$& 
				$\begin{array}[t]{lcl} \\\\	\langle \gra,\grb,\grg\;|\; \gra^2=\text{central},\gra\grb\grg=1\rangle \end{array}$&
				$\begin{array}[t]{lcl} \\
				s&\mapsto &
				\grb\\t&\mapsto& \grg
				\end{array}$ &
				$\begin{array}[t]{lcl} \\
				\gra&\mapsto& (st)^{-1}\\\grb&\mapsto& s\\\grg&\mapsto&t
				\end{array}$\\
				
				\hline\hline
				
				$\begin{array}[t]{lcl} G_6\\ \\
				G_9\\ \\
				G_{17}
				\end{array}$&
				$\begin{array}[t]{lcl} \\\\	\langle s,t\;|\; (st)^3=(ts)^3\rangle\\
				\end{array}$& 
				$\begin{array}[t]{lcl} \\	\\\langle \gra,\grb,\grg\;|\; \grb^3=\text{central},\gra\grb\grg=1\rangle \end{array}$&
				$\begin{array}[t]{lcl}\\
				s&\mapsto &
				\gra\\t&\mapsto& \grg 
				\end{array}$ &
				$\begin{array}[t]{lcl}\\
				\gra&\mapsto& s\\\grb&\mapsto& (ts)^{-1}\\\grg&\mapsto&t
				\end{array}$\\
				
				\hline\hline
				$\begin{array}[t]{lcl} G_7\\ \\
				G_{11}\\ \\
				G_{19}
				\end{array}$&
				$\begin{array}[t]{lcl} \\\\	\langle s,t,u\;|\; stu=tus=ust\rangle\\
				\end{array}$& 
				$\begin{array}[t]{lcl} \\\\	\langle \gra,\grb,\grg\;|\;\gra\grb\grg=\text{central }\rangle \end{array}$&
				$\begin{array}[t]{lcl}\\
				s&\mapsto &
				\gra\\t&\mapsto& \grb\\
				u&\mapsto& \grg
				\end{array}$ &
				$\begin{array}[t]{lcl}\\
				\gra&\mapsto& s\\\grb&\mapsto& t\\\grg&\mapsto&u
				\end{array}$\\

				\hline\hline

				$\begin{array}[t]{lcl} \\
				G_{12}\\
				\end{array}$&
				$\begin{array}[t]{lcl} \\	\langle s,t,u\;|\; stus=tust=ustu\rangle\\
				\end{array}$& 
				$\begin{array}[t]{lcl} \\	\langle \gra,\grb,\grg\;|\; \grb^3=\grg^{-4}=\text{central},\gra\grb\grg=1\rangle \end{array}$&
				$\begin{array}[t]{lcl}
				s&\mapsto &
				\gra\\t&\mapsto& \gra^{-1}\grg^2\grb\\
				u&\mapsto&(\grg\grb)^{-1} 
				\end{array}$ &
				$\begin{array}[t]{lcl}
				\gra&\mapsto& s\\\grb&\mapsto& (stus)^{-1}\\\grg&\mapsto&stu
				\end{array}$\\ 
				
				\hline\hline
				$\begin{array}[t]{lcl} \\ \\
				G_{13}\\
				\end{array}$&
				
				$\begin{array}[t]{lcl}\\\\
				\left\langle
				\begin{array}{l|cl}
				s,t,u& stust=ustus&\\
				&tust=ustu&
				\end{array}\right \rangle
				\end{array}$&

				$\begin{array}[t]{lcl} \\
				\left\langle
				\begin{array}{l|cl}
				& \grb^3=\grd\grg^{-2}=\text{central}&\\
				\gra,\grb,\grg,\grd&\grg^4=\text{central}\\
				&\gra\grb\grg=1&
				\end{array}\right \rangle\\\\
				\end{array}$&

				$\begin{array}[t]{lcl}
				\\s&\mapsto &
				\grd\\t&\mapsto& \gra\\
				u&\mapsto&\grb(\grd\gra)^{-1} 
				\end{array}$ &
				$\begin{array}[t]{lcl}
				\gra&\mapsto& t\\\grb&\mapsto& ust\\\grg&\mapsto&(tust)^{-1}\\
				\grd&\mapsto&s
				\end{array}$\\
				
				\hline\hline
				
				$\begin{array}[t]{lcl} \\
				G_{14}\\
				\end{array}$&
				$\begin{array}[t]{lcl} \\	\langle s,t\;|\; (st)^4=(ts)^4\rangle\\
				\end{array}$& 
				$\begin{array}[t]{lcl} \\	\langle \gra,\grb,\grg\;|\; \grg^4=\text{central},\gra\grb\grg=1\rangle \end{array}$&
				$\begin{array}[t]{lcl}
				s&\mapsto &
				\gra\\t&\mapsto& \grb
				\end{array}$ &
				$\begin{array}[t]{lcl}
				\gra&\mapsto& s\\\grb&\mapsto& t\\\grg&\mapsto&(st)^{-1}
				\end{array}$\\
				
				\hline\hline
				$\begin{array}[t]{lcl} \\\\
				G_{15}\\
				\end{array}$&
				$\begin{array}[t]{lcl}\\\\
				\left\langle
				\begin{array}{l|cl}
				s,t,u& ustut=stutu&\\
				&tus=stu&
				\end{array}\right \rangle
				\end{array}$& 
				
				$\begin{array}[t]{lcl} \\
				\left\langle
				\begin{array}{l|cl}
				& \grd\grg^{-2}=\text{central}&\\
				\gra,\grb,\grg,\grd&\grg^4=\text{central}\\
				&\gra\grb\grg=1&
				\end{array}\right \rangle\\\\
				\end{array}$&

				$\begin{array}[t]{lcl}
				\\s&\mapsto &
				\grd\\t&\mapsto& \gra\\
				u&\mapsto&\grb
				\end{array}$ &
				$\begin{array}[t]{lcl}
				\gra&\mapsto& t\\\grb&\mapsto& u\\\grg&\mapsto&(tu)^{-1}\\
				\grd&\mapsto&s
				\end{array}$\\
				
				\hline\hline
				$\begin{array}[t]{lcl} \\
				G_{20}\\
				\end{array}$&
				$\begin{array}[t]{lcl} \\	\langle s,t\;|\; ststs=tstst\rangle\\
				\end{array}$& 
				$\begin{array}[t]{lcl} \\	\langle \gra,\grb,\grg\;|\; \gra^2=\text{central},\gra^4=\grg^{-5},\gra\grb\grg=1\rangle \end{array}$&
				$\begin{array}[t]{lcl}
				s&\mapsto &
				\grb\\t&\mapsto& (\gra\grg)^{-1}
				\end{array}$ &
				$\begin{array}[t]{lcl}
				\gra&\mapsto& (ststs)^{-1}\\\grb&\mapsto& s\\\grg&\mapsto&tsts
				\end{array}$\\
				\hline\hline
				$\begin{array}[t]{lcl} \\
				G_{21}\\
				\end{array}$&
				$\begin{array}[t]{lcl} \\	\langle s,t\;|\; (st)^5=(ts)^5\rangle\\
				\end{array}$& 
				$\begin{array}[t]{lcl} \\	\langle \gra,\grb,\grg\;|\; \grg^5=\text{central},\gra\grb\grg=1\rangle \end{array}$&
				$\begin{array}[t]{lcl}
				s&\mapsto &
				\gra\\t&\mapsto& \grb
				\end{array}$ &
				$\begin{array}[t]{lcl}
				\gra&\mapsto& s\\\grb&\mapsto& t\\\grg&\mapsto&(st)^{-1}
				\end{array}$\\
				\hline\hline
				$\begin{array}[t]{lcl} \\
				G_{22}\\
				\end{array}$&
				$\begin{array}[t]{lcl} \\	\langle s,t,u\;|\; stust=tustu=ustus\rangle\\
				\end{array}$& 
				$\begin{array}[t]{lcl} \\	\langle \gra,\grb,\grg\;|\; \grb^3=\text{central},\grg^5=\text{central},\gra\grb\grg=1\rangle \end{array}$&
				$\begin{array}[t]{lcl}
				s&\mapsto &
				\gra\\t&\mapsto& (\grb^{-1}\grg\grb^2)^{-1}\\
				u&\mapsto&(\grg\grb)^{-1}
				\end{array}$ &
				$\begin{array}[t]{lcl}
				\gra&\mapsto&s\\\grb&\mapsto& tustu\\\grg&\mapsto&(stustu)^{-1}
				\end{array}$\\
				\hline
			\end{tabular} }
		\end{center}
		
	\end{table}

\begin{table}[htp]
	\begin{center}
		\small
		\caption{
			\bf{BMR and ER presentations for the exceptional groups of rank 2}}
		\label{t3}
		
		\scalebox{0.73}
		{\begin{tabular}{|c|c|c|c|c|c|}
				\hline
				Group & BMR presentation & ER presentation & $\grf_1$: BMR $\leadsto$ ER& $\grf_2$: ER $\leadsto$ BMR \\
				\hline
				$\begin{array}[t]{lcl} \\G_4\\ \\
				G_8\\ \\
				G_{16} 
				\end{array}$&
				$\begin{array}[t]{lcl}\\	\langle s,t\;|\;s^3=t^3=1, sts=tst\rangle\\ \\
					\langle s,t\;|\;s^4=t^4=1, sts=tst\rangle\\ \\
						\langle s,t\;|\;s^5=t^5=1, sts=tst\rangle\\
				\end{array}$& 
				$\begin{array}[t]{lcl}\\	\langle a,b,c\;|\; a^2=b^{-3}=\text{central}, c^3=1, abc=1\rangle\\ \\
				\langle a,b,c\;|\; a^2=b^{-3}=\text{central}, c^4=1, abc=1\rangle\\ \\
				\langle a,b,c\;|\; a^2=b^{-3}=\text{central}, c^5=1, abc=1\rangle\\ \\
				 \end{array}$&
				$\begin{array}[t]{lcl} \\ \\
				
				s&\mapsto &c\\
				\\t&\mapsto& c^{-1}b 
				\end{array}$ &
				$\begin{array}[t]{lcl}\\
				a&\mapsto& (sts)^{-1}\\\\b&\mapsto& st\\\\c&\mapsto&s
				\end{array}$\\
				
				\hline \hline
				$\begin{array}[t]{lcl} \\G_5\\ \\
				G_{10}\\ \\
				G_{18}
				\end{array}$&
				$\begin{array}[t]{lcl}\\ \langle s,t\;|\; s^3=t^3=1,stst=tsts\rangle\\ \\
				\langle s,t\;|\; s^3=t^4=1,stst=tsts\rangle\\ \\
				\langle s,t\;|\; s^3=t^5=1,stst=tsts\rangle\\ \\
				\end{array}$& 
				$\begin{array}[t]{lcl}\\	\langle a,b,c\;|\; a^2=\text{central},b^3=c^3=1, abc=1\rangle \\ \\
					\langle a,b,c\;|\; a^2=\text{central},b^3=c^4=1, abc=1\rangle\\ \\
						\langle a,b,c\;|\; a^2=\text{central},b^3=c^5=1, abc=1\rangle \\ \\
						\end{array}$&
				$\begin{array}[t]{lcl} \\\\
				s&\mapsto &
				b\\\\t&\mapsto& c 
				\end{array}$ &
				$\begin{array}[t]{lcl} \\
				a&\mapsto& (st)^{-1}\\\\b&\mapsto& s\\\\c&\mapsto&t
				\end{array}$\\
				
				\hline \hline
				
				$\begin{array}[t]{lcl}\\ G_6\\ \\
				G_9\\ \\
				G_{17}
				\end{array}$&
				$\begin{array}[t]{lcl} \\\langle s,t\;|\; s^2=t^3=1,(st)^3=(ts)^3\rangle\\ \\
				\langle s,t\;|\; s^2=t^4=1,(st)^3=(ts)^3\rangle\\ \\
				\langle s,t\;|\; s^2=t^5=1,(st)^3=(ts)^3\rangle\\ \\
				\end{array}$& 
				$\begin{array}[t]{lcl} \\\langle a,b,c\;|\; a^2=c^3=1, b^3=\text{central},abc=1\rangle\\ \\
				\langle a,b,c\;|\; a^2=c^4=1, b^3=\text{central},abc=1\rangle\\ \\
				\langle a,b,c\;|\; a^2=c^5=1, b^3=\text{central},abc=1\rangle\\ \\
				 \end{array}
				$&
				$\begin{array}[t]{lcl}\\\\
				s&\mapsto &
				a\\\\t&\mapsto& c 
				\end{array}$ &
				$\begin{array}[t]{lcl}\\
				a&\mapsto& s\\\\b&\mapsto& (ts)^{-1}\\\\c&\mapsto&t
				\end{array}$\\
				
				\hline \hline
				$\begin{array}[t]{lcl} \\G_7\\ \\\\
				G_{11}\\ \\\\
				G_{19}
				\end{array}$&
				$\begin{array}[t]{lcl} \\	\left\langle
				\begin{array}{l|cl}
				\;s,t,u& s^2=t^3=u^3=1&\\
				&stu=tus=ust&
				\end{array}\right \rangle\\\\
		
					\left\langle
					\begin{array}{l|cl}
					\;s,t,u& s^2=t^3=u^4=1&\\
					&stu=tus=ust&
					\end{array}\right \rangle\\\\
					\left\langle
					\begin{array}{l|cl}
					\;s,t,u& s^2=t^3=u^5=1&\\
					&stu=tus=ust&
					\end{array}\right \rangle\\\\
				\end{array}$& 
				$\begin{array}[t]{lcl} \\	\langle a,b,c\;|\;a^2=b^3=c^3=1, abc=\text{central }\rangle\\\\\\
				\langle a,b,c\;|\;a^2=b^3=c^4=1, abc=\text{central }\rangle\\\\\\
				\langle a,b,c\;|\;a^2=b^3=c^5=1, abc=\text{central }\rangle\\\\
				 \end{array}$&
				$\begin{array}[t]{lcl}\\\\
				s&\mapsto &
				a\\\\t&\mapsto& b\\\\
				u&\mapsto& c
				\end{array}$ &
				$\begin{array}[t]{lcl}\\\\
				a&\mapsto& s\\\\b&\mapsto& t\\\\c&\mapsto&u
				\end{array}$\\

				\hline \hline

				$\begin{array}[t]{lcl} \\
				G_{12}\\
				\end{array}$&
				$\begin{array}[t]{lcl} \\	\left\langle
				\begin{array}{l|cl}
				\;s,t,u& s^2=t^2=u^2=1&\\
				&stus=tust=ustu&
				\end{array}\right \rangle\\\\
				\end{array}$& 
				$\begin{array}[t]{lcl} \\ \left	\langle \begin{array}{l|cl} &a^2=1&\\ a,b,c &b^3=c^{-4}=\text{central}&\\ &abc=1&
				\end{array}\right\rangle \\\\\end{array}$&
				$\begin{array}[t]{lcl}
				\\s&\mapsto &
				a\\t&\mapsto& a^{-1}c^2b\\
				u&\mapsto&(cb)^{-1} 
				\end{array}$ &
				$\begin{array}[t]{lcl}
				\\a&\mapsto& s\\b&\mapsto& (stus)^{-1}\\c&\mapsto&stu
				\end{array}$\\ 
				
				\hline\hline
				$\begin{array}[t]{lcl} \\ \\
				G_{13}\\\\
				\end{array}$&
				
				$\begin{array}[t]{lcl}\\
				\left\langle
				\begin{array}{l|cl}
				&s^2=t^2=u^2=1 &\\s,t,u&stust=ustus&\\
				&tust=ustu&
				\end{array}\right \rangle\\ \\
				\end{array}$& 
				
				$\begin{array}[t]{lcl} \\ \left	\langle \begin{array}{l|cl} &a^2=d^2=1&\\a,b,c,d &b^3=dc^{-2}=\text{central}&\\ &abc=1&
				\end{array}\right\rangle \\\\\end{array}$&

				$\begin{array}[t]{lcl}\\
				s&\mapsto &
				d\\t&\mapsto& a\\
				u&\mapsto&b(da)^{-1} 
				\end{array}$ &
				$\begin{array}[t]{lcl}
				\\a&\mapsto& t\\b&\mapsto& ust\\c&\mapsto&(tust)^{-1}\\
				d&\mapsto&s\\ 
				\end{array}$\\
				
				\hline\hline
				
				$\begin{array}[t]{lcl} \\
				G_{14}\\
				\end{array}$&
				$\begin{array}[t]{lcl} \\	\langle s,t\;|\; s^2=t^3=1,(st)^4=(ts)^4\rangle\\
				\end{array}$& 
				$\begin{array}[t]{lcl} \\	\langle a,b,c\;|\;a^2=b^3=1, c^4=\text{central},abc=1\rangle \end{array}$&
				$\begin{array}[t]{lcl}
				s&\mapsto &
				a\\t&\mapsto& b
				\end{array}$ &
				$\begin{array}[t]{lcl}
				a&\mapsto& s\\b&\mapsto& t\\c&\mapsto&(st)^{-1}
				\end{array}$\\
				
				\hline\hline
				$\begin{array}[t]{lcl} \\\\
				G_{15}\\
				\end{array}$&
				$\begin{array}[t]{lcl}\\
				\left\langle
				\begin{array}{l|cl}
				& s^2=t^2=u^3=1&\\s,t,u&ustut=stutu&\\
				&tus=stu&
				\end{array}\right \rangle\\\\
				\end{array}$& 
				
				$\begin{array}[t]{lcl} \\ \left	\langle \begin{array}{l|cl} &a^2=b^3=d^2=1&\\ a,b,c,d &dc^{-2}=\text{central}&\\ &abc=1&
				\end{array}\right\rangle \\\\\end{array}$&
				
				$\begin{array}[t]{lcl}
				\\s&\mapsto &
				d\\t&\mapsto& a\\
				u&\mapsto&b
				\end{array}$ &
				$\begin{array}[t]{lcl}
				a&\mapsto& t\\b&\mapsto& u\\c&\mapsto&(tu)^{-1}\\
				d&\mapsto&s
				\end{array}$\\
				
				\hline \hline
				$\begin{array}[t]{lcl} \\\\
				G_{20}\\
				\end{array}$&
				$\begin{array}[t]{lcl} \\\\	\langle s,t\;|\;s^3=t^3=1, ststs=tstst\rangle\\
				\end{array}$& 
				
				$\begin{array}[t]{lcl} \\ \left	\langle \begin{array}{l|cl} &a^2=\text{central}&\\
				a,b,c&b^3=1,a^4=c^{-5}&\\ &abc=1&
				\end{array}\right\rangle \\\\\end{array}$&

				$\begin{array}[t]{lcl}
				\\s&\mapsto &
				b\\\\t&\mapsto& (ac)^{-1}
				\end{array}$ &
				$\begin{array}[t]{lcl}
				\\a&\mapsto& (ststs)^{-1}\\b&\mapsto& s\\c&\mapsto&tsts
				\end{array}$\\
				\hline\hline
				$\begin{array}[t]{lcl} \\ 
				G_{21}\\
				\end{array}$&
				$\begin{array}[t]{lcl} \\	\langle s,t\;|\; s^2=t^3=1,(st)^5=(ts)^5\rangle\\
				\end{array}$& 
				
					$\begin{array}[t]{lcl} \\	\langle a,b,c\;|\;a^2=b^3=1, c^5=\text{central},abc=1\rangle \end{array}$&
				
				$\begin{array}[t]{lcl}
				s&\mapsto &
				a\\t&\mapsto& b
				\end{array}$ &
				$\begin{array}[t]{lcl}
				a&\mapsto& s\\b&\mapsto& t\\c&\mapsto&(st)^{-1}
				\end{array}$\\
				\hline \hline
				$\begin{array}[t]{lcl} \\
				G_{22}\\
				\end{array}$&
				$\begin{array}[t]{lcl} \\	\left\langle \begin{array}{l|cl}s,t,u&s^2=t^2=u^2=1&\\ &stust=tustu=ustus&
				\end{array}
				
				\right\rangle\\
				\end{array}$& 
				$\begin{array}[t]{lcl} \\\left\langle 
				\begin{array}{l|cl}
				 &b^3=\text{central}&\\a,b,c&a^2=1,c^5=b^{-6}&\\&abc=1&\end{array}\right\rangle\\\\ \end{array}$&
				$\begin{array}[t]{lcl}
				\\s&\mapsto &
				a\\t&\mapsto& (b^{-1}cb^2)^{-1}\\
				u&\mapsto&(cb)^{-1}
				\end{array}$ &
				$\begin{array}[t]{lcl}
				\\a&\mapsto&s\\b&\mapsto& tustu\\c&\mapsto&(stustu)^{-1}
				\end{array}$\\
				\hline
			\end{tabular} }
		\end{center}
		
	\end{table}


\end{appendices}
\renewcommand\bibname{Bibliography}

\addcontentsline{toc}{chapter}{Bibliography}

	\newpage
		~
		\clearpage
		\thispagestyle{empty}
		\newpage
		~
		\clearpage	

\begin{thebibliography}{9999}

\bibitem{ariki}
Susumu Ariki.
\newblock Representation theory of a {H}ecke algebra of ${G}(r, p, n)$.
\newblock {\em Journal of algebra}, 177(1):164--185, 1995.

\bibitem{arikii}
Susumu Ariki and Kazuhiko Koike.
\newblock A {H}ecke algebra of ($\mathbb{Z}/r\mathbb{Z})\wr {S}_n$ and
  construction of its irreducible representations.
\newblock {\em Advances in Mathematics}, 106(2):216--243, 1994.

\bibitem{bfunar}
Paolo Bellingeri and Louis Funar.
\newblock Polynomial invariants of links satisfying cubic skein relations.
\newblock {\em Asian Journal of Mathematics}, 8(3):475--510, 2004.

\bibitem{berenstein}
Carlos~A Berenstein and Roger Gay.
\newblock {\em Complex variables: an introduction}, volume 125.
\newblock Springer Science \& Business Media, 2012.

\bibitem{bessiszariski}
David Bessis.
\newblock Zariski theorems and diagrams for braid groups.
\newblock {\em Inventiones mathematicae}, 145(3):487--507, 2001.

\bibitem{bessiscenter}
David Bessis.
\newblock Finite complex reflection arrangements are ${K}(\pi,1)$.
\newblock {\em Annals of Mathematics}, 181(3):809--904, 2015.

\bibitem{Brieskorn}
Egbert Brieskorn.
\newblock Die {F}undamentalgruppe des {R}aumes der regul{\"a}ren {O}rbits einer
  endlichen komplexen {S}piegelungsgruppe.
\newblock {\em Inventiones mathematicae}, 12(1):57--61, 1971.

\bibitem{broue2000}
Michel Brou{\'e}.
\newblock Reflection groups, {B}raid groups, {H}ecke algebras, finite reductive
  groups.
\newblock In {\em Current developments in mathematics}. Citeseer, 2000.

\bibitem{brouebook}
Michel Brou{\'e}.
\newblock {\em Introduction to complex reflection groups and their braid
  groups}.
\newblock Springer, 2010.

\bibitem{brouem}
Michel Brou{\'e} and Gunter Malle.
\newblock Zyklotomische {H}eckealgebren in {R}epr{\'e}sentations unipotentes
  g{\'e}n{\'e}riques et blocs des groupes r{\'e}ductifs finis.
\newblock {\em Ast{\'e}risque}, 212:119--189, 1993.

\bibitem{bmm}
Michel Brou{\'e}, Gunter Malle, and Jean Michel.
\newblock Towards {S}petses {I}.
\newblock {\em Transformation groups}, 4(2-3):157--218, 1999.

\bibitem{bmr}
Michel Brou{\'e}, Gunter Malle, and Rapha{\"e}l Rouquier.
\newblock Complex reflection groups, braid groups, {H}ecke algebras.
\newblock {\em Journal fur die Reine und Angewandte Mathematik}, 500:127--190,
  1998.

\bibitem{chavli}
Eirini Chavli.
\newblock Universal deformations of the finite quotients of the braid group on
  3 strands.
\newblock {\em To appear in Journal of Algebra}.

\bibitem{chlouverakibook}
Maria Chlouveraki.
\newblock {\em Blocks and families for cyclotomic {H}ecke algebras}.
\newblock Springer, 2009.

\bibitem{chlouveraki}
Maria Chlouveraki and Hyohe Miyachi.
\newblock Decomposition matrices for d-{H}arish-{C}handra series: the
  exceptional rank two cases.
\newblock {\em LMS Journal of Computation and Mathematics}, 14:271--290, 2011.

\bibitem{cohen}
Arjeh~M Cohen.
\newblock Finite complex reflection groups.
\newblock In {\em Annales Scientifiques de L'{\'e}cole Normale Sup{\'e}rieure},
  volume~9, pages 379--436, 1976.

\bibitem{Coxeter}
HSM Coxeter.
\newblock Factor groups of the braid group.
\newblock In {\em Proc. Fourth Canadian Math. Congress, Banff}, pages 95--122,
  1957.

\bibitem{dignecenter}
Fran{\c{c}}ois Digne, Ivan Marin, and Jean Michel.
\newblock The center of pure complex braid groups.
\newblock {\em Journal of Algebra}, 347(1):206--213, 2011.

\bibitem{ERcoxeter}
Pavel Etingof and Eric Rains.
\newblock New deformations of group algebras of {C}oxeter groups.
\newblock {\em International Mathematics Research Notices}, 2005(10):635--646,
  2005.

\bibitem{ERrank2}
Pavel Etingof and Eric Rains.
\newblock Central extensions of preprojective algebras, the quantum
  {H}eisenberg algebra, and 2-dimensional complex reflection groups.
\newblock {\em Journal of Algebra}, 299(2):570--588, 2006.

\bibitem{funar1995}
Louis Funar.
\newblock On the quotients of cubic {H}ecke algebras.
\newblock {\em Communications in mathematical physics}, 173(3):513--558, 1995.

\bibitem{geck}
Meinolf Geck and G{\"o}tz Pfeiffer.
\newblock {\em Characters of finite {C}oxeter groups and {I}wahori-{H}ecke
  algebras}.
\newblock Number~21. Oxford University Press, 2000.

\bibitem{Ocat}
Victor Ginzburg, Nicolas Guay, Eric Opdam, and Rapha{\"e}l Rouquier.
\newblock On the category $\mathcal{O}$ for rational {C}herednik algebras.
\newblock {\em Inventiones mathematicae}, 154(3):617--651, 2003.

\bibitem{turaev}
Christian Kassel and Vladimir Turaev.
\newblock Braid groups, volume 247.
\newblock {\em Graduate Texts in Mathematics}, 2008.

\bibitem{lehrer}
Gustav~I. Lehrer and Donald~E. Taylor.
\newblock {\em Unitary reflection groups}, volume~20.
\newblock Cambridge University Press, 2009.

\bibitem{malle2}
Gunter Malle.
\newblock Degr{\'e}s relatifs des alg\`ebres cyclotomiques associ{\'e}es aux
  groupes de r{\'e}flexions complexes de dimension deux.
\newblock {\em Progress in {M}ath.}, 141:311--332, 1996.

\bibitem{malle1}
Gunter Malle.
\newblock On the {R}ationality and {F}ake {D}egrees of {C}haracters of
  {C}yclotomic {A}lgebras.
\newblock {\em Journal of Mathematical Sciences-University of Tokyo},
  6(4):647--678, 1999.

\bibitem{mallem}
Gunter Malle and Jean Michel.
\newblock Constructing representations of {H}ecke algebras for complex
  reflection groups.
\newblock {\em LMS Journal of Computation and Mathematics}, 13:426--450, 2010.

\bibitem{maller}
Gunter Malle and Rapha{\"e}l Rouquier.
\newblock Familles de caract{\`e}res de groupes de r{\'e}flexions complexes.
\newblock {\em Representation Theory of the American Mathematical Society},
  7(23):610--640, 2003.

\bibitem{marinreport}
Ivan Marin.
\newblock Reprot on the {B}rou\'e-{M}alle-{R}ouquier conjectures.
\newblock {\em To appear in proceedings of the Indam intensive period
  Perspectives in Lie theory}.

\bibitem{marincubic}
Ivan Marin.
\newblock The cubic {H}ecke algebra on at most 5 strands.
\newblock {\em Journal of Pure and Applied Algebra}, 216(12):2754--2782, 2012.

\bibitem{ivancourse}
Ivan Marin.
\newblock Groupes de tresses et groupes de r\'eflexions complexes.
\newblock University Lectures, Universit\'e Paris Diderot, 2012.

\bibitem{marinkrammer}
Ivan Marin.
\newblock Krammer representations for complex braid groups.
\newblock {\em Journal of Algebra}, 371:175--206, 2012.

\bibitem{marinG26}
Ivan Marin.
\newblock The freeness conjecture for {H}ecke algebras of complex reflection
  groups and the case of the {H}essian group {$G_{26}$}.
\newblock {\em Journal of Pure and Applied Algebra}, 218(4):704--720, 2014.

\bibitem{marinpfeiffer}
Ivan Marin and G\"otz Pfeiffer.
\newblock The {BMR} freeness conjecture for the 2-reflection groups.
\newblock {\em arXiv preprint: 1411.4760v2}, 2015.

\bibitem{massey}
William~S Massey.
\newblock {\em A basic course in algebraic topology}, volume 127.
\newblock Springer Science \& Business Media, 1991.

\bibitem{Matsumoto}
Hideya Matsumoto.
\newblock G{\'e}n{\'e}rateurs et relations des groupes de {W}eyl
  g{\'e}n{\'e}ralis{\'e}s.
\newblock {\em Comptes {R}endus hebdomadaires des s\'eances de l'{A}cad\'emie
  des sciences}, 258(13):3419, 1964.

\bibitem{Michelcour}
Jean Michel.
\newblock Groupes finis de r\'eflexion.
\newblock University notes, Universit\'e Paris Diderot, 2004.
\newblock \url{http://webusers.imj-prg.fr/~jean.michel/papiers/cours2004.pdf}.

\bibitem{michelgap}
Jean Michel.
\newblock The development version of the {CHEVIE} package of {GAP3}.
\newblock {\em Journal of Algebra}, 435:308--336, 2015.

\bibitem{shephard}
Geoffrey~C Shephard and John~A Todd.
\newblock Finite unitary reflection groups.
\newblock {\em Canad. J. Math}, 6(2):274--301, 1954.

\bibitem{steinberg}
Robert Steinberg.
\newblock Differential equations invariant under finite reflection groups.
\newblock {\em Transactions of the American Mathematical Society}, pages
  392--400, 1964.

\bibitem{tuba}
Imre Tuba and Hans Wenzl.
\newblock Representations of the braid group {$B_3$} and of {$SL
  (2,\mathbb{Z})$}.
\newblock {\em Pacific J. Math}, 197(2):491--510, 2001.

\end{thebibliography}
	\end{document}